\documentclass[11pt]{article}

\usepackage[margin=1in]{geometry}
\usepackage[dvipsnames]{xcolor}
\usepackage[breakable]{tcolorbox}
\usepackage{mathtools}
\usepackage{enumitem}
\usepackage{amsmath}
\usepackage{amssymb}
\usepackage{amsthm}
\usepackage{tikz}
\usepackage{bbm}
\usepackage{subcaption}
\usepackage{placeins}
\usepackage[expansion=false,final]{microtype}
\usepackage[longnamesfirst]{natbib}
\usepackage[hang,flushmargin]{footmisc}
\usepackage[hypertexnames=false]{hyperref}
\usepackage[norefs,nocites]{refcheck}
\usepackage{multirow}
\usepackage{arydshln}
\usepackage[T1]{fontenc}
\hypersetup{colorlinks=true,
linkcolor=blue,urlcolor=blue,citecolor=blue}
\DeclareTextSymbolDefault{\ae}{T1}
\usetikzlibrary{arrows.meta}
\DeclareMathOperator*{\argmin}{arg\,min}
\DeclareMathOperator{\diag}{diag}
\DeclareMathOperator{\rank}{rank}
\DeclareMathOperator{\tr}{tr}
\DeclareMathOperator{\esssup}{ess\,sup}
\DeclareMathOperator{\essinf}{ess\,inf}
\DeclareMathOperator{\Exp}{Exp}
\DeclareMathOperator{\vecmat}{vec}
\newcommand{\equallaw}{%
\mathrel{\smash{\overset{\mathcal L}{=}}\vphantom{\rule{1pt}{2.5ex}}}}
\renewcommand{\P}{\ensuremath{\mathbb{P}}}
\newcommand{\E}{\ensuremath{\mathbb{E}}}
\newcommand{\cF}{\ensuremath{\mathcal{F}}}
\newcommand{\cP}{\ensuremath{\mathcal{P}}}
\newcommand{\Fr}{\ensuremath{\mathrm{F}}}
\DeclareMathOperator{\Cov}{Cov}
\DeclareMathOperator{\Var}{Var}
\DeclareMathOperator{\Poisson}{Poisson}
\newcommand{\cV}{\ensuremath{\mathcal{V}}}
\newcommand{\cA}{\ensuremath{\mathcal{A}}}
\newcommand{\cB}{\ensuremath{\mathcal{B}}}
\newcommand{\bX}{\ensuremath{\mathbf{X}}}
\newcommand{\bY}{\ensuremath{\mathbf{Y}}}
\newcommand{\bW}{\ensuremath{\mathbf{W}}}
\newcommand{\cH}{\ensuremath{\mathcal{H}}}
\newcommand{\N}{\ensuremath{\mathbb{N}}}
\newcommand{\cN}{\ensuremath{\mathcal{N}}}
\newcommand{\I}{\ensuremath{\mathbbm{1}}}
\newcommand{\cG}{\ensuremath{\mathcal{G}}}
\newcommand{\cX}{\ensuremath{\mathcal{X}}}
\newcommand{\cC}{\ensuremath{\mathcal{C}}}
\newcommand{\cM}{\ensuremath{\mathcal{M}}}
\newcommand{\cW}{\ensuremath{\mathcal{W}}}
\newcommand{\cQ}{\ensuremath{\mathcal{Q}}}
\newcommand{\cJ}{\ensuremath{\mathcal{J}}}
\newcommand{\R}{\ensuremath{\mathbb{R}}}
\newcommand{\T}{\ensuremath{\mathsf{T}}}
\newcommand{\Z}{\ensuremath{\mathbb{Z}}}
\newcommand{\logplus}{\log^+}
\newcommand{\Clp}{\ensuremath{C_{\mathrm{lp}}}}
\renewcommand{\c}{\ensuremath{\mathsf{c}}}
\newcommand{\dprime}{\prime\prime}
\newcommand{\tprime}{\prime\prime\prime}
\newlist{inlineroman}{enumerate*}{1}
\setlist[inlineroman]{afterlabel=~,label=(\roman*)}

\setlist[enumerate,1]{label=(\alph*)}

\newcommand{\diffi}{\,\mathrm{d}}

\newcommand{\vvvert}{{\vert\kern-0.25ex\vert\kern-0.25ex\vert}}
\newcommand{\bigvvvert}{{\big\vert\kern-0.35ex\big\vert\kern-0.35ex\big\vert}}
\newcommand{\Bigvvvert}{{\Big\vert\kern-0.3ex\Big\vert\kern-0.3ex\Big\vert}}
\newcommand{\biggvvvert}{{\bigg\vert\kern-0.3ex\bigg\vert\kern-0.3ex\bigg\vert}}
\newcommand{\Biggvvvert}{{\Bigg\vert\kern-0.3ex\Bigg\vert\kern-0.3ex\Bigg\vert}}

\renewcommand{\arraystretch}{1.25}

\theoremstyle{plain}
\newtheorem{theorem}{Theorem}[section]
\newtheorem{lemma}[theorem]{Lemma}
\newtheorem{corollary}[theorem]{Corollary}
\newtheorem{definition}[theorem]{Definition}
\newtheorem{example}[theorem]{Example}
\newtheorem{assumption}[theorem]{Assumption}

\newcommand{\myfigscale}{0.9}
 
\newtoggle{journal}
\togglefalse{journal}

\makeatletter%
\begin{document}

\title{Higher-order Yurinskii coupling}

\author{
  Matias D.\ Cattaneo\textsuperscript{1}
  \and Ricardo P.\ Masini\textsuperscript{2}
  \and William G.\ Underwood\textsuperscript{3*}
}

\maketitle

\footnotetext[1]{
  Department of Economics,
  Princeton University.
}
\footnotetext[2]{
  Department of Statistics and Data Science,
  University of California, Davis.
}
\footnotetext[3]{
  Department of Statistics,
  University of Warwick.
}
\let\thefootnote\relax
\footnotetext[1]{
  \textsuperscript{*}Corresponding author:
  \href{mailto:william.underwood@warwick.ac.uk}{%
  \texttt{william.underwood@warwick.ac.uk}}
}
\newcommand{\thefootnote}{\arabic{footnote}}

\setcounter{page}{0}\thispagestyle{empty}

\begin{abstract}
  {\normalsize %
Yurinskii's coupling is a popular and widely applicable result
establishing Gaussian strong approximations for sums
of independent random vectors.
We present a higher-order generalization of Yurinskii's method,
approximately matching all moments up to order $k \in \N$ rather than
solely the mean and covariance.
The resulting couplings permit non-Gaussian target distributions
and sharpen classical bounds on the $\ell_p$-coupling error.
A principal motivation for our work is to study
nonparametric and multiplier bootstrap
procedures for statistical inference.
Our approximate moment matching framework is well suited
to such analyses because both of these bootstrap schemes reproduce
certain empirical moments, thereby enabling refined coupling bounds.
We use our main results to construct strong approximations
for empirical processes,
obtaining explicit bootstrap coupling rates in VC-type classes.
As an application, we establish the asymptotic validity
of uniform bootstrap confidence bands for local polynomial
estimators under weaker bandwidth conditions than those
required by existing second-order Yurinskii results.
A detailed literature survey compares a wide range of couplings,
and we identify settings in which our bounds are the best available.
Further contributions include Prokhorov metric bounds,
distributional approximations for maximum statistics
and anti-concentration inequalities for Gaussian process suprema.
}
\end{abstract}

\vspace*{10mm}
\noindent\textbf{Keywords}:
Coupling;
strong approximation;
bootstrapping;
distribution theory;
empirical processes;
nonparametric inference;
local polynomial methods.

\vspace*{4mm}
\noindent\textbf{MSC}:
Primary
62E17, %
62F40; %
Secondary
60F15. %

\clearpage
\pagebreak

\begingroup
\pagestyle{empty}

\linespread{0.95}\selectfont
\makeatletter
\patchcmd{\l@section}{\addvspace{1.0em \@plus\p@}}%
{\addvspace{0.85em \@plus\p@}}{}{}
\makeatother
\tableofcontents

\clearpage
\endgroup
\pagebreak
\setcounter{page}{1}

\section{Introduction}
\label{sec:introduction}

Let $X_1, \ldots, X_n$ be independent
random vectors taking values in $\R^d$ and
write $S_n \coloneqq \sum_{i=1}^n X_i$.
Yurinskii's coupling constructs,
on an extension of the probability space, a second sum
$T_n \coloneqq \sum_{i=1}^n Y_i$ with independent summands such that
$S_n$ and $T_n$ are close with high probability
in Euclidean norm \citep{yurinskii1978error,pollard2002user}
or, more generally, in $\ell_p$-norm
\iftoggle{journal}{(see \citep{le1988prokhorov} and
\citep{belloni2019conditional}, Supplemental Appendix~I.1)}%
{\cite[Supplemental Appendix~I.1]{le1988prokhorov,belloni2019conditional}},
where $p \in [1, \infty]$.
In these classical formulations,
each approximating summand $Y_i$ is Gaussian
with the same mean vector and covariance matrix as $X_i$.
Recently, \citet{cattaneo2025yurinskii} generalized these results
to approximate martingale data and showed how
vanishing conditional third moments may lead to an improved bound.
Yurinskii's coupling stands out for its flexibility and robustness,
with classical versions requiring only independent
summands with finite third moments.
As such, it finds applications ranging from uniform
inference in nonparametric regression
to the analysis of stochastic gradient descent
(see Section~\ref{sec:literature} and references therein).
However, Yurinskii's method often yields coupling error bounds
which decay rather slowly as $n \to \infty$.
This paper studies the role of higher-order moment matching
(beyond the mean and covariance)
in Yurinskii-type coupling constructions.
To this end, for $\alpha \in \N_0^d$ and $i \in [n]$,
define the moment discrepancy
\begin{equation*}
  \delta_{\alpha,i}
  \coloneqq
  \E(X_i^\alpha)-\E(Y_i^\alpha),
\end{equation*}
where $x^\alpha \coloneqq \prod_{j=1}^d x_j^{\alpha_j}$,
assuming the expectations exist and are finite.
Write $|\alpha| \coloneqq \sum_{j=1}^d \alpha_j$.

In Section~\ref{sec:yurinskii}, Theorem~\ref{thm:yurinskii}
gives a non-asymptotic bound on the $\ell_p$-coupling error
between the laws of $S_n$ and $T_n$
in terms of moment bounds, the dimension $d$,
the discrepancies $\delta_{\alpha,i}$, and their averages
$\bar\delta_\alpha \coloneqq \sum_{i=1}^n \delta_{\alpha,i} / n$.
In classical formulations, $\delta_{\alpha,i} = 0$ whenever
$|\alpha| \in \{1, 2\}$,
corresponding to matching the mean vector and covariance matrix;
our result allows for
approximate (average) moment matching up to arbitrary order,
considering $|\alpha| \in [k]$ where $k \in \N$.
For skewed summands, taking $k \geq 3$
generally calls for a non-Gaussian approximation,
in order to replicate the non-vanishing third moments.
In a canonical setting with coordinatewise
bounded moments of order $k + \gamma$, where $\gamma \in (0, 1]$,
and assuming exact moment matching up to order $k$,
we obtain a coupling rate of $n^{-1/2 + 1/(k+\gamma)}$
for the normalized sums $S_n / \sqrt n$ and $T_n / \sqrt n$,
up to an arbitrarily slowly diverging factor
and terms involving the dimension $d$.
When $k + \gamma > 3$,
this improves upon the classical Yurinskii rate of $n^{-1/6}$,
which corresponds to taking $k = 2$ and $\gamma = 1$
under bounded third moments.

In Section~\ref{sec:bootstrap}, we apply our higher-order
Yurinskii coupling to nonparametric and multiplier bootstrap sums.
Under these resampling schemes,
the marginal distribution of the centered sum
$S_n \coloneqq \sum_{i=1}^n \{X_i - \E(X_i)\}$ is compared to the
conditional law of the bootstrap statistic
$S_n^* \coloneqq \sum_{i=1}^n Y_i^*$ given the data.
Under the nonparametric bootstrap, we define
$Y_i^* \coloneqq X_i^* - \bar X_n$,
with $\bar X_n \coloneqq \sum_{i=1}^n X_i / n$
and $X_1^*, \ldots, X_n^*$
drawn conditionally i.i.d.\ from the empirical distribution
of $X_1, \ldots, X_n$.
For the multiplier bootstrap, we take $Y_i^* \coloneqq \chi_i (X_i - \bar X_n)$
with independent multiplier weights $\chi_i$ satisfying
$\E(\chi_i) = 0$ and $\E(\chi_i^\ell) = 1$ for each $\ell \in \{2, k\}$,
where $k \in \{2, 3\}$.
The resulting moment discrepancies $\hat\delta_{\alpha,i}$ are
zero when $|\alpha| = 1$ and are generally random when $|\alpha| \geq 2$.
For $|\alpha| \in \{2, k\}$, the averaged discrepancy may be written as
\begin{equation}
  \label{eq:delta_hat}
  \hat \delta_\alpha
  \coloneqq
  \frac{1}{n}
  \sum_{i=1}^n
  \hat \delta_{\alpha,i}
  = \frac{1}{n}
  \sum_{i=1}^n
  \Bigl(
    \E\bigl\{\bigl(X_i - \E(X_i)\bigr)^\alpha\bigr\}
    - (X_i - \bar X_n)^\alpha
  \Bigr).
\end{equation}
Arguments based on laws of large numbers show that
$\hat\delta_\alpha \to_\P 0$ under appropriate conditions.
Our Theorem~\ref{thm:bootstrap} proves that the nonparametric
and third-order multiplier bootstraps both achieve a
coupling rate of $n^{-1/4}$, up to an arbitrarily slowly diverging factor,
for the normalized sums $S_n / \sqrt n$ and $S_n^* / \sqrt n$
when the dimension is fixed, assuming a common mean vector
and bounded fourth moments.
Here, our formulation based on average moment discrepancies is crucial;
there is no guarantee that $\hat\delta_{\alpha,i} \to_\P 0$
for each $i \in [n]$
when using the multiplier bootstrap or the nonparametric bootstrap
with heterogeneous data.
Theorem~\ref{thm:refined_bootstrap} shows how the coupling error
for the nonparametric bootstrap
may be further reduced under additional homoskedasticity,
covariance spectrum, and sixth moment assumptions,
attaining a coupling rate of $n^{-1/3}$
for the normalized sums, again up to a slowly diverging factor.
All of our bootstrap coupling targets are fully feasible
in the sense that sampling from the bootstrap law
does not require knowledge of any population moments.

Section~\ref{sec:literature} gives a detailed comparison
of our main results with three strands of the literature:
\begin{inlineroman}
  \item Yurinskii-based couplings;
  \item the Koml\'os--Major--Tusn\'ady (KMT) method; and
  \item Wasserstein bounds.
\end{inlineroman}
In Table~\ref{tab:literature}, we quantify the precise dependence of
various bounds on the sample size $n$ and dimension $d$
in several canonical settings.
For example, when $\|X_i\|_2 + \|Y_i\|_2$ is uniformly bounded,
matching more moments improves the sample size exponent
of our $\ell_\infty$-coupling bound,
while its dependence on $d$ remains logarithmic.
We identify settings in which our coupling bounds
are the best available
and provide an in-depth discussion of some recent papers
especially close to our work.

In Section~\ref{sec:empirical_process}, we turn our attention
to strong approximation of empirical processes indexed by function classes
\citep{van1996weak,gine2016mathematical}.
Theorem~\ref{thm:empirical_process}
presents a bound on the supremum norm coupling error between
two empirical processes in an analogous approximate (average)
moment matching framework.
Taking $\gamma = 1$ in Corollary~\ref{cor:empirical_process_vc},
we attain a coupling rate of $n^{(1-k) / (2k + 2 + u k)}$,
up to a logarithmic factor, for the normalized processes
$S_n(f) / \sqrt n$ and $T_n(f) / \sqrt n$
when exactly matching centered moments up to order $k$
for each observation
in a VC-type class with covering exponent $u \geq 0$
and sub-Weibull envelope.
We then provide coupling results for bootstrapped empirical processes:
Theorem~\ref{thm:bootstrap_vc} shows that for both the
nonparametric and third-order multiplier bootstraps,
the normalized processes $S_n(f) / \sqrt n$ and $S_n^*(f) / \sqrt n$
admit a coupling rate of $n^{-1/(4 + 2u)}$,
up to a logarithmic factor, under the same envelope and VC-type conditions,
and assuming a common mean function.
This improves upon the $n^{-1/(6 + 2u)}$ rate under classical second-order
moment matching, without requiring vanishing third moments
or knowledge of any population parameters.
To streamline the presentation,
Appendix~\ref{sec:literature_empirical_process}%
\iftoggle{journal}{ in the supplementary material
\citep{cattaneo2026highsupplement}}{}
compares our empirical
process results with earlier strong approximations in the literature.

Section~\ref{sec:locpol} presents a substantive
statistical application of our bootstrap-based coupling theory.
We consider local polynomial estimators
for nonparametric regression or conditional distribution
function estimation \citep{eggermont2009maximum,samworth2026modern},
constructing confidence bands for the estimand
uniformly over a selected subset of the covariates.
Nonparametric and third-order multiplier bootstrap procedures
both provide asymptotically nominal uniform coverage
whenever the bias is negligible and $n h^{D_W + D_w} \to \infty$
(up to logarithmic factors and under regularity conditions),
where $h$ is the smoothing bandwidth,
$D_W$ is the full covariate dimension,
and $D_w$ is the number of covariates in the selected subset.
Under second-order moment matching,
the corresponding bandwidth condition is
$n h^{D_W + 2 D_w} \to \infty$.
When $D_w \geq 1$,
third-order matching therefore relaxes our sufficient condition on the
H{\"o}lder smoothness order~$\zeta$ from $\zeta > D_w$ to $\zeta > D_w/2$.
Our approach permits non-identically distributed
data points and unbounded sub-Weibull regression errors.

Finally, the
\iftoggle{journal}{supplementary material
\citep{cattaneo2026highsupplement} collects}{%
appendices collect}
proof sketches, additional corollaries,
further comparisons with the literature,
and full derivations of our main results.
Our Prokhorov metric bounds,
distributional approximations for maximum statistics,
empirical process maximal inequalities
and anti-concentration inequalities for Gaussian suprema
may be of independent interest.

\iftoggle{journal}{}{\pagebreak}

\subsection*{Notation}

We write $\N \coloneqq \{1, 2, \ldots\}$ for the natural numbers
and $\N_0 \coloneqq \{0\} \cup \N$.
For $n \in \N$ we write $[n] \coloneqq \{1, 2, \ldots, n\}$,
with $[0] \coloneqq \emptyset$.
For $x, y \in \R$, write $x \lor y \coloneqq \max\{x, y\}$
and $x \land y \coloneqq \min\{x, y\}$.
Indicator functions are denoted by $\I$.
For $d \in \N$, the $\ell_p$-norm of $x \in \R^d$ is
$\|x\|_p \coloneqq \bigl(\sum_{j=1}^d |x_j|^p\bigr)^{1/p}$
for $p \in [1, \infty)$ and
$\|x\|_\infty \coloneqq \max_{j \in [d]} |x_j|$.
We write $0_d$ for the $d$-dimensional zero vector
and $0_{d \times d}$ for the $d \times d$ zero matrix;
the $d \times d$ identity matrix is denoted by $I_d$.
For $A \in \R^{d \times d}$ and $p, q \in [1, \infty]$, write
$\|A\|_{p \to q} \coloneqq
\sup\bigl\{\|A x\|_q : x \in \R^d, \ \|x\|_p \leq 1\bigr\}$.
The $j$th standard basis vector is denoted by $e_j$
when the ambient space is clear from context.
For a real-valued random variable $X$ and $r \in (0, \infty)$,
write $\vvvert X \vvvert_r \coloneqq \{\E(|X|^r)\}^{1/r}$
and $\vvvert X \vvvert_\infty \coloneqq \esssup |X|$
for the essential supremum.
For a non-decreasing convex function
$\omega: [0, \infty) \to \R$ with $\omega(0) = 0$
and $\omega \not\equiv 0$,
the $\omega$-Orlicz norm of
$X$ is $\vvvert X \vvvert_\omega \coloneqq
\inf\{t > 0: \E[\omega(|X|/t)] \leq 1\}$.
For $r \in (0, \infty)$, define $\omega_r: [0, \infty) \to \R$
by $\omega_r(x) \coloneqq \exp(x^r) - 1$ and write
$\vvvert X \vvvert_{\omega_\infty} \coloneqq \vvvert X \vvvert_{\infty}$.
If $r \in (0, 1)$, although $\omega_r$ is non-convex,
we may use the same infimum formula
to define the quasi-norm $\vvvert X \vvvert_{\omega_r}$
\citep{kuchibhotla2022moving}.
If $\cH$ is a $\sigma$-algebra, we write
$\vvvert X \vvvert_{\omega \mid \cH} \coloneqq
\essinf\bigl\{T > 0: \E\{\omega(|X|/T) \mid \cH\} \leq 1
\text{ almost surely}\bigr\}$
and $\vvvert X \vvvert_{\omega_\infty \mid \cH}
\coloneqq \essinf\{T > 0: |X| \leq T \text{ almost surely}\}$,
where the essential infima are over $\cH$-measurable
$(0,\infty]$-valued random variables $T$.
For random variables $X$ and $Y$, we write
$X \equallaw Y$ if they have the same distribution.
The positive logarithm function is
$\logplus x \coloneqq 0 \lor \log x$
and we take $0^0 \coloneqq 1$.
For $d \in \N$ and a multi-index $\alpha \in \N_0^d$, define
$|\alpha| \coloneqq \sum_{j=1}^d \alpha_j$ and
$\alpha! \coloneqq \prod_{j=1}^d \alpha_j!$.
For $x \in \R^d$, we write
$x^\alpha \coloneqq \prod_{j=1}^d {x_j}^{\alpha_j}$,
and for $f: \R^d \to \R$, we denote its $\alpha$-partial derivative by
$\partial_\alpha f \coloneqq \partial^{|\alpha|} f
\big/ \bigl(\partial x_1^{\alpha_1} \cdots \partial x_d^{\alpha_d}\bigr)$.
The Borel-measurable subsets of $\R^d$ are denoted by $\cB(\R^d)$.
With $(x_n)$ and $(y_n)$ sequences of non-negative
real numbers, we write
$x_n \lesssim y_n$ if there exists $C > 0$
with $x_n \leq C y_n$ for all $n \in \N$.
For $(X_n)$ and $(Y_n)$ sequences of non-negative
random variables, we write $X_n \lesssim_\P Y_n$ if
for each $\varepsilon > 0$ there exist $N \in \N$
and $M > 0$
such that $\P(X_n > M Y_n) \leq \varepsilon$
for all $n \geq N$.
For $\mu \in \R^d$ and $\Sigma \in \R^{d \times d}$ positive semidefinite,
the multivariate Gaussian distribution with mean $\mu$ and variance $\Sigma$
is $N(\mu, \Sigma)$.

\section{Higher-order Yurinskii coupling}
\label{sec:yurinskii}

We establish a non-asymptotic higher-order
version of Yurinskii's coupling between sums
$S_n\coloneqq\sum_{i=1}^n X_i$ and $T_n\coloneqq\sum_{i=1}^n Y_i$
of independent random vectors in $\R^d$,
with explicit dependence on the dimension $d$.
The resulting bound depends on the distributions of
the summands $X_i$ and $Y_i$ only via
their absolute moments (written as $\beta_{a,b,i}$)
and their moment discrepancies (denoted by $\delta_{\alpha,i}$).
We state our assumptions and present the main result
as Theorem~\ref{thm:yurinskii}, then discuss its
implications and outline the proof.

\begin{assumption}
  \label{ass:data}
  For $n, d \in \N$,
  let $X_1, \ldots, X_n, Y_1, \ldots, Y_n$ be random vectors
  taking values in $\R^d$
  such that $X_1, \ldots, X_n$ are independent and
  $Y_1, \ldots, Y_n$ are independent.
  Fix $k \in \N_0$ and $\gamma \in (0, 1]$,
  and assume that
  $\sum_{i=1}^n \E\bigl(\|X_i\|_2^{k+\gamma}
  + \|Y_i\|_2^{k+\gamma}\bigr) < \infty$.
  Take $m \in \{0, \ldots, k\}$ and suppose that
  for each $\alpha \in \N_0^d$
  with $|\alpha| \in [m]$,
  there exist $\mu_{X,\alpha}, \mu_{Y,\alpha} \in \R$
  such that $\E(X_i^\alpha) = \mu_{X,\alpha}$
  and $\E(Y_i^\alpha) = \mu_{Y,\alpha}$
  for every $i \in [n]$.
  For $\alpha \in \N_0^d$ with $|\alpha| \in [k]$
  and $i \in [n]$, let
  $\delta_{\alpha,i} \coloneqq \E(X_i^\alpha) - \E(Y_i^\alpha)$
  and write $\bar \delta_\alpha \coloneqq
  \sum_{i=1}^n \delta_{\alpha,i} / n$.
  Fix $p \in [1, \infty]$ and,
  for $a \in \N_0$, $b \geq 0$ and $i \in [n]$, define
  \begin{align*}
    \beta_{a,b,i}
    &\coloneqq
    \frac{1}{\sqrt{a!}} \,
    \E \Bigl(
      \|X_i\|_2^a \| X_i\|_p^{b}
      + \|Y_i\|_2^a \| Y_i\|_p^{b}
    \Bigr),
    &\bar \beta_{a,b}
    &\coloneqq
    \frac{1}{n}
    \sum_{i=1}^n
    \beta_{a,b,i}.
  \end{align*}
  For each $\ell \in \{m+1, \ldots, k\}$,
  fix $k_\ell \in \{0, \ldots, m\}$.
  For $\eta, s > 0$, define
  \[
    \begin{aligned}
      \Delta(s)
      &\coloneqq
      \frac{n}{2}
      \Biggl(
        \sum_{\ell = 1}^k
        \sum_{|\alpha| = \ell}
        \frac{\bar \delta_\alpha^2 s^{2 \ell}}{\alpha!}
      \Biggr)^{1/2},
      \iftoggle{journal}{\\ \Gamma_\eta(s) &\coloneqq}{\qquad
      \Gamma_\eta(s) \coloneqq}
      \sum_{i=1}^n
      \sum_{\ell = m+1}^k
      \frac{(2 s)^{\ell + k_\ell}}
      {2 \eta^{\gamma}}
      \bigl(\bar \beta_{k_\ell,\gamma}
      + \beta_{k_\ell,\gamma,i}\bigr)
      \Biggl(
        \sum_{|\alpha| = \ell}
        \frac{\delta_{\alpha,i}^2}{\alpha!}
      \Biggr)^{1/2}.
    \end{aligned}
  \]
\end{assumption}

We also define the following function,
for $d \in \N$, $p \in [1, \infty]$ and $t \geq 0$:
\begin{equation*}
  \phi_{p,d}(t)
  \coloneqq
  \begin{cases}
    \sqrt{p} d^{1/p} + d^{1/p - 1/2} t
    &\quad \text{ if } p \in [1,2),\\
    \sqrt{p} d^{1/p} + t
    &\quad \text{ if } p \in [2,\infty),\\
    \sqrt{2 \log 2d} + t
    &\quad \text{ if } p =\infty.
  \end{cases}
\end{equation*}
Set $\psi_{p,d} \coloneqq \phi_{p,d}(0)$.
When no confusion is possible, write
$\phi(t) \coloneqq \phi_{p,d}(t)$
and $\psi \coloneqq \psi_{p,d}$.
We show in Lemma~\ref{lem:gaussian_norm_tail}%
\iftoggle{journal}{ in the supplementary material
\citep{cattaneo2026highsupplement}}{}
that if $Z \sim N(0_d, I_d)$, then
$\P\bigl(\|Z\|_p > \phi(t)\bigr) \leq e^{-t^2 / 2}$.

\begin{theorem}[Higher-order Yurinskii coupling]
  \label{thm:yurinskii}
  Grant Assumption~\ref{ass:data}.
  For each $t, \eta > 0$,
  there exists a probability space supporting
  random vectors $S_n'$ and $T_n'$,
  where $S_n' \equallaw S_n$ and $T_n' \equallaw T_n$,
  with
  \begin{equation*}
    \P\bigl(\|S_n' - T_n'\|_p > 3\eta\bigr)
    \leq
    2 e^{-t^2 / 2}
    + \frac{n \bar\beta_{k,\gamma} \phi(t)^k}{\eta^{k+\gamma}}
    + \Delta\biggl(\frac{\phi(t)}{\eta}\biggr)
    + \Gamma_\eta\biggl(\frac{\phi(t)}{\eta}\biggr).
  \end{equation*}
\end{theorem}

Theorem~\ref{thm:yurinskii} is proved via Lindeberg's
method and Taylor's theorem, thereby
identifying three sources of coupling error.
The first is the Taylor remainder,
$n \bar\beta_{k,\gamma} \phi(t)^k/\eta^{k+\gamma}$,
which is the higher-order analog of
the term appearing in classical versions of Yurinskii's theorem
\citep[Chapter~10, Section~4]{pollard2002user}.
These presentations use $\gamma = 1$,
but we generalize to $\gamma \in (0, 1]$,
allowing for heavier-tailed summands.
The second source of error is the moment discrepancy term $\Delta$,
which stems from the Taylor polynomial and depends
on the average discrepancies $\bar \delta_\alpha$
between the $\alpha$-moments of $X_i$ and $Y_i$.
Traditionally, the approximation order is fixed at $k = 2$
and the first two moment tensors are assumed to match exactly;
that is, $\E(X_i) = \E(Y_i) = \mu_i \in \R^d$ and
$\Var(X_i) = \Var(Y_i) = \Sigma_i \in \R^{d \times d}$
for all $i \in [n]$,
yielding $\delta_{\alpha,i} = 0$ whenever
$|\alpha| \in \{1, 2\}$.
For our purposes, this classical framework
is restrictive in two ways:
\begin{inlineroman}
  \item it employs moment matching only up to second order,
    neglecting higher-order effects; and
  \item it requires exact moment matching, and hence does not quantify
    errors due to moment discrepancies,
    for example when using bootstrap methods.
\end{inlineroman}
The third source of error is the homogenization term $\Gamma_\eta$,
which arises from replacing the individual moment discrepancies
$\delta_{\alpha,i}$ by their averages $\bar\delta_\alpha$.
This replacement is important when the moment discrepancies are small
only after averaging over $i \in [n]$
(e.g., under the multiplier bootstrap or the
  nonparametric bootstrap with non-identically distributed data;
see Section~\ref{sec:bootstrap}).
We show in Lemma~\ref{lem:bound_Gamma}%
\iftoggle{journal}{ in the supplementary material
\citep{cattaneo2026highsupplement}}{}
that $\sum_{|\alpha| = \ell} \delta_{\alpha,i}^2 / \alpha!
\leq \beta_{\ell,0,i}^2$ for each $\ell \in [k]$ and $i \in [n]$,
allowing $\Gamma_\eta$ to be controlled in terms of absolute moments.
Although apparently a crude estimate,
this bound suffices for all of our results
involving bootstrap schemes.

If all moments up to order $k$ of $X_i$ and $Y_i$ agree
for every $i \in [n]$, that is, if $\E(X_i^\alpha) = \E(Y_i^\alpha)$
for each $\alpha \in \N_0^d$ with $|\alpha| \in [k]$,
then $\Delta(s) = \Gamma_\eta(s) = 0$ for every $s > 0$.
Therefore, provided that $\bar\beta_{k,\gamma} > 0$,
choosing an approximately optimal value of $t$
in Theorem~\ref{thm:yurinskii} yields
\begin{equation}
  \label{eq:yurinskii_matching}
  \P\bigl(\|S_n' - T_n'\|_p > 3\eta\bigr)
  \leq
  \frac{3 n \bar \beta_{k,\gamma}}{\eta^{k + \gamma}}
  \phi
  \Biggl(
    \sqrt{2 \logplus \frac{\eta^{k + \gamma}}
    {n \bar \beta_{k,\gamma} \psi^k}}
  \Biggr)^k.
\end{equation}
With $(p, k, \gamma) = (2, 2, 1)$ and Gaussian summands $Y_i$,
this recovers Yurinskii's classical coupling
\citep[Chapter~10, Theorem~10]{pollard2002user},
although the bound itself does not require Gaussianity.

We now turn to rates, maintaining the assumption of
$\delta_{\alpha,i} = 0$ for all $i \in [n]$, whenever $|\alpha| \in [k]$.
Suppose that $k \in \N$,
$\gamma \in (0, 1]$, $p \in [1, \infty]$,
$r \in [2 \lor (k+\gamma), \infty]$ and $M > 0$ are fixed,
with $\vvvert e_j^\T X_i \vvvert_r
\lor \vvvert e_j^\T Y_i \vvvert_r \leq M$
for each $i \in [n]$ and $j \in [d]$.
For every sequence $R_n \to \infty$,
Theorem~\ref{thm:yurinskii} gives couplings satisfying
\begin{equation}
  \label{eq:yurinskii_rate_general}
  \frac{\|S_n' - T_n'\|_p}{\sqrt n}
  \lesssim_\P
  \frac{R_n d^{k/(2k + 2\gamma)
    + \gamma / \{(p \land r) (k + \gamma)\}}
  \psi^{k/(k + \gamma)}}
  {n^{1/2 - 1/(k + \gamma)}}.
\end{equation}
See \citet[Remark~1]{cattaneo2025yurinskii}
for further discussion of the role of $R_n$.
From \eqref{eq:yurinskii_rate_general},
for $\ell_p$, with $p\in[1,\infty)$,
and for $\ell_\infty$,
we deduce the following simpler
(but slightly looser) bounds:
\begin{align}
  \label{eq:yurinskii_rate_p}
  \frac{\|S_n' - T_n'\|_p}{\sqrt n}
  &\lesssim_\P
  \frac{R_n d^{1/2 + 1 / p}}
  {n^{1/2 - 1/(k + \gamma)}},
  &\frac{\|S_n' - T_n'\|_\infty}{\sqrt n}
  &\lesssim_\P
  \frac{R_n \sqrt{d \log 2d}}
  {n^{1/2 - 1/(k + \gamma)}}.
\end{align}
If we instead assume that
$\bigvvvert \|X_i\|_2 \bigvvvert_r
\lor \bigvvvert \|Y_i\|_2 \bigvvvert_r \leq M$ for each $i \in [n]$,
where $r \in [k + \gamma, \infty]$, the coupling rate is
$R_n n^{-1/2+1/(k+\gamma)}
d^{\{(1/p-1/2)\lor0\}\gamma/(k+\gamma)} \psi^{k/(k+\gamma)}$
for each $p \in [1, \infty]$.

Detailed derivations of \eqref{eq:yurinskii_rate_general},
\eqref{eq:yurinskii_rate_p} and the preceding rate are given in
\iftoggle{journal}{the supplementary material
\citep{cattaneo2026highsupplement}}{%
Appendix~\ref{sec:proofs_yurinskii}}.
These coupling rates
illustrate the improvement obtained from higher-order
moment matching: the sample size exponent $1/2-1/(k+\gamma)$ approaches
$1/2$ as $k$ diverges.
\iftoggle{journal}{The supplementary material
\citep{cattaneo2026highsupplement}}{Appendix~\ref{sec:prokhorov}}
presents corresponding Prokhorov bounds in
Theorem~\ref{thm:prokhorov},
and a lower bound in Lemma~\ref{lem:lower_bound}
demonstrates near-optimality of the Prokhorov rate as $k \to \infty$.
For fixed dimension $d$, matching moments beyond second order improves upon the
rate of $n^{-1/6}$
provided by the classical Yurinskii coupling.
When $d$ grows with $n$, increasing the matching order $k$
improves the sample size exponent but worsens the dimension
dependence in~\eqref{eq:yurinskii_rate_general};
imposing stronger moment assumptions by increasing $r$
reduces the dimension exponent, up to $r = p$.

When $k=2$, one may take $Y_i \sim N(\mu_i, \Sigma_i)$ for
$i \in [n]$, where $\mu_i = \E(X_i)$ and $\Sigma_i = \Var(X_i)$,
yielding a Gaussian coupling variable $T_n$.
If $k \geq 3$, one must generally take $T_n$
non-Gaussian to ensure third-order moment matching.
One important exception is if
$X_i - \mu_i$ and $\mu_i - X_i$ have the same distribution
for each $i \in [n]$;
in this case, all odd centered moments of each
$X_i$ are zero (when they exist),
and a Gaussian coupling suffices for $k = 3$.
See \citet{cattaneo2025yurinskii} for further discussion
of third-order approximations under symmetry.
Section~\ref{sec:literature} presents a detailed
comparison of Theorem~\ref{thm:yurinskii}
with existing results from the literature.

\begin{figure}[t]
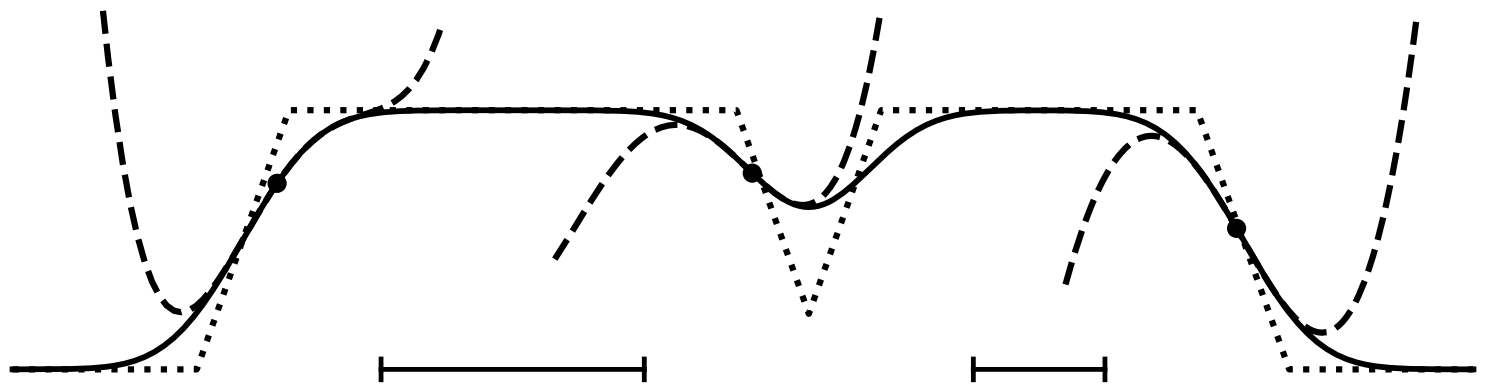

  \centering
  \scalebox{\myfigscale}{%
  \csname input\endcsname{smoothing_exact.pgf}}
  \caption{Smoothing the indicator function
    of a Borel set, here a union of two closed intervals.
    We use linear interpolation
    to construct $g(x)$ (dotted line),
    then apply Gaussian convolution to
    obtain a smooth function $f(x)$ (solid line).
    We depict Taylor polynomials for $f(x)$ of order four
  (dashed lines) at three marked points.}
  \label{fig:smoothing}
\end{figure}

\iftoggle{journal}{%
  The supplementary material \citep{cattaneo2026highsupplement}
  provides a proof sketch and a
full derivation for Theorem~\ref{thm:yurinskii}.}{%
  Appendix~\ref{sec:proof_sketch_yurinskii}
  provides a proof sketch for Theorem~\ref{thm:yurinskii};
a full derivation is given in Appendix~\ref{sec:proofs_yurinskii}.}
By applying Strassen's theorem
\citep[Section~10.3, Theorem~8]{pollard2002user},
it suffices to bound
$\E\bigl(\I\{S_n \in A\} - \I\{T_n \in A_p^{3\eta}\}\bigr)$
uniformly over Borel sets $A \subseteq \R^d$, where
$A_p^{\eta} \coloneqq
\bigl\{ x \in \R^d : \exists y \in A : \|x-y\|_p < \eta \bigr\}$.
To do this, we construct a linear interpolating function
$g(x)$ between the indicators of $A_p^{\eta}$ and $A_p^{2\eta}$
and smooth it by Gaussian convolution to obtain $f(x)$;
see Figure~\ref{fig:smoothing}.
We then use a randomized Lindeberg argument to compare
$\E\{f(S_n)\}$ and $\E\{f(T_n)\}$,
expanding the increments with Taylor's theorem
and controlling the derivatives of $f(x)$ using Hermite polynomials.

\section{Bootstrap-based couplings}
\label{sec:bootstrap}

In statistical applications,
coupling-based inference typically proceeds by computing properties
(such as quantiles of an appropriate transformation)
of the law of the target variable $T_n$,
either analytically or by sampling.
However, the population moments of the data are typically unknown,
so exact moment matching is statistically infeasible.
In this section, we establish coupling results when the target
variable is constructed directly using a bootstrap procedure.
We show how such resampling schemes provide
fully feasible conditional target distributions
with automatic approximate average moment matching, under suitable conditions.
Specifically, we consider both the classical nonparametric (or Efron/empirical)
bootstrap and the multiplier (or wild) bootstrap.
We define each of these methods below;
see \citet{van1996weak} for further details.

For the nonparametric bootstrap, conditional on a realization of
the sample $X_1, \ldots, X_n$,
we draw $X^*_i$ independently for each $i \in [n]$
from the empirical measure associated with the data.
We present here a general version of this definition
in an arbitrary measurable space
for use with the empirical process bootstrap
in the upcoming Section~\ref{sec:empirical_process}.

\begin{assumption}[Nonparametric bootstrap]%
  \label{ass:bootstrap}
  Let $(\cX, \cH)$ be a measurable space and
  $X_1, \ldots, X_n$ be $\cX$-valued random variables;
  write $\bX_n \coloneqq (X_{i}: i \in [n])$.
  Conditional on $\bX_n$,
  let $X^*_1, \ldots, X^*_n$ be i.i.d.\ $\cX$-valued random variables with
  $\P\bigl(X_1^* \in A \mid \bX_n\bigr)
  = \sum_{i=1}^n \I\{X_i \in A\} / n$ for each $A \in \cH$.
  Write $\bX_n^* \coloneqq (X^*_i: i \in [n])$.
\end{assumption}

For the multiplier bootstrap, the auxiliary randomness comes
instead from multiplying the $\R^d$-valued data points
by independent random weights satisfying the following conditions.

\begin{assumption}[Multiplier bootstrap]%
  \label{ass:multiplier}
  Fix $k \in \{2, 3\}$, called the \emph{order}.
  Let $\chi_1, \ldots, \chi_n$ be independent real-valued random variables
  satisfying $\E(\chi_i) = 0$ and $\E\bigl(\chi_i^\ell\bigr) = 1$ for
  each $i \in [n]$
  and $\ell \in \{2, k\}$.
  Suppose that $\vvvert \chi_i \vvvert_4 \leq 2$
  for each $i \in [n]$.
\end{assumption}

We show in Lemma~\ref{lem:multiplier_weights}%
\iftoggle{journal}{ in the supplementary material
\citep{cattaneo2026highsupplement}}{}
that Assumption~\ref{ass:multiplier} holds
with $k=2$ for Rademacher, standard Gaussian
and centered standard exponential weights.
It also holds with $k=3$ for
Mammen and centered Poisson weights.
The value $2$ in the $L_4$-norm bound may be
replaced with any fixed constant.
Throughout the paper, we write
$\P^*$ and $\E^*$ for probability and expectation
respectively under the relevant bootstrap measure,
that is, conditional on the data $X_1, \ldots, X_n$.
In Theorem~\ref{thm:bootstrap}, we derive a coupling result
for the nonparametric and multiplier bootstraps,
assuming independent but not necessarily
identically distributed data.

\begin{theorem}[Bootstrap coupling]
  \label{thm:bootstrap}
  For $n, d \in \N$,
  let $X_1, \ldots, X_n$ be independent $\R^d$-valued random variables
  and write $\bX_n \coloneqq (X_1, \ldots, X_n)$.
  For each $i \in [n]$, assume that $\E(X_i) = \mu \in \R^d$
  and $\E\bigl(\|X_i - \mu\|_2^4\bigr) < \infty$.
  Fix $k \in \{2, 3\}$, $\gamma \in (0, 1]$ and $p \in [1, \infty]$.
  For $q \in [1, \infty]$, let $M_q \geq 0$ satisfy
  $M_q^4 = \sum_{i=1}^n \E\bigl(\|X_i - \mu\|_q^4\bigr) / n$.
  Write $\bar X_n \coloneqq \sum_{i=1}^n X_i/n$ and
  $S_n \coloneqq \sum_{i=1}^n (X_i - \mu)$,
  and grant one of the following.
  \begin{enumerate}
    \item
      \label{it:thm_bootstrap_nonparametric}
      Let $X_1^*, \ldots, X_n^*$ satisfy Assumption~\ref{ass:bootstrap}
      and set $Y_i^* \coloneqq X^*_i - \bar X_n$ for each $i \in [n]$.

    \item
      \label{it:thm_bootstrap_multiplier}
      Let $\chi_1, \ldots, \chi_n$ be random variables
      satisfying Assumption~\ref{ass:multiplier}
      with order $k$, independent of $\bX_n$.
      Set $Y_i^* \coloneqq \chi_i(X_i-\bar X_n)$
      for each $i \in [n]$.
  \end{enumerate}
  Let $S_n^* \coloneqq \sum_{i=1}^n Y_i^*$.
  For each $t, \eta > 0$, there is a probability space
  supporting $\bX_n'$, $S_n'$ and $S_n^{*\prime}$, where
  $S_n' \equallaw S_n$, $S_n'$ is independent of $\bX_n'$,
  and $(\bX_n', S_n^{*\prime}) \equallaw (\bX_n, S_n^*)$, with
  \begin{equation}
    \label{eq:bootstrap}
    \P\bigl(
      \|S_n' - S_n^{* \prime}\|_p
      > \eta
    \bigr)
    \leq
    2 e^{-t^2 / 2}
    + \frac{C \sqrt n \phi(t)^2 M_2^2}{\eta^2}
    + \frac{C n \phi(t)^k M_2^k M_p^\gamma}
    {\eta^{k+\gamma}},
  \end{equation}
  where $C > 0$ is a universal constant.
\end{theorem}

In Theorem~\ref{thm:bootstrap},
the summands in $S_n$ satisfy
$\E(X_i - \mu) = 0_d$, where $\mu$ is the common population mean.
Likewise, under either bootstrap scheme,
the summands in $S_n^*$ satisfy $\E^*(Y_i^*) = 0_d$.
Preserving the law of
$(\bX_n, S_n^*)$ ensures that,
conditional on the copied data $\bX_n'$,
the coupled bootstrap statistic
$S_n^{*\prime}$ follows the correct bootstrap resampling distribution.
Moreover, independence of $S_n'$ and $\bX_n'$ guarantees that
we compare the bootstrap law of $S_n^*$ to the marginal law of $S_n$.
The errors appearing in \eqref{eq:bootstrap}
are attributable to two sources.
Firstly, the term involving $\eta^2$ comes from the $\Delta$ contribution
in Theorem~\ref{thm:yurinskii}; this is due to errors in
the bootstrap replications of the population average moments.
Secondly, the term including $\eta^{k+\gamma}$ bounds those involving
both $\bar\beta_{k,\gamma}$ and $\Gamma_\eta$.

We now consider bootstrap approximation rates
in a setting similar to that of \eqref{eq:yurinskii_rate_general}.
Take $k=3$ and $\gamma=1$, fix $p \in [1, \infty]$,
and suppose that $\max_{i \in [n]} \max_{j \in [d]}
\vvvert e_j^\T (X_i - \mu) \vvvert_r \leq M$,
where $r \in [4, \infty]$ and $M > 0$ are fixed.
For every sequence $R_n\to\infty$,
we obtain couplings for $\ell_p$, with $p\in[1,\infty)$,
and for $\ell_\infty$ satisfying, respectively,
\begin{align}
  \label{eq:bootstrap_rate}
  \frac{\|S_n' - S_n^{*\prime}\|_p}{\sqrt n}
  &\lesssim_\P
  \frac{R_n d^{1/2 + 1/p}}{n^{1/4}},
  &\frac{\|S_n' - S_n^{*\prime}\|_\infty}{\sqrt n}
  &\lesssim_\P
  \frac{R_n \sqrt{d \log 2 d}}{n^{1/4}}.
\end{align}
If instead we assume that $M_2 \leq M$
in Theorem~\ref{thm:bootstrap}, these bounds become
$R_n n^{-1/4} d^{1/p}$ and $R_n n^{-1/4} \sqrt{\log 2d}$,
respectively.
For fixed $d$, the bootstrap achieves a coupling
rate of $R_n n^{-1/4}$
by approximating the average second- and third-order population moments.
For multiplier weights satisfying Assumption~\ref{ass:multiplier}
for $k=2$ but not $k=3$
(e.g., Rademacher, Gaussian or centered exponential weights),
taking $k=2$ and $\gamma=1$ in
Theorem~\ref{thm:bootstrap}
yields the slower rate of $R_n n^{-1/6}$,
the same as that given by the classical Gaussian Yurinskii coupling
with exact second-order moment matching.

Drawing samples from conditional bootstrap distributions
is straightforward, and our coupling bound automatically
accounts for any errors resulting from implicit estimation
of population moments.
As such, Theorem~\ref{thm:bootstrap} is well suited for
direct application to statistical inference procedures;
see, for example, the local polynomial methods in
the upcoming Section~\ref{sec:locpol}
and the distributional approximation for maxima
presented as Theorem~\ref{thm:bootstrap_ks}%
\iftoggle{journal}{ in the supplementary material
\citep{cattaneo2026highsupplement}}{}.
When $\mu=0_d$, the proof of Theorem~\ref{thm:bootstrap}
may easily be modified to allow also for the non-centered
multiplier bootstrap that uses $Y_i^* \coloneqq \chi_i X_i$%
\iftoggle{journal}{}{; see Appendix~\ref{sec:proofs_bootstrap}}.
We compare our bootstrap coupling theory
with other existing bootstrap approximations
in Section~\ref{sec:literature_bootstrap}.

\iftoggle{journal}{%
  The supplementary material \citep{cattaneo2026highsupplement}
  provides a proof sketch and a
full derivation for Theorem~\ref{thm:bootstrap}.}{%
  Appendix~\ref{sec:proof_sketch_bootstrap}
  provides a proof sketch for Theorem~\ref{thm:bootstrap},
  and a full derivation is presented
in Appendix~\ref{sec:proofs_bootstrap}.}
We rely on a Borel set approximation bound
identical to that underlying the proof of Theorem~\ref{thm:yurinskii},
now comparing the marginal law of $S_n$
with the bootstrap law of $S_n^*$ conditional on $\bX_n$.
For $\alpha \in \N_0^d$ with $|\alpha| \in [k]$ and $i \in [n]$, define
$\hat \delta_{\alpha,i}
\coloneqq \E\bigl\{(X_i - \mu)^\alpha\bigr\}
- \E^*\bigl\{(Y_i^*)^\alpha\bigr\}$;
under either bootstrap scheme,
the corresponding average moment discrepancies
are given by \eqref{eq:delta_hat}.
We use symmetrization to bound $\E(\hat\Delta)$,
where $\hat\Delta$ is obtained by replacing
$\bar\delta_\alpha$ with $\hat\delta_\alpha$
in the definition of $\Delta$,
and a conditional version of Strassen's
theorem yields the coupling.

We now present a refined coupling rate
for the nonparametric bootstrap.

\begin{theorem}[Refined nonparametric bootstrap coupling]%
  \label{thm:refined_bootstrap}
  For $n, d \in \N$,
  let $X_1,\ldots,X_n$ be independent $\R^d$-valued
  random variables and
  write $\bX_n \coloneqq (X_1,\ldots,X_n)$.
  For each $i \in [n]$, assume that
  $\E(X_i) = \mu$, $\Var(X_i) = \Sigma \neq 0_{d \times d}$,
  and $\E\bigl(\|X_i-\mu\|_2^6\bigr)<\infty$.
  Fix $k\in\{2,\ldots,5\}$, $\gamma\in(0,1]$ and $p \in [1, \infty]$.
  Let $M_2 > 0$ satisfy
  $M_2^6 = \sum_{i=1}^n \E\bigl(\|X_i-\mu\|_2^6\bigr)/n$.
  Let $X_1^*,\ldots,X_n^*$ satisfy
  Assumption~\ref{ass:bootstrap}, and let
  $\bar X_n\coloneqq\sum_{i=1}^n X_i/n$,
  $S_n\coloneqq\sum_{i=1}^n(X_i-\mu)$ and
  $S_n^*\coloneqq\sum_{i=1}^n(X_i^*-\bar X_n)$.
  Let $\lambda>0$ be the smallest non-zero eigenvalue
  of $\Sigma$ and define $a_p \coloneqq (1/p - 1/2) \lor 0$
  and $B_p \coloneqq \|\Sigma^{1/2}\|_{p \to p}$.
  For each $t, \eta > 0$, there is a probability space
  supporting $\bX_n'$, $S_n'$ and $S_n^{*\prime}$, where
  $S_n' \equallaw S_n$, $S_n'$ is independent of $\bX_n'$,
  and $(\bX_n', S_n^{*\prime}) \equallaw (\bX_n, S_n^*)$, with
  \begin{align*}
    \P\bigl(
      \|S_n' - S_n^{*\prime}\|_p > \eta
    \bigr)
    &\leq
    2 e^{-t^2 / 2}
    + \frac{
      C \sqrt n M_2^5 B_p^3 \phi(t)^3
    }{\lambda^{5/2} \eta^3}
    + \frac{
      C n d^{\gamma a_p}
      (M_2 B_p)^{k+\gamma}
      \phi(t)^k
    }{\lambda^{(k+\gamma)/2} \eta^{k+\gamma}} \\
    &\qquad+
    \frac{C M_2^6}{n^{3/2} \lambda^3}
    + \frac{C d^{2a_p} M_2^6}{\lambda^2 \eta^2},
  \end{align*}
  where $C > 0$ is a universal constant.
\end{theorem}

To obtain coupling rates from Theorem~\ref{thm:refined_bootstrap},
take $k=5$ and $\gamma=1$ and fix $p \in [1, \infty]$.
Suppose that $\max_{i \in [n]} \max_{j \in [d]}
\vvvert e_j^\T (X_i - \mu) \vvvert_r \leq M$,
where $r \in [6, \infty]$ and $M > 0$ are fixed.
For every sequence $R_n\to\infty$,
we obtain couplings for $\ell_p$, with $p\in[1,\infty)$,
and for $\ell_\infty$ satisfying, respectively,
\begin{align}
  \label{eq:refined_rate}
  \frac{\|S_n' - S_n^{*\prime}\|_p}{\sqrt n}
  &\lesssim_\P
  \frac{R_n d^{4/3 + 1/p}}{n^{1/3} \lambda^{5/6}},
  &\frac{\|S_n' - S_n^{*\prime}\|_\infty}{\sqrt n}
  &\lesssim_\P
  \frac{R_n d^{4/3} \sqrt{\log 2d}}{n^{1/3} \lambda^{5/6}}.
\end{align}
If instead we assume that $M_2 \leq M$ in Theorem~\ref{thm:refined_bootstrap},
these bounds become
$R_n n^{-1/3} d^{a_p} \lambda^{-4/3}$
and $R_n n^{-1/3} \lambda^{-4/3}$, respectively.
For fixed $d$ and $\lambda$, this result attains a sample size rate
of $R_n n^{-1/3}$; compare this to the slower
$R_n n^{-1/4}$ rate given in Theorem~\ref{thm:bootstrap}.
This improvement requires bounded sixth moments
(rather than fourth moments)
and assumes a common covariance matrix $\Sigma$;
moreover, it applies only to the nonparametric bootstrap.
The dependence on the dimension is rather worse than
in Theorem~\ref{thm:bootstrap}: comparing the rates
in \eqref{eq:bootstrap_rate} and \eqref{eq:refined_rate}
shows that Theorem~\ref{thm:refined_bootstrap} offers a superior bound
whenever $d / (\lambda n^{1/10}) \to 0$, assuming all other parameters
are fixed.

To prove Theorem~\ref{thm:refined_bootstrap},
we replace the sum $S_n$ from Theorem~\ref{thm:bootstrap}
with its standardized version $\hat S_n \coloneqq \Sigma^{+/2} S_n$,
where $A^{+/2}$ denotes the positive semidefinite square root
of the Moore--Penrose pseudoinverse of $A$.
Likewise, we replace the bootstrapped sum $S_n^*$ with
$\hat S_n^* \coloneqq \hat \Sigma_n^{+/2} S_n^*$, where
$\hat \Sigma_n \coloneqq \sum_{i=1}^n X_i X_i^\T / n - \bar X_n \bar X_n^\T$.
We then compare the population law of $\hat S_n$ with the conditional
bootstrap law of $\hat S_n^*$ given $\bX_n$.
On the event where the matrices $\Sigma$ and
$\hat\Sigma_n$ have the same range,
for each $i \in [n]$, we have
\begin{align*}
  \E\bigl\{\Sigma^{+/2}(X_i-\mu)\bigr\}
  &= \E^*\bigl\{\hat\Sigma_n^{+/2}(X_i^*-\bar X_n)\bigr\}
  = 0_d, \\
  \Var \bigl\{\Sigma^{+/2}(X_i-\mu)\bigr\}
  &= \Var^* \bigl\{\hat\Sigma_n^{+/2}(X_i^*-\bar X_n)\bigr\},
\end{align*}
by the homoskedasticity assumption,
where $\Var^*$ denotes covariance conditional on $\bX_n$.
We may therefore take $m=2$ in Theorem~\ref{thm:yurinskii}.
Moreover, the term $\hat\Delta(s)$
contains only the average moment discrepancy terms
$\hat\delta_{\alpha}$ where
$|\alpha| = \ell \in \{3, \ldots, k\}$.
The absence of second-order terms in $\hat\Delta(s)$ and, as $m = 2$,
in $\hat\Gamma_\eta(s)$, drives the rate refinement.
Finally, multiplying the coupled standardized sums
$\hat S_n'$ and $\hat S_n^{*\prime}$
by $\Sigma^{1/2}$ and $\hat\Sigma_n^{\prime 1/2}$, respectively,
where $\hat\Sigma_n'$ is the sample covariance matrix of $\bX_n'$,
yields a coupling for $S_n$ and $S_n^*$.
This last step introduces a term depending
on the covariance estimation error $\hat\Sigma_n - \Sigma$.
This construction does not extend in general
to the multiplier bootstrap with $k \in \{4, 5\}$
as it is impossible to find
weights $\chi_i$ with $\E(\chi_i) = 0$ and $\E(\chi_i^\ell) = 1$
for $\ell \in \{2,3,4\}$.

\section{Comparison with the literature}
\label{sec:literature}

We compare our results with existing coupling and bootstrap approximations,
focusing on dimension dependence, moment conditions,
and the role of moment matching.
For clarity, we divide the literature into three broad strands,
containing couplings based on Yurinskii's approach,
the Koml\'os--Major--Tusn\'ady (KMT) method,
and Wasserstein bounds, respectively.
We then discuss distributional
approximation methods and the role of moment matching in the bootstrap.

We assume throughout this section
that the summands $X_i$ are independent with $\E(X_i) = 0_d$;
further, we suppose that they satisfy one of three bounds,
with fixed parameters $M > 0$ and $r \in [1, \infty)$:
\begin{inlineroman}
  \item
    $\|X_i\|_\infty \leq M$ almost surely;
  \item
    $\|X_i\|_2 \leq M$ almost surely; or
  \item
    $\max_{j \in [d]} \vvvert e_j^\T X_i \vvvert_r \leq M$.
\end{inlineroman}
Condition~(i) implies that $\|X_i\|_2 \leq \sqrt d M$
almost surely,
while condition~(iii) gives
$\bigl\{\E(\|X_i\|_2^r)\bigr\}{}^{1/r} \leq \sqrt d M$
for $r \geq 2$.
Many coupling results in the literature take
$Y_i \sim N(0_d, \Sigma_i)$ with $\Sigma_i = \Var(X_i)$,
thereby matching
the first $k=2$ moments of each summand ($k=3$ if third moments vanish);
our Theorem~\ref{thm:yurinskii} matches the first $k \in \N$
moments and permits non-Gaussian targets.
Unless stated otherwise, all presented rates are
bounds in probability for the normalized coupling error
$\|S_n' - T_n'\|_p / \sqrt n$, where $p\in[1,\infty]$ is fixed.
Yurinskii-type constructions, including Theorem~\ref{thm:yurinskii},
depend on the required coupling error threshold $\eta$,
and their rates hold up to any deterministic sequence
$R_n \to \infty$; see \citet[Remark~1]{cattaneo2025yurinskii};
KMT and Wasserstein bounds provide a
single coupling valid for every threshold so
their rates need no such factor.

Table~\ref{tab:literature} collects coupling rates from the literature,
with Panels~A--C covering the approaches discussed in
Sections~\ref{sec:literature_yurinskii}--\ref{sec:literature_wasserstein}
respectively.
Our Theorems~\ref{thm:yurinskii}, \ref{thm:bootstrap}
and~\ref{thm:refined_bootstrap}
are presented in Panel~D.
Results restricted to univariate data are omitted.
In the sampling column, i.n.i.d.\ denotes independent
but not necessarily identically distributed summands.
Each rate is formulated as the product
of a \emph{sample size factor} and a \emph{dimension factor}.
These factors are given in separate columns;
the dimension factor depends on the assumed moment bound.
For example, the first row states that
if $\|X_i\|_\infty \leq M$ almost surely,
then Yurinskii's coupling gives
(for every sequence $R_n \to \infty$)
a series of couplings satisfying
$\|S_n' - T_n'\|_2 / \sqrt n \lesssim_\P R_n n^{-1/6} d^{5/6}$.
This bound tends to zero (for a suitable sequence~$R_n$)
whenever $d / n^{1/5} \to 0$.
If instead $\|X_i\|_2 \leq M$ almost surely,
the corresponding dimension factor is
$d^{1/3}$ and $d / n^{1/2} \to 0$ suffices.
Table~\ref{tab:literature} omits $p \in [1, 2)$;
in these cases, the corresponding rate may be calculated by
taking the dimension factor associated with
$p = 2$ and multiplying it by $d^{1/p - 1/2}$.
The column labeled $r$ gives the moment order under which each result holds,
with $\infty$ indicating almost surely bounded summands.
Finite exponential moment assumptions are indicated by exp.
In the last column, N.A.\
indicates that the result requires
bounded summands or exponential moments.
Throughout, $\vartheta > 0$ is arbitrary.
We define $\varrho_r(d) \coloneqq
\max\bigl\{d^{21/2+\vartheta},\, d^{(r+2)/4} (\log 2d)^{(r+1)/2}\bigr\}$
and, in Panel~C, we take $\lambda \coloneqq \lambda_{\min}
\bigl(\sum_{i=1}^n \Var(X_i)/n\bigr)$.
Under $\|X_i\|_2 \leq M$, we have $\lambda \leq M^2/d$ in Panel~C
and $\lambda \leq M^2/ (\rank \Sigma)$
for Theorem~\ref{thm:refined_bootstrap},
so the corresponding dimension factors still grow with $d$ in Panel~C
and with $\rank \Sigma$ for Theorem~\ref{thm:refined_bootstrap}.
For our Theorem~\ref{thm:yurinskii},
we assume exact moment matching for each summand up to order $k$
and take $\gamma \in (0, 1]$ fixed;
the summands $Y_i$ satisfy the same assumptions as $X_i$.
For Theorems~\ref{thm:bootstrap} and~\ref{thm:refined_bootstrap},
the coupling target is the conditional law (given the data)
of the nonparametric or third-order multiplier bootstrap sum $S_n^*$.
Theorem~\ref{thm:refined_bootstrap} applies only to the nonparametric
bootstrap and assumes a common covariance matrix with smallest
non-zero eigenvalue~$\lambda$.

\begin{table}[p]
  \shortcites{fang2026sharp}
  \centering
  \rotatebox{90}{%
    \begin{minipage}{\textheight}
      \centering
      \caption{Coupling rates for random vectors ($d>1$)
      under various moment bounds.}
      \label{tab:literature}
      \iftoggle{journal}{\scriptsize}{\footnotesize}
      \setlength{\tabcolsep}{2.9pt}%
      \renewcommand{\arraystretch}{1.09}%
      \vspace*{-1mm}
      \begin{tabular*}{\linewidth}{@{\extracolsep{\fill}}lcccccccc@{}}
        \hline
        & & & & & \raisebox{-0.6mm}{Sample}
        & \multicolumn{3}{c}{Dimension factor} \\
        \cline{7-9}
        & \raisebox{0.4mm}{Sampling}
        & \raisebox{0.4mm}{$k$}
        & \raisebox{0.4mm}{$p$}
        & \raisebox{0.4mm}{$r$}
        & \raisebox{0.4mm}{size factor}
        & $\|X_i\|_\infty \leq M$ & $\|X_i\|_2 \leq M$
        & \rule{0pt}{2.7ex}$\vvvert e_j^\T X_i \vvvert_r \leq M$ \\
        \hline
        \noalign{\vskip 0.2mm}
        \multicolumn{9}{@{}l}{\textit{Panel A:
        Yurinskii-type couplings}} \\[-0.3mm]
        \quad\citet{yurinskii1978error}
        & i.n.i.d. & $2$ & $2$ & $[3,\infty)$ & $n^{-1/6}$
        & $d^{5/6}$ & $d^{1/3}$ & $d^{5/6}$ \\
        \quad\citet{zaitsev1987estimates}
        & i.n.i.d. & $2$ & $[2,\infty]$ & exp.
        & $n^{-1/2}$
        & $d^{5/2} \log 2d$ & $d^{2} \log 2d$ & N.A. \\
        \quad\citet{le1988prokhorov}
        & i.n.i.d. & $2$ & $[2,\infty]$ & $[3,\infty)$ & $n^{-1/6}$
        & $d^{5/6}$ & $d^{1/3}$ & $d^{5/6}$ \\
        \rule{0pt}{10.5pt}\quad\citet{belloni2019conditional}
        & i.n.i.d. & $2$ & $[2,\infty)$ & $[3,\infty)$ & $n^{-1/6}$
        & $d^{1/3+1/p}$ & $d^{2/(3p)}$
        & $d^{1/3+\{2/p+1/(p \land r)\}/3}$ \\
        & & & $\infty$ & & $n^{-1/6}$
        & $d^{1/3} \sqrt{\log 2d}$ & $(\log 2d)^{1/3}$
        & $d^{1/(3r)} (d \log 2d)^{1/3}$ \\
        \rule{0pt}{10.5pt}\quad\iftoggle{journal}{%
        Cattaneo et al.\ \citep{cattaneo2025yurinskii}}{%
        \citet{cattaneo2025yurinskii}}
        & i.n.i.d. & $3$ & $[2,\infty)$ & $[4,\infty)$ & $n^{-1/4}$
        & $d^{3/8+1/p}$ & $d^{3/(4p)}$
        & $d^{3/8+\{3/p+1/(p \land r)\}/4}$ \\
        & & & $\infty$ & & $n^{-1/4}$
        & $d^{3/8} \sqrt{\log 2d}$ & $(\log 2d)^{3/8}$
        & $d^{1/(4r)} (d \log 2d)^{3/8}$ \\
        \hline
        \noalign{\vskip 0.2mm}
        \multicolumn{9}{@{}l}{\textit{Panel B:
        KMT-type couplings}} \\[-0.3mm]
        \quad\citet{zaitsev1998multidimensional}
        & i.n.i.d. & $2$ & $[2,\infty]$ & exp. & $n^{-1/2} \log 2n$
        & $d^{25/4+\vartheta} \log 2d$
        & $d^{23/4+\vartheta} \log 2d$ & N.A. \\
        \quad\citet{zaitsev2001multidimensional,zaitsev2002multidimensionalII,%
        zaitsev2002multidimensionalIII}
        & i.n.i.d. & $2$ & $[2,\infty]$ & exp. & $n^{-1/2} \log 2n$
        & $d^{8+\vartheta} \log 2d$
        & $d^{15/2+\vartheta} \log 2d$ & N.A. \\
        \quad\citet{zaitsev2007estimates}
        & i.n.i.d. & $2$ & $[2,\infty]$ & $(2,\infty)$ & $n^{-1/2+1/r} \log 2n$
        & $d^{8+\vartheta}$ & $d^{15/2+\vartheta}$
        & $d^{8+\vartheta}$ \\
        \quad\citet{gotze2009bounds}
        & i.n.i.d. & $2$ & $[2,\infty]$ & $(2,\infty)$ & $n^{-1/2+1/r}$
        & $d^{11+\vartheta}$
        & $d^{21/2+\vartheta}$
        & $d^{11+\vartheta}$ \\
        \quad\citet{gotze2011estimates}
        & i.i.d. & $2$ & $[2,\infty]$ & $(2,\infty)$ & $n^{-1/2+1/r}$
        & $d^{1/2} \varrho_r(d)$ & $\varrho_r(d)$
        & $d^{1/2} \varrho_r(d)$ \\
        \hline
        \noalign{\vskip 0.2mm}
        \multicolumn{9}{@{}l}{\textit{Panel C:
        Wasserstein-based couplings}} \\[-0.3mm]
        \quad\citet{zhai2018high}
        & i.n.i.d. & $2$ & $[2,\infty]$ & $\infty$ & $n^{-1/2} \log 2n$
        & $d$ & $d^{1/2}$ & N.A. \\
        \quad\citet{bonis2020stein}
        & i.i.d. & $2$ & $[2,\infty]$ & $[4,\infty)$ & $n^{-1/2}$
        & $d$ & $d^{1/2}$ & $d^{5/4}$ \\
        \rule{0pt}{10.5pt}\quad\iftoggle{journal}{%
        Eldan et al.\ \citep{eldan2020clt}}{%
        \citet{eldan2020clt}}
        & i.i.d. & $2$ & $[2,\infty]$ & $\infty$ & $n^{-1/2} \sqrt{\log 2n}$
        & $d$ & $d^{1/2}$ & N.A. \\
        \quad\citet{mies2023sequential}
        & i.n.i.d. & $2$ & $[2,\infty]$ & $(2,\infty]$
        & $n^{-1/2+1/r} \sqrt{\log 2n}$
        & $d^{1-1/r}$ & $d^{1/2-1/r}$ & $d^{1-1/r}$ \\
        \quad\iftoggle{journal}{%
        Fang et al.\ \citep[Corollary~2.1]{fang2026martingale}}{%
        \citet[Corollary~2.1]{fang2026martingale}}
        & i.n.i.d. & $2$ & $[2,\infty]$ & $[3,\infty)$ & $n^{-1/2} \log 2n$
        & $\lambda^{-1} d^{1+1/p}$ & $\lambda^{-1}$
        & $\lambda^{-1} d^{1+1/(p \land r)}$ \\
        \rule{0pt}{10.5pt}\quad\iftoggle{journal}{%
        Fang et al.\ \citep[Corollary~2.2]{fang2026martingale}}{%
        \citet[Corollary~2.2]{fang2026martingale}}
        & i.n.i.d. & $2$ & $[2,\infty)$ & $[3,\infty)$ & $n^{-1/6}$
        & $d^{1/3+1/p}$ & $d^{2/(3p)}$
        & $d^{1/3+\{2/p+1/(p \land r)\}/3}$ \\
        & & & $\infty$ & & $n^{-1/6}$
        & $(d \log 2d)^{1/3}$ & $(\log 2d)^{1/3}$
        & $d^{1/(3r)} (d \log 2d)^{1/3}$ \\
        \quad\iftoggle{journal}{%
        Fang et al.\ \citep[Theorem~2.1]{fang2026sharp}}{%
        \citet[Theorem~2.1]{fang2026sharp}}
        & i.n.i.d. & $2$ & $[2,\infty]$
        & $[4,\infty)$ & $n^{-1/2}$
        & $d$ & $d^{1/2}$ & $\lambda^{-1/2}d$ \\
        \quad\iftoggle{journal}{%
        Fang et al.\ \citep[Theorem~2.3]{fang2026sharp}}{%
        \citet[Theorem~2.3]{fang2026sharp}}
        & i.i.d. & $2$ & $[2,\infty]$
        & $[4,\infty)$ & $n^{-1/4} \log 2n$
        & $d^{3/4}$ & $d^{1/4}$ & $d^{3/4}$ \\
        \hline
        \noalign{\vskip 0.2mm}
        \multicolumn{9}{@{}l}{\textit{Panel D:
        Results of this paper}} \\[-0.3mm]
        \quad Theorem~\ref{thm:yurinskii}
        & i.n.i.d. & $\N$ & $[2,\infty)$ & $\bigl[2\lor(k+\gamma),\infty\bigr)$
        & $n^{-1/2+1/(k+\gamma)}$
        & $d^{k/(2k+2\gamma)+1/p}$
        & $d^{k/\{p(k+\gamma)\}}$
        & $d^{\{k(1/2 + 1/p)+\gamma/(p \land r)\}/(k+\gamma)}$ \\
        & & & $\infty$ & & $n^{-1/2+1/(k+\gamma)}$
        & $(d \log 2d)^{k/(2k+2\gamma)}$
        & $(\log 2d)^{k/(2k+2\gamma)}$
        & $d^{\gamma/\{r(k+\gamma)\}}
        (d \log 2d)^{k/(2k+2\gamma)}$ \\
        \rule{0pt}{10.5pt}\quad Theorem~\ref{thm:bootstrap}
        & i.n.i.d. & $3$ & $[2,\infty)$ & $[4,\infty)$ & $n^{-1/4}$
        & $d^{1/2+1/p}$ & $d^{1/p}$ & $d^{1/2+1/p}$ \\
        & & & $\infty$ & & $n^{-1/4}$
        & $\sqrt{d \log 2d}$ & $\sqrt{\log 2d}$
        & $\sqrt{d \log 2d}$ \\
        \rule{0pt}{10.5pt}\quad Theorem~\ref{thm:refined_bootstrap}
        & i.n.i.d. & $5$ & $[2,\infty)$ & $[6,\infty)$ & $n^{-1/3}$
        & $\lambda^{-5/6} d^{4/3+1/p}$ & $\lambda^{-4/3}$
        & $\lambda^{-5/6} d^{4/3+1/p}$ \\
        & & & $\infty$ & & $n^{-1/3}$
        & $\lambda^{-5/6} d^{4/3} \sqrt{\log 2d}$ & $\lambda^{-4/3}$
        & $\lambda^{-5/6} d^{4/3} \sqrt{\log 2d}$ \\
        \hline
      \end{tabular*}

      \vspace{2pt}
      \begin{minipage}{\linewidth}
        \scriptsize\raggedright
        \emph{Note:} Rates in Panels~A and~D omit arbitrarily slowly
        diverging factors. Additional assumptions and precise
        statements are given in
        Sections~\ref{sec:literature_yurinskii}%
        --\ref{sec:literature_wasserstein}
        and Appendix~\ref{sec:literature_details}%
        \iftoggle{journal}{ in the supplementary material
        \citep{cattaneo2026highsupplement}}.
      \end{minipage}
  \end{minipage}}
\end{table}

\subsection{Yurinskii-type couplings}
\label{sec:literature_yurinskii}

The classical Yurinskii coupling combines Strassen's theorem
with the Borel set approximation and
Lindeberg method of \citet{yurinskii1978error},
who used these tools to bound the Prokhorov distance
between $S_n$ and $T_n$.
This coupling corresponds to our Theorem~\ref{thm:yurinskii}
with $(p, k, \gamma) = (2, 2, 1)$ and
Gaussian target summands $Y_i$, assuming exactly matching
mean vectors and covariance matrices.
See \citet[Chapter~10]{pollard2002user} for a textbook
treatment of Yurinskii's coupling, and
\citet{chatterjee2006generalization} and \citet{ravckauskas2023asymptotic}
for a broader discussion of Lindeberg's method and its generalizations.
\citet[Theorem~1, part~2]{le1988prokhorov} and
\citet[Supplemental Appendix~I.1, Lemma~38]{belloni2019conditional}
extended Yurinskii's approach to arbitrary norms, including the
$\ell_p$-norms with $p \in [1, \infty]$; the maximum norm
with $p = \infty$ is especially useful in empirical process
or high-dimensional settings \citep{kim2024absolute}.
Under finite exponential moments,
\citet[Theorem~1.1]{zaitsev1987estimates} bounded the
Prokhorov distance at a faster rate
in the sample size;
combined with Strassen's theorem, his bound gives a Gaussian coupling
with sample size factor $n^{-1/2}$ and no conditions on the covariance
matrix, but with larger dimension factors
(Panel~A of Table~\ref{tab:literature}).
\citet{li2020uniform} extended Yurinskii's Euclidean
($p = 2$) coupling
to martingales and mixingales, assuming covariance eigenvalues
bounded away from zero.
\citet{cattaneo2025yurinskii}
further generalized Yurinskii's coupling
to approximate martingale difference sequences,
allowing for Gaussian mixture coupling distributions
and imposing no restrictions on covariance degeneracy.
Moreover, they showed that the coupling bound improves
when the conditional third moments of every summand vanish.
Subsequent work includes strong approximations
for dependent random fields \citep{bonnerjee2026strong}.

Relative to this literature, our Theorem~\ref{thm:yurinskii}
improves Yurinskii's coupling for independent
random vectors in three ways:
(i) the moments of $X_i$ and $Y_i$ may be matched
up to arbitrary order $k \in \N$,
and the approximating sum $T_n$ need not be Gaussian;
(ii) the moment matching need only hold approximately
and on average
over $i \in [n]$, with the resulting errors made explicit through
$\Delta$ and $\Gamma_\eta$; and
(iii) both families of summands need finite moments only of order
$k + \gamma$, where $\gamma \in (0, 1]$.
When the moment discrepancies are small, higher-order matching yields
faster approximation rates
and allows the dimension to grow more quickly with the sample size;
see Section~\ref{sec:yurinskii}.
Panel~A of Table~\ref{tab:literature} compares
Yurinskii coupling rates under bounds
on $\|X_i\|_\infty$, $\|X_i\|_2$ or $\vvvert e_j^\T X_i \vvvert_r$.
Appendix~\ref{sec:literature_details_yurinskii}
presents these rates more precisely in terms
of the quantities $\bar\beta_{a,b}$ defined in
Assumption~\ref{ass:data};
such bounds can be considerably sharper than the
worst-case rates in Table~\ref{tab:literature}.

Yurinskii's coupling also provides a basis for bootstrap approximation.
For example,
\citet[Theorem~7]{belloni2019conditional}
used Gaussian strong approximations to justify the weighted bootstrap
for quantile regression coefficient processes, and
\citet[Proposition~A.2]{cattaneo2025yurinskii}
obtained a Gaussian multiplier bootstrap approximation
by replacing the population covariance matrix with an estimator.
Our Theorems~\ref{thm:bootstrap} and~\ref{thm:refined_bootstrap}
provide bootstrap couplings based on approximate average moment matching,
with matching orders up to $k = 3$ and $k = 5$, respectively.
Table~\ref{tab:literature} presents
Theorem~\ref{thm:bootstrap} with $(k, \gamma) = (3, 1)$
for the nonparametric and third-order multiplier bootstraps.
Under bounded fourth moments, this result allows heterogeneous
covariances and attains the sample size factor $n^{-1/4}$ of
\citet{cattaneo2025yurinskii} without requiring knowledge of the
population moments or assuming vanishing third moments.
When $\|X_i\|_\infty \leq M$, estimating the average population
moments introduces an additional dimension factor of $d^{1/8}$.
For the nonparametric bootstrap, Theorem~\ref{thm:refined_bootstrap}
improves the sample size factor to $n^{-1/3}$ under bounded sixth
moments and a common covariance matrix $\Sigma$.
With $\|X_i\|_2 \leq M$ and $p \geq 2$,
the resulting rate depends on $d$ only via
the smallest non-zero eigenvalue $\lambda$ of $\Sigma$;
since $\lambda \rank \Sigma \leq \tr \Sigma \leq M^2$,
this dependence is of order
$\lambda^{-4/3} \gtrsim (\rank \Sigma)^{4/3}$.
See Section~\ref{sec:literature_bootstrap}
for further discussion of
contemporary bootstrapping methods.

\subsection{KMT-type couplings}

The approach of \citet{komlos1975approximation,komlos1976approximation}
provides an alternative construction,
yielding couplings with Gaussian summands and
controlling the maximum deviation over the first $n$ partial sums.
For real-valued i.i.d.\ summands possessing finite exponential moments,
the resulting coupling rate is $(\log 2n)/\sqrt n$.
\citet{sakhanenko1984rate} extended this method to non-identically
distributed summands with finite exponential moments, and
\citet{sakhanenko1985estimates} established a bound under
polynomial moments.
\citet{einmahl1989extensions} further generalized
these constructions to $\R^d$-valued random vectors,
obtaining a rate of $n^{-1/2+1/r}$ for i.i.d.\ summands
possessing finite $r$th moments, where $r>3$,
and a rate of $(\log 2n)/\sqrt n$ assuming finite exponential
moments alongside an additional smoothness condition.
\citet{zaitsev1998multidimensional} showed how to
remove the latter condition and gave a coupling with
explicit dependence on the dimension $d$, and
\citet{zaitsev2001multidimensional,
  zaitsev2002multidimensionalII,
zaitsev2002multidimensionalIII}
permitted the summands to have different covariance matrices,
under a block covariance condition.
These constructions also use moment matching arguments;
\citet[Theorem~1.5]{zaitsev1998multidimensional}
employed auxiliary random vectors sharing the first three moments
of the summands, and he used this result in the proofs of his
Theorems~1.3 and~1.6.
Moreover, \citet[Theorem~1.6]{zaitsev1998multidimensional} and
\citet[Theorem~1.2]{zaitsev2001multidimensional}
showed that vanishing third moments permit the removal of
a factor of $\log 2d$
from the corresponding exponential moment bounds in
Panel~B of Table~\ref{tab:literature}.
Couplings based on finite moment bounds for heterogeneous arrays are given by
\citet{zaitsev2007estimates} and \citet{gotze2009bounds}, with i.i.d.\
refinements in \citet{gotze2010rates,gotze2011estimates};
see \citet{Zaitsev_2013_RMSurveys}
for a review and further references.
More recently, \citet{ye2025computable} derived computable scalar KMT thresholds
for bounded i.i.d.\ summands at the cost of an additional logarithmic factor,
and \citet{waudbysmith2026nonasymptotic} obtained non-asymptotic and
distribution-uniform versions of scalar KMT approximations.
Empirical process constructions using KMT ideas are discussed in
Appendix~\ref{sec:literature_empirical_process}.

Panel~B of Table~\ref{tab:literature} reports KMT coupling rates;
Appendix~\ref{sec:literature_details_hungarian} gives their
derivations and some more precise bounds.
These constructions control all partial sums simultaneously;
in contrast, Theorem~\ref{thm:yurinskii} targets only the full sum
and its coupling depends on the threshold $\eta$.
With finite $r$th moments,
these constructions attain a sample size factor of $n^{-1/2+1/r}$ without
requiring higher-order matching, and under
finite exponential moments, they attain $n^{-1/2}$,
each up to logarithmic factors.
They also require additional conditions such as
a common covariance matrix
\citep{zaitsev1998multidimensional},
i.i.d.\ summands with a nonsingular covariance matrix
\citep{gotze2011estimates},
or block covariance criteria
\citep{zaitsev2001multidimensional,
  zaitsev2002multidimensionalII,zaitsev2002multidimensionalIII,
zaitsev2007estimates,gotze2009bounds}.
We assume in addition that
the relevant covariance matrices have bounded condition numbers.

Our Theorem~\ref{thm:yurinskii} gives polynomially smaller
dimension factors than these methods
and imposes no covariance restrictions, but
its sample size factor $n^{-1/2+1/(k+\gamma)}$ approaches
$n^{-1/2}$ only as the matching order $k$ diverges.
Under finite $r$th moments, our Theorem~\ref{thm:yurinskii}
attains the same sample size factor $n^{-1/2+1/r}$
by taking $k + \gamma = r$,
with Gaussian targets when $r \in (2, 3]$
and with exact matching up to order $k = \lceil r \rceil - 1$.
For example, with $p = 2$ and $r \in (2, 3]$,
we attain a dimension factor of
$d^{1/2+1/r}$, compared with $d^{8+\vartheta}$
for \citet{zaitsev2007estimates}.
With bounded summands, the rates in Panel~B
tend to zero only if $d / n^{2/23} \to 0$,
while our Theorem~\ref{thm:yurinskii} with $\|X_i\|_2 \leq M$,
$p = 2$, $k \geq 2$ and $\gamma = 1$
only requires $d / n^{1-1/k} \to 0$.

Applied conditionally on the data, these KMT-type
constructions may also provide bootstrap
couplings when their assumptions hold for the bootstrap law.
In that case, the population and bootstrap covariance
matrices must also be compared;
see Section~\ref{sec:literature_bootstrap}.
Our Theorems~\ref{thm:bootstrap} and~\ref{thm:refined_bootstrap}
couple the conditional bootstrap law directly to the marginal law of $S_n$
and show how approximate matching of higher-order average moments
can improve upon the classical second-order Yurinskii coupling rate.

\subsection{Wasserstein-based couplings}
\label{sec:literature_wasserstein}

The $\ell_p$-Wasserstein metric of order $s \geq 1$
between two probability measures $P$ and $Q$ on $\R^d$
with finite $s$th moments is
\begin{equation*}
  W_{s,p}(P, Q)
  \coloneqq
  \inf
  \Bigl\{
    \bigl( \E\|X-Y\|_p^s \bigr)^{1/s}:\,
    X \sim P, \ Y \sim Q
  \Bigr\},
\end{equation*}
where the infimum ranges over couplings of $P$ and $Q$.
Any coupling $(X, Y)$ attaining this infimum therefore satisfies
$\P(\|X - Y\|_p > \eta) \leq W_{s,p}(P, Q)^s / \eta^s$
for every $\eta>0$ by Markov's inequality.

Assuming bounded summands and a common covariance matrix,
\citet{zhai2018high} controlled the quadratic Wasserstein
distance $W_{2,2}$ between
$S_n$ and its Gaussian approximation using Talagrand's
transportation inequality.
\citet{bonis2020stein} used Stein's method to remove
logarithmic factors for i.i.d.\ summands
with nonsingular covariance, also deriving
bounds under finite moments and for i.n.i.d.\ summands;
see also \citet{buzun2022strong} and \citet{fang2023wasserstein}
for higher transport orders $s$.
An improvement of the logarithmic factor for bounded
i.i.d.\ summands with arbitrary covariance
was obtained by \citet{eldan2020clt} using martingale embeddings.
\citet{mies2023sequential} extended their result
to i.n.i.d.\ unbounded summands and used blocking
to construct a coupling for all partial sums,
albeit at a slower rate.

Other transport bounds exploit additional regularity.
With the dimension $d$ fixed,
\citet{davie2017polynomial} and \citet{zhang2019multidimensional}
obtained faster rates assuming higher-order
matching and Cram{\'e}r's condition.
For isotropic and absolutely continuous summands with uniformly bounded
Poincar{\'e} constants, \citet{courtade2019existence} used Stein kernels
to obtain a $W_{2,2}$ bound of order $(d/n)^{1/2}$.
For i.i.d.\ isotropic summands with vanishing third moments
and satisfying a Poincar\'e inequality,
\citet{fathi2021higher}
used higher-order Stein kernels to obtain
a $W_{2,2}$ bound of order $\sqrt{d} / n$
up to a logarithmic factor,
under additional regularity conditions.
\citet{fang2024logconcave} constructed a
coupling with bounds for every one-dimensional projection,
giving much weaker dimension dependence under
$\|\cdot\|_\infty$, assuming i.i.d.\ log-concave summands.
For i.i.d.\ isotropic summands with finite fourth moments and an
absolutely continuous component, \citet{bonis2024improved} obtained
a $W_{2,2}$ bound with a leading $n^{-1/2}$ term
determined by third moments.
\citet{davie2026asymptotic} obtained an asymptotic expansion
of the $W_{s,p}$ distance
assuming i.i.d.\ scalar random variables with all moments finite
and satisfying Cram{\'e}r's condition; the leading term gains
a factor of $n^{-1/2}$ for each additional matching moment.
Several of these works rely on Poincar{\'e} inequalities or
Cram{\'e}r conditions, which generally fail for empirical
distributions appearing in the nonparametric bootstrap.
\citet{wu2026uncertainty} obtained $W_{1,2}$ bounds
for vector-valued martingales with non-random terminal
predictable quadratic variation.
\citet{zhang2026martingale} obtained $W_{s,2}$ bounds,
where $s \geq 2$,
for weighted Poisson martingales generated by uniformly ergodic
Markov chains, under moment conditions of order exceeding $2s$.

In concurrent work, \citet{fang2026martingale} obtain $W_{s,p}$ bounds
for martingales, where $p\in[1,\infty]$ and $s \geq 1$,
by combining Lindeberg's method with a Gaussian convolution inequality.
For independent summands, their Corollary~2.1 gives a
rate of $n^{-1/2}$ up to logarithmic factors,
under a uniform ellipticity condition
relating to lower bounds on the covariance matrices of the summands.
With $s = 1$ and $p \geq 2$, their Corollary~2.2 requires no such
covariance conditions and gives
the same powers
of $n$ and $d$ as our Theorem~\ref{thm:yurinskii}
with $k=2$ and $\gamma=1$,
using a single Gaussian coupling valid at every threshold.
Their bounds on $W_{1,p}$ require finite third moments.
\citet{fang2026sharp} remove the logarithmic factor for
independent summands with bounded Euclidean norms and
arbitrary covariance
matrices, obtaining a rate of $\sqrt{d/n}$ for $W_{2,2}$;
this rate is optimal in general \citep{zhai2018high}.
They also sharpen the dimension dependence for
isotropic summands with finite fourth moments,
and further combine these Wasserstein bounds with a KMT-type dyadic
construction to couple all partial sums without covariance restrictions,
assuming i.i.d.\ summands.

Panel~C of Table~\ref{tab:literature} reports selected rates
associated with these results;
Appendix~\ref{sec:literature_details_wasserstein}
gives further details.
Most rely on bounding $W_{s,2}$ and applying
$\|x\|_p \leq \|x\|_2$ for $p \geq 2$.
The rates for \citet{bonis2020stein}
assume i.i.d.\ summands whose
positive definite covariance matrix $\Sigma$
has bounded condition number;
under $\vvvert e_j^\T X_i \vvvert_r \leq M$,
we also assume $\lambda_{\min}(\Sigma) \gtrsim M^2$.

In Panel~D, the target summands $Y_i$ for our Theorem~\ref{thm:yurinskii}
match the moments of $X_i$ up to order $k$
and satisfy the same bounds as $X_i$.
They are therefore non-Gaussian
under $\|X_i\|_\infty \leq M$ or $\|X_i\|_2 \leq M$,
and are generally non-Gaussian for $k \geq 3$;
in contrast, the results in Panel~C use Gaussian targets.
Firstly, if $\|X_i\|_\infty \leq M$,
\citet[Theorem~2.1 and Corollary~2.1]{fang2026sharp} attain a rate of
$d / \sqrt n$ while Theorem~\ref{thm:yurinskii} achieves
$n^{-1/2 + 1/(k+\gamma)} (d \log 2d)^{k/(2k + 2\gamma)}$,
which is better whenever $n^{2/(k+2\gamma)} / d \to 0$,
up to logarithmic factors.
Taking $k=3$ and $\gamma = 1$ shows
that higher-order moment matching gives a coupling error
bound that converges to zero and is
smaller than that of \citet{fang2026sharp} whenever
$n^{2/5}/d \to 0$ and $n^{2/3}/d \to \infty$, up to logs.
Under the same assumptions,
\citet[Corollary~2.2]{fang2026martingale} attain a rate of
$n^{-1/6} (d\log 2d)^{1/3}$.
Taking $\gamma=1$, our rate is faster whenever
$k \geq 3$ and $(d \log 2d) / n^2 \to 0$;
this includes every regime in which our bound converges to zero.
Secondly, assuming $\|X_i\|_2 \leq M$,
\citet{bonis2020stein}
attains a rate of $\sqrt{d/n}$.
Here, taking $\gamma = 1$, Theorem~\ref{thm:yurinskii} achieves
$n^{-1/2 + 1/(k+1)} (\log 2d)^{k/(2k+2)}$,
and hence our rate is faster whenever $n^{2/(k+1)} / d \to 0$,
up to logarithmic terms.
Thirdly, with $\vvvert e_j^\T X_i \vvvert_r \leq M$
where $r > 2$,
Theorem~\ref{thm:yurinskii} attains a rate of
$n^{-1/2 + 1/r} \sqrt{d\log 2d}$ by taking $k + \gamma = r$,
with Gaussian targets when $2 < r \leq 3$
and with exact matching up to order
$\lceil r \rceil - 1$ when $r > 3$.
\citet{mies2023sequential} obtained the same sample size factor
(up to logs) with a larger dimension factor of $d^{1-1/r}$.

\subsection{Distributional approximation and bootstrapping}
\label{sec:literature_bootstrap}

Another strand of the literature studies distributional approximations
over classes of sets, rather than couplings.
For example, if $T_n$ is Gaussian with the same mean
and covariance as $S_n$, a typical objective is to bound
\begin{equation}
  \label{eq:distributional_approximation}
  \sup_{A \in \cA}
  \bigl|\P\bigl(S_n/\sqrt n\in A\bigr)
  - \P\bigl(T_n/\sqrt n\in A\bigr)\bigr|,
\end{equation}
where $\cA$ is a structured class, such as
the collection of convex sets, Euclidean balls, or
hyperrectangles in $\R^d$.
Classical multivariate Berry--Esseen theorems bound
\eqref{eq:distributional_approximation} when
$\cA$ is the class of convex sets; see, for example,
\citet{bentkus2003dependence} and \citet{fang2024large}.
Approximations over Euclidean balls include results based on
higher-order moment matching and Edgeworth expansions
\citep{zhilova2020nonclassical,zhilova2022new}.
For hyperrectangles, approaches include smoothed Lindeberg interpolations
or Stein's method, along with Gaussian anti-concentration inequalities
\citep{chernozhukov2013gaussian,chernozhukov2017central,
  chernozhuokov2022improved,fang2021high,lopes2022central,
chernozhukov2023nearly}.
These results often use the geometry of $\cA$ to obtain better dimension
dependence than coupling-based approaches.

Bootstrap-based methods for distributional approximation
are of particular relevance to our work; see Appendix~\ref{sec:ks}
for applications of our main results to Kolmogorov--Smirnov approximations
under moment matching conditions and in bootstrap settings.
We refer the reader to \citet{chernozhukov2023high}
for a detailed review.
For approximating $\|S_n\|_\infty$, or with $\cA$ a class of
hyperrectangles,
\citet{chernozhukov2013gaussian,chernozhukov2017central}
established bootstrap approximations with logarithmic
dependence on the dimension, allowing $d$ to be much larger than $n$ under
appropriate moment and tail conditions.
\citet{chernozhuokov2022improved} obtained sharper bounds using an
iterative randomized Lindeberg method, including results for empirical
and multiplier bootstraps and allowing for
singular covariance matrices.
\citet{lopes2022central} and \citet{chernozhukov2023nearly} obtained bootstrap
approximation rates of $n^{-1/2}$ up to logarithmic factors,
under tail and moment conditions and a lower bound
on the eigenvalues of the correlation matrix.
More recently, \citet{lopes2026general} obtained Gaussian and
Gaussian multiplier bootstrap approximation rates of $n^{-1/2}$,
up to logarithmic factors, for i.i.d.\ summands without covariance
eigenvalue lower bounds, under suitable sub-Gaussian conditions.
\citet{chernozhukov2026sharp} obtained the same Gaussian approximation rate
for independent but not necessarily identically distributed summands
with bounded or sub-exponential coordinates.
Such logarithmic dependence on the dimension requires light tails:
with only $r > 2$ finite moments,
\citet{kock2024remark} showed that Gaussian approximation of the maximum,
even at a fixed critical value, can fail when
$d / n^{r/2 - 1 + \varepsilon} \to \infty$
for some $\varepsilon > 0$.
For Euclidean balls, \citet{fang2024large}
achieved an approximation rate of $d / \sqrt n$, up to logs,
assuming an identity covariance matrix and bounded fourth moments.
\citet{choudhury2026optimal} showed that $d / \sqrt n \to 0$ is necessary
in general for Gaussian approximation over all Euclidean balls,
even for i.i.d.\ summands with identity covariance and uniformly
bounded coordinatewise fourth moments.

The use of higher moment matching in bootstrap
theory is familiar from classical asymptotic results,
where bootstraps can reproduce Edgeworth terms
that are absent in Gaussian approximations
\citep[Section~3.3]{shao1995jackknife}.
In high-dimensional settings, a closely related contribution is
that of \citet{deng2020beyond}.
In their Theorems~1 and 2,
they compare the distribution of the maximum
statistic $\max_{j \in [d]} \bigl(e_j^\T S_n\bigr)$
with the bootstrap law of $\max_{j \in [d]} \bigl(e_j^\T S_n^*\bigr)$,
under a nonparametric or third-order multiplier bootstrap.
Their analysis used similar techniques to our main results,
employing smooth approximation, a related randomized Lindeberg method,
approximate average moment matching up to order three,
and anti-concentration inequalities.
However, their approach is explicitly tailored to the maximum
statistic and their results do not include
couplings for the full random vectors $S_n$ and $S_n^*$.
They also considered higher-order approximation under some strong
assumptions; for example, their Theorem~3 requires vanishing
third and fifth moments.
\citet{deng2020conservative} further refined their
randomized Lindeberg method,
obtaining bootstrap coverage error bounds of order $n^{-1/4}$
up to logarithmic factors.

Another closely related work is that of
\citet{zhilova2020nonclassical}, who developed
higher-order bootstrap approximations over Euclidean balls.
Her Theorems~2.1 and~2.2 use Lindeberg's method
under exact moment matching assumptions to derive
Berry--Esseen-type inequalities;
she then constructs distributions with non-degenerate Gaussian
components to provide anti-concentration bounds,
similar to the proof of our Theorem~\ref{thm:ks}%
\iftoggle{journal}{ in the supplementary material
\citep{cattaneo2026highsupplement}}{}.
However, the resulting distributional approximation results
do not quantify the dependence on the dimension
when matching more than three moments; see her Remark~2.2.
Her Theorems~4.1 and~4.3 then turn to nonparametric and
multiplier bootstrap schemes, appealing to approximate
moment matching arguments.
However, some steps in the proof of Theorem~4.3 require
additional justification.
To be precise, Equation~(B.17) in
\citet[supplementary material]{zhilova2020nonclassical} relies on
(i) applying Lemma~B.1 to terms arising from the
bootstrap Lindeberg decomposition as though they were independent; and
(ii) using a pointwise probability bound to directly control
a supremum over Euclidean balls.

\citet{zhilova2022new} derived Edgeworth-type inequalities over
Euclidean balls and half spaces, with bounds depending
on moment discrepancies up to order three.
Her main application was to the nonparametric bootstrap,
and her Theorem~5.1 states a high-probability bound on the bootstrap error,
assuming i.i.d.\ summands with sub-Gaussian coordinates and
a non-degenerate covariance matrix.
However, the proof of Theorem~5.1 on page~10 of
\citet[supplementary material]{zhilova2022new}
uses bounds on empirical third and fourth moments that
appear to require further justification.

More recently, \citet{koike2026high} studied the accuracy of the
bootstrap for maxima through high-dimensional Edgeworth and
Cornish--Fisher expansions.
He showed that, under Stein kernel conditions and other regularity
assumptions,
third moment matching can allow for bootstrap refinements
not provided by Gaussian approximations,
improving coverage error rates from $n^{-1/2}$ to $n^{-1}$,
up to logarithmic factors.
In particular, his Corollary~2.2
establishes higher-order accuracy of the multiplier bootstrap
under additional covariance assumptions and
in sufficiently high-dimensional regimes.
He further shows in Theorem~2.4 that the
double multiplier bootstrap can outperform Gaussian approximations
without these covariance or high-dimensional restrictions.
\citet{feng2026bootstrap} obtained analogous second-order refinements
for coordinate order statistics of fixed rank, under Stein kernel,
tail, variance, and weak dependence conditions.

\section{Strong approximation of empirical processes}
\label{sec:empirical_process}

We use the vector couplings of Sections~\ref{sec:yurinskii}
and~\ref{sec:bootstrap} to obtain non-asymptotic strong approximation
bounds for empirical processes indexed by classes of functions.
For background material on empirical process theory, see
\citet{van1996weak} and \citet{gine2016mathematical}.

The aim of this section is to construct couplings for
empirical processes $\bigl(S_n(f): f \in \cF\bigr)$, where
$S_n(f) \coloneqq \sum_{i=1}^n \bigl(f(X_i) - \E\{f(X_i)\}\bigr)$
with $X_1, \ldots, X_n$ independent, and $\cF$ is a class
of real-valued functions.
Our approach consists of two steps.
Firstly, we discretize the function class by defining
a finite subclass $\tilde\cF \subseteq \cF$ (of cardinality $d$).
Our previous coupling theorems can then be applied
to the $\R^d$-valued vector of evaluations
$\bigl(S_n(f): f \in \tilde\cF\bigr)$.
Secondly, we control the discretization errors
of each process uniformly over
the class $\check\cF$ of differences $f-\pi(f)$,
where $\pi: \cF \to \tilde\cF$.
Our results allow for arbitrary choices of the subclass $\tilde\cF$
and explicitly account for the associated dimension $d$.

We establish a general coupling result between two empirical processes
$S_n(f)$ and $T_n(f)$, analogous to Theorem~\ref{thm:yurinskii}.
Under a VC-type condition
on the uniform $L_2$-covering numbers of $\cF$
and suitable moment matching conditions,
we obtain explicit coupling rates as a corollary.
Following this,
we turn to strong approximations based on
nonparametric or multiplier bootstrap schemes.
Our main result on bootstrapped empirical processes is stated for
VC-type classes; we present a more general version as
Theorem~\ref{thm:empirical_process_bootstrap}%
\iftoggle{journal}{ in the supplementary material
\citep{cattaneo2026highsupplement}}{}.
Appendix~\ref{sec:literature_empirical_process} compares
our empirical process results with the literature.
We begin by giving our assumptions on the indexing
and approximating function classes.

\begin{assumption}[Function class]
  \label{ass:function_class}
  \phantom{a}
  \begin{enumerate}
    \item
      \label{it:ass_function_class}
      Let $\cX$ be a Polish space with Borel $\sigma$-algebra $\cH$
      and take $\cF$ to be a class of measurable functions from $\cX$ to $\R$.
      Suppose there is a countable $\cF' \subseteq \cF$ such that
      each $f \in \cF$ is the pointwise limit of a sequence in $\cF'$.
      Let $F: \cX \to \R$ be a nonzero measurable function satisfying
      $F(x) \geq |f(x)|$ for all $x \in \cX$ and $f \in \cF$.
      For a measurable function $f: \cX \to \R$
      and a measure $Q$ on $(\cX, \cH)$,
      define $\|f\|_{Q,2}^2 \coloneqq \int_\cX |f(x)|^2 \diffi Q(x)$.
      For a class $\cG$ of measurable functions from $\cX$ to $\R$ and
      $\varepsilon \in (0, \infty]$,
      and with $\min \emptyset \coloneqq \infty$,
      define the covering number
      \begin{equation*}
        N\bigl(\varepsilon, \cG, \|\cdot\|_{Q,2}\bigr)
        \coloneqq
        \min
        \Bigl\{
          N \in \N:
          \exists
          f_1, \ldots, f_N \in \cG
          \text{ such that }
          \sup_{f \in \cG}
          \min_{j \in [N]}
          \|f - f_j\|_{Q,2} \leq \varepsilon
        \Bigr\}.
      \end{equation*}
      If $\cG$ is a finite class, write
      $\cG(x) \coloneqq \bigl(f(x): f \in \cG\bigr)$
      for each $x \in \cX$
      and $Q(\cG) \coloneqq \bigl(Q(f): f \in \cG\bigr)$,
      where $Q(f) \coloneqq \int_\cX f(x) \diffi Q(x)$,
      whenever these integrals exist.

    \item
      \label{it:ass_approximating_class}
      Let $\tilde\cF \subseteq \cF$ be a nonempty
      finite subclass of cardinality
      $d \coloneqq |\tilde \cF|$.
      Let $\pi: \cF \to \tilde \cF$ and define
      $\check\cF \coloneqq \{f - \pi(f): f \in \cF\}$.
      Suppose there is a countable $\check\cF' \subseteq \check\cF$ such that
      each $f \in \check\cF$ is the pointwise limit
      of a sequence in $\check\cF'$.
      Let $\check F: \cX \to \R$ be a nonzero measurable function with
      $\check F(x) \geq |f(x)|$
      for all $x \in \cX$ and $f \in \check \cF$.
      Define the uniform $L_2$-entropy integral
      \begin{equation*}
        \cJ(\varepsilon)
        \coloneqq
        \int_0^\varepsilon
        \sup_{Q \in \cQ}
        \sqrt{1 + \log N\bigl(
            \xi \|\check F\|_{Q,2},
            \check\cF, \|\cdot\|_{Q,2}
        \bigr)}
        \diffi \xi,
      \end{equation*}
      where $\cQ$ is the set of probability measures $Q$ on $(\cX, \cH)$ with
      $\|\check F\|_{Q,2} > 0$.

  \end{enumerate}
\end{assumption}

Next, we introduce the random variables underlying the empirical process
and impose some moment conditions.
We use the sub-Weibull Orlicz quasi-norm $\vvvert \cdot \vvvert_{\omega_r}$
as defined in Section~\ref{sec:introduction}.

\begin{assumption}[Empirical process]%
  \label{ass:empirical_process}
  \phantom{a}
  \begin{enumerate}
    \item
      \label{it:ass_empirical_process}
      Grant Assumption~\ref{ass:function_class}\ref{it:ass_function_class}.
      For $n \in \N$, suppose $X_1, \ldots, X_n$
      are independent $\cX$-valued random variables.
      For $i \in [n]$,
      denote the law of $X_i$ by $P_{X,i}$
      and write $\bar P_X \coloneqq \sum_{i=1}^n P_{X,i}/n$.
      Suppose there exist $r \in (0, \infty]$ and $M > 0$ such that
      $\vvvert F(X_i) \vvvert_{\omega_r} \leq M$
      for each $i \in [n]$.

    \item
      \label{it:ass_approximating_process}
      Grant Assumption~\ref{ass:function_class}\ref{it:ass_approximating_class}.
      Suppose that there exists $\check M > 0$ such that
      $\vvvert \check F(X_i) \vvvert_{\omega_r} \leq \check M$
      for each $i \in [n]$.
      Fix $k \in \N$ with $k \geq 2$ and $\gamma \in (0, 1]$.
      Take $m \in [k]$ and suppose that
      for each $\alpha \in \N_0^d$ with $|\alpha| \in [m]$,
      there exists $\mu_{X,\alpha} \in \R$ such that
      $\E\bigl\{\bigl(\tilde\cF(X_i) - P_{X,i}(\tilde\cF)\bigr)^\alpha\bigr\}
      = \mu_{X,\alpha}$ for all $i \in [n]$.
      For $a \in \N_0$ and $b \geq 0$, let
      $\bar\beta_{X,a,b} \coloneqq
      \sum_{i=1}^n \beta_{X,a,b,i} / n$,
      where for each $i \in [n]$,
      \begin{equation*}
        \beta_{X,a,b,i}
        \coloneqq
        \frac{1}{\sqrt{a!}}
        \E \Bigl\{
          \bigl\|\tilde\cF(X_i) - P_{X,i}(\tilde\cF)\bigr\|_2^a
          \bigl\|\tilde\cF(X_i) - P_{X,i}(\tilde\cF)\bigr\|_\infty^b
        \Bigr\}.
      \end{equation*}
      Suppose that $\|\check F\|_{\bar P_X,2} > 0$
      and let $\sigma > 0$ satisfy
      $\|\check f\|_{\bar P_X,2} \leq \sigma$
      for all $\check f \in \check\cF$.
      Take $\nu \in (0, 1]$ satisfying
      $\nu \leq \sigma / \|\check F\|_{\bar P_X,2}$.
      Fix $\kappa \geq 1$ and, for $t > 0$, define
      \begin{equation*}
        R(t)
        \coloneqq
        \biggl\{
          \sigma \sqrt n
          + \check M (\kappa \log 2n)^{1/r}
          \biggl(\frac{\cJ(\nu)}{\nu} + t\biggr)
        \biggr\}
        \biggl(\frac{\cJ(\nu)}{\nu} + t\biggr).
      \end{equation*}

  \end{enumerate}
\end{assumption}

Since the summands in $S_n(f)$ are centered, the condition in
Assumption~\ref{ass:empirical_process}\ref{it:ass_approximating_process}
involving $\mu_{X,\alpha}$ is automatic for $m = 1$.
We now introduce a second empirical process $T_n(f)$ and define
the resulting moment discrepancies.

\begin{assumption}[Strong approximation]%
  \label{ass:strong_approx}
  \phantom{a}
  \begin{enumerate}
    \item
      \label{it:ass_strong_approx_class}
      Suppose that
      $\bX_n \coloneqq (X_1, \ldots, X_n)$ and
      $\bY_n \coloneqq (Y_1, \ldots, Y_n)$
      both satisfy
      Assumption~\ref{ass:empirical_process}\ref{it:ass_empirical_process}
      with the same values of $r$ and $M$.
      For $f \in \cF$,
      let $S_n(f) \coloneqq \sum_{i=1}^n \bigl(f(X_i) - P_{X,i}(f)\bigr)$
      and $T_n(f) \coloneqq \sum_{i=1}^n \bigl(f(Y_i) - P_{Y,i}(f)\bigr)$.

    \item
      Suppose further that $\bX_n$ and $\bY_n$ satisfy
      Assumption~\ref{ass:empirical_process}\ref{it:ass_approximating_process}
      with the same values of
      $\check M$, $\kappa$, $\sigma$, $\nu$,
      $k$, $\gamma$ and $m$.
      Write $\beta_{a,b,i} \coloneqq \beta_{X,a,b,i} + \beta_{Y,a,b,i}$
      and $\bar\beta_{a,b} \coloneqq \bar\beta_{X,a,b} + \bar\beta_{Y,a,b}$.
      For each $\alpha \in \N_0^d$ with $|\alpha| \in [k]$ and $i \in [n]$, let
      \begin{equation*}
        \delta_{\alpha,i}
        \coloneqq
        \E\bigl\{
          \bigl(\tilde\cF(X_i) - P_{X,i}(\tilde\cF)\bigr)^\alpha
        \bigr\}
        -
        \E\bigl\{
          \bigl(\tilde\cF(Y_i) - P_{Y,i}(\tilde\cF)\bigr)^\alpha
        \bigr\},
      \end{equation*}
      and write $\bar\delta_\alpha \coloneqq \sum_{i=1}^n \delta_{\alpha,i}/n$.
      For $\ell \in \{m+1, \ldots, k\}$, take $k_\ell \in \{0, \ldots, m\}$
      and, for $\eta, s > 0$, define
      $\Delta(s)$ and $\Gamma_\eta(s)$
      as in Assumption~\ref{ass:data}.

  \end{enumerate}
\end{assumption}

The next theorem gives a coupling between
the empirical processes $S_n(f)$ and $T_n(f)$.
We write $S_n' \equallaw S_n$ to indicate
equality of all finite-dimensional distributions.

\begin{theorem}[Strong approximation of empirical processes]%
  \label{thm:empirical_process}
  Grant Assumption~\ref{ass:strong_approx}
  and write $\phi(t) \coloneqq \phi_{\infty,d}(t)$.
  There exists $C > 0$ depending only on $r$
  such that for each $\eta, t > 0$,
  there is a probability space supporting
  processes $S_n'$ and $T_n'$, where $S_n' \equallaw S_n$
  and $T_n' \equallaw T_n$, with
  \begin{align*}
    &\P\biggl\{
      \sup_{f \in \cF}
      \bigl|S_n'(f) - T_n'(f)\bigr|
      > 3 \eta
      + C R(t)
    \biggr\} \\
    &\qquad\leq
    4 e^{-t^2 / 2}
    + \frac{n \bar\beta_{k,\gamma} \phi(t)^k}{\eta^{k+\gamma}}
    + \Delta\biggl(\frac{\phi(t)}{\eta}\biggr)
    + \Gamma_\eta\biggl(\frac{\phi(t)}{\eta}\biggr)
    + \frac{1}{n^\kappa}.
  \end{align*}
\end{theorem}

The terms involving $\bar\beta_{k,\gamma}$,
$\Delta$ and $\Gamma_\eta$ are familiar from
Section~\ref{sec:yurinskii} and come from a
coupling construction on the finite subclass $\tilde\cF$
with $p = \infty$.
The resulting random vectors
$S_n'(\tilde \cF)$ and $T_n'(\tilde \cF)$
are then extended to
empirical processes $S_n'$ and $T_n'$
over the full class $\cF$.
The entropy term $R(t)$ controls the magnitude of the
errors in $S_n'$ and $T_n'$
due to replacing $\cF$ with $\tilde\cF$ using the map $\pi$.
More precisely, this term simultaneously bounds
$\sup_{f \in \check\cF} |S_n'(f)|$
and $\sup_{f \in \check\cF} |T_n'(f)|$,
where $\check\cF \coloneqq \{f - \pi(f): f \in \cF\}$,
with probability at least
$1 - 2 e^{-t^2/2} - 1/n^\kappa$.
This maximal inequality relies on Orlicz norm control of the envelope
function $\check F$ and depends on the
uniform $L_2$-entropy integral $\cJ(\varepsilon)$
of the difference class $\check\cF$
\citep[see][]{boucheron2013concentration,chernozhukov2014gaussian}.
Theorem~\ref{thm:empirical_tail}%
\iftoggle{journal}{ in the supplementary material
\citep{cattaneo2026highsupplement}}{}
gives a general version of our empirical process supremum tail bound,
and Figure~\ref{fig:empirical_process} illustrates the two stages
of our strong approximation construction.

\iftoggle{journal}{\begin{figure}[b]
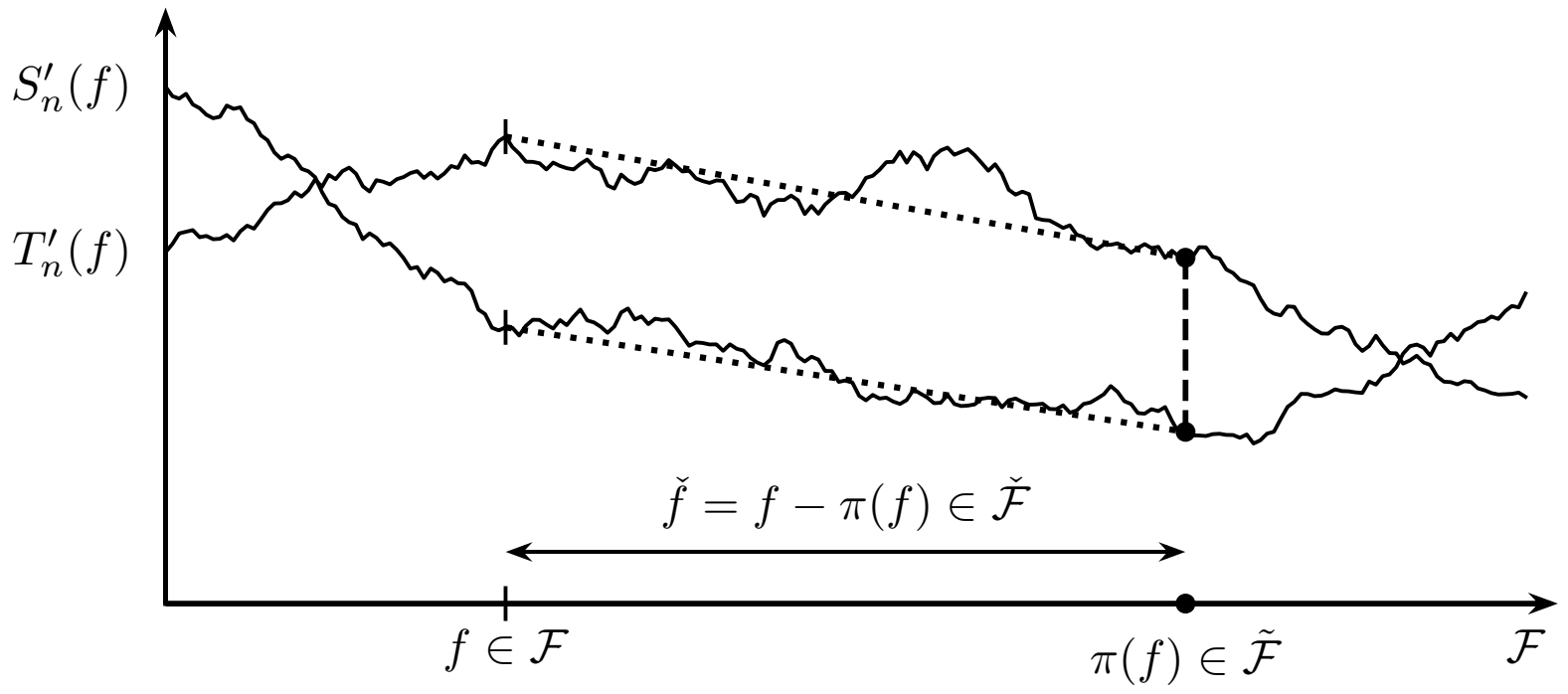
}{\begin{figure}[t]}
  \centering
  \scalebox{\myfigscale}{%
  \csname input\endcsname{strong_approximation.pgf}}
  \caption{Illustration of our empirical process strong approximation
    construction.
    The dashed line represents the coupling error on $\tilde\cF$,
    while the dotted lines represent the discretization errors
    $S_n'(\check f)$
  and $T_n'(\check f)$, where $\check f = f - \pi(f) \in \check\cF$.}
  \label{fig:empirical_process}
\end{figure}

The variables $X_1, \ldots, X_n$ are assumed to be independent
but need not be identically distributed; likewise for $Y_1, \ldots, Y_n$.
This allows applications with grouped or heterogeneous data, for example.
Further, Theorem~\ref{thm:empirical_process} is fully non-asymptotic,
with a constant depending only on the Orlicz parameter $r$,
allowing for the function class
and distributions to vary freely with the sample size $n$.
This theorem is stated as a coupling between
two empirical processes as this is most useful for our upcoming
results on bootstrapping empirical processes.
However, it is also possible to use alternative coupling processes
$T_n$; we provide an analogous Gaussian process strong approximation result
as Theorem~\ref{thm:gaussian_process}%
\iftoggle{journal}{ in the supplementary material
\citep{cattaneo2026highsupplement}}{}.

We now impose a Vapnik--Chervonenkis (VC)-type condition
on the uniform $L_2$-covering numbers of the class $\cF$,
thus allowing us to bound the entropy integral $\cJ(\varepsilon)$
when taking $\check F = 2 F$.

\begin{assumption}[VC-type entropy]%
  \label{ass:vc}
  Grant Assumption~\ref{ass:function_class}\ref{it:ass_function_class}.
  Suppose that there exist $a \geq 1$ and $u \geq 0$ such that
  $\sup_Q N\bigl(\varepsilon \|F\|_{Q,2}, \cF, \|\cdot\|_{Q,2}\bigr)
  \leq a \varepsilon^{-u}$ for all $\varepsilon \in (0, 1]$,
  where the supremum is over
  the set of probability measures $Q$ on $(\cX, \cH)$ with
  $\|F\|_{Q,2} > 0$.
\end{assumption}

The next result combines this assumption with
an exact moment matching condition for each observation,
up to order $k$ and over all finite subclasses.
The proof selects $\tilde\cF$
to balance the resulting coupling and discretization errors;
in general, enlarging $\tilde\cF$ decreases the
discretization error but increases the coupling error
via the dimension $d = |\tilde\cF|$.

\begin{corollary}
  \label{cor:empirical_process_vc}
  Grant Assumptions~\ref{ass:strong_approx}\ref{it:ass_strong_approx_class}
  and~\ref{ass:vc},
  and suppose that $\|F\|_{\bar P_X,2} \land \|F\|_{\bar P_Y,2} > 0$.
  Fix $k \in \N$ with $k \geq 2$ and $\gamma\in(0,1]$.
  Assume further that
  for every nonempty finite $\cG \subseteq \cF$
  and each $\alpha \in \N_0^{|\cG|}$
  with $|\alpha| \in [k]$, we have
  \begin{equation*}
    \E\bigl\{
      \bigl(\cG(X_i) - P_{X,i}(\cG)\bigr)^\alpha
    \bigr\}
    = \E\bigl\{
      \bigl(\cG(Y_i) - P_{Y,i}(\cG)\bigr)^\alpha
    \bigr\},
  \end{equation*}
  for all $i \in [n]$.
  Then there exists $C > 0$
  depending only on $r$, $k$, $a$ and $u$ such that
  for each $\eta > 0$,
  there is a probability space supporting
  processes $S_n'$ and $T_n'$, where $S_n' \equallaw S_n$
  and $T_n' \equallaw T_n$, with
  \begin{equation*}
    \P\biggl\{
      \sup_{f \in \cF} \bigl|S_n'(f) - T_n'(f)\bigr|
      >
      C M \eta
      n^{\frac{4 + u k}{4k + 4\gamma + 2u k}}
      \sqrt{\log 2n}
    \biggr\}
    \leq
    \frac{1}{\eta^{k + \gamma}}.
  \end{equation*}
\end{corollary}

To interpret this bound, take
$k$, $r$, $M$, $a$ and $u$ to be fixed,
set $\gamma = 1$, and let $R_n \to \infty$, so that
\begin{equation*}
  \sup_{f \in \cF}
  \frac{\bigl|S_n'(f) - T_n'(f)\bigr|}{\sqrt n}
  \lesssim_\P
  R_n
  n^{\frac{1 - k}{2k + 2 + u k}}
  \sqrt{\log 2n}.
\end{equation*}
For fixed VC exponent $u$, increasing the
number of matching moments $k$ improves the rate.
As $k\to\infty$, the sample size rate exponent approaches $-1/(2+u)$.
Corresponding Gaussian process bounds are given in
Corollary~\ref{cor:gaussian_vc}%
\iftoggle{journal}{ in the supplementary material
\citep{cattaneo2026highsupplement}}{}.

Coupling to a Gaussian target generally permits
moment matching only up to order two, so consider $k = 2$.
Corollary~\ref{cor:empirical_process_vc} then gives a rate of
$n^{-1/(6+2u)}$, also obtained by \citet{cattaneo2025yurinskii} in a similar
setting; this is faster than the rate $n^{-1/(2+5u)}$ of
\citet{berthet2006revisiting} whenever $u > 4/3$.
If all third moments vanish,
Corollary~\ref{cor:empirical_process_vc}
gives an empirical process coupling with
$k = 3$ and rate $n^{-2/(8+3u)}$, which is
faster than that of \citet{berthet2006revisiting} whenever $u > 4/7$.
These rates are all given up to logarithmic factors;
Appendix~\ref{sec:literature_empirical_process} provides further discussion
of empirical process coupling rates in the literature.

Next, we use the empirical moment approximation ideas
of Section~\ref{sec:bootstrap}
to construct bootstrap-based strong approximations
for empirical processes.
Specifically, we apply arguments from
Theorem~\ref{thm:bootstrap} to establish
a coupling on the finite class $\tilde\cF$
under a nonparametric or multiplier bootstrap scheme.
It remains to control the discretization errors
of the relevant bootstrap process over $\check\cF$;
this is accomplished via a uniform bound for empirical $L_2$-norms
(see Lemmas~\ref{lem:empirical_norm} and~\ref{lem:bootstrap_tail}%
  \iftoggle{journal}{ in the supplementary material
\citep{cattaneo2026highsupplement}}{}).
The general result is given
as Theorem~\ref{thm:empirical_process_bootstrap}%
\iftoggle{journal}{ in the supplementary material
\citep{cattaneo2026highsupplement}}{};
here, we present a VC-type specialization
with an explicit coupling bound.
We show that using a nonparametric or third-order multiplier
bootstrap generally recovers better rates than
those obtained by second-order Yurinskii-based Gaussian approximations.

\begin{theorem}[VC-type bootstrap empirical process coupling]%
  \label{thm:bootstrap_vc}
  Grant Assumption~\ref{ass:vc}.
  Let $\bX_n \coloneqq (X_i: i \in [n])$ satisfy
  Assumption~\ref{ass:empirical_process}\ref{it:ass_empirical_process}
  and suppose that $\|F\|_{\bar P_X,2}>0$.
  For each $i \in [n]$ and $f \in \cF$,
  assume that $P_{X,i}(f) = P_X(f)$ does not depend on $i$.
  Let $k \in \{2, 3\}$,
  write $S_n(f) \coloneqq \sum_{i=1}^n \bigl(f(X_i) - P_X(f)\bigr)$
  and $P_n(f) \coloneqq \sum_{i=1}^{n} f(X_i) / n$,
  and grant one of the following.
  \begin{enumerate}
    \item
      \label{it:thm_ep_bootstrap_vc_nonparametric}
      Let $X_1^*, \ldots, X_n^*$ satisfy Assumption~\ref{ass:bootstrap}.
      For $f \in \cF$,
      take $Y_i^*(f) \coloneqq f(X_i^*) - P_n(f)$ for each $i \in [n]$.
      Set $r' \coloneqq \infty$.

    \item
      \label{it:thm_ep_bootstrap_vc_multiplier}
      Let $\chi_1, \ldots, \chi_n$ be random variables
      satisfying Assumption~\ref{ass:multiplier}
      with order $k$, independent of $\bX_n$.
      Suppose further that
      $\vvvert \chi_i \vvvert_{\omega_{r'}} \leq 2$
      for each $i \in [n]$, where $r' \in (0, \infty]$.
      Set $Y_i^*(f) \coloneqq \chi_i \{f(X_i) - P_n(f)\}$
      for each $i \in [n]$.
  \end{enumerate}
  Write $S_n^*(f) \coloneqq \sum_{i=1}^n Y_i^*(f)$.
  For each $\eta > 0$,
  there is a probability space supporting
  $\bX_n'$, $S_n'$ and $S_n^{*\prime}$, where
  $S_n' \equallaw S_n$, $S_n'$ is independent of $\bX_n'$,
  and $(\bX_n', S_n^{*\prime}) \equallaw (\bX_n, S_n^*)$, with
  \begin{equation*}
    \P\biggl(
      \sup_{f \in \cF}
      \bigl|S_n'(f) - S_n^{*\prime}(f)\bigr|
      > C M \eta
      n^{\frac{4-k+u}{10 - 2k + 2u}} \sqrt{\log 2 n}
    \biggr)
    \leq
    \frac{1}{\eta^2},
  \end{equation*}
  where $C > 0$ depends only on $r$, $r'$, $a$ and $u$.
\end{theorem}

Lemma~\ref{lem:multiplier_weights} in
\iftoggle{journal}{the supplementary material
\citep{cattaneo2026highsupplement}}{%
Appendix~\ref{sec:bootstrap_ep_general}}
verifies the multiplier conditions with $k=2$ for Rademacher ($r'=\infty$),
standard Gaussian ($r'=2$)
and centered standard exponential ($r'=1$) weights,
and with $k \in \{2, 3\}$ for
Mammen $(r'=\infty)$ and centered Poisson
($r'=1$) weights.

Fix $r$, $r'$, $M$, $a$ and $u$, and let $R_n\to\infty$.
For $k = 2$, Theorem~\ref{thm:bootstrap_vc} gives
\begin{equation*}
  \sup_{f \in \cF}
  \frac{\bigl|S_n'(f) - S_n^{*\prime}(f)\bigr|}{\sqrt n}
  \lesssim_\P
  R_n n^{-\frac{1}{6+2u}} \sqrt{\log 2n},
\end{equation*}
matching the rate attained in Corollary~\ref{cor:empirical_process_vc}
with $k=2$, even when the finite-dimensional
population moments of the empirical process
$S_n$ are unknown.
For $k=3$, under the nonparametric or third-order multiplier bootstrap,
we have the improved rate of
\begin{equation*}
  \sup_{f \in \cF}
  \frac{\bigl|S_n'(f) - S_n^{*\prime}(f)\bigr|}{\sqrt n}
  \lesssim_\P
  R_n n^{-\frac{1}{4+2u}} \sqrt{\log 2n}.
\end{equation*}
Up to logarithmic factors and assuming a common mean $P_X(f)$,
we therefore improve upon the Gaussian approximation rate of
$n^{-1/(6 + 2u)}$ given by \citet{cattaneo2025yurinskii}
for all $u \geq 0$.
Our rate is also superior to that of
$n^{-1/(2+5u)}$ achieved by \citet{berthet2006revisiting}
whenever $u > 2/3$.
The more general bootstrap coupling given as
Theorem~\ref{thm:empirical_process_bootstrap}%
\iftoggle{journal}{ in the supplementary material
\citep{cattaneo2026highsupplement}}{}
can give sharper bounds for localized function classes;
see Section~\ref{sec:locpol} for an example.

\section{Application to local polynomial estimation}
\label{sec:locpol}

Local polynomial methods form a popular class
of techniques for nonparametric statistical estimation,
often attaining minimax-optimal rates in smooth settings
by virtue of their bias reduction properties.
We demonstrate how our general results
on bootstrapping empirical processes can be used
to construct strong approximations for such estimators,
focusing on nonparametric regression
and conditional distribution estimation.
In combination with careful anti-concentration arguments,
this approach permits valid and feasible nonparametric statistical inference
based on the construction of uniform bootstrap confidence bands.
We provide analogous results for Gaussian-based inference in
\iftoggle{journal}{the supplementary material
\citep{cattaneo2026highsupplement}}{%
Appendix~\ref{sec:gaussian_approx_locpol}}.
For background material on local polynomial techniques, see
\citet{eggermont2009maximum,samworth2026modern}.
We now present our local polynomial setup.

Let $D_W, D_V \in \N$ and
suppose that $(W_1, V_1), \ldots, (W_n, V_n)$
are independent random variables taking values in
$\R^{D_W} \times \R^{D_V}$.
For $i \in [n]$, write $X_i \coloneqq (W_i, V_i)$,
and set $\bX_n \coloneqq (X_i: i \in [n])$.
Let $g: \R^{D_V} \times \cV \to \R$,
where either $\cV = [0, 1]^{D_V}$ or $\cV = \emptyset$.
When $\cV = \emptyset$, we interpret $g(V, v)$ to mean $g(V)$,
and read conditions stated for all $v \in \cV$ as conditions on $g(V)$,
with a slight abuse of notation.
Write $D_v \coloneqq 0$ if $\cV = \emptyset$
and $D_v \coloneqq D_V$ if $\cV = [0, 1]^{D_V}$.
Suppose that $K: \R^{D_W} \to [0, \infty)$ is Borel measurable
and supported on $[-1, 1]^{D_W}$.
Let $h > 0$ and define $K_h(w) \coloneqq K(w/h)/h^{D_W}$
for $w \in \R^{D_W}$.
For $q \in \N_0$,
write $D_q \coloneqq \binom{D_W + q}{q}$
and let $Q_h: \R^{D_W} \to \R^{D_q}$ be defined by
$Q_h(w) \coloneqq \bigl(w^\alpha / (\alpha! h^{|\alpha|}): \,
\alpha \in \N_0^{D_W}, \ 0 \leq |\alpha| \leq q\bigr)$,
with components ordered in increasing lexicographic fashion in $\alpha$,
so that $Q_h(w)^\T e_1 = 1$.
For $w \in [0, 1]^{D_W}$, define the diagonal matrix
$\hat K_h(w) \in \R^{n \times n}$ with
$\hat K_h(w)_{i i} \coloneqq K_h(W_i - w)$ for $i \in [n]$.
Let $\hat Q_h(w) \in \R^{n \times D_q}$ be given by
$\hat Q_h(w)_{i j} \coloneqq Q_h(W_i - w)^\T e_j$
for $i \in [n]$ and $j \in [D_q]$.
Define $\hat H_h(w) \coloneqq \hat Q_h(w)^\T \hat K_h(w) \hat Q_h(w) / n
\in \R^{D_q \times D_q}$.
With $\hat g(v) \coloneqq \bigl(g(V_1, v), \ldots, g(V_n, v)\bigr)^\T$
for $v \in \cV$ and
$N_{n,h}(w, v) \coloneqq \hat Q_h(w)^\T \hat K_h(w) \hat g(v) / n$,
define the degree-$q$ local polynomial estimator
\begin{equation}
  \label{eq:locpol}
  \hat \theta_{n,h}(w, v)
  \coloneqq
  e_1^\T \hat H_h(w)^{-1} N_{n,h}(w, v),
\end{equation}
whenever the inverse matrix exists,
and $\hat \theta_{n,h}(w, v) \coloneqq 0$ otherwise.
Define the estimand
$\theta(w, v) \coloneqq \E\bigl\{g(V_i, v) \mid W_i = w\bigr\}$
whenever the conditional expectation exists,
and suppose that this quantity does not depend on $i$.

The general form of \eqref{eq:locpol} covers two
common settings.
Firstly, if $\cV = \emptyset$, then taking
$g(V) \coloneqq V$ yields the local polynomial estimator
for the regression function $\theta(w)$,
where $V_i = \theta(W_i) + \varepsilon_i$
with $\E(\varepsilon_i \mid W_i) = 0$ for all $i \in [n]$.
Secondly, with $\cV = [0, 1]^{D_V}$ and
$g(V, v) \coloneqq \I\{V \leq v\}$,
with the inequality defined entrywise,
we obtain a local polynomial estimator
for the multivariate conditional cumulative
distribution function
$\theta(w, v) = \P(V_i \leq v \mid W_i = w)$.
We now make some further assumptions on the
data generating process and estimator.

\begin{assumption}[Local polynomial estimator]%
  \label{ass:locpol}
  Fix $\Clp \geq 1$.
  Suppose that the kernel satisfies
  $K(w) \leq \Clp$ for $w \in \R^{D_W}$
  and $K(w) \geq 1/\Clp$ for
  $w \in [-1/\Clp, 1/\Clp]^{D_W}$.
  Further, assume that $K$ is Lipschitz, with
  $|K(w) - K(w')| \leq \Clp \|w' - w\|_\infty$
  for all $w, w' \in \R^{D_W}$.
  For each $i \in [n]$, suppose that
  $W_i$ takes values in $[0, 1]^{D_W}$
  and has Lebesgue density function $f_{W,i}$,
  with $1/\Clp \leq f_{W,i}(w) \leq \Clp$
  for all $w \in [0, 1]^{D_W}$.
  Assume that $|\theta(w, v)| \leq \Clp$
  for all $w$ and $v$ and, for some $\zeta > 0$,
  that $\theta(w, v)$ is $\lceil\zeta\rceil-1$ times
  differentiable in $w$ for all $v$.
  Suppose that $\theta(w, v)$ satisfies the
  H{\"o}lder condition
  \begin{equation*}
    \biggl|
    \frac{\partial^{|\alpha|} \theta(w, v)}{\partial w^\alpha}
    - \frac{\partial^{|\alpha|} \theta(w', v)}{\partial w^\alpha}
    \biggr|
    \leq
    \Clp \|w - w'\|_{2}^{\zeta - |\alpha|},
  \end{equation*}
  for all $w$, $w'$ and $v$, and for
  every $\alpha \in \N_0^{D_W}$ with
  $|\alpha|= \lceil\zeta\rceil-1$.
  Assume that, for some $r \in (0, \infty]$, we have
  $\vvvert g(V_i, v) \vvvert_{\omega_r \mid W_i} \leq \Clp$
  almost surely for all $v$ and $i \in [n]$.
  Suppose that either
  \begin{inlineroman}
    \item $\cV = \emptyset$, $D_V = 1$ and $g(V) = V$ or
    \item $\cV = [0, 1]^{D_V}$ and $g(V, v) = \I\{V \leq v\}$.
  \end{inlineroman}
  Let $D_w \in \{0\} \cup [D_W]$ and define
  $\cW \coloneqq [0, 1]^{D_w} \times \prod_{s=D_w+1}^{D_W} \{w_s\}
  \subseteq [0, 1]^{D_W}$, where $w_s \in [0, 1]$.
  When $\cV = \emptyset$, we interpret $\cW \times \cV$ to mean $\cW$,
  with a slight abuse of notation.
  Finally, suppose that the bandwidth $h$ satisfies
  $n h^{D_W} \geq n^{1/\Clp}$
  and $h \leq n^{-1/\Clp}$, and that the polynomial degree
  is chosen so that $q \geq \lceil \zeta \rceil - 1$.
\end{assumption}

By permitting $D_w < D_W$, we allow for inference
which is uniformly valid only over a strict subset
of the design space.
This is useful, for example, if one wants to conduct
a statistical test uniformly over the response variable $v \in \cV$
at a fixed covariate value $w \in \cW$.

\begin{theorem}[Bootstrap approximation for local polynomial estimators]%
  \label{thm:locpol_bootstrap}
  Grant the setup of Section~\ref{sec:locpol} and
  Assumption~\ref{ass:locpol}.
  Fix $k \in \{2, 3\}$ and grant one of the following.
  \begin{enumerate}
    \item
      \label{it:thm_locpol_bootstrap_nonparametric}
      Let $(W_1^*, V_1^*), \ldots, (W_n^*, V_n^*)$
      satisfy Assumption~\ref{ass:bootstrap}.
      Write $\hat H_h^*(w)$ and $N_{n,h}^*(w, v)$
      for the bootstrap analogs of
      $\hat H_h(w)$ and $N_{n,h}(w, v)$
      respectively, computed using
      $(W_i^*, V_i^*)$ for $i \in [n]$.
      Set $r' \coloneqq \infty$ and define
      $\tilde \theta_{n,h}^*(w, v) \coloneqq e_1^\T \hat H_h^*(w)^{-1}
      N_{n,h}^*(w, v) - \hat \theta_{n,h}(w, v)$.

    \item
      \label{it:thm_locpol_bootstrap_multiplier}
      Let $\chi_1, \ldots, \chi_n$ be random variables
      satisfying Assumption~\ref{ass:multiplier}
      with order $k$, independent of $\bX_n$.
      Suppose that
      $\vvvert \chi_i \vvvert_{\omega_{r'}} \leq 2$
      for each $i \in [n]$, where $r' \in (0, \infty]$.
      Let $\hat g^\circ(v) \in \R^n$ have entries given by
      $\hat g^\circ(v)_i \coloneqq \chi_i
      \bigl(g(V_i, v) - \hat\theta_{n,h}(W_i, v)\bigr)$,
      and define $N_{n,h}^\circ(w, v) \coloneqq \hat Q_h(w)^\T
      \hat K_h(w) \hat g^\circ(v) / n$
      and $\tilde \theta_{n,h}^*(w, v)
      \coloneqq e_1^\T \hat H_h(w)^{-1} N_{n,h}^\circ(w, v)$.

  \end{enumerate}
  In either case, if the relevant inverse matrix does not exist,
  we set $\tilde \theta_{n,h}^*(w, v) \coloneqq 0$.
  For each $\eta > 0$, there is a probability space supporting
  $\bX_n'$, $\hat\theta'_{n,h}$
  and $\tilde\theta^{*\prime}_{n,h}$, where
  $\hat\theta'_{n,h} \equallaw \hat\theta_{n,h}$,
  $\hat\theta'_{n,h}$ is independent of $\bX_n'$,
  and $(\bX_n', \tilde\theta^{*\prime}_{n,h})
  \equallaw (\bX_n, \tilde\theta^*_{n,h})$, with
  \begin{align*}
    &\P\Biggl\{
      \sup_{(w, v) \in \cW \times \cV}
      \sqrt{n h^{D_W}}
      \bigl|
      \hat \theta_{n,h}'(w, v) - \theta(w, v)
      - \tilde \theta_{n,h}^{*\prime}(w, v)
      \bigr| \\
      &\qquad> \frac{C \eta \sqrt{\log 2 n} }
      {(n h^{D_W + (4-k)D_w})^{1/(10-2k+2D_w+4D_v)}}
      + C \sqrt{n h^{D_W}} h^\zeta
    \Biggr\}
    \leq
    2 e^{-n h^{D_W}/C} + \frac{1}{\eta^2},
  \end{align*}
  where $C > 0$ depends only on
  $\Clp$, $D_W$, $D_V$, $r$, $r'$ and $q$.
\end{theorem}

Theorem~\ref{thm:locpol_bootstrap} demonstrates that,
under the bandwidth restrictions in Assumption~\ref{ass:locpol},
a bootstrap process approximation for the scaled
local polynomial estimator holds whenever
\begin{inlineroman}
  \item the bias is negligible in the sense that
    $\sqrt{n h^{D_W}} h^\zeta \to 0$, and
  \item the bandwidth is not too small;
    i.e., $n h^{D_W + 2 D_w} \to \infty$ when $k = 2$
    or $n h^{D_W + D_w} \to \infty$ when $k = 3$,
    up to logarithmic factors.
\end{inlineroman}
These conditions impose undersmoothing
with respect to the usual pointwise mean squared error-optimal bandwidth
$h^* \asymp n^{-1/(2 \zeta + D_W)}$.
A valid bandwidth sequence for bootstrap coupling exists
whenever $\zeta > D_w$ for $k = 2$,
or $\zeta > D_w / 2$ for $k = 3$;
in particular, the required smoothness depends only on the
inference dimension $D_w$ and not on the full covariate dimension $D_W$.

In the regression setting,
when $D_w = D_W$ (that is, uniform inference is required over
the entire support of the covariates)
and $k = 2$, we recover a rate comparable to
the Gaussian approximation rate achieved in
\citet[Proposition~4.3]{cattaneo2025yurinskii}.
Moreover, since we rely on bootstrapping,
we do not require knowledge of any population moments.
Further, \citet{cattaneo2025yurinskii} note that a third-order
coupling improves their rate when
$\E(\varepsilon_i^3 \mid W_i) = 0$;
our Theorem~\ref{thm:locpol_bootstrap} with $k = 3$
gives an explicit rate without this additional condition.
For completeness, we present a Gaussian process approximation
version of this result
as Theorem~\ref{thm:locpol_gaussian}%
\iftoggle{journal}{ in the supplementary material
\citep{cattaneo2026highsupplement}}{}.

In practice, for statistical power reasons,
it is desirable to conduct
inference based on a Studentized process.
To this end, define $\tilde g(v) \in \R^n$ with
$\tilde g(v)_i \coloneqq g(V_i, v) - \theta(W_i, v)$ and write
$S_{n,h}(w, v) \coloneqq \hat Q_h(w)^\T \hat K_h(w) \tilde g(v) / n$.
Let $H_h(w) \coloneqq \E\bigl\{\hat H_h(w)\bigr\}$ and define
\begin{equation*}
  \Sigma_h(w, v, w', v') \coloneqq
  e_1^\T
  H_h(w)^{-1}
  \Cov\bigl\{
    S_{n,h}(w, v),
    S_{n,h}(w', v')
  \bigr\}
  H_h(w')^{-1}
  e_1
\end{equation*}
for $(w, v), (w', v') \in \cW \times \cV$,
whenever $H_h(w)$ and $H_h(w')$ are invertible.
The Studentized process is given by
$\hat\Sigma_{n,h}(w, v, w, v)^{-1/2}
\bigl\{\hat\theta_{n,h}(w, v) - \theta(w, v)\bigr\}$,
whenever $\hat\Sigma_{n,h}(w, v, w, v) > 0$,
where $\hat\Sigma_{n,h}(w, v, w', v')$ is an
estimator of $\Sigma_h(w, v, w', v')$.
We propose the following estimator.
Let $\check g(v) \in \R^n$ with
$\check g(v)_i \coloneqq g(V_i, v) - \hat\theta_{n,h}(W_i, v)$,
take $\hat D_h(v,v') \in \R^{n \times n}$ to be the diagonal matrix
with entries $\hat D_h(v, v')_{ii}
\coloneqq \check g(v)_i \check g(v')_i$,
and define
\begin{equation*}
  \hat\Sigma_{n,h}(w, v, w', v')
  \coloneqq
  \frac{1}{n^2}
  e_1^\T \hat H_h(w)^{-1}
  \hat Q_h(w)^\T \hat K_h(w) \hat D_h(v, v')
  \hat K_h(w') \hat Q_h(w')
  \hat H_h(w')^{-1} e_1,
\end{equation*}
when the relevant matrices are invertible, and
$\hat\Sigma_{n,h}(w, v, w', v') \coloneqq 0$ otherwise.

By combining Theorem~\ref{thm:locpol_bootstrap} with
an auxiliary conditional Gaussian coupling
and a Gaussian anti-concentration inequality,
we deduce a Kolmogorov--Smirnov approximation result
(see Lemmas~\ref{lem:anti-concentration_gp} and~\ref{lem:locpol_bootstrap_ks}%
  \iftoggle{journal}{ in the supplementary material
\citep{cattaneo2026highsupplement}}{}).
From this, we deduce the validity of
the following bootstrap uniform confidence bands
for $\theta(w, v)$, under some additional assumptions.

\begin{theorem}[Bootstrap confidence bands for local polynomial estimators]%
  \label{thm:locpol_bootstrap_cb}
  Grant the setup of Theorem~\ref{thm:locpol_bootstrap}
  and assume $h^\zeta \sqrt{n h^{D_W} \log 2n} \to 0$ and
  $\Var\bigl\{g(V_i, v) \mid W_i\bigr\} \geq 1 /\Clp$
  almost surely for $i \in [n]$ and $v \in \cV$.
  Suppose that
  $n h^{D_W + (4-k)D_w} / (\log 2n)^{10-2k+2D_w+4D_v}
  \allowbreak \to \infty$.
  Under either condition~\ref{it:thm_locpol_bootstrap_nonparametric}
  or condition~\ref{it:thm_locpol_bootstrap_multiplier},
  for $\alpha \in (0, 1)$ and $B \in \N$, define
  the empirical quantile
  \begin{equation*}
    \hat \xi_{\alpha,B}
    \coloneqq
    \inf
    \Biggl\{
      \xi:
      \frac{1}{B}
      \sum_{b=1}^B
      \I\Biggl\{
        \sup_{(w, v) \in \cW \times \cV}
        \,
        \Biggl|
        \frac{\tilde\theta_{n,h}^{*(b)}(w, v)}
        {\sqrt{\smash[b]{\hat \Sigma_{n,h}(w,v,w,v)}}}
        \Biggr| \leq \xi
      \Biggr\}
      \geq 1 - \alpha
    \Biggr\},
  \end{equation*}
  where $\tilde\theta_{n,h}^{*(1)},\ldots,\tilde\theta_{n,h}^{*(B)}$
  are i.i.d.\ draws of $\tilde\theta_{n,h}^{*}$, conditional on $\bX_n$.
  We interpret the ratio as zero when
  $\hat\Sigma_{n,h}(w,v,w,v) = 0$.
  Then, as $n, B \to \infty$,
  \begin{equation*}
    \P\biggl(
      \theta(w, v)
      \in
      \biggl[
        \hat\theta_{n,h}(w, v)
        \pm
        \hat\xi_{\alpha,B}
        \sqrt{\smash[b]{\hat \Sigma_{n,h}(w,v,w,v)}}
      \biggr]
      \text{ for all }
      (w, v) \in \cW \times \cV
    \biggr)
    \to 1 - \alpha.
  \end{equation*}
\end{theorem}

The lower bound on the conditional variance of the responses
ensures that the width of the confidence band
is positive everywhere with high probability.
In the case of conditional distribution function estimation,
this corresponds to assuming that $\theta(w, v)$ is bounded
away from zero and one.
Theorem~\ref{thm:locpol_gaussian_ks}%
\iftoggle{journal}{ in the supplementary material}{}
provides a corresponding
Kolmogorov--Smirnov result for Gaussian-based inference;
this approach usually requires estimating the full covariance
function $\Sigma_h(w,v,w',v')$, whereas our bootstrap version
estimates only the diagonal section $\Sigma_h(w,v,w,v)$
for Studentization.
 
\section{Acknowledgments and funding}

The authors used large language models,
including the OpenAI GPT,
Anthropic Claude and Google Gemini series,
to assist with mathematics and manuscript preparation.
This included checking proofs,
strengthening results and
broadening the literature review.
The authors carefully reviewed all suggestions incorporated into
the manuscript and accept full responsibility for its content.

Cattaneo gratefully acknowledges financial support from the National
Science Foundation via grants DMS-2210561, SES-2241575, and
SES-2342226; the National Institute of Food and Agriculture (NIFA)
through grant 2024-67023-42704; and the John Simon Guggenheim
Memorial Foundation through a 2026 Guggenheim Fellowship.
 
\appendix

\section{Proof sketches}

We provide sketch proofs of our main results,
Theorem~\ref{thm:yurinskii} and Theorem~\ref{thm:bootstrap}.

\subsection{Proof sketch for Theorem~\ref{thm:yurinskii}}
\label{sec:proof_sketch_yurinskii}

We outline the proof of Theorem~\ref{thm:yurinskii};
a full derivation is provided in
Appendix~\ref{sec:proofs_yurinskii}.
By Strassen's theorem
\citep[Section~10.3, Theorem~8]{pollard2002user},
it suffices to bound
\begin{equation}
  \label{eq:borel_approx}
  \P(S_n \in A) - \P\big( T_n \in A_p^{3\eta} \big)
  = \E\bigl(\I\{S_n \in A\} - \I\{T_n \in A_p^{3\eta}\}\bigr),
\end{equation}
uniformly over Borel sets $A \subseteq \R^d$, where
\begin{equation}
  \label{eq:enlargement}
  A_p^{\eta}
  \coloneqq \bigl\{
    x \in \R^d :
    \exists y \in A :
    \|x-y\|_p < \eta
  \bigr\}
\end{equation}
is the $\ell_p$-norm $\eta$-enlargement of $A$.
To do this, we fix $A$, $p$ and $\eta$,
and apply a two-step smoothing procedure to the
indicator functions in \eqref{eq:borel_approx} as follows.

Firstly, define the linear interpolant
$g(x) \coloneqq \bigl(1 - \|x-A_p^\eta\|_p / \eta \bigr) \lor 0$,
where $\|x - B\|_p \coloneqq \inf_{y \in B} \|x - y\|_p$,
and observe that $g = 1$ on $A_p^\eta$ and $g = 0$ outside $A_p^{2\eta}$.
Secondly, take $\sigma > 0$ and set
$f(x) \coloneqq \E\bigl\{ g(x + \sigma Z) \bigr\}$,
where $Z \sim N(0_d, I_d)$.
Then $f$ is infinitely differentiable and satisfies
\begin{equation*}
  (1 - \varepsilon) \I\{x \in A\}
  \leq f(x)
  \leq \varepsilon + (1 - \varepsilon)
  \I\big\{x \in A_p^{3\eta}\big\},
\end{equation*}
for all $x \in \R^d$, where
$\varepsilon \coloneqq \P\bigl(\|Z\|_p > \eta / \sigma\bigr)$.
It follows that \eqref{eq:borel_approx} may be bounded as
\begin{equation*}
  \P(S_n \in A)
  - \P(T_n \in A_p^{3\eta})
  \leq
  2\varepsilon + \E\big\{f(S_n) - f(T_n)\big\}.
\end{equation*}

Next, we randomize the order of the summands in $S_n$ and $T_n$
and apply Lindeberg's method;
note that we may assume without loss of generality
that $(X_1, \ldots, X_n)$ and $(Y_1, \ldots, Y_n)$ are independent.
Specifically, let $\pi$ be a uniformly random permutation on $[n]$,
independent of the data, and for $r \in [n]$, define
$W_r \coloneqq \sum_{s=1}^{r-1} X_{\pi(s)} + \sum_{s=r+1}^n Y_{\pi(s)}$.
Let $W_{r,-i}$ be a random variable independent of $(X_i, Y_i)$
with distribution given by the conditional law of $W_r$ given $\pi(r) = i$.
That is, $W_{r,-i}$ includes $r - 1$ summands from $S_n$ and
$n - r$ from $T_n$, while avoiding index $i$.
By a telescoping sum and applying Taylor's theorem to order $k$,
then bounding the remainder using properties of Hermite polynomials
and the multinomial theorem, we obtain
\begin{equation*}
  \E\bigl\{
    f(S_n) - f(T_n)
  \bigr\}
  \leq
  \frac{1}{n}
  \sum_{i=1}^n
  \sum_{r=1}^{n}
  \sum_{|\alpha| = 1}^k
  \frac{\delta_{\alpha,i}}{\alpha!}
  \E\bigl\{
    \partial_\alpha f(W_{r,-i})
  \bigr\}
  + \frac{n \bar\beta_{k,\gamma}}{\sigma^k \eta^{\gamma}}.
\end{equation*}

For the Taylor polynomial terms,
for each $r$ and $\alpha$, we may write
\begin{align*}
  \sum_{i=1}^n
  \delta_{\alpha,i}
  \,
  \E\bigl\{
    \partial_\alpha f(W_{r,-i})
  \bigr\}
  &=
  \sum_{i=1}^{n}
  \bar\delta_{\alpha}
  \,
  \E\bigl\{
    \partial_\alpha f(W_{r,-i})
  \bigr\}
  \iftoggle{journal}{\\ &\qquad}{}
  + \frac{1}{n}
  \sum_{i=1}^n
  \sum_{j=1}^{n}
  \delta_{\alpha,i}
  \,
  \E\bigl\{
    \partial_\alpha f(W_{r,-i})
    - \partial_\alpha f(W_{r,-j})
  \bigr\}.
\end{align*}
By properties of Hermite polynomials,
the first term above is bounded by
$\Delta(1/\sigma)$ after averaging over $r$,
dividing by $\alpha!$ and summing over $\alpha$.
For the second term,
note that there is no contribution
if $|\alpha| \leq m$,
since then $\delta_{\alpha,i} = \mu_{X,\alpha} - \mu_{Y,\alpha}$
does not depend on $i$.
Otherwise, if $|\alpha| \in \{m+1, \ldots, k\}$,
we replace $W_{r,-j}$ (where $i \neq j$)
with an auxiliary coupling
$W_{r,-j}' \equallaw W_{r,-j}$, constructed as follows:
if $W_{r,-i}$ contains $X_j$ as a summand
then replace $X_j$ with $X_i$; conversely,
if it contains $Y_j$ then replace $Y_j$ with $Y_i$.
We then apply Taylor's theorem again,
this time to order $k_\ell$, where $\ell = |\alpha|$,
noting that the resulting Taylor polynomial vanishes
in expectation as $\E(X_i^\nu) = \mu_{X,\nu}$ and
$\E(Y_i^\nu) = \mu_{Y,\nu}$ whenever $|\nu| \in [k_\ell] \subseteq [m]$.
After the same averaging over $r$ and weighted summation over $\alpha$,
the Taylor remainder is bounded by $\Gamma_\eta(1/\sigma)$.
To conclude, we set
$\sigma \coloneqq \eta / t$ and replace $t$ with $\phi(t)$
so that $\varepsilon \leq e^{-t^2/2}$.

\subsection{Proof sketch for Theorem~\ref{thm:bootstrap}}
\label{sec:proof_sketch_bootstrap}

We outline the proof of Theorem~\ref{thm:bootstrap};
a full derivation is provided in
Appendix~\ref{sec:proofs_bootstrap}.
The proof begins with a Borel set approximation argument
identical to that underlying the proof of Theorem~\ref{thm:yurinskii},
taking $m = 1$ and comparing the marginal law of $S_n$
with the bootstrap law of $S_n^*$ conditional on $\bX_n$.
The remaining work is to control the average moment discrepancies
under the relevant bootstrap distributions and
to construct a conditional coupling.

Under either bootstrap scheme,
for $\alpha \in \N_0^d$ with $|\alpha| \in [k]$ and $i \in [n]$, we have
\begin{equation*}
  \hat \delta_{\alpha,i}
  =
  \E\bigl\{\bigl(X_i - \E(X_i)\bigr)^\alpha\bigr\}
  - \E^*\bigl\{(Y_i^*)^\alpha\bigr\}.
\end{equation*}
Noting that $\hat \delta_{\alpha,i} = 0$ if $|\alpha| = 1$
and since $k \in \{2, 3\}$, suppose now that $|\alpha| = \ell \in \{2, k\}$.
For the nonparametric bootstrap, since $Y_i^* = X_i^* - \bar X_n$,
we deduce that
\begin{equation*}
  \E^*\bigl\{(Y_i^*)^\alpha\bigr\}
  =
  \frac{1}{n}
  \sum_{j=1}^n (X_j - \bar X_n)^\alpha.
\end{equation*}
Alternatively, for the multiplier bootstrap,
as $Y_i^* = \chi_i (X_i - \bar X_n)$ and $\chi_i$ is an
independent multiplier weight of order $k$, we have
\begin{equation*}
  \E^*\bigl\{(Y_i^*)^\alpha\bigr\}
  =
  \E\bigl(\chi_i^\ell\bigr)
  (X_i - \bar X_n)^\alpha
  =
  (X_i - \bar X_n)^\alpha.
\end{equation*}
If the data are not identically distributed, then
$\hat\delta_{\alpha,i}$ may not be small under the nonparametric bootstrap.
Similarly, even with i.i.d.\ data,
$\hat\delta_{\alpha,i}$ is not generally small
when using the multiplier bootstrap.
However, under either bootstrap scheme, the average moment
discrepancies are
\begin{equation*}
  \hat \delta_\alpha
  \coloneqq
  \frac{1}{n}
  \sum_{i=1}^n
  \hat \delta_{\alpha,i}
  =
  \frac{1}{n}
  \sum_{i=1}^n
  \Bigl(
    \E\bigl\{\bigl(X_i - \E(X_i)\bigr)^\alpha\bigr\}
    - (X_i - \bar X_n)^\alpha
  \Bigr).
\end{equation*}
Here, the averaging induces cancellation,
and a symmetrization argument yields
\begin{equation*}
  \E\Biggl\{
    \Biggl(
      \sum_{|\alpha|=\ell}
      \frac{\hat\delta_\alpha^2}{\alpha!}
    \Biggr)^{1/2}
  \Biggr\}
  \leq
  \frac{C M_2^\ell}{n^{1-\ell/4}}.
\end{equation*}
Let $\hat\Delta$ denote the bootstrap version of $\Delta$,
defined using $\hat\delta_\alpha$ in place of $\bar\delta_\alpha$.
It follows that
$\E\bigl\{\hat\Delta(s)\bigr\} \land 1 \leq C \sqrt n M_2^2 s^2$
for all $s > 0$.
The terms involving bootstrap versions of
$\bar\beta_{k,\gamma}$ and $\Gamma_\eta$
are controlled using $M_2$ and $M_p$.
To conclude the proof, we apply the aforementioned
Borel set approximation result conditional on the
data $X_1, \ldots, X_n$ and take an expectation.
Finally, a coupling is deduced via a conditional
version of Strassen's theorem; see
\citet[Theorem~B.2]{chen2020jackknife}
and \citet[Lemma~SA.1]{cattaneo2025yurinskii}.
This gives $S_n'\equallaw S_n$ with $S_n'$ independent of $\bX_n'$,
while preserving the bootstrap relationship
$(\bX_n',S_n^{*\prime})\equallaw(\bX_n,S_n^*)$.

The nonparametric bootstrap is a special case of the exchangeably weighted
bootstrap \citep{shao1995jackknife},
and our multiplier bootstrap weights are exchangeable whenever
they are identically distributed.
However, our proof does not apply directly in the general
exchangeable setting
because we require a representation of the bootstrap sum
with independent summands, conditional on the original data.

\section{Further applications}
\label{sec:consequences}

We derive bounds on the Prokhorov metric between
sums of random vectors and obtain distributional
approximations for maximum statistics.
We also briefly discuss alternatives to bootstrapping.

\subsection{Bounds on the Prokhorov metric}
\label{sec:prokhorov}

Under the setup of Section~\ref{sec:yurinskii},
we derive a bound on the Prokhorov $\ell_p$-metric between
$S_n$ and $T_n$ in terms of the moment discrepancies.
We also derive a lower bound on this metric,
thus demonstrating near optimality of our results
when arbitrarily many moments can be matched.
\begin{definition}[Prokhorov $\ell_p$-metric]%
  For $p \in [1,\infty]$ and $A \in \cB(\R^d)$,
  define $A_p^\eta$ as in \eqref{eq:enlargement}.
  The Prokhorov $\ell_p$-metric between
  two probability distributions
  $P_1$ and $P_2$ on $\cB(\R^d)$ is defined as
  \begin{equation*}
    \rho_p(P_1, P_2)
    \coloneqq
    \inf \Bigl\{
      \eta > 0 : \,
      \forall A \in \cB\bigl(\R^d\bigr), \,
      P_1(A) - P_2(A_p^{\eta}) \leq \eta
    \Bigr\}.
  \end{equation*}
\end{definition}

The Prokhorov metric is symmetric and
bounded as $0 \leq \rho_p(P_1, P_2) \leq 1$,
and it metrizes weak convergence
\citep[see][Section~2, p.~423]{gibbs2002choosing}.
If $X \sim P_1$ and $Y \sim P_2$, we write
$\rho_p(X, Y) \coloneqq \rho_p(P_1, P_2)$.
When $d=1$, we may omit the $p$ subscript as
all $\ell_p$-norms are equal.

\begin{theorem}[Bound on the Prokhorov metric]%
  \label{thm:prokhorov}
  Grant Assumption~\ref{ass:data}, let $t, \eta > 0$, and suppose
  \begin{equation*}
    2 e^{-t^2 / 2}
    + \frac{n \bar\beta_{k,\gamma} \phi(t)^k}{\eta^{k+\gamma}}
    + \Delta\biggl(\frac{\phi(t)}{\eta}\biggr)
    + \Gamma_\eta\biggl(\frac{\phi(t)}{\eta}\biggr)
    \leq 3 \eta.
  \end{equation*}
  Then $\rho_p(S_n, T_n) \leq 3\eta$.
\end{theorem}

If $\E(X_i^\alpha) = \E(Y_i^\alpha)$ for each $i \in [n]$
and $\alpha \in \N_0^d$ with $|\alpha| \in [k]$, then
Theorem~\ref{thm:prokhorov} gives
\begin{equation}
  \label{eq:prokhorov_matching}
  \rho_p(S_n, T_n)
  \leq
  3 \bigl(n \bar\beta_{k,\gamma}\bigr)^{\frac{1}{k+\gamma+1}}
  \phi
  \Biggl(
    \sqrt{
      \frac{2}{k+\gamma+1}
      \logplus
      \frac{1}
      {n \bar\beta_{k,\gamma} \psi^k}
    }
    \,
  \Biggr)^{\frac{k}{k + \gamma + 1}},
\end{equation}
whenever $\bar\beta_{k,\gamma} > 0$.
To interpret this,
consider again the setting of \eqref{eq:yurinskii_rate_general}.
For $p \in [1, \infty)$, we have
\begin{align}
  \label{eq:prokhorov_rate_p}
  \rho_p\biggl(
    \frac{S_n}{\sqrt n}, \frac{T_n}{\sqrt n}
  \biggr)
  &\lesssim
  n^{\frac{2 - k - \gamma}{2k + 2\gamma + 2}}
  d^{\frac{k/2 + k/p + \gamma/(p \land r)}{k + \gamma + 1}}
  \biggl(
    p + \frac{\log 2n}{d^{1 \land (2/p)}}
  \biggr)^{\frac{k}{2k + 2\gamma + 2}}, \\
  \label{eq:prokhorov_rate_infty}
  \rho_\infty\biggl(
    \frac{S_n}{\sqrt n}, \frac{T_n}{\sqrt n}
  \biggr)
  &\lesssim
  n^{\frac{2 - k - \gamma}{2k + 2\gamma + 2}}
  d^{\frac{k/2 + \gamma/r}{k + \gamma + 1}}
  \bigl( \log 2 n d \bigr)^{\frac{k}{2k + 2\gamma + 2}}.
\end{align}

The bounds given in \eqref{eq:prokhorov_rate_p} and
\eqref{eq:prokhorov_rate_infty}
show that the Prokhorov distance between
$S_n/\sqrt n$ and $T_n/\sqrt n$ converges to zero whenever
the dimension $d$ satisfies
$d \lesssim n^u$ with
$u < (k + \gamma - 2) / \{k + 2k/p + 2\gamma/(p \land r)\}$.
If, in addition, $d$ is fixed, then the rate of approximation is given by
$n^{(2 - k - \gamma)/(2k + 2\gamma + 2)}$,
up to logarithmic terms.
For example, with $k = 2$, we obtain a rate of
$n^{-\gamma/(6 + 2\gamma)}$ under finite moments of order $2+\gamma$.
When $\gamma = 1$, this recovers the $n^{-1/8}$ rate of
\citet[Theorem~1]{yurinskii1978error},
which cannot be improved in general for Gaussian approximation
under bounded third absolute moments;
see \citet{bentkus1985estimates} and
\citet[Equation~1.2.13]{senatov1998normal}.
\citet{senatov1985estimate}
obtained the same rate with a dimension factor of $d^{5/8}$
for i.i.d.\ summands.
The sample size exponent of $(2 - k - \gamma)/(2k + 2\gamma + 2)$
approaches $-1/2$ as $k$ increases.
We now show that when $d = 1$,
the rate $1/\sqrt n$ is a lower bound on the Prokhorov distance
between $S_n/\sqrt n$ and $T_n/\sqrt n$,
regardless of how many moments are matched,
thereby establishing near optimality of Theorem~\ref{thm:prokhorov}
when $k$ is large.

\begin{lemma}[Lower bound on the Prokhorov metric with exact
  moment matching]%
  \label{lem:lower_bound}
  Let $Z \sim N(0, 1)$. For each $k \in \N_0$,
  there exists a finitely supported distribution $P_k$ on $\R$
  such that if $X \sim P_k$ then $\E(Z^r) = \E(X^r)$
  for each $r \in [k]$. Further, if
  $X_1, \ldots, X_n \sim P_k$ are i.i.d.,
  then with $S_n \coloneqq \sum_{i=1}^{n} X_i$,
  there exists a quantity $C > 0$ depending only on $k$ such that
  \begin{equation*}
    \rho\biggl(\frac{S_n}{\sqrt n}, Z\biggr)
    \geq \frac{1}{C \sqrt n}.
  \end{equation*}
\end{lemma}

We remark that a (nonparametric or multiplier)
bootstrap version of Theorem~\ref{thm:prokhorov}
may be easily formulated, following the setting and derivation
of Theorem~\ref{thm:bootstrap} to
bound the expected Prokhorov distance between
the law of $S_n/\sqrt n$ and the conditional law of
the bootstrap statistic $S_n^* / \sqrt n$.
We omit a precise statement for brevity;
under the conditions of Theorem~\ref{thm:bootstrap},
with fixed $d,p,M_2$ and $M_p$, setting $\gamma=1$,
the expected Prokhorov distance is bounded by $n^{-1/8}$ for $k=2$
and $n^{-1/6}$ for $k=3$, up to logarithmic factors.
Consistency of the bootstrap then follows by Markov's inequality.

\subsection{Approximation of maxima}
\label{sec:ks}

We show how our coupling results may be combined with
appropriate anti-concentration inequalities to yield
Kolmogorov--Smirnov bounds for transformations
of the sums $S_n$ and $T_n$.
To begin, we consider a function $f: \R^d \to \R$
satisfying the H{\"o}lder condition
$|f(x) - f(y)| \leq L \|x - y\|_p^\zeta$,
where $L > 0$ and $\zeta \in (0, 1]$.
For any $\eta > 0$, we may write
\begin{align}
  \label{eq:holder_ks}
  &\sup_{\xi \in \R}
  \bigl|
  \P\bigl(f(S_n) \leq \xi\bigr)
  -\P\bigl(f(T_n) \leq \xi\bigr)
  \bigr|
  \iftoggle{journal}{\\ \nonumber &\qquad}{}
  \leq
  \P\bigl(\|S_n' - T_n'\|_p > \eta\bigr)
  + \sup_{\xi\in\R}
  \P\bigl(\xi < f(T_n) \leq \xi + L \eta^\zeta\bigr),
\end{align}
where $S_n' \equallaw S_n$ and $T_n' \equallaw T_n$.
The first term on the right-hand side of
\eqref{eq:holder_ks} is controlled by our coupling results,
while the second depends on anti-concentration properties of $f(T_n)$.
Example applications include norms,
order statistics of the coordinates,
and maxima of linear contrasts.

Deriving anti-concentration inequalities for general
distributions is a challenging problem, so
for concreteness, we now specialize to the case where
$f(x) = \|x\|_\infty$, taking $\zeta = 1$.
Gaussian anti-concentration results for such maximum statistics
are readily available in the literature
\citep{chernozhukov2014anti,chernozhukov2015comparison,
belloni2024anti,cattaneo2026sharp}.
We may therefore derive an explicit
bound on the Kolmogorov--Smirnov error
by employing an auxiliary conditional Gaussian coupling.
From now on, we impose a coordinatewise moment bound
on the summands for simplicity;
more general versions of the upcoming results could
be formulated without much difficulty.

\begin{theorem}[Kolmogorov--Smirnov approximation of
  the maximum statistic]%
  \label{thm:ks}
  Grant the setup of Theorem~\ref{thm:yurinskii}
  with $p = \infty$ and suppose that
  $\sigma^2 \coloneqq \min_{j \in [d]} \Var\bigl(e_j^\T S_n\bigr) > 0$.
  Take $k \in \{2, 3\}$,
  let $r \geq k + \gamma$ and suppose that
  $\vvvert e_j^\T X_i \vvvert_r
  \lor \vvvert e_j^\T Y_i \vvvert_r \leq M$
  for each $i \in [n]$ and $j \in [d]$.
  Then, for any $\eta, t > 0$,
  \begin{align*}
    &\sup_{\xi \in \R} \,
    \Bigl|
    \P\bigl(\|S_n\|_\infty \leq \xi\bigr)
    - \P\bigl(\|T_n\|_\infty \leq \xi\bigr)
    \Bigr|
    \leq
    4 e^{-t^2 / 2}
    + \Delta\biggl(\frac{\sqrt{\log 2d} + t}{\eta}\biggr) \\
    &\qquad+
    \Gamma_\eta\biggl(\frac{\sqrt{\log 2d} + t}{\eta}\biggr)
    + \frac{C \eta \sqrt{\log 2 d}}{\sigma}
    + \frac{C d^{k/2 + \gamma / r}
    \bigl(\sqrt{\log 2 d} + t\bigr)^{k + \gamma} n M^{k + \gamma}}
    {\eta^{k + \gamma}},
  \end{align*}
  where $C > 0$ is a universal constant.
\end{theorem}

For interpretation, suppose additionally that
all moments of $X_i$ and $Y_i$ match up to order $k$.
Take $M$ to be fixed and assume $\sigma^2 \gtrsim n$
(this holds, for example, if the variance of each
component of $X_i$ is uniformly bounded away from zero).
Then the Kolmogorov--Smirnov error is bounded by
$n^{(2-k-\gamma)/(2k+2\gamma+2)}
d^{k/(2k+2\gamma+2) + \gamma/(r k + r \gamma + r)}$,
up to logarithmic factors in $n$ and $d$.
Thus, taking $\gamma = 1$,
the rate of approximation in Kolmogorov--Smirnov distance is
$n^{-1/8} d^{1/4 + 1/(4 r)}$
with $k = 2$ and
$n^{-1/5} d^{3/10 + 1/(5 r)}$
with $k = 3$, up to logarithmic factors.
Theorem~\ref{thm:ks} can therefore be used for
valid statistical inference whenever the corresponding
bound converges to zero.
An approximation for the one-sided maximum statistic
$\max_{j \in [d]} \bigl(e_j^\T S_n\bigr)$
may be derived in exactly the same manner.
Such results are important in statistical settings
where the test statistic may be represented as the maximum
or maximum absolute value of several univariate statistics.

The proof of Theorem~\ref{thm:ks} employs an auxiliary coupling
to an intermediate distribution with matching moments
up to order $k \in \{2, 3\}$.
This distribution is constructed to contain a Gaussian component,
thereby supplying the required anti-concentration bound
via standard results
without sacrificing coupling accuracy \citep[see][]{zhilova2022new}.

We now show how nonparametric and multiplier
bootstrap schemes can be used to obtain
a fully feasible procedure for inference based on the
maximum statistic $\|S_n\|_\infty$.
We fix $k = 3$ and $\gamma = 1$ to fully exploit the ability
of the bootstrap to match higher-order moments, as in the discussion
following Theorem~\ref{thm:bootstrap}.

\begin{theorem}[Bootstrap-based inference using the maximum statistic]%
  \label{thm:bootstrap_ks}
  Grant the setup of Theorem~\ref{thm:bootstrap},
  including either condition \ref{it:thm_bootstrap_nonparametric}
  or condition \ref{it:thm_bootstrap_multiplier},
  taking $p = \infty$, $k = 3$ and $\gamma = 1$.
  Assume that $\vvvert e_j^\T(X_i-\mu)\vvvert_r \leq M$
  for each $i\in[n]$ and $j\in[d]$,
  where $r\in[4,\infty]$ and $M>0$.
  Suppose further that
  $\sigma^2 \coloneqq \min_{j \in [d]} \Var\bigl(e_j^\T S_n\bigr) > 0$.
  Then for any $\eta, t > 0$,
  \begin{align*}
    &\E \biggl\{
      \sup_{\xi \in \R} \,
      \Bigl|
      \P\bigl(\|S_n\|_\infty \leq \xi\bigr)
      - \P^*\bigl(\|S_n^*\|_\infty \leq \xi\bigr)
      \Bigr|
    \biggr\}
    \iftoggle{journal}{\\ &\qquad}{}
    \leq
    4 e^{-t^2 / 2}
    + \frac{C \eta \sqrt{\log 2 d}}{\sigma}
    + \frac{C \sqrt{n} d (\log 2d + t^2) M^2}{\eta^2},
  \end{align*}
  where $C > 0$ is a universal constant.
\end{theorem}

To interpret this result, suppose
that $M$ is fixed and $\sigma^2 \asymp n$.
Then the expectation in Theorem~\ref{thm:bootstrap_ks}
may be bounded by $n^{-1/6} d^{1/3}$,
up to logarithmic factors in $n$ and $d$.
As such, we recover a valid and feasible procedure for
inference whenever $d^2 / n \to 0$,
up to logarithmic factors.
An analogous result can be given for
the one-sided maximum statistic
$\max_{j \in [d]} \bigl(e_j^\T S_n\bigr)$.

\subsection{Alternatives to bootstrapping}

Throughout the paper,
we have taken advantage of the fact that
bootstrap procedures produce conditional distributions
with certain averaged conditional moments exactly equal to the
sample moments of the data;
these then offer high-probability approximations of the
true population moments.
However, our main coupling results can be applied to
any sampling scheme which produces a law with moments
suitably close to their population counterparts.

In this section, we briefly discuss an alternative
method for constructing a distribution
whose moments up to a fixed order
$k \in \N_0$ are equal to some specified values.
For example, in high dimensions one may wish to use
a regularized covariance estimator in place of the sample covariance matrix
\citep[see][for a review]{lam2020high}.
Such constructions are more broadly known as solutions to
the \emph{moment problem};
see \citet{schmudgen2017moment} for a contemporary review.
We assume that a sequence of target moments
$\mu_\alpha$ for $\alpha \in \N_0^d$ with $|\alpha| \in [k]$ is provided;
in practice, these may be based on moment estimators.
We will further take as given that there exists a
probability measure $P$ on
$\R^d$ such that if $X \sim P$, then $\E(X^\alpha) = \mu_\alpha$
for each $\alpha \in \N_0^d$ with $|\alpha| \in [k]$;
our aim is to construct an explicit instance of such a measure.
Determining the existence of these measures in general is known as the
\emph{multidimensional truncated Hamburger moment problem}
\citep[Part~IV]{schmudgen2017moment}.
As a simple example, we present Prony's method
in the univariate setting ($d = 1$), following the presentation of
\citet[Section~3.5.5]{blekherman2012semidefinite}.

\begin{example}[Prony's method]
  \label{ex:prony}
  Let $k \in \N$ and suppose $X$ is a real-valued random variable
  with $\E\bigl(|X|^{2k-1}\bigr) < \infty$.
  For $i \in [2k-1]$, define
  $\mu_i \coloneqq \E(X^i)$ and set $\mu_0 \coloneqq 1$.
  Define the Hankel matrix
  $H \in \R^{k \times k}$ by $H_{i j} = \mu_{i + j}$
  for $0 \leq i, j \leq k-1$ and assume it is positive definite.
  Set $(v_0, \ldots, v_{k-1})^\T \coloneqq
  - H^{-1}(\mu_k, \ldots, \mu_{2 k - 1})^\T$ and
  let $x_1, \ldots, x_k$ be the roots of the polynomial
  $p(x) \coloneqq x^{k} + \sum_{j=0}^{k-1} v_j x^j$.
  These can be computed to arbitrary precision
  using the method of \citet{aberth1973iteration}.
  Define $(w_1, \ldots, w_k)$ as the solution to the
  nonsingular Vandermonde system $\sum_{i=1}^k x_i^j w_i = \mu_j$
  for $0 \leq j \leq k-1$.
  Let $Y$ be a random variable with
  $\P(Y = x_i) = w_i$ for $1 \leq i \leq k$.
  Then $\E(Y^i) = \mu_i$ for $0 \leq i \leq 2k-1$.
\end{example}

The multivariate setting ($d \geq 2$) is rather more complicated.
Given a valid target moment sequence
$\mu_\alpha$ for $\alpha \in \N_0^d$ with $|\alpha| \in [k]$,
the Richter--Tchakaloff theorem
\citep[Theorem~1.24]{schmudgen2017moment} asserts
the existence of a discrete probability measure
$P$ on $\R^d$ with at most
$(d+k)! / (d! k!)$ atoms such that if $Y \sim P$ then
$\E(Y^\alpha) = \mu_\alpha$ for each
$\alpha \in \N_0^d$ with $|\alpha| \in [k]$.
However, in general, constructing such a measure is known to be
computationally hard \citep[Section~3.5.6]{blekherman2012semidefinite}.
Nonetheless, under certain \emph{flat extension} conditions,
multivariate generalizations of Prony's method can be applied
\citep[see][]{curto1996solution}.

\section{Additional empirical process results}
\label{sec:further_empirical}

We present some further results relating to
the empirical process strong approximations of
Section~\ref{sec:empirical_process}.

\subsection{Gaussian approximation}
\label{sec:gaussian_process}

While Theorem~\ref{thm:empirical_process} gives an
``empirical process-to-empirical process'' approximation,
we provide here an
``empirical process-to-Gaussian process'' result.
Instead of assuming the existence of variables
$Y_1, \ldots, Y_n$ such that the induced empirical process
$T_n$ satisfies an approximate moment matching condition
relative to $S_n$,
we instead consider a separable Gaussian process which
matches the first two finite-dimensional
moments of the given empirical process.

\begin{assumption}[Gaussian approximation]%
  \label{ass:gaussian_process}
  \phantom{a}
  \begin{enumerate}
    \item
      \label{it:ass_gaussian_process}
      Let $\bX_n \coloneqq (X_1, \ldots, X_n)$ satisfy
      Assumption~\ref{ass:empirical_process}\ref{it:ass_empirical_process}.
      Let $\bY_n \coloneqq (Y_1, \ldots, Y_n)$ be independent
      separable Gaussian processes on $\cF$ with $\E\{Y_i(f)\} = 0$ and
      $\E\bigl\{Y_i(f) Y_i(f')\bigr\} =
      \E\bigl\{f(X_i) f'(X_i)\bigr\} - P_{X,i}(f) P_{X,i}(f')$
      for each $i \in [n]$ and $f, f' \in \cF$.
      For $f \in \cF$, let
      $S_n(f) \coloneqq \sum_{i=1}^n \bigl(f(X_i) - P_{X,i}(f)\bigr)$
      and $T_n(f) \coloneqq \sum_{i=1}^n Y_i(f)$.

    \item
      Suppose that $\bX_n$ further satisfies
      Assumption~\ref{ass:empirical_process}\ref{it:ass_approximating_process}.
      For $i \in [n]$, define the $\R^d$-valued variable
      $\tilde\cF(Y_i) \coloneqq \bigl(Y_i(f): f \in \tilde \cF\bigr)$.
      For $a \in \N_0$, $b \geq 0$ and $i \in [n]$, define
      \begin{align*}
        \beta_{Y,a,b,i}
        &\coloneqq
        \frac{1}{\sqrt{a!}}
        \E \Bigl(
          \bigl\|\tilde\cF(Y_i)\bigr\|_2^a
          \bigl\|\tilde\cF(Y_i)\bigr\|_\infty^{b}
        \Bigr),
        &\bar\beta_{Y,a,b}
        &\coloneqq
        \frac{1}{n}
        \sum_{i=1}^n
        \beta_{Y,a,b,i}.
      \end{align*}
      Let $\beta_{a,b,i} \coloneqq \beta_{X,a,b,i} + \beta_{Y,a,b,i}$ and
      $\bar\beta_{a,b} \coloneqq \bar\beta_{X,a,b} + \bar\beta_{Y,a,b}$.
      For $\alpha \in \N_0^d$ with $|\alpha| \in [k]$, define
      \begin{equation*}
        \delta_{\alpha,i}
        \coloneqq
        \E\bigl\{\bigl(\tilde\cF(X_i) - P_{X,i}(\tilde\cF)\bigr)^\alpha\bigr\}
        - \E\bigl\{\tilde\cF(Y_i)^\alpha\bigr\}
      \end{equation*}
      for $i \in [n]$, and write $\bar \delta_\alpha \coloneqq
      \sum_{i=1}^n \delta_{\alpha,i} / n$.
      For $\ell \in \{m+1, \ldots, k\}$, take $k_\ell \in \{0, \ldots, m\}$
      and, for $\eta, s > 0$, define
      $\Delta(s)$ and $\Gamma_\eta(s)$
      as in Assumption~\ref{ass:data}.

  \end{enumerate}
\end{assumption}

\begin{theorem}[Gaussian approximation of empirical processes]%
  \label{thm:gaussian_process}
  Grant Assumption~\ref{ass:gaussian_process}
  and write $\phi(t) \coloneqq \phi_{\infty,d}(t)$.
  There exists $C > 0$ depending only on $r$
  such that for each $\eta, t > 0$,
  there exists a probability space supporting
  processes $S_n'$ and $T_n'$,
  where $S_n' \equallaw S_n$ and $T_n' \equallaw T_n$, with
  \begin{align*}
    &\P\biggl\{
      \sup_{f \in \cF}
      \bigl|S_n'(f) - T_n'(f)\bigr|
      > 3 \eta
      + C R(t)
    \biggr\} \\
    &\qquad\leq
    4 e^{-t^2 / 2}
    + \frac{n \bar\beta_{k,\gamma} \phi(t)^k}{\eta^{k+\gamma}}
    + \Delta\biggl(\frac{\phi(t)}{\eta}\biggr)
    + \Gamma_\eta\biggl(\frac{\phi(t)}{\eta}\biggr)
    + \frac{1}{n^\kappa}.
  \end{align*}
\end{theorem}

By definition of the Gaussian processes $Y_i$,
we have $\delta_{\alpha,i} = 0$ automatically whenever
$|\alpha| \in \{1, 2\}$.
If $\tilde\cF(X_i) - P_{X,i}(\tilde\cF)$ has a law which is
symmetric about zero for every $i$,
then also $\delta_{\alpha,i} = 0$ when $|\alpha| = 3$,
leading to refined approximations.
However, in the general case,
since the laws of separable Gaussian processes
are defined by their first two moments,
approximation rates are limited by
nonvanishing discrepancies in the third moments.
Unlike the bootstrap approach of
Theorem~\ref{thm:bootstrap_vc},
Theorem~\ref{thm:gaussian_process} does not immediately
yield a feasible procedure for sampling from the coupling process,
as the covariance function of the Gaussian process $T_n$
is typically unknown.

In the next result, we establish
a Gaussian process version of the VC-type result given as
Corollary~\ref{cor:empirical_process_vc}.

\begin{corollary}
  \label{cor:gaussian_vc}
  Grant Assumptions~\ref{ass:gaussian_process}\ref{%
  it:ass_gaussian_process}
  and Assumption~\ref{ass:vc}, and suppose that $\|F\|_{\bar P_X,2} > 0$.
  Fix $k \in \N$ with $k \geq 2$ and $\gamma \in (0,1]$.
  Assume further that for every nonempty finite $\cG \subseteq \cF$
  and each $\alpha \in \N_0^{|\cG|}$
  with $3 \leq |\alpha| \leq k$, we have
  \begin{equation*}
    \E\bigl\{
      \bigl(\cG(X_i) - P_{X,i}(\cG)\bigr)^\alpha
    \bigr\}
    = \E\bigl\{ \cG(Y_i)^\alpha \bigr\},
  \end{equation*}
  for all $i \in [n]$, where $\cG(Y_i) \coloneqq \bigl(Y_i(f): f \in \cG\bigr)$.
  Then there exists $C > 0$
  depending only on $r$, $k$, $a$ and $u$ such that
  for each $\eta > 0$,
  there is a probability space supporting
  processes $S_n'$ and $T_n'$, where $S_n' \equallaw S_n$
  and $T_n' \equallaw T_n$, with
  \begin{equation*}
    \P\biggl\{
      \sup_{f \in \cF} \bigl|S_n'(f) - T_n'(f)\bigr|
      >
      C M \eta
      n^{\frac{4 + u k}{4k + 4\gamma + 2u k}}
      \sqrt{\log 2n}
    \biggr\}
    \leq
    \frac{1}{\eta^{k + \gamma}}.
  \end{equation*}
\end{corollary}

\subsection{Bootstrap approximation}
\label{sec:bootstrap_ep_general}

We present a general version of our
bootstrap approximation result for empirical processes,
permitting both nonparametric and multiplier bootstrap schemes.
This general version is used to establish
Theorem~\ref{thm:locpol_bootstrap} in Section~\ref{sec:locpol}
of the main paper.
We also verify that some commonly used multiplier bootstrap
weights satisfy our high-level assumptions.

\begin{theorem}[Bootstrap empirical process coupling]%
  \label{thm:empirical_process_bootstrap}
  Let $\bX_n \coloneqq (X_i: i \in [n])$
  satisfy Assumption~\ref{ass:empirical_process} with $m = 1$
  and $k \geq 2$.
  Write $S_n(f) \coloneqq \sum_{i=1}^n \bigl(f(X_i) - P_{X,i}(f)\bigr)$
  and $P_n(f) \coloneqq \sum_{i=1}^{n} f(X_i) / n$,
  and grant one of the following.
  \begin{enumerate}
    \item
      \label{it:thm_ep_bootstrap_nonparametric}
      Let $X_1^*, \ldots, X_n^*$ satisfy Assumption~\ref{ass:bootstrap}.
      For $f \in \cF$,
      take $Y_i^*(f) \coloneqq f(X_i^*) - P_n(f)$ for each $i \in [n]$.
      Set $r' \coloneqq \infty$.

    \item
      \label{it:thm_ep_bootstrap_multiplier}
      Let $\chi_1, \ldots, \chi_n$ be
      independent real-valued random variables, independent of $\bX_n$,
      satisfying $\E(\chi_i) = 0$, $\E(\chi_i^2) = 1$
      and $\vvvert \chi_i \vvvert_{\omega_{r'}} \leq 2$
      for each $i \in [n]$, where $r' \in (0, \infty]$.
      Either
      \begin{inlineroman}
        \item
          set $Y_i^*(f) \coloneqq \chi_i \{f(X_i) - P_n(f)\}$
          for each $i \in [n]$; or
        \item
          assume $P_{X,i}(f) = 0$ for all $f \in \cF$
          and $i \in [n]$, and
          define $Y_i^*(f) \coloneqq \chi_i f(X_i)$.
      \end{inlineroman}
  \end{enumerate}
  Write $S_n^*(f) \coloneqq \sum_{i=1}^n Y_i^*(f)$ and
  $Y_i^*(\tilde\cF) \coloneqq \bigl(Y_i^*(f): f \in \tilde\cF\bigr)$.
  For $a \in \N_0$, $b \geq 0$ and $i \in [n]$, define
  \begin{equation*}
    \hat\beta_{a,b,i}
    \coloneqq
    \beta_{X,a,b,i}
    + \frac{1}{\sqrt{a!}}
    \E^* \Bigl(
      \bigl\|Y_i^*(\tilde\cF)\bigr\|_2^a
      \bigl\|Y_i^*(\tilde\cF)\bigr\|_\infty^b
    \Bigr),
  \end{equation*}
  and write $\hat\beta_{a,b} \coloneqq \sum_{i=1}^n \hat\beta_{a,b,i} / n$.
  For $\alpha \in \N_0^d$ with $|\alpha| \in [k]$ and $i \in [n]$, let
  \begin{equation*}
    \hat\delta_{\alpha,i}
    \coloneqq
    \E\bigl\{\bigl(\tilde\cF(X_i) - P_{X,i}(\tilde\cF)\bigr)^\alpha\bigr\}
    - \E^*\bigl\{Y_i^*(\tilde\cF)^\alpha\bigr\},
  \end{equation*}
  and write $\hat\delta_\alpha \coloneqq
  \sum_{i=1}^n \hat\delta_{\alpha,i}/n$.
  For $\ell \in \{2, \ldots, k\}$, let $k_\ell \in \{0, 1\}$ and define
  $\hat\Delta(s)$ and $\hat\Gamma_\eta(s)$ as in
  Assumption~\ref{ass:data}, with
  $\beta_{a,b,i}$, $\bar\beta_{a,b}$, $\delta_{\alpha,i}$
  and $\bar\delta_\alpha$ replaced by
  $\hat\beta_{a,b,i}$, $\hat\beta_{a,b}$, $\hat\delta_{\alpha,i}$
  and $\hat\delta_\alpha$, respectively.
  For $\kappa \geq 1$ and $t > 0$, let
  \begin{equation*}
    R^*(t)
    \coloneqq
    \biggl\{
      \sigma \sqrt n
      + \check M (\kappa \log 2n)^{1/r + 1/r'}
      \biggl(\frac{\cJ(\nu)}{\nu} + t\biggr)
    \biggr\}
    \biggl(\frac{\cJ(\nu)}{\nu} + t\biggr),
  \end{equation*}
  and write $\phi(t) \coloneqq \phi_{\infty,d}(t)$.
  Then, under condition~\ref{it:thm_ep_bootstrap_nonparametric},
  for each $\eta, t > 0$, there is a probability space supporting
  $\bX_n'$, $\bX_n''$, $\bX_n^{*\prime}$,
  $S_n'$ and $S_n^{*\prime}$, where
  $(\bX_n'', S_n') \equallaw (\bX_n, S_n)$,
  $\bX_n''$ is independent of $\bX_n'$, and
  $(\bX_n', \bX_n^{*\prime}, S_n^{*\prime}) \equallaw (\bX_n, \bX_n^*, S_n^*)$,
  with
  \begin{align*}
    &
    \P\biggl\{
      \sup_{f \in \cF}
      \bigl|S_n'(f) - S_n^{*\prime}(f)\bigr|
      > 3 \eta
      + C R^*(t)
    \biggr\} \\
    &\qquad\leq
    4 e^{-t^2 / 2}
    + \E\Biggl(
      \frac{n \hat\beta_{k,\gamma} \phi(t)^k}{\eta^{k+\gamma}}
      + \hat\Delta\biggl(\frac{\phi(t)}{\eta}\biggr)
      + \hat\Gamma_\eta\biggl(\frac{\phi(t)}{\eta}\biggr)
    \Biggr)
    + \frac{1}{n^\kappa},
  \end{align*}
  where $C > 0$ depends only on $r$ and $r'$.
  Under condition~\ref{it:thm_ep_bootstrap_multiplier},
  the same result holds after replacing
  $\bX_n^*$ with $(\chi_i: i \in [n])$
  and $\bX_n^{*\prime}$ with $(\chi_i': i \in [n])$
  in the coupling construction.
\end{theorem}

\begin{lemma}[Multiplier weights]%
  \label{lem:multiplier_weights}
  Consider the following distributions and values of $k$ and $r'$.
  \begin{enumerate}
    \item Rademacher: $k = 2$, $r' = \infty$,
      $\P(\chi = 1) = \P(\chi = -1) = 1/2$.
    \item Gaussian: $k = 2$, $r' = 2$,
      $\chi \sim N(0, 1)$.
    \item Centered exponential: $k = 2$, $r' = 1$,
      $\chi \sim \Exp(1) - 1$.
    \item Mammen: $k = 3$, $r' = \infty$, \\
      $\P\bigl(\chi = 1/2 + \sqrt{5}/2\bigr)
      = 1/2 - 1/\bigl(2\sqrt{5}\bigr)$,
      $\P\bigl(\chi = 1/2 - \sqrt{5}/2\bigr)
      = 1/2 + 1/\bigl(2\sqrt{5}\bigr)$.
    \item Centered Poisson: $k = 3$, $r' = 1$,
      $\chi \sim \Poisson(1) - 1$.
  \end{enumerate}
  Each of these satisfies
  $\E(\chi) = 0$, $\E\bigl(\chi^\ell\bigr) = 1$ for
  $\ell \in \{2, k\}$,
  $\vvvert \chi \vvvert_4 \leq 2$ and
  $\vvvert \chi \vvvert_{\omega_{r'}} \leq 2$.
\end{lemma}

\subsection{Maximal inequalities}

We present a general-purpose tail bound for the supremum of an
empirical process.
We describe the setup in Assumption~\ref{ass:general_process};
we do not need the function class to be
centered or uniformly bounded,
requiring only a sub-Weibull envelope bound involving
the Orlicz quasi-norm
$\vvvert \cdot \vvvert_{\omega_r}$ where $r \in (0, \infty]$.

\begin{assumption}[General empirical process]%
  \label{ass:general_process}
  Let $X_1, \dots, X_n$
  be independent random variables taking values in a
  measurable space $(\cX, \cH)$.
  For $i \in [n]$, write $P_{X,i}$ for the law of $X_i$,
  and define $\bar P_X \coloneqq \sum_{i=1}^n P_{X,i}/n$.
  Let $\cF$ be a class of measurable functions from $\cX$ to $\R$
  and suppose that there is a countable subclass $\cF' \subseteq \cF$ such that
  for every $f \in \cF$, there exists a sequence $(f_m)$ in $\cF'$
  with $f_m(x) \to f(x)$ for each $x \in \cX$.
  For a measurable function $f: \cX \to \R$
  and a measure $Q$ on $(\cX, \cH)$,
  define $\|f\|_{Q,2}^2 \coloneqq \int_\cX |f(x)|^2 \diffi Q(x)$.
  Let $F: \cX \to \R$ be a measurable function satisfying $F(x) \geq |f(x)|$
  for all $x \in \cX$ and $f \in \cF$.
  Suppose that $\|F\|_{\bar P_X, 2} > 0$ and that
  there exist $r \in (0, \infty]$ and $M > 0$ such that
  $\vvvert F(X_i) \vvvert_{\omega_r} \leq M$
  for each $i \in [n]$.
  For $\varepsilon \in (0, \infty]$ and $Q$ a measure on $(\cX, \cH)$,
  define the covering number
  \begin{equation*}
    N\bigl(\varepsilon, \cF, \|\cdot\|_{Q,2}\bigr)
    \coloneqq
    \min
    \Bigl\{
      N \in \N:
      \exists
      f_1, \ldots, f_N \in \cF
      \text{ such that }
      \sup_{f \in \cF}
      \min_{j \in [N]}
      \|f - f_j\|_{Q,2} \leq \varepsilon
    \Bigr\},
  \end{equation*}
  with $\min \emptyset \coloneqq \infty$.
  Let $\kappa \geq 1$ and, for $\varepsilon > 0$,
  define the uniform $L_2$-entropy integral
  \begin{equation*}
    \cJ(\varepsilon)
    \coloneqq
    \int_0^\varepsilon
    \sup_{Q \in \cQ}
    \sqrt{1 + \log N\bigl(
        \xi \|F\|_{Q,2}, \cF, \|\cdot\|_{Q,2}
    \bigr)}
    \diffi \xi,
  \end{equation*}
  where $\cQ$ is the set of all
  probability measures $Q$ on $\cX$ with $\|F\|_{Q,2} > 0$.
  Let $\sigma > 0$ satisfy
  $\sup_{f \in \cF} \|f\|_{\bar P_X,2} \leq \sigma$
  and $\nu \in (0, 1]$ satisfy
  $\nu \leq \sigma / \|F\|_{\bar P_X,2}$.
  Define $S_n(f) \coloneqq
  \sum_{i=1}^n \bigl(f(X_i) - P_{X,i}(f)\bigr)$
  for each $f \in \cF$.
\end{assumption}

\begin{theorem}[Tail bound for empirical process suprema]%
  \label{thm:empirical_tail}
  Grant Assumption~\ref{ass:general_process}.
  There exists a quantity $C > 0$
  depending only on $r$ such that for all
  $t > 0$,
  \begin{equation}
    \label{eq:empirical_tail}
    \sup_{f \in \cF}
    \bigl| S_n(f) \bigr|
    \leq
    C \biggl\{
      \sigma \sqrt n
      + M (\log 2n)^{1/r}
      \biggl(
        \frac{\cJ(\nu)}{\nu} + \kappa^{1/r} t
      \biggr)
    \biggr\}
    \biggl(\frac{\cJ(\nu)}{\nu} + t\biggr),
  \end{equation}
  with probability at least $1 - e^{-t^2/2} - 1/n^\kappa$.
  If $r = \infty$, then
  \eqref{eq:empirical_tail} holds
  with probability at least $1 - e^{-t^2/2}$.
\end{theorem}

The bound in Theorem~\ref{thm:empirical_tail}
is localized in the sense that the leading term
$\sigma \sqrt n$ depends only on the pointwise
$L_2$ bound $\sigma$, and not on the envelope function $F$.
The remaining terms depend on the envelope Orlicz bound $M$
and the uniform $L_2$-metric entropy of the class $\cF$.
We now present a uniform bound on the empirical $L_2$-norms
over $\cF$.

\begin{lemma}[Empirical norms]%
  \label{lem:empirical_norm}
  Grant Assumption~\ref{ass:general_process} and,
  for a function $f: \cX \to \R$, write
  $\|f\|_{P_n,2}^2 \coloneqq \sum_{i=1}^n f(X_i)^2 / n$.
  There exists a quantity $C > 0$
  depending only on $r$ such that for all
  $t > 0$, with probability at least
  $1 - e^{-t^2/2} - 1/n^\kappa$,
  \begin{equation*}
    \sup_{f \in \cF}
    \|f\|_{P_n,2}^2
    \leq
    2 \sigma^2
    + \frac{C M^2 (\kappa \log 2n)^{2/r}}{n}
    \biggl(
      \frac{\cJ(\nu)}{\nu} + t
    \biggr)^2.
  \end{equation*}
\end{lemma}

As the final result of this section,
we use the previous lemmas to derive a bound on the
supremum of a nonparametric or multiplier bootstrap process.

\begin{lemma}[Tail bound for bootstrap process suprema]%
  \label{lem:bootstrap_tail}
  Grant Assumption~\ref{ass:general_process} and write
  $\bX_n \coloneqq (X_1, \ldots, X_n)$.
  Write $P_n(f) \coloneqq \sum_{i=1}^{n} f(X_i) / n$,
  and suppose one of the following holds.
  \begin{enumerate}
    \item
      \label{it:bootstrap_tail_nonparametric}
      Let $X_1^*, \ldots, X_n^*$ satisfy Assumption~\ref{ass:bootstrap}.
      For $f \in \cF$,
      take $Y_i^*(f) \coloneqq f(X_i^*) - P_n(f)$ for each $i \in [n]$.
      Set $r' \coloneqq \infty$.

    \item
      \label{it:bootstrap_tail_multiplier}
      Let $\chi_1, \ldots, \chi_n$ be independent
      real-valued random variables, independent of $\bX_n$,
      satisfying $\E(\chi_i) = 0$, $\E(\chi_i^2) = 1$ and
      $\vvvert \chi_i \vvvert_{\omega_{r'}} \leq 2$
      for each $i \in [n]$, where $r' \in (0, \infty]$.
      Either
      \begin{inlineroman}
        \item
          set $Y_i^*(f) \coloneqq \chi_i \{f(X_i) - P_n(f)\}$
          for each $i \in [n]$; or
        \item
          define $Y_i^*(f) \coloneqq \chi_i f(X_i)$
          for each $i \in [n]$.
      \end{inlineroman}
  \end{enumerate}
  Let $S_n^*(f) \coloneqq \sum_{i=1}^n Y_i^*(f)$
  for each $f \in \cF$.
  There exists a quantity $C > 0$
  depending only on $r$ and $r'$ such that for all
  $t > 0$,
  \begin{equation}
    \label{eq:bootstrap_tail}
    \sup_{f \in \cF}
    |S_n^*(f)|
    \leq
    C \biggl\{
      \sigma \sqrt n
      + M (\kappa \log 2n)^{1/r + 1/r'}
      \biggl(
        \frac{\cJ(\nu)}{\nu} + t
      \biggr)
    \biggr\}
    \biggl( \frac{\cJ(\nu)}{\nu} + t \biggr),
  \end{equation}
  with probability at least
  $1 - e^{-t^2/2} - 1/n^\kappa$.
  If $r = r' = \infty$,
  \eqref{eq:bootstrap_tail} holds
  with probability at least $1 - e^{-t^2/2}$.
\end{lemma}

\subsection{Anti-concentration of Gaussian process suprema}

We give a general bound on the concentration function
of the supremum of a separable Gaussian process.
Unlike several existing results in the literature
\citep{nazarov2004maximal,chernozhukov2014anti,
  chernozhukov2015comparison, chernozhukov2017detailed,belloni2024anti,
cattaneo2026sharp},
we explicitly allow for non-centered, heteroskedastic and
infinite-dimensional Gaussian processes.

\begin{lemma}[Anti-concentration for Gaussian process suprema]%
  \label{lem:anti-concentration_gp}
  Let $X$ be a Gaussian process indexed by a set $T$.
  Write $\mu_t \coloneqq \E(X_t)$ and
  assume that $\sigma^2_t \coloneqq \Var(X_t) \geq \sigma^2 > 0$
  for all $t \in T$.
  Suppose that there exist a countable subset $T_0\subseteq T$
  and maps $t_n: T \to T_0$ for $n \in \N$, such that,
  on an event with probability one, we have
  $X_{t_n(t)} \to X_t$ for every $t\in T$.
  Assume that $\sup_{t \in T} |X_t - \mu_t| / \sigma_t < \infty$ almost surely.
  Then $E \coloneqq \E\bigl(\sup_{t \in T}
  |X_t - \mu_t| / \sigma_t\bigr)
  \in \bigl[\sqrt{2 / \pi}, \infty\bigr)$ and,
  for all $\eta > 0$, we have
  \begin{equation*}
    \sup_{\xi \in \R}
    \P\Bigl(\xi <
      \sup_{t \in T} X_t
      \leq \xi + \eta
    \Bigr)
    \leq
    \frac{\eta E}{\sigma}.
  \end{equation*}
  The same bound holds with $\sup_{t \in T} X_t$
  replaced by $\sup_{t \in T} |X_t|$.
\end{lemma}

\section{\texorpdfstring{Gaussian approximation of
  \iftoggle{journal}{\ifx\protect\relax\\\fi}{}local polynomial estimators}
{Gaussian approximation of local polynomial estimators}}
\label{sec:gaussian_approx_locpol}

We present a Gaussian approximation version of
Theorem~\ref{thm:locpol_bootstrap},
assuming exact moment matching up to order $k \in \{2, 3\}$.

\begin{theorem}[Gaussian approximation for local polynomial estimators]%
  \label{thm:locpol_gaussian}
  Grant the setup given in Section~\ref{sec:locpol}
  and Assumption~\ref{ass:locpol}.
  Take $T_{n,h}$ to be a separable zero-mean $D_q$-dimensional
  Gaussian process on $\cW \times \cV$ with
  $\Cov\bigl\{T_{n,h}(w, v), T_{n,h}(w', v')\bigr\}
  = \Cov\bigl\{S_{n,h}(w, v), S_{n,h}(w', v')\bigr\}$
  for each $(w, v), (w', v') \in \cW \times \cV$
  and let $Z_{n,h}(w, v) \coloneqq e_1^\T H_h(w)^{-1} T_{n,h}(w, v)$.
  With $k = 2$, there exists a quantity $C > 0$ depending only on
  $\Clp, D_W, D_V, r$ and $q$
  such that for each $\eta > 0$, there is a
  probability space supporting $\hat\theta'_{n,h}$
  and $Z'_{n,h}$, where
  $\hat\theta'_{n,h} \equallaw \hat\theta_{n,h}$
  and $Z'_{n,h} \equallaw Z_{n,h}$, with
  \begin{align}
    \iftoggle{journal}{\nonumber}{\label{eq:locpol_gaussian}}
    &\P\Biggl\{
      \sup_{(w, v) \in \cW \times \cV}
      \sqrt{n h^{D_W}}
      \bigl|
      \hat\theta'_{n,h}(w, v) - \theta(w, v) - Z'_{n,h}(w, v)
      \bigr| \\
      \nonumber
      &
      \iftoggle{journal}{\qquad\quad}{\qquad}
      >
      C \eta
      \Bigl\{
        \bigl(n h^{D_W}\bigr)^{1 - k}
        h^{- 2 D_w}
      \Bigr\}^{\frac{1}{2k + 2 + k D_w + 2 k D_v}}
      \!\!\sqrt{\log 2n}
      + C \sqrt{n h^{D_W}} h^\zeta
    \Biggr\}
    \iftoggle{journal}{\\ \label{eq:locpol_gaussian} &\qquad}{}
    \leq
    2 e^{-{n h^{D_W} / C}}
    + \frac{1}{\eta^{k + 1}}.
  \end{align}
  If, in addition, we have
  $\E\bigl\{ \tilde g(v)_i \tilde g(v')_i
  \tilde g(v'')_i \mid W_i\bigr\} = 0$
  almost surely for all $v, v', v'' \in \cV$
  and each $i \in [n]$,
  then \eqref{eq:locpol_gaussian} also holds with $k = 3$.
\end{theorem}

We also give a Kolmogorov--Smirnov approximation bound
for the Studentized process under this Gaussian approximation.
Conditional on $\bX_n$, let $\hat Z_{n,h}$ be a separable
zero-mean Gaussian process on $\cW \times \cV$ with covariance function
$\hat\Sigma_{n,h}(w,v,w',v')$.

\begin{theorem}[Gaussian inference for local polynomial estimators]%
  \label{thm:locpol_gaussian_ks}
  Grant the setup of Theorem~\ref{thm:locpol_gaussian}
  and assume that $h^\zeta \sqrt{n h^{D_W} \log 2n} \to 0$ and
  $\Var\bigl\{g(V_i, v) \mid W_i\bigr\} \geq 1 /\Clp$
  almost surely for each $i \in [n]$ and $v \in \cV$.
  If $k = 2$, suppose further that
  $n h^{2 D_w + D_W} / (\log 2n)^{6 + 2 D_w + 4 D_v} \to \infty$,
  while if $k = 3$, assume instead that we have
  $\E\bigl\{ \tilde g(v)_i \tilde g(v')_i
  \tilde g(v'')_i \mid W_i\bigr\} = 0$
  almost surely for all $v, v', v'' \in \cV$
  and each $i \in [n]$, and that
  $n h^{D_w + D_W} / (\log 2n)^{4 + 3 D_w / 2 + 3 D_v} \to \infty$.
  Then
  \begin{align*}
    &\E\Biggl\{
      \sup_{\xi \in \R}
      \,
      \Biggl|
      \P\Biggl(
        \sup_{(w, v) \in \cW \times \cV}
        \,
        \Biggl|
        \frac{\hat\theta_{n,h}(w,v)-\theta(w,v)}
        {\sqrt{\smash[b]{\hat \Sigma_{n,h}(w,v,w,v)}}}
        \Biggr| \leq \xi
      \Biggr) \\
      &\qquad\qquad\qquad-
      \P\Biggl(
        \sup_{(w, v) \in \cW \times \cV}
        \,
        \Biggl|
        \frac{\hat Z_{n,h}(w,v)}
        {\sqrt{\smash[b]{\hat \Sigma_{n,h}(w,v,w,v)}}}
        \Biggr| \leq \xi
        \Bigm|
        \bX_n
      \Biggr)
      \Biggr|
    \Biggr\}
    \to 0.
  \end{align*}
  We interpret both ratios as zero when
  $\hat\Sigma_{n,h}(w,v,w,v) = 0$.
\end{theorem}

Theorem~\ref{thm:locpol_gaussian_ks} may be used
to construct uniform confidence bands for $\theta(w, v)$
in the same manner as in Theorem~\ref{thm:locpol_bootstrap_cb};
we omit the details for brevity.
For conditional distribution function estimation in the $k = 3$ case,
the vanishing conditional third moments assumption
requires $\theta(W_i,v) \in \{0,1/2,1\}$ almost surely;
the variance lower bound then forces $\theta(W_i,v) = 1/2$.
As such, the third-order Gaussian approximation is primarily
of interest in the regression setting;
compare this with the bootstrap of Theorem~\ref{thm:locpol_bootstrap_cb},
which imposes no such third-order constraint on the responses.

\section{Comparison with the empirical process literature}
\label{sec:literature_empirical_process}

Strong approximations for empirical processes are often constructed
via the decomposition
\begin{equation}
  \label{eq:ep_lit_decomp}
  \sup_{f \in \cF} \bigl|S_n'(f) - T_n'(f)\bigr|
  \leq
  \max_{f \in \tilde\cF} \bigl|S_n'(f) - T_n'(f)\bigr|
  + \sup_{f \in \check\cF} \bigl|S_n'(f)\bigr|
  + \sup_{f \in \check\cF} \bigl|T_n'(f)\bigr|,
\end{equation}
where $\tilde\cF$ is a finite subclass of $\cF$ and
$\check\cF$ is the class
of differences $\bigl\{f - \pi(f): f \in \cF\bigr\}$,
as in Assumption~\ref{ass:function_class}
and Figure~\ref{fig:empirical_process}.
The last two terms in \eqref{eq:ep_lit_decomp}
are discretization errors and may be controlled by maximal
inequalities for empirical and Gaussian processes.
The first term is the coupling error for the random vector
$\bigl(S_n(f): f \in \tilde\cF\bigr)$, which has dimension
$d = |\tilde\cF|$.
The class $\tilde \cF$ is chosen to balance the
coupling and discretization errors;
typically, $d \to \infty$ as $n \to \infty$.
Following Section~\ref{sec:literature}, we classify
various strong approximation constructions according to the
coupling method.
Yurinskii-type couplings use a Prokhorov-type
bound and Strassen's theorem, with resulting rates
generally depending on entropy quantities associated with $\cF$.
Koml\'os--Major--Tusn\'ady (KMT)
couplings use recursive quantile transformations to couple
conditional binomial cell counts to Gaussian variables.
In structured classes, the resulting bounds may
depend on the ambient dimension $D$ and function regularity.
Wasserstein-based constructions obtain couplings from transport bounds.

Table~\ref{tab:empirical_process_literature}
collects selected strong approximation results
from the literature in Panels~A--C and our results in Panel~D.
We consider couplings of the entire process for fixed values of $n$,
with errors given by $\sup_{f \in \cF} |S_n'(f) - T_n'(f)| / \sqrt n$,
where $T_n'$ is a Gaussian process or another empirical process.
Unless stated otherwise, rates are in probability
for fixed function classes with bounded envelope functions.
We hold the data dimension $D$, matching order $k$
and entropy exponent $u$ fixed.
VC-type class ($u$) refers to Assumption~\ref{ass:vc};
elsewhere, $u$ denotes an entropy exponent whose
definition varies between rows.
The displayed rates omit logarithmic and arbitrarily
slowly diverging factors.
For our results in Panel~D, we take $\gamma = 1$.
Corollaries~\ref{cor:gaussian_vc} and~\ref{cor:empirical_process_vc}
require the mixed moments of every centered evaluation vector to match
through order $k$ for each observation.
Matching a Gaussian target through order $k \geq 3$ requires all mixed
third moments to vanish, and for $k \geq 4$ all fourth-order joint
cumulants must also vanish.
Theorem~\ref{thm:bootstrap_vc} uses the bootstrap
for automatic approximate matching.
Appendix~\ref{sec:literature_details_empirical_process}
gives precise versions of the results in
Table~\ref{tab:empirical_process_literature},
as well as details for all rates presented in this section.

\subsection{Yurinskii-type couplings}

This approach goes back to \citet{dudley1983invariance},
who obtained relatively slow Gaussian approximation rates
with i.i.d.\ data for VC classes of sets using an extension
of Yurinskii's method \citep{yurinskii1978error}
due to \citet{dehling1983limit}.
With the same approach, \citet[Theorem~5.10]{massart1986rates}
obtained a rate of $n^{-1/(8+5u)}$ for function classes with
uniform $L_2$ entropy exponent $u$ and $L_3$ envelope function.
He further obtained slower rates under the assumption that
the envelope possesses only
$2 + \delta$ finite moments, where $\delta \in (0, 1)$.
\citet[Theorem~2 and Section~6]{le1988prokhorov} also achieved
a rate of $n^{-1/(8+5u)}$ for VC classes of sets with i.n.i.d.\ data.
\citet[Section~6]{le1988prokhorov} showed how
replacing Yurinskii's coupling with a Gaussian approximation theorem of
\citet{zaitsev1987estimates} may improve these rates
by providing a better dependence on the sample size when assuming finite
exponential moments (see Section~\ref{sec:literature_yurinskii}).
For axis-aligned boxes and i.i.d.\ data with a
bounded density on $[0,1]^D$,
\citet[Theorem~1]{rio1996vitesses} obtained the rate $n^{-1/12}$
using Zaitsev's result; here, the sample size exponent
does not depend on $D$.
\citet[Fact~2.2 and Proposition~2.2]{einmahl1997gaussian} combined
Zaitsev's coupling with a discretization error bound based on
Talagrand's inequality.
Building on this combination,
\citet[Proposition~1]{berthet2006revisiting} obtained a rate of
$n^{-1/(2+5u)}$ for i.i.d.\ data.
\citet[Theorem~3.4]{settati2009gaussian} obtained the same rate of
$n^{-1/(2+5u)}$
under population and empirical entropy conditions,
with $u$ denoting the population $L_2$ covering exponent.
The discussion following
Proposition~3.1 of \citet{cattaneo2025yurinskii} gives the rate
$n^{-1/(6+2u)}$ for i.i.d.\ data with a bounded envelope function,
improving on that of \citet{berthet2006revisiting} when $u>4/3$;
Appendix~\ref{sec:literature_details_empirical_process}
explains why the same rate holds for i.n.i.d.\ data,
and \citet{cattaneo2025yurinskii} further allow for martingale dependence.
The same discussion describes a
third-order refinement under vanishing conditional mixed third moments,
but omits the corresponding empirical process rate.

Our Yurinskii-type empirical process coupling
in Theorem~\ref{thm:empirical_process}
allows approximate and averaged moment matching,
similar to Theorem~\ref{thm:yurinskii}.
It permits i.n.i.d.\ data points taking values in a Polish space
and assumes a sub-Weibull envelope function;
cf.\ \citet{massart1986rates}, who require only finite
envelope moments of order $2+\delta$.
Under exact individual-level moment matching up to order $k \geq 2$,
Corollary~\ref{cor:empirical_process_vc}
gives an explicit rate of $n^{-(k-1)/(2k+2+u k)}$ when $\gamma = 1$.
This equals $n^{-1/(6+2u)}$ when $k=2$ and improves on the rate of
\citet{berthet2006revisiting} whenever $u > 4/(4k-5)$;
the sample size exponent approaches $1/(2+u)$ as $k$ diverges.
We observe that for the class of lower orthants,
exact moment matching up to order three forces
the centered empirical process laws to be identical;
in this case, a coupling exists with zero error
(see Appendix~\ref{sec:literature_details_empirical_process}).
Theorem~\ref{thm:gaussian_process}
and Corollary~\ref{cor:gaussian_vc}
in Appendix~\ref{sec:gaussian_process}
provide Gaussian process versions of Theorem~\ref{thm:empirical_process}
and Corollary~\ref{cor:empirical_process_vc}, respectively.
With $k = 3$, the third centered moments must vanish,
and the resulting Gaussian approximation rate of $n^{-2/(8+3u)}$
then improves on that of \citet{berthet2006revisiting}
when $u > 4/7$ and on the second-order rate of
\citet{cattaneo2025yurinskii} for every $u \geq 0$.
All these rates improve on the $n^{-1/(8+5u)}$ rate
attained by \citet{massart1986rates}.

\begin{table}[p]
  \centering
  \rotatebox{90}{%
    \begin{minipage}{\textheight}
      \centering
      \caption{Strong approximation of empirical processes.}
      \label{tab:empirical_process_literature}
      \footnotesize
      \setlength{\tabcolsep}{2.9pt}%
      \renewcommand{\arraystretch}{1.09}%
      \vspace*{-1mm}
      \begin{tabular*}{\linewidth}{%
        @{\extracolsep{\fill}}lccp{0.40\linewidth}cc@{}}
        \hline
        & Sampling & $k$ & Class and assumptions
        & Target & Rate \\
        \hline
        \noalign{\vskip 0.2mm}
        \multicolumn{6}{@{}l}{\textit{Panel A:
        Yurinskii-type couplings}} \\[-0.3mm]
        \quad\citet{massart1986rates} & i.i.d. & $2$
        & VC-type class ($u$); arbitrary law
        & Gaussian & $n^{-1/(8+5u)}$ \\
        \quad\citet{le1988prokhorov} & i.n.i.d. & $2$
        & VC sets ($u$); arbitrary law
        & Gaussian & $n^{-1/(8+5u)}$ \\
        \quad\citet{rio1996vitesses} & i.i.d. & $2$
        & Axis-aligned boxes; bounded density on $[0,1]^D$
        & Gaussian & $n^{-1/12}$ \\
        \quad\citet{berthet2006revisiting} & i.i.d. & $2$
        & VC-type class ($u$); arbitrary law
        & Gaussian & $n^{-1/(2+5u)}$ \\
        \quad\citet{settati2009gaussian} & i.i.d. & $2$
        & VC-type random entropy ($u$); arbitrary law
        & Gaussian & $n^{-1/(2+5u)}$ \\
        \quad\citet{cattaneo2025yurinskii} & i.n.i.d. & $2$
        & VC-type class ($u$); arbitrary law
        & Gaussian & $n^{-1/(6+2u)}$ \\
        \hline
        \noalign{\vskip 0.2mm}
        \multicolumn{6}{@{}l}{\textit{Panel B: KMT-type couplings}} \\[-0.3mm]
        \quad\citet{komlos1975approximation} & i.i.d. & $2$
        & Half-lines; arbitrary law on $\R$
        & Gaussian & $n^{-1/2}$ \\
        \quad\citet{borisov1982approximation} & i.i.d. & $2$
        & Lower orthants; arbitrary law on $\R^D$
        & Gaussian & $n^{-1/(4D-2)}$ \\
        \quad\citet{massart1989strong} & i.i.d. & $2$
        & VC sets, regular boundaries;
        uniform on $[0,1]^D$, $D\geq2$
        & Gaussian & $n^{-1/(2D)}$ \\
        \quad\citet{koltchinskii1994komlos} & i.i.d. & $2$
        & VC sets, regular boundaries;
        uniform on $[0,1]^D$, $D\geq2$
        & Gaussian & $n^{-1/(2D)}$ \\
        \quad\citet{rio1994local} & i.i.d. & $2$
        & VC-type ($u$), bounded variation;
        pos.\ $C^1$ density on $[0,1]^D$, $D\geq2$
        & Gaussian & $n^{-1/(2D)}$ \\
        \quad\citet{castelle1998strong}
        & i.i.d. & $2$
        & Lower rectangles; uniform on $[0,1]^2$
        & Gaussian & $n^{-1/2}$ \\
        \quad\citet{gine2004kernel} & i.i.d. & $2$
        & Bounded variation functions; density on $\R$
        & Gaussian & $n^{-1/2}$ \\
        \quad\citet{cattaneo2024uniform} & i.i.d. & $2$
        & Bounded variation functions; arbitrary law on $\R$
        & Gaussian & $n^{-1/2}$ \\
        \quad\citet{cattaneo2025strong} & i.i.d. & $2$
        & VC-type Lipschitz functions;
        transformation assumptions on $\R^D$
        & Gaussian & $n^{-1/(D\lor2)}$ \\
        \hline
        \noalign{\vskip 0.2mm}
        \multicolumn{6}{@{}l}{\textit{Panel C: Wasserstein-based
        couplings}} \\[-0.3mm]
        \quad\citet{dedecker2014strong} & i.i.d. & $2$
        & Lower orthants; arbitrary law on $\R^D$
        & Gaussian & $n^{-1/6}$ \\
        \hline
        \noalign{\vskip 0.2mm}
        \multicolumn{6}{@{}l}{\textit{Panel D: Results of this
        paper}} \\[-0.3mm]
        \quad Corollary~\ref{cor:gaussian_vc} & i.n.i.d. & $\geq2$
        & VC-type class ($u$), exact moment matching;
        arbitrary law
        & Gaussian & $n^{-(k-1)/(2k+2+u k)}$ \\
        \quad Corollary~\ref{cor:empirical_process_vc} & i.n.i.d. & $\geq2$
        & VC-type class ($u$), exact moment matching;
        arbitrary law
        & Empirical & $n^{-(k-1)/(2k+2+u k)}$ \\
        \quad Theorem~\ref{thm:bootstrap_vc} & i.n.i.d. & $3$
        & VC-type class ($u$), common means; arbitrary law
        & Bootstrap & $n^{-1/(4+2u)}$ \\
        \hline
      \end{tabular*}

      \vspace{2pt}
      \begin{minipage}{\linewidth}
        \scriptsize\raggedright
        \emph{Note:} Rates omit logarithmic and
        arbitrarily slowly diverging factors. Precise assumptions
        and results are given in
        Appendices~\ref{sec:literature_empirical_process} and
        \ref{sec:literature_details_empirical_process}.
      \end{minipage}
    \end{minipage}%
  }
\end{table}
\subsection{KMT-type couplings}

For i.i.d.\ data uniform on $[0, 1]$,
\citet[Theorem~3]{komlos1975approximation} constructed a Brownian bridge
approximation of the normalized empirical distribution function
with the optimal rate $n^{-1/2} \log n$.
For $D=1$ and i.i.d.\ data with a Lebesgue density, the same rate holds
for classes of functions with uniformly bounded pointwise total variation;
see \citet[Section~2]{gine2004kernel} and
\citet[proof of Proposition~5]{gine2010confidence}.
\citet[Lemma~SA26]{cattaneo2024uniform} give a self-contained proof
that allows for arbitrary distributions.
\citet[Corollary~1 and Remark~1]{cattaneo2025strong} recover this rate
using a weaker notion of total variation, under additional entropy,
local variation, and normalizing transformation conditions.

For the normalized empirical distribution function of i.i.d.\ uniform data on
$[0,1]^2$, \citet{tusnady1977remark} stated a coupling rate of
$n^{-1/2}(\log n)^2$ with an exponential tail bound,
omitting the details of the proof;
\citet[Theorem~2.3]{castelle1998strong} established a
localized refinement that recovers this rate.
For lower orthants and i.i.d.\ data with an arbitrary distribution on $\R^D$,
\citet[Theorem~3]{borisov1982approximation} obtained a rate of
$n^{-1/(4D-2)}$ by a construction based on the univariate KMT approximation.
For i.i.d.\ data uniform on $[0, 1]^D$ with $D \geq 2$,
\citet[Theorem~1]{massart1989strong} achieved a rate of $n^{-1/(2D)}$
for VC classes of sets with regular boundaries.
\citet[Theorem~1]{beck1985lower} showed that
no Brownian bridge approximation may attain a faster rate
for the class of Euclidean balls; for orthants,
his Theorem~2 gave only the lower bound $n^{-1/2}(\log n)^{(D-1)/2}$.
\citet{koltchinskii1994komlos} extended these constructions to
function classes with rates determined by
Haar approximations and empirical entropy conditions.
For i.i.d.\ data admitting an appropriately regular density
on $[0,1]^D$, with $D\geq2$,
\citet[Theorem~1.1]{rio1994local} provided a rate of $n^{-1/(2D)}$ for
uniformly bounded VC-type function classes with uniformly bounded
total variation.
For VC-type classes of uniformly Lipschitz functions,
\citet{cattaneo2025strong} improved this rate to
$n^{-1/2}$ when $D = 2$ and to $n^{-1/D}$ when $D \geq 3$,
under regularity conditions on a normalizing transformation;
these hold for a continuously differentiable, strictly positive
Lebesgue density on the unit cube.

To compare the Gaussian coupling in Corollary~\ref{cor:gaussian_vc}
with the KMT-type coupling of \citet[Theorem~1.1]{rio1994local},
suppose that the data are i.i.d.\ uniform on $[0,1]^D$, with $D \geq 2$.
Consider a pointwise measurable class
$\cF = \{f_\theta: \theta \in \Theta\}$ with a bounded envelope $F$
and uniformly bounded variation in the data argument.
If $\Theta$ is a bounded subset of $\R^p$ and
$|f_\theta - f_{\theta'}| \leq \|\theta - \theta'\|_2 F$,
we may take $u = p$
\citep[Theorem~2.7.11]{van1996weak}.
Ignoring logarithmic factors, our rate with $\gamma = 1$ and $k = 2$
is $n^{-1/(6+2p)}$, which is faster than Rio's $n^{-1/(2D)}$
when $D > 3 + p$.
If all mixed third moments of the centered evaluations vanish,
then $k = 3$ gives $n^{-2/(8+3p)}$, which is faster
when $D > 2 + 3p/4$.
Alternatively, $u$ may be taken to be twice
the VC dimension of certain classes:
$u=2D$ for lower orthants, $u=2D+2$ for half-spaces and Euclidean balls,
and $u=4D$ for rectangles \citep[Theorem~2.6.4]{van1996weak}.
Under the corresponding distributional assumptions,
the KMT rate $n^{-1/(2D)}$ for these classes
is faster than our general bound (with $k = 2$) in every dimension.

\subsection{Wasserstein-based couplings}

For lower orthants and i.i.d.\ data on $\R^D$ with an arbitrary
distribution, \citet[Theorem~3.2]{dedecker2014strong} obtained
a rate of $n^{-1/6}$ using the Kantorovich--Rubinstein theorem
and the Lindeberg method;
their results also allow for strictly stationary $\beta$-mixing data.
For $D=2$, this rate has the same sample size exponent as the bound of
\citet{borisov1982approximation}, and for $D > 2$ it provides a faster rate.
\citet[Theorem~2.1]{dedecker2013strong} used the same method for
univariate strictly stationary sequences with dependence coefficients
of order $n^{-1-\delta}$ where $\delta > 0$,
obtaining a rate of $n^{-\delta^2/\{22(\delta+2)^2\}}$.

\subsection{Bootstrap-based couplings}

Existing bootstrap couplings often exploit specific structure
in the function class or data distribution.
For i.i.d.\ data, \citet{praestgaard1993exchangeably} established
conditional weak convergence of exchangeably weighted bootstrap
empirical processes indexed by Donsker classes.
\citet[Theorem~3.2]{neumann1998strong} showed how to
couple a centered kernel density estimator computed from i.i.d.\ or
geometrically $\beta$-mixing data on $\R^D$
with its empirical bootstrap counterpart
\citep[see also][Remark~3.1ii]{chernozhukov2014gaussian};
after normalization by $\sqrt{n h^D}$, where $h$ is the bandwidth,
his approach yields a rate of order
$(n h^D)^{-1/(D+4)}$ for sufficiently small $h$,
up to logarithmic factors (see Section~\ref{sec:locpol}
and Appendix~\ref{sec:literature_details_empirical_process}).
For univariate i.i.d.\ data, \citet[Theorem~1.2]{csorgo2000approximation}
approximated the bootstrap empirical process based on $n$ draws using
Brownian bridges in an unconditional coupling at the rate $n^{-1/2}\log n$.

Our Theorem~\ref{thm:bootstrap_vc}
gives a feasible approximation for empirical processes indexed by
general VC-type classes by coupling the process to a
nonparametric or multiplier bootstrap process.
This result is derived using the more general
Theorem~\ref{thm:empirical_process_bootstrap}
in Appendix~\ref{sec:bootstrap_ep_general},
which in turn relies on Theorem~\ref{thm:empirical_process}.
Under the common mean assumption that $P_{X,i}(f) = P_X(f)$
for every $i \in [n]$ and
$f\in\cF$, the nonparametric and third-order multiplier bootstraps attain
a rate of $n^{-1/(4+2u)}$, without requiring vanishing third moments.
This improves on the second-order rate
displayed by \citet{cattaneo2025yurinskii} for every $u \geq 0$
and on that of \citet{berthet2006revisiting} when $u > 2/3$.
The Gaussian approximation rate of $n^{-2/(8+3u)}$
attained in Corollary~\ref{cor:gaussian_vc}
in Appendix~\ref{sec:gaussian_process}
remains faster for every $u>0$, but requires vanishing third moments.
Further, the bootstrap process may be simulated,
conditional on the data, without knowledge of the population
empirical process moments.
In contrast, Gaussian process approximations
depend on the population covariance function, which is generally unknown
in applications.
For lower orthants, the common mean condition requires all observations
to have the same distribution function, implying i.i.d.\ sampling;
for other function classes, it generally permits
heterogeneous distributions.
Our bootstrap result does not require
any density or transformation conditions such as those imposed by
\citet{rio1994local} or \citet{cattaneo2025strong}.

\subsection{Other related approximations}

Another approach to approximation of stochastic processes
uses finite-dimensional projections of Hilbert space-valued random variables.
\citet{gotze2011estimates} constructed strong Gaussian approximations for
sums of i.i.d.\ Hilbert space-valued random variables with finite moments of
order exceeding two, applying vector coupling results to leading
covariance eigenspaces; the resulting rates depend on the
covariance spectrum and moment bounds.
\citet[Theorem~2.2]{fang2026sharp} also gave a
Wasserstein $W_{4,2}$ bound for sums of independent centered random
variables in a separable Hilbert space with finite fourth moments;
their bound depends on the trace of the square root of a fourth moment operator.
These results provide a uniform approximation to the stochastic process
whenever the supremum norm may be controlled
in terms of the Hilbert norm.

Other results approximate only the distribution of the supremum of the process.
Under an anti-concentration condition on the Gaussian supremum,
\citet[Theorem~3]{norvaisa1991rate} approximated the law of
$\sup_{f\in\cF} |S_n(f)|/\sqrt n$ by that of a Gaussian supremum
at a rate of $n^{-1/6}$, up to logarithmic factors,
for bounded VC subgraph classes.
\citet{chernozhukov2014gaussian} showed how to couple
$\sup_{f \in \cF} S_n(f) / \sqrt n$ to the supremum of a Gaussian process,
permitting $n$-varying function classes,
and \citet{chernozhukov2016empirical} established related results
for nonparametric and multiplier bootstraps.
These results are well suited to inference based on the supremum,
but they do not couple the two processes;
see Appendix~\ref{sec:ks} for our results on maximum statistics.
Section~\ref{sec:literature_bootstrap} discusses further distributional
and bootstrap approximations for high-dimensional vectors and their maxima.
Finally, \citet{dudley1983invariance}, \citet{komlos1975approximation},
\citet{massart1986rates}, \citet{berthet2006revisiting}
and \citet{dedecker2014strong}
also give sequential approximations by Kiefer processes,
thus controlling the partial sums for all $n$ simultaneously.

\section{Proofs and technical results}

We provide complete proofs and supporting lemmas for all of the
results in the paper,
beginning by defining some additional notation.
For a matrix $A \in \R^{d \times d}$,
the entrywise (Frobenius-type) $\ell_p$-norm is
$\|A\|_{\Fr,p} \coloneqq
\bigl(\sum_{i=1}^d \sum_{j=1}^d |A_{i j}|^p\bigr)^{1/p}$
for $p \in [1, \infty)$.
The kernel (null space) of $A$
is written $\ker A$ and the trace of $A$ is $\tr A$.

\subsection{Proofs for Section~\ref{sec:yurinskii}}
\label{sec:proofs_yurinskii}

We begin with a short series of lemmas required
to establish Theorem~\ref{thm:yurinskii}.

\begin{lemma}[Borel set smoothing]%
  \label{lem:smoothing}
  Let $A \in \cB(\R^d)$,
  $p \in [1,\infty]$ and $\eta > 0$.
  Define $A_p^{\eta}$ as in \eqref{eq:enlargement}
  and, for $x \in \R^d$ and $B \subseteq \R^d$,
  write $\|x - B\|_p \coloneqq \inf_{y \in B} \|x - y\|_p$
  with $\inf\emptyset \coloneqq \infty$.
  For $\sigma > 0$ and with $Z \sim N(0_d, I_d)$,
  define the functions
  \begin{align*}
    g(x)
    \coloneqq
    g_{A,p,\eta}(x)
    &\coloneqq
    \biggl( 1 - \frac{\|x-A_p^\eta\|_p}{\eta} \biggr) \lor 0,
    &f(x)
    \coloneqq
    f_{A,p,\eta,\sigma}(x)
    &\coloneqq
    \E\bigl\{ g(x + \sigma Z) \bigr\}.
  \end{align*}
  Writing $\varphi(x)$ for the Lebesgue density of
  $N(0, 1)$, define
  $H_r(x) \coloneqq (-1)^r \partial_r \varphi(x) / \varphi(x)$,
  the probabilist's Hermite polynomial of degree $r \in \N_0$.
  Then $f$ is infinitely differentiable and,
  for $\alpha \in \N_0^d$,
  \begin{equation*}
    \partial_\alpha f(x)
    =
    \frac{1}{\sigma^{|\alpha|}}
    \E\bigl\{
      g(x + \sigma Z)
      H_\alpha(Z)
    \bigr\},
  \end{equation*}
  where $H_\alpha(z) \coloneqq \prod_{j=1}^d H_{\alpha_j}(z_j)$
  for $z \in \R^d$.
  Further,
  $|\partial_\alpha f(x)| \leq \sigma^{-|\alpha|} \sqrt{\alpha!}$
  for all $x \in \R^d$.
  Finally, with
  $\varepsilon \coloneqq \varepsilon_{d,p,\eta,\sigma}
  \coloneqq \P\bigl(\|Z\|_p > \eta / \sigma\bigr)$,
  we have for all $x \in \R^d$ that
  \begin{equation*}
    (1 - \varepsilon) \I\{x \in A\}
    \leq f(x)
    \leq \varepsilon + (1 - \varepsilon)
    \I\big\{x \in A_p^{3\eta}\big\}.
  \end{equation*}
\end{lemma}

\begin{proof}[Proof of Lemma~\ref{lem:smoothing}]
  Take $x \in \R^d$. Changing variables
  by taking $u \coloneqq x + \sigma v$, we have
  \begin{equation*}
    f(x)
    =
    \int_{\R^d}
    g(x + \sigma v)
    \prod_{j=1}^d \varphi(v_j) \diffi v
    =
    \frac{1}{\sigma^d}
    \int_{\R^d}
    g(u)
    \prod_{j=1}^d \varphi\biggl(\frac{u_j - x_j}{\sigma}\biggr)
    \diffi u.
  \end{equation*}
  As $\partial_r \varphi(x) = (-1)^r H_r(x) \varphi(x)$
  for each $r \in \N_0$ and by the dominated convergence theorem,
  \begin{align*}
    \partial_\alpha f(x)
    &=
    \frac{1}{\sigma^{d}}
    \int_{\R^d}
    g(u)
    \prod_{j=1}^d
    \biggl(- \frac{1}{\sigma}\biggr)^{\alpha_j}
    \partial_{\alpha_j}
    \varphi\biggl(\frac{u_j - x_j}{\sigma}\biggr)
    \diffi u \\
    &=
    \frac{1}{\sigma^{d + |\alpha|}}
    \int_{\R^d}
    g(u)
    \prod_{j=1}^d
    H_{\alpha_j}
    \biggl(\frac{u_j - x_j}{\sigma}\biggr)
    \varphi\biggl(\frac{u_j - x_j}{\sigma}\biggr)
    \diffi u \\
    &=
    \frac{1}{\sigma^{|\alpha|}}
    \! \int_{\R^d}
    g(x + \sigma v)
    \prod_{j=1}^d
    H_{\alpha_j}(v_j)
    \varphi(v_j)
    \diffi v \\
    &=
    \frac{1}{\sigma^{|\alpha|}}
    \E\bigl\{
      g(x + \sigma Z)
      H_{\alpha}(Z)
    \bigr\}.
  \end{align*}
  By orthogonality of Hermite polynomials and integration by parts
  \citep[see][Exercise~4.1]{samworth2026modern}, we have
  for $\alpha, \alpha' \in \N_0^d$ that
  $\E\bigl\{H_\alpha(Z) H_{\alpha'}(Z)\bigr\}
  = \alpha! \, \I\{\alpha = \alpha'\}$.
  By Jensen's inequality, it follows that
  $|\partial_\alpha f(x)| \leq
  \sigma^{-|\alpha|} \bigl\{\E\bigl(H_\alpha(Z)^2\bigr)\bigr\}^{1/2}
  = \sigma^{-|\alpha|} \sqrt{\alpha!}$.
  For the final claim, as
  $\I\big\{x \in A_p^{\eta}\big\}
  \leq g(x) \leq \I\big\{x \in A_p^{2\eta}\big\}$,
  we have
  \begin{align*}
    f(x)
    &\geq \P\bigl(x + \sigma Z \in A_p^{\eta}\bigr)
    \geq \P\bigl(x \in A, \|\sigma Z\|_p < \eta\bigr)
    = (1 - \varepsilon) \I\{x \in A\}, \\
    f(x)
    &\leq
    \P\bigl(x + \sigma Z \in A_p^{2\eta}\bigr)
    \iftoggle{journal}{\\ &}{}
    \leq
    \P\bigl(\|\sigma Z\|_p > \eta\bigr)
    + \I\bigl\{x \in A_p^{3\eta}\bigr\}
    \P\bigl(\|\sigma Z\|_p \leq \eta\bigr)
    \leq
    \varepsilon + (1-\varepsilon)\I\big\{x \in A_p^{3\eta}\big\},
  \end{align*}
  as required.
\end{proof}

\begin{lemma}[Higher-order Lindeberg interpolation]
  \label{lem:lindeberg}
  Grant Assumption~\ref{ass:data}.
  Let $A \in \cB(\R^d)$ and $\eta, \sigma > 0$.
  With $f \coloneqq f_{A,p,\eta,\sigma}$
  the function defined in Lemma~\ref{lem:smoothing},
  we have
  \begin{equation*}
    \E\bigl\{
      f(S_n) - f(T_n)
    \bigr\}
    \leq
    \frac{n \bar\beta_{k,\gamma}}{\sigma^k \eta^\gamma}
    + \Delta\biggl(\frac{1}{\sigma}\biggr)
    + \Gamma_\eta\biggl(\frac{1}{\sigma}\biggr).
  \end{equation*}
\end{lemma}

\begin{proof}[Proof of Lemma~\ref{lem:lindeberg}]
  Since the result does not depend on the copula between
  $\bX_n \coloneqq (X_1, \ldots, X_n)$ and
  $\bY_n \coloneqq (Y_1, \ldots, Y_n)$,
  we assume that they are independent,
  enlarging the probability space if necessary.
  Let $\pi$ be a uniformly random permutation on $[n]$,
  independent of $(\bX_n, \bY_n)$.
  For $r \in [n]$, define
  $W_r \coloneqq \sum_{s=1}^{r-1} X_{\pi(s)} + \sum_{s=r+1}^n Y_{\pi(s)}$.
  For $r,i \in [n]$, let $\pi_{r,i}$ be uniform over permutations
  with $\pi_{r,i}(r)=i$, independently of $(\bX_n,\bY_n)$,
  and set $W_{r,-i}\coloneqq\sum_{s<r}X_{\pi_{r,i}(s)}
  +\sum_{s>r}Y_{\pi_{r,i}(s)}$.
  Thus $W_{r,-i}$ follows the conditional law of $W_r$
  given $\pi(r)=i$ and is independent of $(X_i,Y_i)$.
  By a telescoping sum and the argument of
  \citet[Theorem~18.18]{aliprantis2006infinite},
  applied to the random functions $t \mapsto f(W_{r,-i} + t X_i)$
  and $t \mapsto f(W_{r,-i} + t Y_i)$, for each $r, i \in [n]$,
  there exist random variables
  $\tilde X_{r,i}$ on the line segment between $0_d$ and $X_i$, and
  $\tilde Y_{r,i}$ on the line segment between $0_d$ and $Y_i$, such that
  \begin{align}
    \nonumber
    \E\bigl\{
      f(S_n) - f(T_n)
    \bigr\}
    &=
    \sum_{r=1}^n
    \E\bigl\{
      f(W_r + X_{\pi(r)})
      - f(W_r + Y_{\pi(r)})
    \bigr\} \\
    \nonumber
    &=
    \frac{1}{n}
    \sum_{i=1}^n
    \sum_{r=1}^n
    \E\bigl\{
      f(W_{r,-i} + X_i)
      - f(W_{r,-i} + Y_i)
    \bigr\} \\
    \label{eq:new_lindeberg_matching}
    &=
    \frac{1}{n}
    \sum_{i=1}^n
    \sum_{r=1}^{n}
    \sum_{|\alpha| = 1}^k
    \frac{1}{\alpha!}
    \E\bigl\{
      \partial_\alpha f(W_{r,-i})
      (X_i^\alpha - Y_i^\alpha)
    \bigr\} \\
    \label{eq:new_lindeberg_error_1}
    &\qquad+
    \frac{1}{n}
    \sum_{i=1}^n
    \sum_{r=1}^{n}
    \sum_{|\alpha| = k}
    \E \biggl\{
      \frac{X_i^\alpha}{\alpha!}
      \bigl(\partial_\alpha f(W_{r,-i} + \tilde X_{r,i})
      - \partial_\alpha f(W_{r,-i})\bigr)
    \biggr\} \\
    \label{eq:new_lindeberg_error_2}
    &\qquad-
    \frac{1}{n}
    \sum_{i=1}^n
    \sum_{r=1}^{n}
    \sum_{|\alpha| = k}
    \E \biggl\{
      \frac{Y_i^\alpha}{\alpha!}
      \bigl(\partial_\alpha f(W_{r,-i} + \tilde Y_{r,i})
      - \partial_\alpha f(W_{r,-i})\bigr)
    \biggr\}.
  \end{align}

  For the term in \eqref{eq:new_lindeberg_matching},
  since $X_i$ and $Y_i$ are each independent of $W_{r,-i}$
  and $f$ has bounded derivatives by Lemma~\ref{lem:smoothing},
  for each $r \in [n]$ and $\alpha \in \N_0^d$ with
  $|\alpha| \in [k]$, we may write
  \begin{align}
    \nonumber
    &\sum_{i=1}^n
    \E\bigl\{
      \partial_\alpha f(W_{r,-i})
      (X_i^\alpha - Y_i^\alpha)
    \bigr\}
    = \sum_{i=1}^n
    \delta_{\alpha,i}
    \E\bigl\{
      \partial_\alpha f(W_{r,-i})
    \bigr\} \\
    \nonumber
    &\qquad=
    \sum_{i=1}^n
    \delta_{\alpha,i}
    \Biggl(
      \frac{1}{n}
      \sum_{j=1}^{n}
      \E\bigl\{
        \partial_\alpha f(W_{r,-j})
      \bigr\}
    \Biggr)
    \iftoggle{journal}{\\ \nonumber &\qquad\qquad}{}
    +\sum_{i=1}^n
    \delta_{\alpha,i}
    \Biggl(
      \E\bigl\{
        \partial_\alpha f(W_{r,-i})
      \bigr\}
      - \frac{1}{n}
      \sum_{j=1}^{n}
      \E\bigl\{
        \partial_\alpha f(W_{r,-j})
      \bigr\}
    \Biggr) \\
    \label{eq:new_lindeberg_matching_average_remainder}
    &\qquad=
    \sum_{i=1}^{n}
    \bar\delta_{\alpha}
    \E\bigl\{
      \partial_\alpha f(W_{r,-i})
    \bigr\}
    + \frac{1}{n}
    \sum_{i=1}^n
    \sum_{j=1}^{n}
    \delta_{\alpha,i}
    \E\bigl\{
      \partial_\alpha f(W_{r,-i})
      - \partial_\alpha f(W_{r,-j})
    \bigr\}.
  \end{align}
  Consider the first term in
  \eqref{eq:new_lindeberg_matching_average_remainder}.
  Throughout, $Z \sim N(0_d, I_d)$ is an auxiliary
  Gaussian random vector independent of
  all the other quantities.
  By Lemma~\ref{lem:smoothing}, Jensen's inequality
  and orthogonality of Hermite polynomials
  \citep[see][Exercise~4.1]{samworth2026modern},
  \begin{align}
    \nonumber
    &\sum_{|\alpha| = 1}^k
    \frac{\bar\delta_{\alpha}}{\alpha!}
    \E\bigl\{
      \partial_\alpha f(W_{r,-i})
    \bigr\}
    =
    \sum_{|\alpha| = 1}^k
    \frac{\bar\delta_{\alpha}}{\alpha! \sigma^{|\alpha|}}
    \E\bigl\{
      g(W_{r,-i} + \sigma Z)
      H_\alpha(Z)
    \bigr\} \\
    \nonumber
    &\qquad=
    \sum_{|\alpha| = 1}^k
    \frac{\bar\delta_{\alpha}}{\alpha! \sigma^{|\alpha|}}
    \E\Biggl\{
      \biggl(g(W_{r,-i} + \sigma Z) - \frac{1}{2}\biggr)
      H_\alpha(Z)
    \Biggr\} \\
    \nonumber
    &\qquad\leq
    \Biggl\{\E\Biggl(
        \sum_{|\alpha| = 1}^k
        \sum_{|\alpha'| = 1}^k
        \frac{\bar\delta_{\alpha}}{\alpha! \sigma^{|\alpha|}}
        \frac{\bar\delta_{\alpha'}}{\alpha'! \sigma^{|\alpha'|}}
        \biggl(g(W_{r,-i} + \sigma Z) - \frac{1}{2}\biggr)^2
        H_\alpha(Z)
        H_{\alpha'}(Z)
    \Biggr)\Biggr\}^{1/2} \\
    \label{eq:lindeberg_Delta}
    &\qquad\leq
    \frac{1}{2}
    \Biggl(
      \sum_{|\alpha| = 1}^k
      \frac{\bar\delta_{\alpha}^2}{\alpha! \sigma^{2|\alpha|}}
    \Biggr)^{1/2}
    = \frac{1}{n}
    \Delta\biggl(\frac{1}{\sigma}\biggr).
  \end{align}
  The second term in
  \eqref{eq:new_lindeberg_matching_average_remainder}
  vanishes if $|\alpha| \leq m$,
  because then $\delta_{\alpha,i} = \mu_{X,\alpha} - \mu_{Y,\alpha}$
  does not depend on $i$.
  It also vanishes if $n = 1$, so suppose that $n \geq 2$.
  For $|\alpha| \in \{m+1, \ldots, k\}$,
  fix $i, j \in [n]$ with $i \neq j$ and
  write $B_{r,i} \coloneqq \{\pi_{r,i}(s): s \in [r - 1]\}$.
  With $B_{r,i}^\c \coloneqq ([n] \setminus \{i\}) \setminus B_{r,i}$,
  we have $W_{r,-i} = \sum_{s \in B_{r,i}} X_s + \sum_{s \in B_{r,i}^\c} Y_s$.
  Define $U_{r,i,j} \coloneqq
  \I\{j \in B_{r,i}\} X_i + \I\{j \in B_{r,i}^\c\} Y_i$
  and $V_{r,i,j} \coloneqq
  \I\{j \in B_{r,i}\} X_j + \I\{j \in B_{r,i}^\c\} Y_j$,
  and set $W_{r,-j}' \coloneqq W_{r,-i} + U_{r,i,j} - V_{r,i,j}$.
  In other words, if $W_{r,-i}$ contains $X_j$ as a summand
  then we replace $X_j$ with $X_i$; conversely,
  if it contains $Y_j$ then we replace $Y_j$ with $Y_i$.
  We claim that $W_{r,-j}'$ has the same law as $W_{r,-j}$.
  Indeed, let $\pi_{r,i}'$ be the permutation obtained by
  first applying $\pi_{r,i}$, then transposing $i \leftrightarrow j$.
  This transposition is a bijection between the
  set of permutations mapping $r \to i$
  and the set sending $r \to j$,
  so $\pi_{r,i}'$ is uniform over permutations
  with $\pi_{r,i}'(r)=j$.
  Hence $W_{r,-j}'$ has the same law as $W_{r,-j}$.
  Write $W_{r,-i,-j} \coloneqq W_{r,-i} - V_{r,i,j} = W_{r,-j}' - U_{r,i,j}$.

  Apply the argument of
  \citet[Theorem~18.18]{aliprantis2006infinite}
  to the random functions
  $t \mapsto \sum_{|\alpha| = \ell} \delta_{\alpha,i}
  \partial_\alpha f(W_{r,-i,-j} + t V_{r,i,j}) / \alpha!$,
  and likewise with $U_{r,i,j}$ in place of $V_{r,i,j}$,
  for each $\ell \in \{m+1, \ldots, k\}$
  to see that there exist random variables
  $\tilde U_{\ell,r,i,j}$ on the line segment between $0_d$ and $U_{r,i,j}$,
  and $\tilde V_{\ell,r,i,j}$ on the line segment between $0_d$ and
  $V_{r,i,j}$, such that
  \begin{align}
    \nonumber
    &\sum_{\ell = m + 1}^k
    \sum_{|\alpha| = \ell}
    \frac{\delta_{\alpha,i}}{\alpha!}
    \bigl\{
      \partial_\alpha f(W_{r,-i})
      - \partial_\alpha f(W_{r,-j}')
    \bigr\} \\
    \nonumber
    &\quad=
    \sum_{\ell = m + 1}^k
    \sum_{|\alpha| = \ell}
    \frac{\delta_{\alpha,i}}{\alpha!}
    \Bigl\{
      \bigl(
        \partial_\alpha f(W_{r,-i})
        - \partial_\alpha f(W_{r,-i,-j})
      \bigr)
      -
      \bigl(
        \partial_\alpha f(W_{r,-j}')
        - \partial_\alpha f(W_{r,-i,-j})
      \bigr)
    \Bigr\} \\
    \label{eq:new_lindeberg_coupling_exact}
    &\quad=
    \sum_{\ell = m + 1}^k
    \sum_{|\alpha| = \ell}
    \frac{\delta_{\alpha,i}}{\alpha!}
    \sum_{|\nu|=1}^{k_\ell}
    \frac{1}{\nu!}
    \partial_{\alpha + \nu} f(W_{r,-i,-j})
    (V_{r,i,j}^{\nu} - U_{r,i,j}^{\nu}) \\
    \label{eq:new_lindeberg_coupling_one}
    &\qquad\qquad+
    \sum_{\ell = m + 1}^k
    \sum_{|\alpha| = \ell}
    \frac{\delta_{\alpha,i}}{\alpha!}
    \sum_{|\nu|={k_\ell}}
    \frac{V_{r,i,j}^{\nu}}{\nu!}
    \bigl(
      \partial_{\alpha + \nu} f(W_{r,-i,-j} + \tilde V_{\ell,r,i,j})
      - \partial_{\alpha + \nu} f(W_{r,-i,-j})
    \bigr) \\
    \label{eq:new_lindeberg_coupling_two}
    &\qquad\qquad-
    \sum_{\ell = m + 1}^k
    \sum_{|\alpha| = \ell}
    \frac{\delta_{\alpha,i}}{\alpha!}
    \sum_{|\nu|={k_\ell}}
    \frac{U_{r,i,j}^{\nu}}{\nu!}
    \bigl(
      \partial_{\alpha + \nu} f(W_{r,-i,-j} + \tilde U_{\ell,r,i,j})
      - \partial_{\alpha + \nu} f(W_{r,-i,-j})
    \bigr).
  \end{align}
  With $\bX_{n,-i,-j} \coloneqq (X_s: s \in [n] \setminus \{i,j\})$
  and $\bY_{n,-i,-j} \coloneqq (Y_s: s \in [n] \setminus \{i,j\})$,
  $W_{r,-i,-j}$ depends only on
  $\pi_{r,i}$, $\bX_{n,-i,-j}$ and $\bY_{n,-i,-j}$.
  Since $U_{r,i,j}$ is a function of $\pi_{r,i}$, $X_i$ and $Y_i$,
  and $V_{r,i,j}$ is a function of $\pi_{r,i}$, $X_j$ and $Y_j$,
  we deduce that $W_{r,-i,-j}$ is independent of
  $(U_{r,i,j}, V_{r,i,j})$ conditional on $\pi_{r,i}$.
  For \eqref{eq:new_lindeberg_coupling_exact},
  taking an expectation yields
  \begin{align*}
    &\E \bigl\{
      \partial_{\alpha + \nu} f(W_{r,-i,-j})
      (V_{r,i,j}^{\nu} - U_{r,i,j}^{\nu})
    \bigr\}
    =
    \E \Bigl(
      \E \bigl\{
        \partial_{\alpha + \nu} f(W_{r,-i,-j})
        (V_{r,i,j}^{\nu} - U_{r,i,j}^{\nu})
        \mid \pi_{r,i}
      \bigr\}
    \Bigr) \\
    &\qquad=
    \E \Bigl(
      \E \bigl\{
        \partial_{\alpha + \nu} f(W_{r,-i,-j})
        \mid \pi_{r,i}
      \bigr\}
      \E \bigl\{
        V_{r,i,j}^{\nu} - U_{r,i,j}^{\nu}
        \mid \pi_{r,i}
      \bigr\}
    \Bigr).
  \end{align*}
  Moreover, as $|\nu| \in [k_\ell] \subseteq [m]$, we have
  \begin{align*}
    \E \bigl\{
      V_{r,i,j}^{\nu} - U_{r,i,j}^{\nu}
      \mid \pi_{r,i}
    \bigr\}
    &=
    \I\{j \in B_{r,i}\}
    \E (X_j^{\nu} - X_i^{\nu})
    + \I\{j \in B_{r,i}^\c\}
    \E (Y_j^{\nu} - Y_i^{\nu}) \\
    &=
    \I\{j \in B_{r,i}\}
    (\mu_{X,\nu} - \mu_{X,\nu})
    + \I\{j \in B_{r,i}^\c\}
    (\mu_{Y,\nu} - \mu_{Y,\nu})
    = 0,
  \end{align*}
  so the expectation of
  \eqref{eq:new_lindeberg_coupling_exact} vanishes.
  Define $G_{\ell,r,i,j} \coloneqq
  g(W_{r,-i,-j} + \tilde V_{\ell,r,i,j} + \sigma Z) - g(W_{r,-i,-j} +
  \sigma Z)$
  and $F_{\ell,r,i,j} \coloneqq
  g(W_{r,-i,-j} + \tilde U_{\ell,r,i,j} + \sigma Z) - g(W_{r,-i,-j} +
  \sigma Z)$.
  Since $|g(x+y) - g(x)| \leq (\|y\|_p / \eta) \land 1
  \leq (\|y\|_p / \eta)^\gamma$
  where $\gamma \in (0, 1]$, we have
  $|G_{\ell,r,i,j}| \leq
  \|V_{r,i,j}\|_p^{\gamma} / \eta^{\gamma}$;
  likewise for $F_{\ell,r,i,j}$.
  Set $\cA_{\ell,\theta} \coloneqq \bigl\{(\alpha, \nu):
  |\alpha| = \ell,\, |\nu| = k_\ell,\, \alpha + \nu = \theta\bigr\}$.
  Hence, using the terms in \eqref{eq:new_lindeberg_coupling_one}
  and \eqref{eq:new_lindeberg_coupling_two},
  along with Lemma~\ref{lem:smoothing}, Jensen's inequality
  and orthogonality of Hermite polynomials, we obtain
  \begin{align*}
    &\sum_{|\alpha| = m + 1}^k
    \frac{\delta_{\alpha,i}}{\alpha!}
    \E \bigl\{
      \partial_\alpha f(W_{r,-i})
      - \partial_\alpha f(W_{r,-j})
    \bigr\} \\
    &\quad=
    \E\Biggl\{
      \sum_{\ell = m + 1}^k
      \! \!
      G_{\ell,r,i,j}
      \!
      \sum_{|\alpha| = \ell}
      \sum_{|\nu|=k_\ell}
      \!
      \frac{\delta_{\alpha,i} V_{r,i,j}^\nu H_{\alpha + \nu}(Z)}
      {\alpha! \nu! \sigma^{\ell + k_\ell}}
    \Biggr\}
    \iftoggle{journal}{\\ &\qquad\qquad}{}
    -
    \E\Biggl\{
      \sum_{\ell = m + 1}^k
      \! \!
      F_{\ell,r,i,j}
      \!
      \sum_{|\alpha| = \ell}
      \sum_{|\nu|=k_\ell}
      \!
      \frac{\delta_{\alpha,i} U_{r,i,j}^\nu H_{\alpha + \nu}(Z)}
      {\alpha! \nu! \sigma^{\ell + k_\ell}}
    \Biggr\} \\
    &\quad\leq
    \sum_{\ell = m + 1}^k
    \E\Biggl\{
      \Biggl|
      \sum_{|\alpha| = \ell}
      \sum_{|\nu|=k_\ell}
      \frac{\delta_{\alpha,i} \|V_{r,i,j}\|_p^{\gamma}
      V_{r,i,j}^\nu H_{\alpha + \nu}(Z)}
      {\alpha! \nu! \sigma^{\ell + k_\ell} \eta^{\gamma}}
      \Biggr|
      \iftoggle{journal}{\\ &\qquad\qquad}{}
      + \Biggl|
      \sum_{|\alpha| = \ell}
      \sum_{|\nu|=k_\ell}
      \frac{\delta_{\alpha,i} \|U_{r,i,j}\|_p^{\gamma}
      U_{r,i,j}^\nu H_{\alpha + \nu}(Z)}
      {\alpha! \nu! \sigma^{\ell + k_\ell} \eta^{\gamma}}
      \Biggr|
    \Biggr\} \\
    &\quad=
    \sum_{\ell = m + 1}^k
    \E\Biggl\{
      \Biggl|
      \sum_{|\theta| = \ell + k_\ell}
      \sum_{(\alpha, \nu) \in \cA_{\ell,\theta}}
      \frac{\delta_{\alpha,i} \|V_{r,i,j}\|_p^{\gamma}
      V_{r,i,j}^\nu H_{\theta}(Z)}
      {\alpha! \nu! \sigma^{\ell + k_\ell} \eta^{\gamma}}
      \Biggr|
    \Biggr\} \\
    &\qquad\qquad+
    \sum_{\ell = m + 1}^k
    \E\Biggl\{
      \Biggl|
      \sum_{|\theta| = \ell + k_\ell}
      \sum_{(\alpha, \nu) \in \cA_{\ell,\theta}}
      \frac{\delta_{\alpha,i} \|U_{r,i,j}\|_p^{\gamma}
      U_{r,i,j}^\nu H_{\theta}(Z)}
      {\alpha! \nu! \sigma^{\ell + k_\ell} \eta^{\gamma}}
      \Biggr|
    \Biggr\} \\
    &\quad\leq
    \sum_{\ell = m + 1}^k
    \E\Biggl[
      \Biggl(
        \sum_{|\theta| = \ell + k_\ell}
        \theta!
        \Biggl\{
          \sum_{(\alpha, \nu) \in \cA_{\ell,\theta}}
          \frac{\delta_{\alpha,i} \|V_{r,i,j}\|_p^{\gamma}
          V_{r,i,j}^\nu}
          {\alpha! \nu! \sigma^{\ell + k_\ell} \eta^{\gamma}}
        \Biggr\}^2
      \Biggr)^{1/2}
    \Biggr] \\
    &\qquad\qquad+
    \sum_{\ell = m + 1}^k
    \E\Biggl[
      \Biggl(
        \sum_{|\theta| = \ell + k_\ell}
        \theta!
        \Biggl\{
          \sum_{(\alpha, \nu) \in \cA_{\ell,\theta}}
          \frac{\delta_{\alpha,i} \|U_{r,i,j}\|_p^{\gamma}
          U_{r,i,j}^\nu}
          {\alpha! \nu! \sigma^{\ell + k_\ell} \eta^{\gamma}}
        \Biggr\}^2
      \Biggr)^{1/2}
    \Biggr].
  \end{align*}
  Define $N_{\ell,\theta} \coloneqq |\cA_{\ell,\theta}|$,
  so that $N_{\ell, \theta} \leq 2^{\ell + k_\ell}$.
  Moreover, $\theta! \leq 2^{|\theta|} \alpha! \nu!$
  whenever $\alpha + \nu = \theta$,
  so by Cauchy--Schwarz and the multinomial theorem,
  \begin{align*}
    &\sum_{|\theta| = \ell + k_\ell}
    \theta!
    \Biggl\{
      \sum_{(\alpha, \nu) \in \cA_{\ell,\theta}}
      \frac{\delta_{\alpha,i} \|V_{r,i,j}\|_p^{\gamma}
      V_{r,i,j}^\nu}
      {\alpha! \nu! \sigma^{\ell + k_\ell} \eta^{\gamma}}
    \Biggr\}^2
    \leq
    \sum_{|\theta| = \ell + k_\ell}
    \theta!
    N_{\ell,\theta}
    \sum_{(\alpha, \nu) \in \cA_{\ell,\theta}}
    \frac{\delta_{\alpha,i}^2 \|V_{r,i,j}\|_p^{2 \gamma}
    V_{r,i,j}^{2\nu}}
    {\alpha!^2 \nu!^2 \sigma^{2\ell + 2k_\ell} \eta^{2\gamma}} \\
    &\qquad\leq
    \frac{\|V_{r,i,j}\|_p^{2 \gamma}}{\eta^{2\gamma}}
    \biggl(\frac{2}{\sigma}\biggr)^{2\ell + 2k_\ell}
    \frac{\|V_{r,i,j}\|_2^{2 k_\ell}}{k_\ell!}
    \sum_{|\alpha| = \ell}
    \frac{\delta_{\alpha,i}^2}{\alpha!}.
  \end{align*}
  A bound for the term involving $U_{r,i,j}$ follows in exactly
  the same fashion. Therefore,
  \begin{align*}
    &\sum_{\ell = m + 1}^k
    \E\Biggl[
      \Biggl(
        \sum_{|\theta| = \ell + k_\ell}
        \theta!
        \Biggl\{
          \sum_{(\alpha,\nu) \in \cA_{\ell,\theta}}
          \frac{\delta_{\alpha,i} \|V_{r,i,j}\|_p^{\gamma} V_{r,i,j}^\nu}
          {\alpha! \nu! \sigma^{\ell + k_\ell} \eta^\gamma}
        \Biggr\}^2
      \Biggr)^{1/2}
    \Biggr] \\
    &\qquad\leq
    \sum_{\ell = m + 1}^k
    \biggl(\frac{2}{\sigma}\biggr)^{\ell + k_\ell}
    \E\Biggl\{
      \frac{\|V_{r,i,j}\|_p^{\gamma}}{\eta^{\gamma}}
      \frac{\|V_{r,i,j}\|_2^{k_\ell}}{\sqrt{k_\ell!}}
    \Biggr\}
    \Biggl(
      \sum_{|\alpha| = \ell}
      \frac{\delta_{\alpha,i}^2}{\alpha!}
    \Biggr)^{1/2},
  \end{align*}
  and similarly for $U_{r,i,j}$.
  Since $\P(j \in B_{r,i})
  = \P\bigl(j \in \pi([r-1]) \mid \pi(r) = i\bigr)
  = (r-1) / (n-1)$ for $i \neq j$,
  \begin{align*}
    &
    \frac{1}{n \sqrt{k_\ell!}}
    \sum_{r=1}^{n}
    \E\Bigl\{
      \|V_{r,i,j}\|_2^{k_\ell} \|V_{r,i,j}\|_p^{\gamma}
      + \|U_{r,i,j}\|_2^{k_\ell} \|U_{r,i,j}\|_p^{\gamma}
    \Bigr\} \\
    &\quad=
    \frac{1}{n \sqrt{k_\ell!}}
    \sum_{r=1}^{n}
    \biggl(
      \frac{r-1}{n-1}
      \E\Bigl\{
        \|X_j\|_2^{k_\ell} \|X_j\|_p^{\gamma}
        + \|X_i\|_2^{k_\ell} \|X_i\|_p^{\gamma}
      \Bigr\} \\
      &\qquad\qquad\qquad\qquad+
      \frac{n-r}{n-1}
      \E\Bigl\{
        \|Y_j\|_2^{k_\ell} \|Y_j\|_p^{\gamma}
        + \|Y_i\|_2^{k_\ell} \|Y_i\|_p^{\gamma}
      \Bigr\}
    \biggr) \\
    &\quad=
    \frac{1}{2 \sqrt{k_\ell!}}
    \E\Bigl\{
      \|X_j\|_2^{k_\ell} \|X_j\|_p^{\gamma}
      + \|X_i\|_2^{k_\ell} \|X_i\|_p^{\gamma}
    \Bigr\}
    + \frac{1}{2 \sqrt{k_\ell!}}
    \E\Bigl\{
      \|Y_j\|_2^{k_\ell} \|Y_j\|_p^{\gamma}
      + \|Y_i\|_2^{k_\ell} \|Y_i\|_p^{\gamma}
    \Bigr\} \\
    &\quad=
    \frac{\beta_{k_\ell,\gamma,j}
    + \beta_{k_\ell,\gamma,i}}{2}.
  \end{align*}
  Therefore,
  \begin{align}
    \nonumber
    &\frac{1}{n^2}
    \sum_{i=1}^n
    \sum_{j=1}^{n}
    \sum_{r=1}^{n}
    \sum_{|\alpha| = 1}^k
    \frac{\delta_{\alpha,i}}{\alpha!}
    \E\bigl\{
      \partial_\alpha f(W_{r,-i})
      - \partial_\alpha f(W_{r,-j})
    \bigr\} \\
    \nonumber
    &\quad\leq
    \frac{1}{n^2}
    \sum_{i=1}^n
    \sum_{j \in [n] \setminus \{i\}}
    \sum_{r=1}^{n}
    \sum_{\ell = m + 1}^k
    \biggl(\frac{2}{\sigma}\biggr)^{\ell + k_\ell}
    \Biggl(
      \sum_{|\alpha| = \ell}
      \frac{\delta_{\alpha,i}^2}{\alpha!}
    \Biggr)^{1/2}
    \iftoggle{journal}{\\ \nonumber &\qquad\qquad\times}{}
    \E\biggl\{
      \frac{\|V_{r,i,j}\|_2^{k_\ell} \|V_{r,i,j}\|_p^{\gamma}
      + \|U_{r,i,j}\|_2^{k_\ell} \|U_{r,i,j}\|_p^{\gamma}}
      {\sqrt{k_\ell!} \eta^{\gamma}}
    \biggr\} \\
    \nonumber
    &\quad\leq
    \frac{1}{n}
    \sum_{i=1}^n
    \sum_{j=1}^{n}
    \sum_{\ell = m + 1}^k
    \biggl(\frac{2}{\sigma}\biggr)^{\ell + k_\ell}
    \Biggl(
      \sum_{|\alpha| = \ell}
      \frac{\delta_{\alpha,i}^2}{\alpha!}
    \Biggr)^{1/2}
    \frac{\beta_{k_\ell,\gamma,j}
    + \beta_{k_\ell,\gamma,i}}{2 \eta^\gamma} \\
    \label{eq:lindeberg_Gamma}
    &\quad=
    \sum_{i=1}^n
    \sum_{\ell = m + 1}^k
    \biggl(\frac{2}{\sigma}\biggr)^{\ell + k_\ell}
    \frac{\bar \beta_{k_\ell,\gamma}
    + \beta_{k_\ell,\gamma,i}}{2 \eta^\gamma}
    \Biggl(
      \sum_{|\alpha| = \ell}
      \frac{\delta_{\alpha,i}^2}{\alpha!}
    \Biggr)^{1/2}
    = \Gamma_\eta\biggl(\frac{1}{\sigma}\biggr).
  \end{align}

  It remains to handle
  \eqref{eq:new_lindeberg_error_1} and
  \eqref{eq:new_lindeberg_error_2}.
  By Lemma~\ref{lem:smoothing},
  Jensen's inequality,
  orthogonality of Hermite polynomials
  and the multinomial theorem, and since
  $|g(x+y) - g(x)| \leq \bigl(\|y\|_p / \eta\bigr) \land 1
  \leq \|y\|_p^{\gamma} / \eta^{\gamma}$, we have that
  \begin{align*}
    &\Biggl|\sum_{|\alpha| = k}
    \E \biggl\{
      \frac{X_i^\alpha}{\alpha!}
      \bigl(\partial_\alpha f(W_{r,-i} + \tilde X_{r,i})
      - \partial_\alpha f(W_{r,-i})\bigr)
    \biggr\}\Biggr| \\
    &\qquad=
    \Biggl| \E \Biggl\{
      \sum_{|\alpha| = k}
      \frac{X_i^\alpha}{\alpha! \sigma^k}
      \bigl(
        g(W_{r,-i} + \tilde X_{r,i} + \sigma Z)
        - g(W_{r,-i} + \sigma Z)
      \bigr)
      H_\alpha(Z)
    \Biggr\}\Biggr| \\
    &\qquad\leq
    \E \Biggl\{
      \frac{\|X_i\|_p^{\gamma}}{\sigma^k \eta^{\gamma}}
      \Biggl(
        \sum_{|\alpha| = k}
        \frac{X_i^{2\alpha}}{\alpha!}
      \Biggr)^{1/2}
    \Biggr\}
    = \frac{\E\bigl\{\|X_i\|_2^{k} \|X_i\|_p^{\gamma}\bigr\}}
    {\sqrt{k!} \sigma^k \eta^{\gamma}}.
  \end{align*}
  Applying the same approach to
  \eqref{eq:new_lindeberg_error_2}, we obtain
  \begin{align}
    \nonumber
    &\frac{1}{n}
    \sum_{i=1}^n
    \sum_{r=1}^{n}
    \sum_{|\alpha| = k}
    \E \biggl\{
      \frac{X_i^\alpha}{\alpha!}
      \bigl(\partial_\alpha f(W_{r,-i} + \tilde X_{r,i})
      - \partial_\alpha f(W_{r,-i})\bigr)
    \biggr\} \\
    \label{eq:lindeberg_beta}
    &\qquad-
    \frac{1}{n}
    \sum_{i=1}^n
    \sum_{r=1}^{n}
    \sum_{|\alpha| = k}
    \E \biggl\{
      \frac{Y_i^\alpha}{\alpha!}
      \bigl(\partial_\alpha f(W_{r,-i} + \tilde Y_{r,i})
      - \partial_\alpha f(W_{r,-i})\bigr)
    \biggr\}
    \leq
    \frac{n \bar\beta_{k,\gamma}}
    {\sigma^k \eta^{\gamma}}.
  \end{align}
  The result now follows by
  \eqref{eq:lindeberg_Delta},
  \eqref{eq:lindeberg_Gamma} and \eqref{eq:lindeberg_beta}.
\end{proof}

\begin{lemma}[Gaussian $\ell_p$-norm tail bound]%
  \label{lem:gaussian_norm_tail}
  Let $d \in \N$ and $Z \sim N(0_d, I_d)$.
  Then for each $p \in [1, \infty]$ and $t \geq 0$,
  with $\phi(t) \coloneqq \phi_{p,d}(t)$
  defined as in Section~\ref{sec:yurinskii}, we have
  \begin{equation*}
    \P\bigl(\|Z\|_p > \phi(t)\bigr)
    \leq e^{-t^2 / 2}.
  \end{equation*}
\end{lemma}

\begin{proof}[Proof of Lemma~\ref{lem:gaussian_norm_tail}]
  By \citet[Lemma~SA.4]{cattaneo2025yurinskii},
  \begin{equation*}
    \E\bigl(\|Z\|_p\bigr)
    \leq
    \begin{cases}
      \sqrt{p} d^{1/p}
      & \text{ if } p \in [1,\infty),\\
      \sqrt{2\log 2d}
      & \text{ if } p =\infty.
    \end{cases}
  \end{equation*}
  For $x \in \R^d$, if $p \in [2, \infty]$
  then $\|x\|_p \leq \|x\|_2$.
  Otherwise, if $p \in [1, 2)$ then by H{\"o}lder's inequality,
  $\|x\|_p \leq d^{1/p-1/2} \|x\|_2$.
  The result follows by the Gaussian concentration inequality
  given by \citet[Theorem~5.6]{boucheron2013concentration}.
\end{proof}

\begin{lemma}[Borel set approximation]%
  \label{lem:borel}
  Grant Assumption~\ref{ass:data}.
  For each $t, \eta > 0$ and $A \in \cB(\R^d)$,
  \begin{equation*}
    \P(S_n \in A)
    - \P\big( T_n \in A_p^{3\eta} \big)
    \leq
    2 e^{-t^2 / 2}
    + \frac{n \bar\beta_{k,\gamma} \phi(t)^k}{\eta^{k+\gamma}}
    + \Delta\biggl(\frac{\phi(t)}{\eta}\biggr)
    + \Gamma_\eta\biggl(\frac{\phi(t)}{\eta}\biggr).
  \end{equation*}
  If $\Delta(s) = \Gamma_\eta(s) = 0$ for every $s > 0$
  and $\bar\beta_{k,\gamma} > 0$, then further
  \begin{equation}
    \label{eq:borel_matching}
    \P(S_n \in A)
    - \P\big( T_n \in A_p^{3\eta} \big)
    \leq
    \frac{3 n \bar\beta_{k,\gamma}}{\eta^{k + \gamma}}
    \phi\Biggl(
      \sqrt{2 \logplus
        \frac{\eta^{k + \gamma}}
        {n \bar\beta_{k,\gamma} \psi^k}
      }
      \,
    \Biggr)^k.
  \end{equation}
  If $\bar\beta_{k,\gamma} = 0$, then $S_n = T_n = 0$ almost surely,
  so the probability difference is nonpositive.
\end{lemma}

\begin{proof}[Proof of Lemma~\ref{lem:borel}]
  Take $\sigma > 0$ and
  $f \coloneqq f_{A,p,\eta,\sigma}$
  as in Lemma~\ref{lem:smoothing}.
  By Lemmas~\ref{lem:smoothing} and~\ref{lem:lindeberg},
  \begin{align*}
    \P(S_n \in A)
    &=
    \E\big\{\I\{S_n \in A\} - f(S_n)\big\}
    + \E\big\{f(S_n) - f(T_n)\big\}
    + \E \big\{f(T_n)\big\} \\
    &\leq
    \varepsilon \, \P(S_n \in A)
    + \frac{n \bar\beta_{k,\gamma}}{\sigma^k \eta^\gamma}
    + \Delta\biggl(\frac{1}{\sigma}\biggr)
    + \Gamma_\eta\biggl(\frac{1}{\sigma}\biggr)
    + \varepsilon
    + (1 - \varepsilon) \, \P\big(T_n \in A_p^{3\eta}\big) \\
    &\leq
    \P\big( T_n \in A_p^{3\eta} \big)
    + 2\varepsilon
    + \frac{n \bar\beta_{k,\gamma}}{\sigma^k \eta^\gamma}
    + \Delta\biggl(\frac{1}{\sigma}\biggr)
    + \Gamma_\eta\biggl(\frac{1}{\sigma}\biggr).
  \end{align*}
  Set $\sigma \coloneqq \eta / t$ for $t>0$ and recall that
  $\varepsilon \coloneqq \P\bigl(\|Z\|_p > \eta / \sigma\bigr)$
  to obtain
  \begin{equation*}
    \P(S_n \in A)
    - \P\big( T_n \in A_p^{3\eta} \big)
    \leq
    2 \, \P\bigl(\|Z\|_p > t\bigr)
    + \frac{n \bar\beta_{k,\gamma} t^k}{\eta^{k+\gamma}}
    + \Delta\biggl(\frac{t}{\eta}\biggr)
    + \Gamma_\eta\biggl(\frac{t}{\eta}\biggr).
  \end{equation*}
  The first result now follows by Lemma~\ref{lem:gaussian_norm_tail},
  replacing $t$ with $\phi(t)$.
  For \eqref{eq:borel_matching}, suppose that
  $\Delta(s) = \Gamma_\eta(s) = 0$ for all $s > 0$
  and $\bar\beta_{k,\gamma} > 0$.
  If $\eta^{k+\gamma} \leq n \bar\beta_{k,\gamma} \psi^k$,
  the bound is trivial because its right-hand side is at least $3$.
  Otherwise, taking
  $t \coloneqq \bigl\{2 \logplus \bigl(n^{-1} \bar \beta_{k,\gamma}^{-1}
  \eta^{k + \gamma} \psi^{-k}\bigr)\bigr\}^{1/2}$
  and noting that $\psi \leq \phi(t)$
  for all $t \geq 0$ yields
  \begin{align*}
    \P(S_n \in A)
    - \P\big( T_n \in A_p^{3\eta} \big)
    &\leq
    \frac{2 n \bar \beta_{k,\gamma} \psi^k}{\eta^{k + \gamma}}
    + \frac{n \bar \beta_{k,\gamma}}{\eta^{k + \gamma}}
    \phi\Biggl(
      \sqrt{2 \logplus \frac{\eta^{k + \gamma}}{n \bar \beta_{k,\gamma}
      \psi^k}}
      \,
    \Biggr)^k \\
    &\leq
    \frac{3 n \bar \beta_{k,\gamma}}{\eta^{k + \gamma}}
    \phi\Biggl(
      \sqrt{2 \logplus \frac{\eta^{k + \gamma}}
      {n \bar \beta_{k,\gamma} \psi^k}}
      \,
    \Biggr)^k.
    \iftoggle{journal}{}{\qedhere}
  \end{align*}
\end{proof}

\begin{proof}[Proof of Theorem~\ref{thm:yurinskii}
  and \eqref{eq:yurinskii_matching}]
  Apply Strassen's theorem \citep[Section~10.3, Theorem~8]{pollard2002user}
  to the bound in Lemma~\ref{lem:borel},
  noting that $A_p^{3\eta}\subseteq\{x:\|x-A\|_p\leq3\eta\}$.
  To obtain \eqref{eq:yurinskii_matching},
  similarly apply Strassen's theorem to \eqref{eq:borel_matching}.
\end{proof}

\begin{lemma}[Bound for $\Gamma_\eta$]%
  \label{lem:bound_Gamma}
  Grant Assumption~\ref{ass:data}.
  For each $i \in [n]$ and $\ell \in [k]$, we have
  \begin{equation*}
    \sum_{|\alpha| = \ell}
    \frac{\delta_{\alpha,i}^2}{\alpha!}
    \leq
    \beta_{\ell,0,i}^2.
  \end{equation*}
  Therefore,
  \begin{equation*}
    \Gamma_\eta(s)
    \leq
    \sum_{i=1}^n
    \sum_{\ell = m+1}^k
    \frac{(2 s)^{\ell + k_\ell}}
    {2 \eta^{\gamma}}
    \bigl(\bar \beta_{k_\ell,\gamma}
    + \beta_{k_\ell,\gamma,i}\bigr) \beta_{\ell,0,i}.
  \end{equation*}
\end{lemma}

\begin{proof}[Proof of Lemma~\ref{lem:bound_Gamma}]
  For $u, v \in \R^d$,
  writing $u \circ v \in \R^d$ for the entrywise product of
  $u$ and $v$, by the multinomial theorem and Cauchy--Schwarz,
  \begin{equation*}
    \sum_{|\alpha| = \ell}
    \frac{u^\alpha v^\alpha}{\alpha!}
    =
    \sum_{|\alpha| = \ell}
    \frac{(u \circ v)^\alpha}{\alpha!}
    =
    \frac{1}{\ell!}
    \Biggl\{
      \sum_{j=1}^d (u \circ v)_j
    \Biggr\}^\ell
    =
    \frac{(u^\T v)^\ell}{\ell!}
    \leq
    \frac{\|u\|_2^\ell \|v\|_2^\ell}{\ell!}.
  \end{equation*}
  Therefore, taking $\tilde X_i$ and $\tilde Y_i$
  to be copies of $X_i$ and $Y_i$, respectively,
  independent of each other and of $(X_i, Y_i)$
  for each $i \in [n]$,
  \begin{align*}
    \sum_{|\alpha| = \ell}
    \frac{\delta_{\alpha,i}^2}{\alpha!}
    &=
    \sum_{|\alpha| = \ell}
    \frac{\{\E(X_i^\alpha) - \E(Y_i^\alpha)\}^2}
    {\alpha!}
    =
    \sum_{|\alpha| = \ell}
    \frac{\{\E(X_i^\alpha)\}^2 - 2 \E(X_i^\alpha) \E(Y_i^\alpha)
    + \{\E(Y_i^\alpha)\}^2}
    {\alpha!} \\
    &=
    \E \Biggl\{
      \sum_{|\alpha| = \ell}
      \frac{X_i^\alpha \tilde X_i^\alpha
        - 2 X_i^\alpha \tilde Y_i^\alpha
      + Y_i^\alpha \tilde Y_i^\alpha}
      {\alpha!}
    \Biggr\} \\
    &\leq
    \E \Biggl\{
      \frac{
        \|X_i\|_2^\ell \|\tilde X_i\|_2^\ell
        + 2 \|X_i\|_2^\ell \|\tilde Y_i\|_2^\ell
        + \|Y_i\|_2^\ell \|\tilde Y_i\|_2^\ell
      }
      {\ell!}
    \Biggr\} \\
    &=
    \frac{
      \bigl\{\E\bigl(\|X_i\|_2^\ell\bigr)\bigr\}^2
      + 2 \E\bigl(\|X_i\|_2^\ell\bigr)
      \E\bigl(\|Y_i\|_2^\ell\bigr)
      + \bigl\{\E\bigl(\|Y_i\|_2^\ell\bigr)\bigr\}^2
    }
    {\ell!} \\
    &=
    \frac{
      \bigl\{
        \E\bigl(\|X_i\|_2^\ell + \|Y_i\|_2^\ell\bigr)
      \bigr\}^2
    }
    {\ell!}
    = \beta_{\ell,0,i}^2,
  \end{align*}
  establishing the first claim.
  The second follows as $\beta_{\ell,0,i} \geq 0$.
\end{proof}

\begin{proof}[Details for \eqref{eq:yurinskii_rate_general} and
  \eqref{eq:yurinskii_rate_p}]
  As $r \geq 2 \lor (k + \gamma)$,
  H{\"o}lder's inequality gives, for all $x \in \R^d$,
  \begin{align*}
    \|x\|_2
    &\leq d^{1/2 - 1/r} \|x\|_r,
    &
    \|x\|_p
    &\leq d^{(1/p - 1/r) \lor 0} \|x\|_r.
  \end{align*}
  Thus, by Jensen's inequality, since $r \geq k + \gamma$,
  for each $i \in [n]$,
  \begin{align*}
    \beta_{k,\gamma,i}
    &=
    \frac{1}{\sqrt{k!}}
    \E \Bigl(
      \|X_i\|_2^k
      \| X_i\|_p^{\gamma}
      + \|Y_i\|_2^k
      \| Y_i\|_p^{\gamma}
    \Bigr) \\
    &\leq
    \frac{ d^{k/2 - k/r} d^{(\gamma/p - \gamma/r) \lor 0} }
    {\sqrt{k!}}
    \E \Bigl(
      \| X_i\|_r^{k+\gamma}
      + \| Y_i\|_r^{k+\gamma}
    \Bigr) \\
    &\leq
    \frac{ d^{k/2 - k/r} d^{(\gamma/p - \gamma/r) \lor 0} }
    {\sqrt{k!}}
    \Bigl(
      \bigvvvert
      \| X_i\|_r
      \bigvvvert_r^{k+\gamma}
      + \bigvvvert
      \| Y_i\|_r
      \bigvvvert_r^{k+\gamma}
    \Bigr) \\
    &\leq
    \frac{2 d^{k/2 - k/r} d^{(\gamma/p - \gamma/r) \lor 0} d^{k/r + \gamma/r}
    M^{k+\gamma}}
    {\sqrt{k!}}
    = \frac{2 d^{k/2 + \gamma / (p \land r)} M^{k+\gamma}}
    {\sqrt{k!}}.
  \end{align*}
  Substituting this into Theorem~\ref{thm:yurinskii},
  then choosing
  \begin{equation*}
    t
    \coloneqq
    \sqrt{2\log^+
      \frac{n^{k/2+\gamma/2-1}\eta^{k+\gamma}}
    {d^{k/2+\gamma/(p\land r)}M^{k+\gamma}\psi^k}},
  \end{equation*}
  yields the following bound,
  which is trivial when $t = 0$:
  \begin{align*}
    &
    \P\biggl(
      \frac{1}{\sqrt n}
      \|S_n' - T_n'\|_p > 3\eta
    \biggr)
    \iftoggle{journal}{\\&\qquad}{}
    \leq
    \frac{6 d^{k/2 + \gamma / (p \land r)}
    M^{k+\gamma} \psi^k}
    {n^{k/2 + \gamma/2 - 1} \eta^{k + \gamma} }
    \Biggl\{
      1 +
      \sqrt{2 \logplus
        \frac{n^{k/2 + \gamma/2 - 1} \eta^{k + \gamma}}
        {d^{k/2 + \gamma / (p \land r)}
        M^{k+\gamma} \psi^k}
      }
      \,
    \Biggr\}^k.
  \end{align*}
  We now find a sequence $\eta_n$ such that
  setting $\eta = \eta_n$ causes the right-hand side
  of the above display to converge to zero;
  we may then deduce that
  $\P\bigl(\|S_n' - T_n'\|_p / \sqrt n > 3\eta_n\bigr) \to 0$
  and so $\|S_n' - T_n'\|_p / \sqrt n \lesssim_\P \eta_n$.
  Since $x \{\log (1/x)\}^{k/2} \to 0$ as $x \to 0$,
  it suffices to let $R_n \to \infty$ and take
  \begin{equation*}
    \eta_n
    \coloneqq
    \frac{R_n d^{k/(2k + 2\gamma)
      + \gamma / \{(p \land r) (k + \gamma)\}}
    \psi^{k/(k + \gamma)} M}
    {n^{1/2 - 1/(k + \gamma)}}.
  \end{equation*}

  If $p \in [1, \infty)$ then
  $\psi \asymp d^{1/p}$; using
  $\gamma / (p \land r) \leq \gamma / p + \gamma / r
  \leq \gamma / p + \gamma / 2$, we obtain
  \begin{equation*}
    \eta_n
    \lesssim
    \frac{R_n d^{1/2 + 1 / p} M}
    {n^{1/2 - 1/(k + \gamma)}},
  \end{equation*}
  giving \eqref{eq:yurinskii_rate_p}.
  For $p = \infty$ we have
  $\psi \asymp \sqrt{\log 2 d}$;
  since $r \geq 2$, we obtain
  \begin{equation*}
    \eta_n
    \lesssim
    \frac{R_n (d \log 2d)^{1/2} M}
    {n^{1/2 - 1/(k + \gamma)}}.
  \end{equation*}

  If $\bigvvvert \|X_i\|_2 \bigvvvert_r
  \lor \bigvvvert \|Y_i\|_2 \bigvvvert_r \leq M$,
  where $r \in [k + \gamma, \infty]$, then
  $\|x\|_p \leq d^{(1/p - 1/2) \lor 0} \|x\|_2$ gives
  \begin{equation*}
    \bar\beta_{k,\gamma}
    \leq \frac{d^{\gamma\{(1/p - 1/2) \lor 0 \}}}{n\sqrt{k!}}
    \sum_{i=1}^n
    \E\bigl(\|X_i\|_2^{k+\gamma}
    + \|Y_i\|_2^{k+\gamma}\bigr)
    \leq
    \frac{2d^{\gamma\{(1/p-1/2)\lor0\}}M^{k+\gamma}}{\sqrt{k!}}.
  \end{equation*}
  Repeating the preceding argument with this bound gives
  \begin{equation*}
    \frac{\|S_n'-T_n'\|_p}{\sqrt n}
    \lesssim_\P
    \frac{R_n d^{\{(1/p - 1/2) \lor 0 \} \gamma/(k+\gamma)}
    \psi^{k/(k+\gamma)}M}
    {n^{1/2-1/(k+\gamma)}}.
    \iftoggle{journal}{}{\qedhere}
  \end{equation*}
\end{proof}

\subsection{Proofs for Section~\ref{sec:bootstrap}}
\label{sec:proofs_bootstrap}

\begin{lemma}[A multi-index Lipschitz bound]%
  \label{lem:multi_index}
  Let $d \in \N$ and $\alpha \in \N_0^d$.
  For each $x, y \in \R^d$,
  \begin{equation*}
    |x^\alpha - y^\alpha|
    \leq
    \sum_{r=1}^d
    \alpha_r
    v^{\alpha - e_r}
    |x_r - y_r|,
  \end{equation*}
  where $v \in \R^d$ is defined by
  $v_r \coloneqq |x_r| \lor |y_r|$
  for $r \in [d]$
  and $e_r$ is the $r$th standard basis vector in $\R^d$.
  Here, summands with $\alpha_r = 0$ are interpreted as zero.
\end{lemma}

\begin{proof}[Proof of Lemma~\ref{lem:multi_index}]
  By a telescoping sum and the mean value theorem,
  \begin{align*}
    |x^\alpha - y^\alpha|
    &=
    \Biggl|
    \prod_{j=1}^d x_j^{\alpha_j}
    - \prod_{j=1}^d y_j^{\alpha_j}
    \Biggr|
    =
    \Biggl|
    \sum_{r=1}^d
    \Biggl(
      \prod_{j=r+1}^d x_j^{\alpha_j}
    \Biggr)
    \Biggl(
      \prod_{\ell=1}^{r-1} y_\ell^{\alpha_\ell}
    \Biggr)
    \bigl(x_r^{\alpha_r} - y_r^{\alpha_r}\bigr)
    \Biggr| \\
    &\leq
    \sum_{r=1}^d
    \Biggl(
      \prod_{j=r+1}^d |x_j|^{\alpha_j}
    \Biggr)
    \Biggl(
      \prod_{\ell=1}^{r-1} |y_\ell|^{\alpha_\ell}
    \Biggr)
    \bigl|x_r^{\alpha_r} - y_r^{\alpha_r}\bigr| \\
    &\leq
    \sum_{r=1}^d
    \Biggl(
      \prod_{j=r+1}^d \bigl(|x_j| \lor |y_j|\bigr)^{\alpha_j}
    \Biggr)
    \Biggl(
      \prod_{\ell=1}^{r-1} \bigl(|x_\ell| \lor |y_\ell|\bigr)^{\alpha_\ell}
    \Biggr)
    \alpha_r
    \bigl(|x_r| \lor |y_r|\bigr)^{\alpha_r - 1}
    |x_r - y_r| \\
    &=
    \sum_{r=1}^d
    \alpha_r
    v^{\alpha - e_r}
    |x_r - y_r|.
    \iftoggle{journal}{}{\qedhere}
  \end{align*}
\end{proof}

\begin{lemma}[A multi-index sum inequality]%
  \label{lem:multi_index_sum}
  Let $d, \ell \in \N$. For each $x, y \in \R^d$,
  \begin{equation*}
    \sum_{|\alpha| = \ell}
    \frac{1}{\alpha!}
    (x^\alpha - y^\alpha)^2
    \leq
    \frac{\ell \cdot 2^{\ell-2}}{(\ell-1)!}
    \|x - y\|_2^2
    \bigl(
      \|x\|_2^{2\ell-2} + \|y\|_2^{2\ell-2}
    \bigr),
  \end{equation*}
  where the sum is over $\alpha \in \N_0^d$
  with $|\alpha| = \ell$.
\end{lemma}

\begin{proof}[Proof of Lemma~\ref{lem:multi_index_sum}]
  By Lemma~\ref{lem:multi_index},
  Cauchy--Schwarz and the multinomial theorem,
  as $\ell \geq 1$ and taking $v \in \R^d$ with
  $v_r \coloneqq |x_r| \lor |y_r|$,
  \begin{align*}
    \sum_{|\alpha| = \ell}
    \frac{1}{\alpha!}
    (x^\alpha - y^\alpha)^2
    &\leq
    \sum_{|\alpha| = \ell}
    \frac{1}{\alpha!}
    \Biggl(
      \sum_{r=1}^d
      \alpha_r
      v^{\alpha - e_r}
      |x_r - y_r|
    \Biggr)^2
    \iftoggle{journal}{\\ &}{}
    \leq
    \sum_{|\alpha| = \ell}
    \frac{\ell}{\alpha!}
    \sum_{r=1}^d
    \alpha_r
    v^{2\alpha - 2e_r}
    (x_r - y_r)^2 \\
    &=
    \sum_{r=1}^d
    (x_r - y_r)^2
    \sum_{|\alpha| = \ell}
    \I\{\alpha_r \geq 1\}
    \frac{\ell}{(\alpha - e_r)!}
    v^{2\alpha - 2e_r} \\
    &=
    \frac{\ell}{(\ell-1)!}
    \sum_{r=1}^d
    (x_r - y_r)^2
    \sum_{|\alpha| = \ell - 1}
    \frac{(\ell-1)!}{\alpha!}
    v^{2\alpha}
    \iftoggle{journal}{\\ &}{}
    =
    \frac{\ell}{(\ell-1)!}
    \|x - y\|_2^2
    \cdot
    \|v\|_2^{2 \ell - 2} \\
    &\leq
    \frac{\ell}{(\ell-1)!}
    \|x - y\|_2^2
    \bigl(
      \|x\|_2^2 + \|y\|_2^2
    \bigr)^{\ell-1}
    \iftoggle{journal}{\\ &}{}
    \leq
    \frac{\ell \cdot 2^{\ell-2}}{(\ell-1)!}
    \|x - y\|_2^2
    \bigl(
      \|x\|_2^{2\ell-2} + \|y\|_2^{2\ell-2}
    \bigr),
  \end{align*}
  as required.
\end{proof}

\begin{lemma}[Bootstrap moment bounds]%
  \label{lem:beta_hat_bootstrap}
  Grant the setup of Theorem~\ref{thm:bootstrap}
  and, for $a \in \N_0$, $b \geq 0$ and $i \in [n]$, define
  \begin{equation*}
    \hat\beta_{a,b,i}
    \coloneqq
    \frac{1}{\sqrt{a!}}
    \E \Bigl( \|X_i - \mu\|_2^a \|X_i - \mu\|_p^b \Bigr)
    + \frac{1}{\sqrt{a!}}
    \E^* \Bigl(
      \|Y_i^*\|_2^a \|Y_i^*\|_p^b
    \Bigr).
  \end{equation*}
  Write $\hat \beta_{a,b} \coloneqq
  \sum_{i=1}^n \hat \beta_{a,b,i} / n$.
  Then, for $a, c \in \N_0$ and $b \geq 0$ with $a + b + c \leq 4$,
  we have
  \begin{align*}
    \E\bigl(\hat\beta_{a,b}\bigr)
    &\leq
    \frac{1 + 4^{a+b}}
    {\sqrt{a!}}
    M_2^a M_p^b, \\
    \E\bigl( \hat\beta_{a,b} \hat\beta_{c,0} \bigr)
    \lor
    \frac{1}{n}
    \sum_{i=1}^n
    \E\bigl(
      \hat\beta_{a,b,i}
      \hat\beta_{c,0,i}
    \bigr)
    &\leq
    \frac{\bigl(1 + 4^{a+b}\bigr) \bigl(1 + 4^c\bigr)
    M_2^{a+c} M_p^b}{\sqrt{a! c!}}.
  \end{align*}
\end{lemma}

\begin{proof}[Proof of Lemma~\ref{lem:beta_hat_bootstrap}]
  Let $I$ be a uniform random variable on $[n]$,
  independent of $(\bX_n, Y_1^*,\ldots,Y_n^*)$,
  and write $U \coloneqq X_I - \mu$ and $V \coloneqq Y_I^*$.
  For $q \in \{2,p\}$, Minkowski's and Jensen's inequalities give
  \begin{equation*}
    \bigvvvert \|X_I-\bar X_n\|_q \bigvvvert_4
    \leq
    \bigvvvert \|X_I-\mu\|_q \bigvvvert_4
    + \bigvvvert \|\bar X_n-\mu\|_q \bigvvvert_4
    \leq 2M_q.
  \end{equation*}
  Under the nonparametric bootstrap,
  $\E\bigl(\|V\|_q^4\bigr)=\E\bigl(\|X_I-\bar X_n\|_q^4\bigr)$.
  Under the multiplier bootstrap,
  independence of the weights and the data gives
  \begin{equation*}
    \E\bigl(\|V\|_q^4\bigr)
    =
    \frac1n \sum_{i=1}^n
    \E\bigl(\chi_i^4\bigr)
    \E\bigl(\|X_i-\bar X_n\|_q^4\bigr)
    \leq 4^4 M_q^4.
  \end{equation*}
  Thus, under either bootstrap scheme,
  $\bigvvvert \|U\|_q \bigvvvert_4 = M_q$ and
  $\bigvvvert \|V\|_q \bigvvvert_4 \leq 4M_q$.
  By definition,
  \begin{equation*}
    \hat\beta_{a,b,I}
    =
    \frac{
      \E\bigl(\|U\|_2^a\|U\|_p^b \mid I\bigr)
      +
      \E\bigl(\|V\|_2^a\|V\|_p^b \mid \bX_n,I\bigr)
    }{\sqrt{a!}}.
  \end{equation*}
  The conditional Jensen inequality, Minkowski's inequality
  and H{\"o}lder's inequality therefore give, for $0<a+b\leq4$,
  \begin{equation*}
    \vvvert \hat\beta_{a,b,I} \vvvert_{4/(a+b)}
    \leq
    \frac{
      \bigvvvert \|U\|_2^a\|U\|_p^b \bigvvvert_{4/(a+b)}
      +
      \bigvvvert \|V\|_2^a\|V\|_p^b \bigvvvert_{4/(a+b)}
    }{\sqrt{a!}}
    \leq
    \frac{1+4^{a+b}}{\sqrt{a!}}M_2^aM_p^b
    \eqqcolon B_{a,b}.
  \end{equation*}
  For $a=b=0$, the same bound holds with $4/0 \coloneqq \infty$,
  since $\hat\beta_{0,0,I}=2$.
  As $\hat\beta_{a,b}
  =\E(\hat\beta_{a,b,I}\mid\bX_n)$,
  the conditional Jensen inequality also gives
  $\vvvert \hat\beta_{a,b} \vvvert_{4/(a+b)} \leq B_{a,b}$.
  In particular,
  $\E(\hat\beta_{a,b})\leq B_{a,b}$.
  Finally, since $(a+b)/4+c/4\leq1$,
  H{\"o}lder's inequality gives
  \begin{align*}
    \frac1n\sum_{i=1}^n
    \E\bigl(\hat\beta_{a,b,i}\hat\beta_{c,0,i}\bigr)
    &=
    \E\bigl(\hat\beta_{a,b,I}\hat\beta_{c,0,I}\bigr)
    \leq B_{a,b}B_{c,0}, \\
    \E\bigl(\hat\beta_{a,b}\hat\beta_{c,0}\bigr)
    &\leq B_{a,b}B_{c,0},
  \end{align*}
  as required.
\end{proof}

\begin{proof}[Proof of Theorem~\ref{thm:bootstrap}]
  We begin with a Borel set approximation
  as in Lemma~\ref{lem:borel},
  comparing the law of $S_n$ with the
  (random) law of $S_n^* \mid \bX_n$.
  For $\alpha \in \N_0^d$ with
  $1 \leq |\alpha| \leq k$ and $i \in [n]$, write
  \begin{align*}
    \hat\delta_{\alpha,i}
    &\coloneqq
    \E\bigl\{(X_i-\mu)^\alpha\bigr\}
    -\E^*\bigl\{(Y_i^*)^\alpha\bigr\},
    &\hat\delta_\alpha
    &\coloneqq \frac1n\sum_{i=1}^n\hat\delta_{\alpha,i}.
  \end{align*}
  If $|\alpha|=1$, then $\hat\delta_{\alpha,i}=0$ for each $i\in[n]$
  because $\E^*(Y_i^*) = 0$ either under
  condition \ref{it:thm_bootstrap_nonparametric}
  or under condition \ref{it:thm_bootstrap_multiplier}.
  Fix $\ell \in \{2, k\}$.
  If $|\alpha|=\ell$, by definition of the bootstrap
  expectation in the nonparametric case and since $\E(\chi_i^\ell)=1$
  for each $i \in [n]$ in the multiplier case, we have
  \begin{align}
    \nonumber
    \hat\delta_\alpha
    &=\frac1n\sum_{i=1}^n
    \Bigl(\E\bigl\{(X_i-\mu)^\alpha\bigr\}
    -(X_i-\bar X_n)^\alpha\Bigr) \\
    \label{eq:bootstrap_delta}
    &=\frac1n\sum_{i=1}^n
    \Bigl(\E\bigl\{(X_i-\mu)^\alpha\bigr\}-(X_i-\mu)^\alpha\Bigr)
    +\frac1n\sum_{i=1}^n
    \Bigl((X_i-\mu)^\alpha-(X_i-\bar X_n)^\alpha\Bigr).
  \end{align}
  When $\mu = 0_d$ and $Y_i^* \coloneqq \chi_i X_i$,
  we have $\E^*\bigl\{(Y_i^*)^\alpha\bigr\} = X_i^\alpha$
  for $|\alpha| \in \{2, k\}$,
  so the second sum in \eqref{eq:bootstrap_delta} vanishes
  and the remaining bounds hold with $\bar X_n$ replaced by $0_d$.
  Let $\varepsilon_1,\ldots,\varepsilon_n$ be independent Rademacher
  variables, independent of the data.
  If $u_\alpha \in \R$ for each $\alpha \in \N_0^d$ with $|\alpha| = \ell$,
  note that $\bigl(\sum_{|\alpha| = \ell}
  u_\alpha^2/\alpha!\bigr)^{1/2} = \|w\|_2
  = \sup_{v: \|v\|_2 = 1} \langle w, v \rangle_2$,
  where $w \coloneqq (u_\alpha/\sqrt{\alpha!}:
  \alpha \in \N_0^d, |\alpha| = \ell)$.
  Therefore, for the first term in
  \eqref{eq:bootstrap_delta},
  by a symmetrization inequality
  \citep[Lemma~2.3.6]{van1996weak},
  Jensen's inequality, the multinomial theorem,
  subadditivity of $x \mapsto x^{2/\ell}$ for $x \geq 0$,
  and Jensen's inequality again,
  \begin{align*}
    &\E\left[
      \Biggl\{
        \sum_{|\alpha|=\ell}\frac1{\alpha!}
        \Biggl(\frac1n\sum_{i=1}^n
          \Bigl\{(X_i-\mu)^\alpha-\E\bigl\{(X_i-\mu)^\alpha\bigr\}\Bigr\}
        \Biggr)^2
      \Biggr\}^{1/2}
    \right]\\
    &\quad\leq
    \frac{2}{n} \,
    \E\left[
      \Biggl\{
        \sum_{|\alpha|=\ell}\frac1{\alpha!}
        \Biggl(\sum_{i=1}^n
          \varepsilon_i (X_i-\mu)^\alpha
        \Biggr)^2
      \Biggr\}^{1/2}
    \right] \\
    &\quad\leq
    \frac{2}{n} \,
    \E\left[
      \Biggl\{
        \sum_{i=1}^n
        \sum_{|\alpha|=\ell}
        \frac{(X_i-\mu)^{2\alpha}}{\alpha!}
      \Biggr\}^{1/2}
    \right]
    =
    \frac{2}{n \sqrt{\ell!}} \,
    \E\left[
      \Biggl\{
        \sum_{i=1}^n
        \|X_i - \mu\|_2^{2 \ell}
      \Biggr\}^{1/2}
    \right] \\
    &\quad\leq
    \frac{2}{n \sqrt{\ell!}} \,
    \E\left[
      \Biggl\{
        \sum_{i=1}^n
        \|X_i - \mu\|_2^{4}
      \Biggr\}^{\ell/4}
    \right]
    \leq
    \frac{2}{n \sqrt{\ell!}} \,
    \Biggl\{
      \sum_{i=1}^n
      \E\bigl\{
        \|X_i - \mu\|_2^{4}
      \bigr\}
    \Biggr\}^{\ell/4}
    \leq
    \frac{2 M_2^\ell}{n^{1-\ell/4} \sqrt{\ell!}}.
  \end{align*}
  For the second term in
  \eqref{eq:bootstrap_delta},
  by the triangle inequality and
  Lemma~\ref{lem:multi_index_sum},
  since $\ell \cdot 2^{\ell-2} / (\ell-1)! \leq 3$,
  as $x \mapsto \sqrt{x}$ is subadditive for $x \geq 0$,
  using Cauchy--Schwarz, the $c_r$-inequality,
  and Jensen's inequality, since $2 \ell - 2 \leq 4$
  and $\E\bigl(\|\bar X_n-\mu\|_2^2\bigr)\leq M_2^2/n$,
  \begin{align*}
    &\E\left[
      \Biggl\{
        \sum_{|\alpha|=\ell}
        \frac{1}{\alpha!}
        \Biggl(
          \frac{1}{n}
          \sum_{i=1}^n
          \Bigl\{(X_i-\mu)^\alpha - (X_i-\bar X_n)^\alpha\Bigr\}
        \Biggr)^2
      \Biggr\}^{1/2}
    \right] \\
    &\quad\leq
    \frac{1}{n}
    \sum_{i=1}^n
    \E\left[
      \Biggl\{
        \sum_{|\alpha|=\ell}
        \frac{1}{\alpha!}
        \Bigl\{(X_i-\mu)^\alpha - (X_i-\bar X_n)^\alpha\Bigr\}^2
      \Biggr\}^{1/2}
    \right] \\
    &\quad\leq
    \frac{1}{n}
    \sum_{i=1}^n
    \E\left[
      \Biggl\{
        \frac{\ell \cdot 2^{\ell-2}}{(\ell-1)!}
        \|\bar X_n - \mu\|_2^2
        \bigl(
          \|X_i - \mu\|_2^{2\ell-2} + \|X_i - \bar X_n\|_2^{2\ell-2}
        \bigr)
      \Biggr\}^{1/2}
    \right] \\
    &\quad\leq
    \frac{\sqrt 3}{n}
    \sum_{i=1}^n
    \E\Bigl\{
      \|\bar X_n - \mu\|_2
      \bigl(
        \|X_i - \mu\|_2^{\ell-1} + \|X_i - \bar X_n\|_2^{\ell-1}
      \bigr)
    \Bigr\} \\
    &\quad\leq
    \bigl\{\E\bigl(
        \|\bar X_n - \mu\|_2^2
    \bigr)\bigr\}^{1/2}
    \frac{\sqrt 3}{n}
    \sum_{i=1}^n
    \biggl(
      \Bigl\{\E\bigl(
          \|X_i - \mu\|_2^{2\ell-2}
      \bigr)\Bigr\}^{1/2}
      + \Bigl\{\E\bigl(
          \|X_i - \bar X_n\|_2^{2\ell-2}
      \bigr)\Bigr\}^{1/2}
    \biggr) \\
    &\quad\leq
    \frac{ \sqrt{3} M_2}{\sqrt n}
    \Biggl\{
      \frac{1 + 2^{\ell-3/2}}{n}
      \sum_{i=1}^n
      \Bigl\{\E\bigl(
          \|X_i - \mu\|_2^{2\ell-2}
      \bigr)\Bigr\}^{1/2}
      \iftoggle{journal}{\\ &\qquad\qquad\qquad}{}
      + 2^{\ell-3/2}
      \Biggl(
        \frac{1}{n}
        \sum_{i=1}^n
        \E\bigl\{
          \|X_i - \mu\|_2^{2\ell-2}
        \bigr\}
      \Biggr)^{1/2}
    \Biggr\} \\
    &\quad\leq
    \frac{\sqrt{3} (1 + 2^{\ell-1/2}) M_2^\ell}{\sqrt n}
    \leq \frac{12 M_2^\ell}{\sqrt n}.
  \end{align*}
  Combining the bounds on the terms in
  \eqref{eq:bootstrap_delta} gives
  \begin{equation*}
    \E\left\{
      \Biggl(
        \sum_{|\alpha|=\ell}
        \frac{\hat\delta_\alpha^2}{\alpha!}
      \Biggr)^{1/2}
    \right\}
    \leq
    \biggl(
      \frac{2}{\sqrt{\ell!} n^{1-\ell/4}}
      +\frac{12}{\sqrt n}
    \biggr)
    M_2^\ell
    \leq
    \frac{14 M_2^\ell}{n^{1-\ell/4}}.
  \end{equation*}
  For $s>0$, define the random variable
  \begin{equation*}
    \hat\Delta(s)
    \coloneqq
    \frac{n}{2}
    \Biggl(
      \sum_{\ell=2}^k s^{2\ell}
      \sum_{|\alpha|=\ell}\frac{\hat\delta_\alpha^2}{\alpha!}
    \Biggr)^{1/2}.
  \end{equation*}
  If $n^{1/4} s M_2 \leq 1$, then
  by subadditivity of $x \mapsto \sqrt x$,
  \begin{equation}
    \label{eq:bootstrap_Delta}
    \E\bigl\{\hat\Delta(s)\bigr\}
    \leq
    \frac{n}{2}
    \sum_{\ell=2}^k
    s^\ell \,
    \E\left\{
      \Biggl(
        \sum_{|\alpha|=\ell}
        \frac{\hat\delta_\alpha^2}{\alpha!}
      \Biggr)^{1/2}
    \right\}
    \leq
    7 \sum_{\ell=2}^k
    \bigl(n^{1/4} s M_2\bigr)^\ell
    \leq
    14 s^2 \sqrt n M_2^2.
  \end{equation}

  For $a \in \N_0$ and $b \geq 0$,
  define $\hat\beta_{a,b,i}$ and $\hat\beta_{a,b}$
  as in Lemma~\ref{lem:beta_hat_bootstrap}.
  By Lemma~\ref{lem:beta_hat_bootstrap},
  for a universal constant $C>0$ and all $s,\eta>0$,
  \begin{equation}
    \label{eq:bootstrap_beta}
    \frac{n s^k}{\eta^\gamma}\E\bigl\{\hat\beta_{k,\gamma}\bigr\}
    \leq
    \frac{n \bigl(1 + 4^{k + \gamma}\bigr)
    M_2^k M_p^\gamma s^k}
    {\sqrt{k!} \eta^\gamma}
    \leq
    \frac{C n M_2^k M_p^\gamma s^k}
    {\eta^\gamma}.
  \end{equation}

  Set $m = 1$, recalling that
  $\E(X_i-\mu)=\E^*(Y_i^*)=0_d$
  for each $i \in [n]$.
  For $\ell \in \{2, k\}$,
  set $k_\ell \coloneqq 1 \land (k - \ell)$,
  so that $k_2 = 1 \land (k - 2) = \I\{k=3\}$
  and, when $k = 3$, $k_3 = 1 \land (k - 3) = 0$.
  For $s > 0$, by Lemma~\ref{lem:bound_Gamma}, and since
  $\sum_{i=1}^n \hat\beta_{\ell,0,i} = n \hat\beta_{\ell,0}$, we have
  \begin{align*}
    \hat\Gamma_\eta(s)
    &\coloneqq
    \sum_{i=1}^n
    \sum_{\ell = m+1}^k
    \frac{(2 s)^{\ell + k_\ell}}
    {2 \eta^{\gamma}}
    \bigl(\hat \beta_{k_\ell,\gamma}
    + \hat \beta_{k_\ell,\gamma,i}\bigr)
    \Biggl(
      \sum_{|\alpha| = \ell}
      \frac{\hat \delta_{\alpha,i}^2}{\alpha!}
    \Biggr)^{1/2} \\
    &\,\leq
    \sum_{i=1}^n
    \sum_{\ell = m+1}^k
    \frac{(2 s)^{\ell + k_\ell}}
    {2 \eta^{\gamma}}
    \bigl(
      \hat \beta_{k_\ell,\gamma}
      \hat \beta_{\ell,0}
      + \hat \beta_{k_\ell,\gamma,i}
      \hat \beta_{\ell,0,i}
    \bigr) \\
    &\,=
    \frac{(2 s)^k}
    {2 \eta^{\gamma}}
    \sum_{i=1}^n
    \biggl\{
      \hat \beta_{0,\gamma} \hat \beta_{k,0}
      + \hat \beta_{0,\gamma,i} \hat \beta_{k,0,i}
      + \I\{k=3\}
      \Bigl(
        \hat \beta_{1,\gamma} \hat \beta_{2,0}
        + \hat \beta_{1,\gamma,i} \hat \beta_{2,0,i}
      \Bigr)
    \biggr\}.
  \end{align*}
  Applying Lemma~\ref{lem:beta_hat_bootstrap},
  since $k \in \{2,3\}$ and $\gamma\in(0,1]$, we have
  \begin{equation}
    \label{eq:bootstrap_Gamma}
    \E\bigl\{\hat\Gamma_\eta(s)\bigr\}
    \leq
    \frac{C n M_2^k M_p^\gamma s^k}
    {\eta^\gamma}.
  \end{equation}
  Now set $s \coloneqq \phi(t) / \eta$ and
  apply Lemma~\ref{lem:borel}
  to the law of $S_n$ and the law of $S_n^*$ conditional
  on~$\bX_n$, using
  \eqref{eq:bootstrap_Delta},
  \eqref{eq:bootstrap_beta} and \eqref{eq:bootstrap_Gamma},
  then take an outer expectation, to see that
  if $n^{1/4} s M_2 \leq 1$,
  \begin{align}
    \nonumber
    &\E^\star \biggl(
      \sup_{A \in \cB(\R^d)}
      \Bigl\{
        \P(S_n \in A)
        - \P\big( S_n^* \in A_p^{3\eta} \mid \bX_n \big)
      \Bigr\}
    \biggr) \\
    \nonumber
    &\qquad\leq
    2 e^{-t^2 / 2}
    + \E\Biggl\{
      \frac{n \hat\beta_{k,\gamma}\phi(t)^k}{\eta^{k+\gamma}}
      + \hat \Delta\biggl(\frac{\phi(t)}{\eta}\biggr)
      + \hat \Gamma_\eta\biggl(\frac{\phi(t)}{\eta}\biggr)
    \Biggr\} \\
    \label{eq:bootstrap_outer_expectation}
    &\qquad\leq
    2 e^{-t^2 / 2}
    + \frac{14 \sqrt n M_2^2 \phi(t)^2}{\eta^2}
    + \frac{C n M_2^k M_p^\gamma \phi(t)^k}
    {\eta^{k+\gamma}}.
  \end{align}
  If $n^{1/4} s M_2 > 1$,
  then the second term on the right-hand side above is
  $14 (n^{1/4} s M_2)^2 > 1$, so
  this final inequality holds trivially.

  We now prove the coupling bound. Let $(\Omega, \cH, \P)$
  be a probability space
  supporting $\bX_n' \coloneqq (X_i': i \in [n])$
  and $S_n'$, where $\bX_n'$
  has the same law as $\bX_n$,
  $S_n'$ is distributed as $S_n$, and
  $S_n'$ is independent of $\bX_n'$.
  We may assume that $(\Omega, \cH, \P)$ is rich enough
  to also support
  a standard uniform random variable independent
  of $(\bX_n', S_n')$.
  Define $\cH'$ to be the $\sigma$-algebra generated by $\bX_n'$.
  Take a regular conditional distribution
  $F(\cdot \mid \bX_n)$ of $S_n^*$ given $\bX_n$.
  Evaluating the same conditional distribution at $\bX_n'$
  gives an $\cH'$-measurable regular conditional distribution
  $F(\cdot \mid \bX_n')$ on $(\Omega,\cH,\P)$.
  Also, $\P(S_n \in A)=\P(S_n' \in A \mid \cH')$
  for each $A \in \cB(\R^d)$ by independence.
  Therefore,
  \eqref{eq:bootstrap_outer_expectation}
  holds with $S_n$ replaced by $S_n'$
  and $\P(S_n^* \in \cdot \mid \bX_n)$ replaced with
  $F(\cdot \mid \bX_n')$.
  The result follows by a conditional
  version of Strassen's theorem
  \citep[Lemma~SA.1]{cattaneo2025yurinskii},
  replacing $\eta$ with $\eta/3$.
\end{proof}

\begin{proof}[Details for \eqref{eq:bootstrap_rate}]
  Under the conditions imposed, Minkowski's and
  H{\"o}lder's inequalities yield
  $M_2 \leq \sqrt d M$ and $M_p \leq d^{1/(p \land r)} M$.
  After rescaling the bound from Theorem~\ref{thm:bootstrap}
  by $\sqrt n$ and taking $k = 3$ and $\gamma = 1$,
  it suffices to choose
  $t_n$ and $\eta_n$ such that
  \begin{equation*}
    e^{-t_n^2 / 2}
    + \frac{d \phi(t_n)^2 M^2}{\sqrt n\eta_n^2}
    + \frac{d^{3/2 + 1 / (p \land r)} \phi(t_n)^3 M^4}{n \eta_n^4}
    \to 0.
  \end{equation*}
  To this end, we take $t_n \coloneqq \sqrt R_n$
  where $R_n \to \infty$ and set
  \begin{align*}
    \eta_n
    &\coloneqq
    \frac{t_n d^{1/2} \phi(t_n) M}{n^{1/4}}
    + \frac{t_n d^{3/8 + 1 / \{4 (p \land r)\}}
    \phi(t_n)^{3 / 4} M}{n^{1/4}} \\
    &\,\lesssim
    \frac{R_n d^{1/2} \psi^{3/4} M}{n^{1/4}}
    \Bigl(
      \psi^{1/4}
      + d^{1 / \{4 (p \land r)\} - 1/8}
    \Bigr),
  \end{align*}
  as $\phi(t) \leq (1 + t) \psi$ for all $t \geq 0$.

  We first consider $p \in [1, \infty)$,
  so that $\psi \asymp d^{1/p}$.
  Observe that if $r > p$, then
  $1 / \{4 (p \land r)\} - 1/8 = 1 / (4 p) - 1/8 \leq 1/(4 p)$.
  Otherwise, if $r \leq p$, then since $r \geq 2$, we have that
  $1 / \{4 (p \land r)\} - 1/8 = 1 / (4 r) - 1/8 \leq 0$.
  We deduce that
  \begin{equation*}
    \eta_n
    \lesssim
    \frac{R_n d^{1/2 + 1/p} M}{n^{1/4}}.
  \end{equation*}
  Now, if $p = \infty$, then $\psi \asymp \sqrt{\log 2 d}$.
  Since $r \geq 2$, we have
  $1 / \{4 (p \land r)\} - 1/8 = 1 / (4 r) - 1/8 \leq 0$, so
  \begin{equation*}
    \eta_n
    \lesssim
    \frac{R_n (d \log 2 d)^{1/2} M}{n^{1/4}}.
  \end{equation*}

  If instead $M_2 \leq M$, then
  $\|x\|_p \leq d^{1/p}\|x\|_2$ gives
  $M_p \leq d^{1/p}M$, with $1/\infty=0$.
  Repeating the previous argument with these bounds gives
  \begin{equation*}
    \eta_n \lesssim
    \frac{R_n \psi^{3/4} M}{n^{1/4}}
    \bigl(\psi^{1/4} + d^{1/(4p)}\bigr)
    \lesssim \frac{R_n M \psi}{n^{1/4}},
  \end{equation*}
  since $\psi \gtrsim d^{1/p}$.
  Substituting the value of $\psi$ gives
  \begin{align*}
    \frac{\|S_n' - S_n^{*\prime}\|_p}{\sqrt n}
    &\lesssim_\P
    \frac{R_n d^{1/p}}{n^{1/4}},
    &\frac{\|S_n' - S_n^{*\prime}\|_\infty}{\sqrt n}
    &\lesssim_\P
    \frac{R_n \sqrt{\log 2 d}}{n^{1/4}}.
    \iftoggle{journal}{}{\qedhere}
  \end{align*}
\end{proof}

\begin{lemma}[Refined bootstrap moment bounds]%
  \label{lem:beta_hat_refined}
  Grant the setup of Theorem~\ref{thm:refined_bootstrap} and
  let $\Sigma^{+/2}$, $\hat\Sigma_n^{+/2}$
  and $\cA_n$ be as defined in the proof
  of Theorem~\ref{thm:refined_bootstrap}.
  For $a \in \N_0$ and $b \geq 0$,
  and for each $i \in [n]$, let
  \begin{align*}
    \hat\beta_{a,b,i}
    &\coloneqq
    \frac{1}{\sqrt{a!}}
    \E \Bigl(
      \bigl\|\Sigma^{+/2} (X_i - \mu)\bigr\|_2^a
      \bigl\|\Sigma^{+/2} (X_i - \mu)\bigr\|_p^b
    \Bigr) \\
    &\qquad+
    \frac{1}{\sqrt{a!}}
    \E^* \Bigl(
      \bigl\|\hat\Sigma_n^{+/2} (X_i^* - \bar X_n)\bigr\|_2^a
      \bigl\|\hat\Sigma_n^{+/2} (X_i^* - \bar X_n)\bigr\|_p^b
    \Bigr).
  \end{align*}
  Write $\hat\beta_{a,b}\coloneqq
  \sum_{i=1}^n\hat\beta_{a,b,i}/n$.
  For $a,c\in\N_0$ and $b\geq0$ with $a+b+c\leq6$,
  we have
  \begin{align*}
    \E\bigl\{\hat\beta_{a,b}\I(\cA_n)\bigr\}
    &\leq
    \frac{1+2^{3(a+b)/2}}{\sqrt{a!}}
    \frac{d^{b a_p}M_2^{a+b}}
    {\lambda^{(a+b)/2}}, \\
    \E\bigl\{
      \hat\beta_{a,b}\hat\beta_{c,0}\I(\cA_n)
    \bigr\}
    \lor
    \frac1n\sum_{i=1}^n
    \E\bigl\{
      \hat\beta_{a,b,i}\hat\beta_{c,0,i}\I(\cA_n)
    \bigr\}
    &\leq
    \frac{
      \bigl(1+2^{3(a+b)/2}\bigr)
      \bigl(1+2^{3c/2}\bigr)
    }{\sqrt{a!c!}}
    \iftoggle{journal}{\\ &\qquad\qquad\times}{}
    \frac{d^{b a_p}M_2^{a+b+c}}
    {\lambda^{(a+b+c)/2}}.
  \end{align*}
\end{lemma}

\begin{proof}[Proof of Lemma~\ref{lem:beta_hat_refined}]
  Let $I$ be uniform on $[n]$, independent of
  $(\bX_n, X_1^*, \ldots, X_n^*)$. Write
  $U\coloneqq X_I-\mu$ and $V\coloneqq X_I^*-\bar X_n$.
  By Minkowski's and Jensen's inequalities,
  $\bigvvvert \|U\|_2 \bigvvvert_6 = M_2$ and
  \begin{equation*}
    \bigvvvert \|V\|_2 \bigvvvert_6
    = \bigvvvert \|X_I-\bar X_n\|_2 \bigvvvert_6
    \leq
    \bigvvvert \|X_I-\mu\|_2 \bigvvvert_6
    + \bigvvvert \|\bar X_n-\mu\|_2 \bigvvvert_6
    \leq 2M_2.
  \end{equation*}
  Further,
  $\|\Sigma^{+/2}\|_{2\to2}=\lambda^{-1/2}$ and,
  almost surely on $\cA_n$,
  $\|\hat\Sigma_n^{+/2}\|_{2\to2}\leq(2/\lambda)^{1/2}$,
  as shown in the proof of Theorem~\ref{thm:refined_bootstrap} below.
  Since $\|x\|_p\leq d^{a_p}\|x\|_2$ for $x\in\R^d$,
  for $0<a+b\leq6$ we have
  \begin{equation*}
    \hat\beta_{a,b,I}\I(\cA_n)
    \leq
    \frac{d^{b a_p}}{\sqrt{a!}\lambda^{(a+b)/2}}
    \Bigl\{
      \E\bigl(\|U\|_2^{a+b}\mid I\bigr)
      + 2^{(a+b)/2}
      \, \E\bigl(\|V\|_2^{a+b}\mid\bX_n,I\bigr)
    \Bigr\}\I(\cA_n).
  \end{equation*}
  By the conditional Jensen and Minkowski inequalities,
  \begin{align*}
    \bigvvvert
    \hat\beta_{a,b,I}\I(\cA_n)
    \bigvvvert_{6/(a+b)}
    &\leq
    \frac{d^{b a_p}}{\sqrt{a!}\lambda^{(a+b)/2}}
    \Bigl\{
      \bigvvvert \|U\|_2 \bigvvvert_6^{a+b}
      + 2^{(a+b)/2}
      \bigvvvert \|V\|_2 \bigvvvert_6^{a+b}
    \Bigr\} \\
    &\leq
    \frac{1+2^{3(a+b)/2}}{\sqrt{a!}}
    \frac{d^{b a_p}M_2^{a+b}}{\lambda^{(a+b)/2}}
    \eqqcolon B_{a,b}.
  \end{align*}
  For $a=b=0$, the same bound also holds with
  $6/0\coloneqq\infty$, because $\hat\beta_{0,0,I}=2$.
  As $\hat\beta_{a,b}\I(\cA_n) =
  \E\bigl( \hat\beta_{a,b,I}\I(\cA_n)\mid\bX_n \bigr)$,
  the conditional Jensen inequality also gives
  \[
    \bigvvvert
    \hat\beta_{a,b}\I(\cA_n)
    \bigvvvert_{6/(a+b)}
    \leq B_{a,b}.
  \]
  In particular,
  $\E\{\hat\beta_{a,b}\I(\cA_n)\}\leq B_{a,b}$.
  Finally, as $(a+b)/6+c/6\leq1$, H{\"o}lder's inequality gives
  \begin{align*}
    \frac1n\sum_{i=1}^n
    \E\bigl\{
      \hat\beta_{a,b,i}\hat\beta_{c,0,i}\I(\cA_n)
    \bigr\}
    &=
    \E\bigl\{
      \hat\beta_{a,b,I}\hat\beta_{c,0,I}\I(\cA_n)
    \bigr\}
    \leq B_{a,b}B_{c,0}, \\
    \E\bigl\{
      \hat\beta_{a,b}\hat\beta_{c,0}\I(\cA_n)
    \bigr\}
    &\leq B_{a,b}B_{c,0},
  \end{align*}
  as required.
\end{proof}

\begin{lemma}[A Marcinkiewicz--Zygmund-type inequality]
  \label{lem:mz}
  Fix $n, d \in \N$ and $p \in [2, \infty)$.
  Let $X_1, \ldots, X_n$ be independent $\R^d$-valued random vectors
  with $\E(X_i) = 0_d$ and $\E\bigl(\|X_i\|_2^p\bigr) < \infty$.
  Then
  \begin{equation*}
    \E \Biggl(
      \biggl\|
      \frac{1}{n}
      \sum_{i=1}^n
      X_i
      \biggr\|_2^p
    \Biggr)
    \leq \left(\frac{p-1}{n}\right)^{p/2}
    \frac{1}{n}
    \sum_{i=1}^n
    \E\bigl(\|X_i\|_2^p\bigr).
  \end{equation*}
\end{lemma}

\begin{proof}
  For each fixed $u\in\R^d$ and with $p > 2$,
  \citet[Theorem~2.1]{rio2009moment}
  and Jensen's inequality give
  \begin{align*}
    \E\Biggl(
      \biggl|
      \sum_{i=1}^n
      u^\T X_i
      \biggr|^p
    \Biggr)
    &\leq
    (p-1)^{p/2}
    \Biggl(
      \sum_{i=1}^n
      \Bigl\{
        \E\bigl(|u^\T X_i|^p\bigr)
      \Bigr\}^{2/p}
    \Biggr)^{p/2} \\
    &\leq
    (p-1)^{p/2}
    n^{p/2-1}
    \sum_{i=1}^n
    \E\bigl(|u^\T X_i|^p\bigr).
  \end{align*}
  This bound also holds for $p = 2$ by independence
  and the zero-mean property of the summands.
  Let $Z \sim N(0_d, I_d)$ be independent of $X_1, \ldots, X_n$
  and define $c_p \coloneqq \E\bigl(|Z_1|^p\bigr) \in (0,\infty)$.
  Since $\E\bigl(|Z^\T x|^p\bigr) = c_p \|x\|_2^p$
  for every $x \in \R^d$,
  the previous bound and Fubini's theorem yield
  \begin{equation*}
    c_p \,
    \E\Biggl(
      \biggl\|
      \sum_{i=1}^n
      X_i
      \biggr\|_2^p
    \Biggr)
    \leq
    c_p (p-1)^{p/2}
    n^{p/2-1}
    \sum_{i=1}^n \E\|X_i\|_2^p.
  \end{equation*}
  Dividing by $c_p n^p$ proves the result.
\end{proof}

\begin{proof}[Proof of Theorem~\ref{thm:refined_bootstrap}]
  Let $\hat S_n \coloneqq \Sigma^{+/2} S_n$.
  Define $\hat \Sigma_n = \sum_{i=1}^n X_i X_i^\T / n - \bar X_n \bar X_n^\T$,
  and set $\hat S_n^* \coloneqq \hat \Sigma_n^{+/2} S_n^*$.
  We establish a coupling between the law of $\hat S_n$ and the conditional
  distribution of $\hat S_n^*$, given $\bX_n$.
  Firstly, we provide a bound for the probability
  that the empirical covariance matrix is close to singular,
  in an appropriate sense.
  To that end,
  write $\lambda_{1} \leq \cdots \leq \lambda_{d}$ and
  $\hat\lambda_{n,1} \leq \cdots \leq \hat\lambda_{n,d}$
  for the ordered eigenvalues of $\Sigma$ and
  $\hat\Sigma_{n}$ respectively,
  and let $j_0 \in [d]$ be the index of the first
  non-zero eigenvalue of
  $\Sigma$; that is, $\lambda_{j_0} = \lambda > 0$ and
  $\lambda_{j'} = 0$ whenever $j' < j_0$.
  Let $\hat\lambda_n \coloneqq \hat\lambda_{n,j_0}$.
  As $\Sigma$ is positive semidefinite,
  by Jensen's inequality,
  $\lambda \leq \tr \Sigma
  = \frac1n\sum_{i=1}^n
  \E\bigl\{\|X_i - \mu\|_2^2\bigr\} \leq M_2^2$.
  Define the event
  \begin{equation*}
    \cA_n
    \coloneqq
    \biggl\{
      \frac{\lambda}{2}
      \leq
      \hat \lambda_{n}
      \leq
      \frac{3\lambda}{2}
    \biggr\}
    =
    \biggl\{
      \bigl|\hat \lambda_{n} - \lambda\bigr|
      \leq
      \frac{\lambda}{2}
    \biggr\}.
  \end{equation*}

  We claim that $\hat\lambda_{n,j'} = 0$ almost surely
  whenever $j' < j_0$.
  To see this, suppose that
  $v \in \R^d$ satisfies $\Sigma v = 0_d$. Then
  $\E\bigl\{\{v^\T (X_i - \mu)\}^2\bigr\}
  = v^\T \Sigma v = 0$,
  so $v^\T (X_i - \mu) = 0$ almost surely
  for each $i \in [n]$.
  Therefore, $v^\T \bar X_{n} = v^\T \mu$ and
  $v^\T \hat\Sigma_{n} v
  = \sum_{i=1}^{n} v^\T X_{i} X_{i}^\T v / n
  - v^\T \mu \mu^\T v = 0$,
  so $\hat\Sigma_n v = 0_d$ almost surely.
  Applying this to a basis of $\ker \Sigma$,
  we deduce that $\ker \Sigma \subseteq \ker \hat\Sigma_n$
  almost surely, and so
  $j' \leq j_0 - 1 = \dim \ker \Sigma
  \leq \dim \ker \hat\Sigma_n$
  and $\hat\lambda_{n,j'} = 0$ almost surely.
  On $\cA_n$ we have
  $\hat\lambda_{n,j_0} \geq \lambda / 2 > 0$,
  so $\dim \ker \hat\Sigma_n \leq j_0 - 1$.
  Hence, almost surely on $\cA_n$, we have
  $\ker \hat\Sigma_n = \ker \Sigma$, and
  $\hat\lambda_n \geq \lambda / 2$ is the smallest
  non-zero eigenvalue of $\hat\Sigma_n$.
  Next, by Weyl's inequality
  \citep[see, for example,][Corollary~III.2.6]{bhatia1997matrix},
  \begin{align*}
    &\bigvvvert
    \hat \lambda_{n} - \lambda
    \bigvvvert_{3}
    =
    \bigvvvert
    \hat \lambda_{n,j_0} - \lambda_{j_0}
    \bigvvvert_{3}
    \leq
    \Bigvvvert
    \bigl\| \hat\Sigma_{n} - \Sigma \bigr\|_{2 \to 2}
    \Bigvvvert_{3} \\
    &\qquad=
    \Biggvvvert
    \biggl\|
    \frac{1}{n} \sum_{i=1}^{n} X_{i} X_{i}^\T
    - \bar X_{n} \bar X_{n}^\T
    - \Sigma
    \biggr\|_{2 \to 2}
    \Biggvvvert_{3} \\
    &\qquad=
    \Biggvvvert
    \biggl\|
    \frac{1}{n} \sum_{i=1}^{n}
    (X_{i} - \mu) (X_{i} - \mu)^\T
    - \Sigma
    - (\bar X_{n} - \mu) (\bar X_{n} - \mu)^\T
    \biggr\|_{2 \to 2}
    \Biggvvvert_{3} \\
    &\qquad\leq
    \Biggvvvert
    \biggl\|
    \frac{1}{n} \sum_{i=1}^{n}
    (X_{i} - \mu) (X_{i} - \mu)^\T
    - \Sigma
    \biggr\|_{\Fr,2}
    \Biggvvvert_{3}
    +
    \Bigvvvert
    \bigl\|
    \bar X_{n} - \mu
    \bigr\|_{2}
    \Bigvvvert_{6}^2.
  \end{align*}
  By Lemma~\ref{lem:mz} with $\|A\|_{\Fr, 2} = \|\vecmat(A)\|_2$,
  then applying Minkowski's and Jensen's inequalities,
  \begin{align*}
    &\Biggvvvert
    \biggl\|
    \frac1n\sum_{i=1}^n
    (X_i-\mu)(X_i-\mu)^\T-\Sigma
    \biggr\|_{\Fr,2}
    \Biggvvvert_3^3 \\
    &\qquad\leq
    \biggl(\frac2n\biggr)^{3/2}
    \frac1n\sum_{i=1}^n
    \E\Bigl\{
      \bigl\|(X_i-\mu)(X_i-\mu)^\T-\Sigma\bigr\|_{\Fr,2}^3
    \Bigr\} \\
    &\qquad\leq
    \frac{2^{9/2}}{n^{3/2}}
    \frac1n\sum_{i=1}^n
    \E\bigl(\|X_i-\mu\|_2^6\bigr)
    =
    \frac{2^{9/2}M_2^6}{n^{3/2}}.
  \end{align*}
  Similarly,
  \begin{equation*}
    \Bigvvvert
    \|\bar X_n - \mu\|_2
    \Bigvvvert_6^6
    \leq
    \biggl(\frac5n\biggr)^3
    \frac1n \sum_{i=1}^n
    \E\bigl(\|X_i - \mu\|_2^6\bigr)
    =
    \frac{125 M_2^6}{n^3}.
  \end{equation*}
  We deduce that
  \begin{align}
    \label{eq:refined_V1}
    \bigvvvert
    \hat \lambda_{n} - \lambda
    \bigvvvert_{3}
    &\leq
    \Bigvvvert
    \bigl\| \hat\Sigma_{n} - \Sigma \bigr\|_{2 \to 2}
    \Bigvvvert_{3}
    \leq
    \frac{2^{3/2} M_2^2}{n^{1/2}}
    + \frac{5 M_2^2}{n}
    \leq
    \frac{8 M_2^2}{\sqrt n}, \\
    \label{eq:refined_V2}
    \Bigvvvert
    \bigl\| \bar X_{n} - \mu \bigr\|_{2}
    \Bigvvvert_{6}
    &\leq
    \frac{5^{1/2} M_2}{\sqrt n}.
  \end{align}
  By Markov's inequality,
  \begin{equation}
    \label{eq:refined_P}
    \P(\cA_n^\c)
    =
    \P\biggl(
      \bigl|\hat \lambda_{n} - \lambda\bigr|^{3}
      > \frac{\lambda^{3}}{8}
    \biggr)
    \leq
    \frac{8}{\lambda^{3}}
    \bigvvvert
    \hat \lambda_{n} - \lambda
    \bigvvvert_{3}^{3}
    \leq
    \frac{4096 M_2^6}{n^{3/2} \lambda^{3}}.
  \end{equation}

  Next, let $\Sigma^+$ and $\hat \Sigma_n^+$
  be the Moore--Penrose pseudoinverses of
  $\Sigma$ and $\hat \Sigma_n$, respectively
  \citep[Section~6.7]{friedberg2003linear}.
  We write $\Sigma^{+/2} \coloneqq (\Sigma^+)^{1/2}$
  for brevity; likewise for $\hat\Sigma_n^{+/2}$.
  For $\alpha \in \N_0^d$ with $|\alpha| \in [k]$
  and $i \in [n]$,
  define
  \begin{align*}
    \hat\delta_{\alpha,i}
    &\coloneqq
    \E \Bigl\{
      \bigl(
        \Sigma^{+/2}
        (X_i - \mu)
      \bigr)^\alpha
    \Bigr\}
    - \E^* \Bigl\{
      \bigl(
        \hat \Sigma_n^{+/2}
        (X_i^* - \bar X_n)
      \bigr)^\alpha
    \Bigr\} \\
    &\,=
    \E \biggl\{
      \Bigl\{
        \Sigma^{+/2}
        (X_i - \mu)
      \Bigr\}^\alpha
    \biggr\}
    - \frac{1}{n}
    \sum_{j=1}^n
    \Bigl\{
      \hat \Sigma_n^{+/2}
      (X_j - \bar X_n)
    \Bigr\}^\alpha, \\
    \hat \delta_{\alpha}
    &\coloneqq
    \frac{1}{n}
    \sum_{i=1}^n
    \hat\delta_{\alpha,i}.
  \end{align*}
  Note that for each $i \in [n]$,
  \begin{equation*}
    \frac{1}{n}
    \sum_{j=1}^n
    \hat \Sigma_n^{+/2}
    (X_j - \bar X_n)
    = 0_d
    = \E \Bigl\{
      \Sigma^{+/2}
      (X_i - \mu)
    \Bigr\},
  \end{equation*}
  so $\hat\delta_{\alpha,i} = 0$ if $|\alpha| = 1$.
  Further,
  \begin{align*}
    &\E\Bigl\{
      \Sigma^{+/2}
      (X_i - \mu)
      (X_i - \mu)^\T
      \Sigma^{+/2}
    \Bigr\}
    - \frac{1}{n}
    \sum_{i=1}^{n}
    \hat\Sigma_n^{+/2}
    (X_i - \bar X_n)
    (X_i - \bar X_n)^\T
    \hat\Sigma_n^{+/2} \\
    &\qquad=
    \Sigma^{+/2}
    \Sigma
    \Sigma^{+/2}
    -
    \hat\Sigma_n^{+/2}
    \hat \Sigma_n
    \hat\Sigma_n^{+/2}.
  \end{align*}
  The linear map
  $v \mapsto \Sigma^{+/2} \Sigma^{1/2} v$
  is equal to the orthogonal projection from $\R^d$ onto the orthogonal
  complement of $\ker \Sigma$
  \citep[Chapter~6]{friedberg2003linear}.
  Likewise, $v \mapsto \hat\Sigma_n^{+/2} \hat\Sigma_n^{1/2} v$
  is the projection onto the complement of $\ker \hat\Sigma_n$.
  Since $\ker \Sigma = \ker \hat \Sigma_n$
  almost surely on $\cA_n$,
  we deduce that if $|\alpha| = 2$, then
  $\hat\delta_{\alpha,i} = 0$ for all $i \in [n]$
  almost surely on $\cA_n$.
  We now consider $|\alpha| = \ell$ where $3 \leq \ell \leq k$.
  Write the decomposition
  \begin{align}
    \label{eq:whitened_bootstrap_1}
    \hat\delta_{\alpha}
    &=
    \frac{1}{n} \sum_{i=1}^{n}
    \Bigl\{
      \Sigma^{+/2}
      (X_i - \bar X_n)
    \Bigr\}^\alpha
    -
    \frac{1}{n} \sum_{i=1}^{n}
    \Bigl\{
      \hat\Sigma_n^{+/2}
      (X_i - \bar X_n)
    \Bigr\}^\alpha \\
    &\qquad+
    \label{eq:whitened_bootstrap_2}
    \frac{1}{n} \sum_{i=1}^{n}
    \Bigl\{
      \Sigma^{+/2}
      (X_i - \mu)
    \Bigr\}^\alpha
    -
    \frac{1}{n} \sum_{i=1}^{n}
    \Bigl\{
      \Sigma^{+/2}
      (X_i - \bar X_n)
    \Bigr\}^\alpha \\
    \label{eq:whitened_bootstrap_3}
    &\qquad+
    \frac{1}{n} \sum_{i=1}^{n}
    \E\biggl\{
      \Bigl\{
        \Sigma^{+/2}
        (X_i - \mu)
      \Bigr\}^\alpha
    \biggr\}
    -
    \frac{1}{n} \sum_{i=1}^{n}
    \Bigl\{
      \Sigma^{+/2}
      (X_i - \mu)
    \Bigr\}^\alpha.
  \end{align}
  We bound the first absolute moments of these three terms.
  By Minkowski's inequality and Jensen's inequality, for
  $q \in (0, 6]$,
  \begin{align}
    \label{eq:refined_norms}
    \frac{1}{n} \sum_{i=1}^n
    \E\bigl(\|X_i - \mu\|_2^q\bigr)
    &\leq
    \frac1n\sum_{i=1}^n
    \bigl\{\E\bigl(\|X_i-\mu\|_2^6\bigr)\bigr\}^{q/6}
    \leq M_2^q, \\
    \nonumber
    \frac{1}{n}
    \sum_{i=1}^n
    \E\bigl(\|X_i - \bar X_n\|_2^q\bigr)
    &\leq
    \frac1n\sum_{i=1}^n
    \bigl\{\E\bigl(\|X_i-\bar X_n\|_2^6\bigr)\bigr\}^{q/6}
    \leq (2 M_2)^q.
  \end{align}
  For \eqref{eq:whitened_bootstrap_1},
  since $\ker \hat\Sigma_n = \ker \Sigma$
  almost surely on $\cA_n$,
  we have
  $\hat\Sigma_n^{+/2} - \Sigma^{+/2}
  = \hat\Sigma_n^{+/2}
  \bigl( \Sigma^{1/2} - \hat\Sigma_n^{1/2} \bigr)
  \Sigma^{+/2}$ almost surely on $\cA_n$.
  Note that
  $\|\Sigma^{1/2} v\|_2 \geq \lambda^{1/2} \|v\|_2$
  for all $v$ in the orthogonal complement of $\ker \Sigma$.
  Thus, by applying
  \citet[Problem~X.5.5]{bhatia1997matrix}
  on the orthogonal complement of $\ker \Sigma$,
  \begin{align*}
    \bigl\|
    \hat\Sigma_n^{+/2} - \Sigma^{+/2}
    \bigr\|_{2 \to 2}
    \I(\cA_n)
    &\leq
    \bigl\| \hat\Sigma_n^{+/2} \bigr\|_{2 \to 2}
    \bigl\| \hat\Sigma_n^{1/2} - \Sigma^{1/2} \bigr\|_{2 \to 2}
    \bigl\| \Sigma^{+/2} \bigr\|_{2 \to 2}
    \I(\cA_n) \\
    &\leq
    \frac{2^{1/2}}{\lambda^{3/2}}
    \bigl\| \hat\Sigma_n - \Sigma \bigr\|_{2 \to 2}.
  \end{align*}
  Also,
  $\bigl\| \hat\Sigma_n^{+/2} - \Sigma^{+/2} \bigr\|_{2 \to 2}
  \I(\cA_n)
  \leq (1 + \sqrt 2) / \lambda^{1/2}
  \leq 2^{3/2} / \lambda^{1/2}$
  by the triangle inequality.
  For $\ell \in \{3,4,5\}$,
  set $\theta_\ell \coloneqq 2 - \ell/3 \in [1/3, 1]$.
  Interpolating these two bounds gives
  \begin{align}
    \nonumber
    \bigl\|\hat\Sigma_n^{+/2}-\Sigma^{+/2}\bigr\|_{2\to2}
    \I(\cA_n)
    &\leq
    \biggl(
      \frac{2^{1/2}}{\lambda^{3/2}}
      \bigl\| \hat\Sigma_n - \Sigma \bigr\|_{2 \to 2}
    \biggr)^{\theta_\ell}
    \biggl(
      \frac{2^{3/2}}{\lambda^{1/2}}
    \biggr)^{1 - \theta_\ell} \\
    \label{eq:bootstrap_refined_interp}
    &\leq
    \frac{2^{3/2-\theta_\ell}}{\lambda^{1/2+\theta_\ell}}
    \bigl\|\hat\Sigma_n-\Sigma\bigr\|_{2\to2}^{\theta_\ell}.
  \end{align}
  Write $v_\ell(x)\coloneqq\bigl(x^\alpha/\sqrt{\alpha!}:
  \alpha \in \N_0^d, |\alpha|=\ell\bigr)$,
  so that $\|v_\ell(x)\|_2=\|x\|_2^\ell/\sqrt{\ell!}$
  by the multinomial theorem.
  By the triangle inequality and
  Lemma~\ref{lem:multi_index_sum},
  since $\ell \cdot 2^{\ell-2} \cdot (1 + 2^{\ell-1})
  / (\ell - 1)! \leq 29$,
  using \eqref{eq:bootstrap_refined_interp},
  as $6 \cdot 2^{3/2 - \theta_\ell} \leq 14$,
  applying the generalized H{\"o}lder inequality,
  by \eqref{eq:refined_V1} and
  \eqref{eq:refined_norms},
  and since $M_2^2 / \lambda \geq 1$,
  \begin{align*}
    &\E\Biggl\{
      \biggl\|
      \frac{1}{n}
      \sum_{i=1}^n
      \Bigl\{
        v_\ell\bigl(\Sigma^{+/2}(X_i-\bar X_n)\bigr)
        - v_\ell\bigl(\hat\Sigma_n^{+/2}(X_i-\bar X_n)\bigr)
      \Bigr\}
      \biggr\|_2
      \I(\cA_n)
    \Biggr\} \\
    &\quad\leq
    \frac{1}{n}
    \sum_{i=1}^n
    \E\biggl\{
      \Bigl\|
      v_\ell\bigl(\Sigma^{+/2}(X_i-\bar X_n)\bigr)
      - v_\ell\bigl(\hat\Sigma_n^{+/2}(X_i-\bar X_n)\bigr)
      \Bigr\|_2
      \I(\cA_n)
    \biggr\} \\
    &\quad\leq
    \frac{1}{n}
    \sum_{i=1}^n
    \E\Biggl\{
      \biggl(
        \frac{\ell \cdot 2^{\ell-2}}{(\ell-1)!}
        \bigl\|
        \bigl(\hat\Sigma_n^{+/2} - \Sigma^{+/2}\bigr)
        (X_i-\bar X_n)
        \bigr\|_2^2 \\
        &\qquad\qquad\qquad\times
        \Bigl(
          \bigl\|\hat\Sigma_n^{+/2}(X_i-\bar X_n)\bigr\|_2^{2\ell-2}
          + \bigl\|\Sigma^{+/2}(X_i-\bar X_n)\bigr\|_2^{2\ell-2}
        \Bigr)
      \biggr)^{1/2}
      \I(\cA_n)
    \Biggr\} \\
    &\quad\leq
    \sqrt{\frac{\ell \cdot 2^{\ell-2} \cdot (1 + 2^{\ell-1})}
    {(\ell-1)!}}
    \frac{\lambda^{1/2 - \ell/2}}{n}
    \sum_{i=1}^n
    \E\Bigl\{
      \bigl\|\hat\Sigma_n^{+/2} - \Sigma^{+/2}\bigr\|_{2 \to 2}
      \|X_i-\bar X_n\|_2^{\ell}
      \I(\cA_n)
    \Bigr\} \\
    &\quad\leq
    \frac{6 \lambda^{1/2 - \ell/2}}{n}
    \frac{2^{3/2-\theta_\ell}}{\lambda^{1/2+\theta_\ell}}
    \sum_{i=1}^n
    \E\Bigl\{
      \bigl\|\hat\Sigma_n-\Sigma\bigr\|_{2\to2}^{\theta_\ell}
      \|X_i-\bar X_n\|_2^{\ell}
      \I(\cA_n)
    \Bigr\} \\
    &\quad\leq
    \frac{14}{\lambda^{\ell/2 + \theta_\ell}}
    \biggl\{\E\Bigl(
        \bigl\|\hat\Sigma_n-\Sigma\bigr\|_{2\to2}^3
    \Bigr)\biggr\}^{\theta_\ell / 3}
    \frac{1}{n}
    \sum_{i=1}^n
    \biggl\{\E\Bigl(
        \|X_i-\bar X_n\|_2^{6}
    \Bigr)\biggr\}^{\ell/6} \\
    &\quad\leq
    \frac{14}{\lambda^{\ell/2 + \theta_\ell}}
    \biggl(
      \frac{8 M_2^2}{\sqrt n}
    \biggr)^{\theta_\ell}
    \bigl(2 M_2\bigr)^\ell
    \leq
    \frac{896}{n^{\theta_\ell/2}}
    \biggl(
      \frac{M_2^2}{\lambda}
    \biggr)^{\ell/2 + \theta_\ell}
    \leq
    \frac{C}{n^{1 - \ell/6}}
    \biggl(
      \frac{M_2^2}{\lambda}
    \biggr)^{\ell/2 + 1},
  \end{align*}
  for a universal constant $C$.
  For \eqref{eq:whitened_bootstrap_2},
  by the triangle inequality and
  Lemma~\ref{lem:multi_index_sum},
  since $\ell \cdot 2^{\ell-2}/(\ell-1)! \leq 3$,
  by subadditivity of $x \mapsto x^{3/5}$ for $x \geq 0$,
  using H{\"o}lder's inequality
  followed by \eqref{eq:refined_V2}
  and \eqref{eq:refined_norms} again,
  and since
  $\sqrt{15} \bigl(1 + 2^{6\ell/5-6/5}\bigr)^{5/6} \leq 64$,
  \begin{align*}
    &\E\Biggl\{
      \biggl\|
      \frac{1}{n}
      \sum_{i=1}^n
      \Bigl\{
        v_\ell\bigl(\Sigma^{+/2}(X_i-\mu)\bigr)
        - v_\ell\bigl(\Sigma^{+/2}(X_i-\bar X_n)\bigr)
      \Bigr\}
      \biggr\|_2
    \Biggr\} \\
    &\quad\leq
    \frac{1}{n}
    \sum_{i=1}^n
    \E\biggl\{
      \Bigl\|
      v_\ell\bigl(\Sigma^{+/2}(X_i-\mu)\bigr)
      - v_\ell\bigl(\Sigma^{+/2}(X_i-\bar X_n)\bigr)
      \Bigr\|_2
    \biggr\} \\
    &\quad\leq
    \frac{1}{n}
    \sum_{i=1}^n
    \E\Biggl\{
      \biggl(
        \frac{\ell \cdot 2^{\ell-2}}{(\ell-1)!}
        \bigl\| \Sigma^{+/2} (\bar X_n - \mu) \bigr\|_2^2
        \iftoggle{journal}{\\ &\qquad\qquad\qquad\times}{}
        \Bigl(
          \bigl\|\Sigma^{+/2}(X_i-\bar X_n)\bigr\|_2^{2\ell-2}
          + \bigl\|\Sigma^{+/2}(X_i-\mu)\bigr\|_2^{2\ell-2}
        \Bigr)
      \biggr)^{1/2}
    \Biggr\} \\
    &\quad\leq
    \frac{\sqrt{3}}{\lambda^{\ell/2}}
    \frac{1}{n}
    \sum_{i=1}^n
    \E\biggl\{
      \|\bar X_n - \mu\|_2
      \Bigl(
        \|X_i-\bar X_n\|_2^{2\ell-2}
        + \|X_i-\mu\|_2^{2\ell-2}
      \Bigr)^{1/2}
    \biggr\} \\
    &\quad\leq
    \frac{\sqrt{3}}{\lambda^{\ell/2}}
    \biggl\{\E\Bigl(
        \|\bar X_n - \mu\|_2^6
    \Bigr)\biggr\}^{1/6}
    \frac{1}{n}
    \sum_{i=1}^n
    \biggl\{\E\Bigl(
        \|X_i-\bar X_n\|_2^{6\ell/5-6/5}
        + \|X_i-\mu\|_2^{6\ell/5-6/5}
    \Bigr)\biggr\}^{5/6} \\
    &\quad\leq
    \sqrt{15} \bigl(1 + 2^{6\ell/5-6/5}\bigr)^{5/6}
    \frac{M_2^{\ell}}{\sqrt n \lambda^{\ell/2}}
    \leq
    \frac{64 M_2^{\ell}}{\sqrt n \lambda^{\ell/2}}.
  \end{align*}
  For \eqref{eq:whitened_bootstrap_3}, let
  $\varepsilon_1,\ldots,\varepsilon_n$ be independent Rademacher
  variables, independent of the data.
  By a symmetrization inequality \citep[Lemma~2.3.6]{van1996weak},
  the conditional Jensen inequality,
  the multinomial theorem,
  subadditivity of $x \mapsto x^{3/\ell}$ and Jensen's inequality,
  \begin{align*}
    &\E\Biggl\{
      \biggl\|
      \frac{1}{n}
      \sum_{i=1}^n
      \Bigl\{
        \E\bigl\{v_\ell\bigl(\Sigma^{+/2}(X_i-\mu)\bigr)\bigr\}
        -v_\ell\bigl(\Sigma^{+/2}(X_i-\mu)\bigr)
      \Bigr\}
      \biggr\|_2
    \Biggr\} \\
    &\quad\leq
    \frac{2}{n}
    \E\Biggl\{
      \biggl\|
      \sum_{i=1}^n
      \varepsilon_i
      v_\ell\bigl(\Sigma^{+/2}(X_i-\mu)\bigr)
      \biggr\|_2
    \Biggr\}
    \leq
    \frac{2}{n \sqrt{\ell!}}
    \E\left[
      \Biggl\{
        \sum_{i=1}^n
        \bigl\|\Sigma^{+/2}(X_i-\mu)\bigr\|_2^{2\ell}
      \Biggr\}^{1/2}
    \right] \\
    &\quad\leq
    \frac{2}{n \sqrt{\ell!} \lambda^{\ell/2}}
    \E\left[
      \Biggl\{
        \sum_{i=1}^n
        \|X_i-\mu\|_2^{6}
      \Biggr\}^{\ell/6}
    \right]
    \leq
    \frac{2}{n \sqrt{\ell!} \lambda^{\ell/2}}
    \Biggl(
      \sum_{i=1}^n
      \E\bigl\{
        \|X_i-\mu\|_2^{6}
      \bigr\}
    \Biggr)^{\ell/6} \\
    &\quad\leq
    \frac{2 M_2^\ell}
    {n^{1 - \ell/6} \sqrt{\ell!} \lambda^{\ell/2}}.
  \end{align*}
  Combining the bounds on
  \eqref{eq:whitened_bootstrap_1},
  \eqref{eq:whitened_bootstrap_2} and
  \eqref{eq:whitened_bootstrap_3}
  by the triangle inequality, as
  $\lambda \leq M_2^2$ and $1-\ell/6\leq1/2$,
  we have for $\ell\in\{3,\ldots,k\}$ that
  \begin{align*}
    \E\Biggl\{
      \Biggl(\sum_{|\alpha|=\ell}
        \frac{\hat\delta_\alpha^2}{\alpha!}
      \Biggr)^{1/2}
      \I(\cA_n)
    \Biggr\}
    &\leq
    \frac{C}{n^{1 - \ell/6}}
    \biggl(
      \frac{M_2^2}{\lambda}
    \biggr)^{\ell/2 + 1}
    + \frac{C M_2^{\ell}}{\sqrt n \lambda^{\ell/2}}
    + \frac{2 M_2^\ell}
    {n^{1 - \ell/6} \sqrt{\ell!} \lambda^{\ell/2}} \\
    &\leq
    \frac{C M_2^{\ell + 2}}
    {n^{1 - \ell/6} \lambda^{\ell/2 + 1}}.
  \end{align*}
  For $s > 0$, define the random variable
  \begin{equation*}
    \hat\Delta(s)
    \coloneqq
    \frac{n}{2}
    \Biggl(
      \sum_{\ell = 1}^k
      \sum_{|\alpha| = \ell}
      \frac{\hat\delta_{\alpha}^2 s^{2 \ell}}{\alpha!}
    \Biggr)^{1/2}
    =
    \frac{n}{2}
    \Biggl(
      \sum_{\ell = 3}^k
      s^{2 \ell}
      \sum_{|\alpha| = \ell}
      \frac{\hat\delta_{\alpha}^2}{\alpha!}
    \Biggr)^{1/2},
  \end{equation*}
  where the second equality holds almost surely on $\cA_n$.
  By subadditivity of $x \mapsto \sqrt{x}$,
  if $n^{1/3} s^2 M_2^2 \leq \lambda$,
  \begin{align}
    \nonumber
    \E\bigl\{
      \hat\Delta(s)
      \I(\cA_n)
    \bigr\}
    &\leq
    \frac{n}{2}
    \sum_{\ell = 3}^k
    s^{\ell} \,
    \E\Biggl\{
      \Biggl(
        \sum_{|\alpha| = \ell}
        \frac{\hat\delta_{\alpha}^2}{\alpha!}
      \Biggr)^{1/2}
      \I(\cA_n)
    \Biggr\}
    \leq
    \frac{C M_2^2}{\lambda}
    \sum_{\ell = 3}^k
    \biggl(
      \frac{n^{1/3} s^{2} M_2^{2}}{\lambda}
    \biggr)^{\ell/2} \\
    \label{eq:refined_Delta}
    &\leq
    \frac{C M_2^2}{\lambda}
    \biggl(
      \frac{n^{1/3} s^{2} M_2^{2}}{\lambda}
    \biggr)^{3/2}
    =
    \frac{C s^{3} n^{1/2} M_2^5}{\lambda^{5/2}},
  \end{align}
  increasing the universal constant $C$ if necessary.

  For $a \in \N_0$ and $b \geq 0$,
  define $\hat\beta_{a,b,i}$ and $\hat\beta_{a,b}$
  as in Lemma~\ref{lem:beta_hat_refined}.
  By Lemma~\ref{lem:beta_hat_refined},
  as $1 + 2^{3x/2+3/2} \leq 59 \sqrt{x!}$ for all $x \in \N_0$,
  we have, for $s,\eta>0$,
  \begin{equation}
    \label{eq:refined_beta}
    \frac{n s^k}{\eta^\gamma}
    \E\bigl\{\hat\beta_{k,\gamma}\I(\cA_n)\bigr\}
    \leq
    \frac{1 + 2^{3k/2 + 3\gamma/2}}{\sqrt{k!}}
    \frac{n d^{\gamma a_p} M_2^{k+\gamma} s^k}
    {\lambda^{(k+\gamma)/2} \eta^\gamma}
    \leq
    \frac{59 n d^{\gamma a_p} M_2^{k+\gamma} s^k}
    {\lambda^{(k+\gamma)/2} \eta^\gamma}.
  \end{equation}

  Set $m \coloneqq 2$, recalling that
  $\hat\delta_{\alpha,i} = 0$ almost surely on $\cA_n$
  for $|\alpha| \in \{1, 2\}$, and observing that
  the first two moments of
  $\Sigma^{+/2} (X_i - \mu)$ do not depend on $i$.
  For $\ell \in \{m+1, \ldots, k\}$,
  set $k_\ell \coloneqq k - \ell$, so that
  $k_\ell\leq m$ and
  $k_\ell + \gamma + \ell = k + \gamma \leq 6$.
  By Lemma~\ref{lem:bound_Gamma} and
  Lemma~\ref{lem:beta_hat_refined}, for $s > 0$,
  \begin{align*}
    \hat\Gamma_\eta(s)
    &\coloneqq
    \sum_{i=1}^n
    \sum_{\ell = m+1}^k
    \frac{(2 s)^{\ell + k_\ell}}
    {2 \eta^{\gamma}}
    \bigl(\hat \beta_{k_\ell,\gamma}
    + \hat \beta_{k_\ell,\gamma,i}\bigr)
    \Biggl(
      \sum_{|\alpha| = \ell}
      \frac{\hat \delta_{\alpha,i}^2}{\alpha!}
    \Biggr)^{1/2} \\
    &\,\leq
    \frac{(2 s)^k}{2 \eta^\gamma}
    \sum_{\ell = 3}^k
    \Biggl\{
      n \hat\beta_{k-\ell,\gamma} \hat\beta_{\ell,0}
      + \sum_{i=1}^n
      \hat\beta_{k-\ell,\gamma,i} \hat\beta_{\ell,0,i}
    \Biggr\}.
  \end{align*}
  Taking expectations on $\cA_n$ and
  applying Lemma~\ref{lem:beta_hat_refined} again, we obtain
  \begin{align}
    \nonumber
    \E\bigl\{\hat\Gamma_\eta(s) \I(\cA_n)\bigr\}
    &\leq
    \frac{n (2 s)^{k}}{\eta^{\gamma}}
    \frac{d^{\gamma a_p} M_2^{k+\gamma}}{\lambda^{(k+\gamma)/2}}
    \sum_{\ell = 3}^k
    \frac{\bigl(1 + 2^{3(k - \ell + \gamma)/2}\bigr)
    \bigl(1 + 2^{3\ell/2}\bigr)}
    {\sqrt{(k-\ell)! \, \ell!}} \\
    \label{eq:refined_Gamma}
    &\leq
    \frac{C n d^{\gamma a_p} M_2^{k+\gamma} s^{k}}
    {\eta^{\gamma} \lambda^{(k+\gamma)/2}},
  \end{align}
  increasing the universal constant $C$ if necessary.
  We now set $s \coloneqq \phi(t) / \eta$ and
  apply Lemma~\ref{lem:borel}, taking an
  outer expectation and using
  \eqref{eq:refined_P},
  \eqref{eq:refined_Delta},
  \eqref{eq:refined_beta} and
  \eqref{eq:refined_Gamma} to obtain,
  for $n^{1/3} s^2 M_2^2 \leq \lambda$,
  \begin{align*}
    &
    \E^\star \biggl(
      \sup_{A \in \cB(\R^d)}
      \Bigl\{
        \P(\hat S_n \in A)
        - \P\big( \hat S_n^* \in A_p^{3\eta} \mid \bX_n \big)
      \Bigr\}
    \biggr) \\
    &\quad\leq
    2 e^{-t^2 / 2}
    + \E\Biggl(
      \biggl\{
        \frac{n \hat\beta_{k,\gamma}\phi(t)^k}{\eta^{k+\gamma}}
        + \hat\Delta\biggl(\frac{\phi(t)}{\eta}\biggr)
        + \hat\Gamma_\eta\biggl(\frac{\phi(t)}{\eta}\biggr)
      \biggr\}
      \I(\cA_n)
    \Biggr)
    + \P(\cA_n^\c) \\
    &\quad\leq
    2 e^{-t^2 / 2}
    + \frac{C n^{1/2} M_2^5 \phi(t)^{3}}
    {\lambda^{5/2} \eta^3}
    + \frac{C n d^{\gamma a_p} M_2^{k+\gamma} \phi(t)^k}
    {\lambda^{(k+\gamma)/2} \eta^{k+\gamma}}
    + \frac{C M_2^6}{n^{3/2} \lambda^3}.
  \end{align*}
  If $n^{1/3} s^2 M_2^2 > \lambda$, then with $C \geq 1$,
  the second term on the right-hand side above exceeds one
  as $\lambda \leq M_2^2$, so this bound holds trivially.

  We now establish a coupling between $\hat S_n$ and $\hat S_n^*$
  by Strassen's theorem, as in the proof of Theorem~\ref{thm:bootstrap}.
  For each $t, \eta > 0$, there is a probability space
  supporting $\bX_n'$, $\hat S_n'$ and $\hat S_n^{*\prime}$, where
  $\hat S_n' \equallaw \hat S_n$, $\hat S_n'$ is independent of $\bX_n'$,
  and $(\bX_n', \hat S_n^{*\prime}) \equallaw (\bX_n, \hat S_n^*)$, with
  \begin{equation*}
    \P\bigl(
      \|\hat S_n' - \hat S_n^{* \prime}\|_p
      > 3\eta
    \bigr)
    \leq
    2 e^{-t^2 / 2}
    + \frac{C n^{1/2} M_2^5 \phi(t)^{3}}{\lambda^{5/2} \eta^3}
    + \frac{C n d^{\gamma a_p} M_2^{k+\gamma} \phi(t)^k}
    {\lambda^{(k+\gamma)/2} \eta^{k+\gamma}}
    + \frac{C M_2^6}{n^{3/2} \lambda^3},
  \end{equation*}
  where $C > 0$ is a universal constant.

  For the remainder of the proof, take this coupling
  with $\eta/B_p$ in place of $\eta$.
  It remains to relate this coupling
  (between $\hat S_n$ and $\hat S_n^*$)
  to a coupling of $S_n$ and $S_n^*$.
  To this end, define
  $S_n' \coloneqq \Sigma^{1/2} \hat S_n'$,
  let $\bar X_n' \coloneqq \sum_{i=1}^n X_i' / n$ and
  $\hat \Sigma_n' \coloneqq
  \sum_{i=1}^n X_i' X_i^{\prime\T} / n - \bar X_n' \bar X_n^{\prime\T}$,
  and take
  $S_n^{*\prime} \coloneqq \hat \Sigma_n^{\prime 1/2} \hat S_n^{*\prime}$.
  Observe that $S_n'$ is a function of $\hat S_n'$ only,
  so $S_n'$ is independent of $\bX_n'$.
  Further, for every $v \in \ker \Sigma$, we have
  $v^\T S_n = 0$ almost surely,
  so $\Sigma^{1/2} \hat S_n = \Sigma^{1/2} \Sigma^{+/2} S_n = S_n$
  and hence $S_n' \equallaw S_n$.
  Similarly, $S_n^*$ lies in the image of $\hat \Sigma_n$ almost surely, so
  $\hat \Sigma_n^{1/2} \hat S_n^*
  = \hat \Sigma_n^{1/2} \hat \Sigma_n^{+/2} S_n^* = S_n^*$
  and hence $(\bX_n', S_n^{*\prime}) \equallaw (\bX_n, S_n^*)$.

  Note that the image of $\Sigma$
  contains the image of $\hat\Sigma_n'$ almost surely.
  Apply \citet[Problem~X.5.5]{bhatia1997matrix} on the image of
  $\Sigma$, on which $\Sigma$ is positive definite
  with smallest eigenvalue
  $\lambda$, along with H{\"o}lder's inequality, to obtain
  \begin{align*}
    \bigl\|
    S_n' - S_n^{*\prime}
    \bigr\|_p
    &=
    \bigl\|
    \Sigma^{1/2} \hat S_n'
    - \hat \Sigma_n^{\prime 1/2} \hat S_n^{*\prime}
    \bigr\|_p \\
    &=
    \bigl\|
    \Sigma^{1/2}
    \bigl( \hat S_n' - \hat S_n^{*\prime} \bigr)
    + \bigl( \Sigma^{1/2} - \hat \Sigma_n^{\prime 1/2} \bigr)
    \hat S_n^{*\prime}
    \bigr\|_p \\
    &\leq
    \bigl\|
    \Sigma^{1/2}
    \bigr\|_{p \to p}
    \bigl\| \hat S_n' - \hat S_n^{*\prime} \bigr\|_p
    + d^{(1/p - 1/2) \lor 0}
    \bigl\| \Sigma^{1/2} - \hat \Sigma_n^{\prime 1/2} \bigr\|_{2 \to 2}
    \bigl\|\hat S_n^{*\prime}\bigr\|_2 \\
    &\leq
    B_p
    \bigl\| \hat S_n' - \hat S_n^{*\prime} \bigr\|_p
    + d^{a_p}
    \lambda^{-1/2}
    \bigl\| \Sigma - \hat \Sigma_n' \bigr\|_{2 \to 2}
    \bigl\|\hat S_n^{*\prime}\bigr\|_2.
  \end{align*}
  Conditional on $\bX_n$, the summands of $S_n^*$ are i.i.d.\
  with mean zero and covariance matrix $\hat\Sigma_n$, so
  $\E\bigl\{ \|S_n^*\|_2^2 \mid \bX_n \bigr\}
  = n \tr \hat\Sigma_n
  = \sum_{i=1}^n \|X_i - \bar X_n\|_2^2
  \leq \sum_{i=1}^n \|X_i - \mu\|_2^2$.
  By H{\"o}lder's inequality,
  \eqref{eq:refined_V1} and \eqref{eq:refined_norms},
  since $\hat \lambda_n \geq \lambda/2$ almost surely on $\cA_n$,
  \begin{align*}
    &\E\Bigl\{
      \bigl\| \hat \Sigma_n - \Sigma \bigr\|_{2 \to 2}^2
      \bigl\|\hat S_n^*\bigr\|_2^2
      \I(\cA_n)
    \Bigr\}
    \leq
    \frac{2}{\lambda}
    \E\Biggl\{
      \bigl\| \hat \Sigma_n - \Sigma \bigr\|_{2 \to 2}^2
      \sum_{i=1}^n \|X_i - \mu\|_2^2
    \Biggr\} \\
    &\qquad\leq
    \frac{2}{\lambda}
    \Bigvvvert
    \bigl\| \hat \Sigma_n - \Sigma \bigr\|_{2 \to 2}
    \Bigvvvert_3^2
    \sum_{i=1}^n
    \Bigvvvert
    \bigl\| X_i - \mu \bigr\|_2
    \Bigvvvert_6^2
    \leq
    \frac{2}{\lambda}
    \frac{64 M_2^4}{n}
    n M_2^2
    = \frac{128 M_2^6}{\lambda}.
  \end{align*}
  By Markov's inequality,
  \begin{align*}
    \P\Bigl(
      \bigl\| \hat \Sigma_n - \Sigma \bigr\|_{2 \to 2}
      \bigl\|\hat S_n^*\bigr\|_2
      > \eta
    \Bigr)
    &\leq
    \P(\cA_n^\c)
    + \frac{1}{\eta^2}
    \E\Bigl\{
      \bigl\| \hat \Sigma_n - \Sigma \bigr\|_{2 \to 2}^2
      \bigl\|\hat S_n^*\bigr\|_2^2
      \I(\cA_n)
    \Bigr\} \\
    &\leq
    \P(\cA_n^\c)
    + \frac{128 M_2^6}{\lambda \eta^2}.
  \end{align*}

  Finally, by a union bound and equality in law,
  \begin{align*}
    &\P\Bigl(
      \bigl\|
      S_n' - S_n^{*\prime}
      \bigr\|_p
      > 4 \eta
    \Bigr)
    \leq
    \P\biggl(
      \bigl\| \hat S_n' - \hat S_n^{*\prime} \bigr\|_p
      > \frac{3 \eta}{B_p}
    \biggr)
    + \P\biggl(
      \bigl\| \hat \Sigma_n - \Sigma \bigr\|_{2 \to 2}
      \bigl\|\hat S_n^*\bigr\|_2
      > \frac{\lambda^{1/2} \eta}{d^{a_p}}
    \biggr) \\
    &\qquad\leq
    2 e^{-t^2 / 2}
    + \frac{C n^{1/2} M_2^5 B_p^3 \phi(t)^{3}}{\lambda^{5/2} \eta^3}
    + \frac{C n d^{\gamma a_p}
    (M_2 B_p)^{k+\gamma} \phi(t)^k}
    {\lambda^{(k+\gamma)/2} \eta^{k+\gamma}}
    + \frac{C M_2^6}{n^{3/2} \lambda^3}
    + \frac{C d^{2 a_p} M_2^6}{\lambda^2 \eta^2}.
  \end{align*}
  The result follows by replacing $\eta$ with $\eta/4$
  and increasing the universal constant $C$.
\end{proof}

\begin{proof}[Details for \eqref{eq:refined_rate}]
  Under the conditions imposed, since $r \geq 6$,
  Minkowski's and H{\"o}lder's inequalities yield
  $M_2 \leq \sqrt d M$,
  $\max_{j \in [d]} \Sigma_{j j} \leq M^2$
  and $\lambda \leq \tr \Sigma \leq d M^2$.
  Since $\Sigma^{1/2}$ is symmetric, Cauchy--Schwarz gives
  \begin{align*}
    \|\Sigma^{1/2}\|_{1 \to 1}
    &=
    \|\Sigma^{1/2}\|_{\infty \to \infty}
    =
    \max_{j \in [d]}
    \sum_{r=1}^d \bigl|(\Sigma^{1/2})_{jr}\bigr| \\
    &\leq
    \sqrt d \max_{j \in [d]}
    \Biggl(
      \sum_{r=1}^d \bigl(\Sigma^{1/2}\bigr)_{jr}^2
    \Biggr)^{1/2}
    =
    \sqrt{d} \max_{j \in [d]} \sqrt{\Sigma_{j j}}
    \leq \sqrt d M.
  \end{align*}
  For $p \in (1, \infty)$ and $x \in \R^d$,
  H{\"o}lder's inequality gives
  \begin{align*}
    \bigl\|\Sigma^{1/2} x\bigr\|_p^p
    &=
    \sum_{j=1}^d
    \Biggl|
    \sum_{r=1}^d (\Sigma^{1/2})_{jr} x_r
    \Biggr|^p
    \leq
    \sum_{j=1}^d
    \Biggl(
      \sum_{r=1}^d \bigl|(\Sigma^{1/2})_{jr}\bigr|
    \Biggr)^{p-1}
    \sum_{s=1}^d
    \bigl|(\Sigma^{1/2})_{j s}\bigr| \cdot |x_s|^p \\
    &\leq
    \bigl(\sqrt d M\bigr)^{p-1}
    \sum_{r=1}^d |x_r|^p
    \sum_{j=1}^d |(\Sigma^{1/2})_{jr}|
    \leq
    \bigl(\sqrt d M\bigr)^p \|x\|_p^p.
  \end{align*}
  Thus, $B_p \leq \sqrt d M$ for every $p \in [1, \infty]$.
  After rescaling the bound from Theorem~\ref{thm:refined_bootstrap}
  by $\sqrt n$, setting $k = 5$ and $\gamma = 1$,
  it suffices to choose $t_n$ and $\eta_n$ such that
  \begin{equation*}
    e^{-t_n^2 / 2}
    + \frac{d^{4} \phi(t_n)^3}{n \lambda^{5/2} \eta_n^3}
    + \frac{d^{11/2 + 1/(p \land 2)} \phi(t_n)^5}
    {n^2 \lambda^3 \eta_n^6}
    + \frac{d^3}{n^{3/2} \lambda^3}
    + \frac{d^{2 + 2/(p \land 2)}}{n \lambda^2 \eta_n^2}
    \to 0.
  \end{equation*}
  To this end, fix $R_n \to \infty$, assuming without
  loss of generality that $R_n \geq 1$, let
  $\cA \coloneqq \{n \in \N: d^3 / (n^{3/2} \lambda^3) \leq 1/R_n\}$
  and take $t_n \coloneqq \sqrt R_n$.
  For $n \in \cA$,
  \begin{equation*}
    \eta_n
    \coloneqq
    \frac{t_n d^{4/3} \phi(t_n)}{n^{1/3} \lambda^{5/6}}
    + \frac{t_n d^{11/12 + 1/\{6 (p \land 2)\}} \phi(t_n)^{5/6}}
    {n^{1/3} \lambda^{1/2}}
    + \frac{t_n d^{1 + 1/(p \land 2)}}{n^{1/2} \lambda}
    \lesssim
    \frac{R_n d^{4/3} \psi}{n^{1/3} \lambda^{5/6}}.
  \end{equation*}
  Here, we used
  $\phi(t)\leq (1 + t)\psi$, $\lambda \leq d M^2$,
  noting that $\lambda^{1/3} d^{-5/12 + 1/\{6(p \land 2)\}}
  \phi(t_n)^{-1/6} \lesssim 1$ and
  $d^{-1/3+1/(p\land2)} n^{-1/6} \lambda^{-1/6} \phi(t_n)^{-1}
  \lesssim d^{-1/2+1/(p\land2)} \psi^{-1} \lesssim 1$.
  Substituting $\psi$ gives \eqref{eq:refined_rate}.
  For $n \in \cA^\c \coloneqq \N \setminus \cA$,
  we have $R_n d^{4/3} \psi n^{-1/3} \lambda^{-5/6}
  \gtrsim d^{1/2} \psi \gtrsim d^{1 / (p \land 2)}$.
  Since $\E\bigl\{\|S_n\|_2^2 + \|S_n^*\|_2^2\bigr\} \leq 2 n d M^2$,
  H{\"o}lder's and Markov's inequalities give
  $\|S_n\|_p / \sqrt n + \|S_n^*\|_p / \sqrt n
  \lesssim_\P d^{1/(p \land 2)}$, yielding the bound.

  Finally, suppose that $M_2 \leq M$.
  Let $m \coloneqq \rank \Sigma$ and choose
  $U\in\R^{d \times m}$ with $U^\T U = I_m$
  and $U U^\T \Sigma = \Sigma$.
  Note that  $U U^\T (X_i - \mu) = X_i - \mu$ almost surely,
  $\bigl\|U^\T (X_i - \mu)\bigr\|_2 = \|X_i - \mu\|_2$ almost surely,
  and $\lambda_{\min}\bigl(U^\T \Sigma U\bigr) = \lambda$.
  Apply Theorem~\ref{thm:refined_bootstrap} to
  $U^\T (X_i - \mu)$ with $p=2$, $k=5$ and $\gamma=1$.
  Denote the resulting coupled variables by
  $\tilde X_i'$, $\tilde S_n'$ and $\tilde S_n^{*\prime}$,
  and set $X_i' \coloneqq \mu + U \tilde X_i'$,
  $S_n' \coloneqq U \tilde S_n'$ and
  $S_n^{*\prime} \coloneqq U \tilde S_n^{*\prime}$.
  It suffices to choose $t_n$ and $\eta_n$ such that
  \begin{equation*}
    e^{-t_n^2 / 2}
    + \frac{M_2^5 B_2^3 \bigl(\sqrt m + t_n\bigr)^3}
    {n \lambda^{5/2} \eta_n^3}
    + \frac{(M_2 B_2)^{k+\gamma} \bigl(\sqrt m + t_n\bigr)^k}
    {n^{(k+\gamma)/2 - 1} \lambda^{(k+\gamma)/2} \eta_n^{k+\gamma}}
    + \frac{M_2^6}{n^{3/2} \lambda^3}
    + \frac{M_2^6}{n \lambda^2 \eta_n^2}
    \to 0.
  \end{equation*}
  Taking $t_n \coloneqq \sqrt{R_n}$,
  using $\lambda m \leq \tr \Sigma$
  and $B_2^2 \leq \tr \Sigma \leq M_2^2 \leq M^2$,
  this holds if
  \begin{equation*}
    \frac{R_n^{3/2}}{n \lambda^4 \eta_n^3}
    + \frac{R_n^{k/2}}
    {n^{(k+\gamma)/2 - 1} \lambda^{k+\gamma/2} \eta_n^{k+\gamma}}
    + \frac{1}{n^{3/2} \lambda^3}
    + \frac{1}{n \lambda^2 \eta_n^2}
    \to 0.
  \end{equation*}
  Take $\eta_n \coloneqq R_n n^{-1/3} \lambda^{-4/3}$.
  For $n \lambda^4 \geq 1$, the above display tends to zero.
  Otherwise, $\eta_n \geq R_n$, and Markov's inequality applies since
  $\E\{\|S_n / \sqrt n\|_2^2 + \|S_n^* / \sqrt n\|_2^2\}
  \leq 2 M^2$.
  Using $\|x\|_p\leq d^{(1/p-1/2) \lor 0} \|x\|_2$
  gives the claimed rates for all $p \in [1, \infty]$.
\end{proof}

\subsection{Details for Section~\ref{sec:literature}}
\label{sec:literature_details}

We provide details for all of the results and rates
presented in Section~\ref{sec:literature},
including Table~\ref{tab:literature}.

\subsubsection{Yurinskii-type couplings}
\label{sec:literature_details_yurinskii}

\citet[Theorem~1.1]{zaitsev1987estimates}
gave the Prokhorov-type bound
\begin{equation*}
  P_1(A) - P_2\bigl(A_2^{\eta}\bigr)
  \leq C d^2
  \exp\biggl(-\frac{\eta}{C d^2\tau}\biggr),
\end{equation*}
for each $\eta > 0$, all Borel sets $A \subseteq \R^d$ and
every law $P_1 \in \cA_d(\tau)$.
Here, membership in
$\cA_d(\tau)$ is given by existence and appropriate regularity
of the cumulant generating function,
$P_2$ is a Gaussian measure with the same mean and covariance as $P_1$,
and $C > 0$ is a universal constant.
Taking $P_1$ to be the law of $S_n / \sqrt n = \sum_{i=1}^n X_i / \sqrt n$
and assuming $\|X_i\|_2 \leq M$, along with other conditions, we have
$P_1 \in \cA_d\bigl(C M / \sqrt n\bigr)$.
Using $A_2^\eta \subseteq A_p^\eta$ and applying
Strassen's theorem gives, for each $\eta \geq 1$,
a coupling with Gaussian target $T_n$ such that
\begin{equation*}
  \P\Biggl(
    \frac{\|S_n' - T_n'\|_p}{\sqrt n}
    > \frac{C M d^2}{\sqrt n}
    \bigl\{ \eta + \log \bigl(C d^2\bigr) \bigr\}
  \Biggr)
  \leq
  e^{-\eta},
\end{equation*}
for $p \in [2, \infty]$.
We therefore obtain a sample size factor of $n^{-1/2}$ and a
dimension factor of $d^2 \log 2d$ if $\|X_i\|_2 \leq M$,
and $d^{5/2} \log 2d$ if $\|X_i\|_\infty \leq M$,
using $\|X_i\|_2 \leq \sqrt d \|X_i\|_\infty$.

\citet[Theorem~1]{le1988prokhorov} gave the Prokhorov bound
\begin{equation*}
  \rho_2\biggl(\frac{S_n}{\sqrt n}, \frac{T_n}{\sqrt n}\biggr)
  \leq
  C
  \Biggl\{
    \frac{\bar\beta_{3,0}}{\sqrt n}
    \biggl(
      d + \biggl|
      \log \frac{\bar\beta_{3,0} d}{\sqrt n}
      \biggr|
    \biggr)
  \Biggr\}^{1/4},
\end{equation*}
with $\bar\beta_{3,0}$ defined as in Assumption~\ref{ass:data}.
Dividing $S_n / \sqrt n$ and $T_n / \sqrt n$ by
$R_n \bigl(\bar \beta_{3,0} d / \sqrt n\bigr)^{1/3}$,
then applying Strassen's theorem, gives an $\ell_p$ approximation
rate of $n^{-1/6} d^{1/3} \bar\beta_{3,0}^{1/3}$,
for $p \in [2, \infty]$.
His Theorem~1 also gives Prokhorov bounds for arbitrary norms and
for the maximum norm;
under the assumptions of Table~\ref{tab:literature}
these rates are no faster.

\citet[Supplemental Appendix~I.1, Lemma~38]{belloni2019conditional}
obtained a result equivalent to our Theorem~\ref{thm:yurinskii}
with $(k, \gamma, p) = (2, 1, p)$, for $p \in [1, \infty]$,
assuming a Gaussian target distribution and exact moment matching.
When considering independent summands and with Gaussian targets
and exact moment matching,
\citet[Corollary~2.3]{cattaneo2025yurinskii}
gives the same coupling rates as our Theorem~\ref{thm:yurinskii}
with $k \in \{2, 3\}$ and $\gamma = 1$.

Table~\ref{tab:rates_moments} summarizes the rates attained
by these Yurinskii-type couplings for $p \in [2, \infty]$
in terms of the average moment bounds $\bar\beta_{k, \gamma}$
and the dimension factor $\psi$, as defined in
Section~\ref{sec:yurinskii}.
Rates for $p \in [1, 2)$ may be derived by multiplying
the corresponding $p = 2$ rate by $d^{1/p - 1/2}$.
\citet{zaitsev1987estimates} requires finite exponential moments,
so is omitted.
The rates in Table~\ref{tab:rates_moments} imply the corresponding
entries in Table~\ref{tab:literature} by the following arguments.
If the target summands $Y_i$ satisfy the same bounds as $X_i$,
then $\|X_i\|_\infty \leq M$ almost surely implies
$\bar\beta_{k,\gamma} \lesssim M^{k+\gamma} d^{k/2 + \gamma/p}$
and $\max_{j \in [d]} \vvvert e_j^\T X_i \vvvert_r \leq M$ gives
$\bar\beta_{k,\gamma} \lesssim M^{k+\gamma} d^{k/2 + \gamma/(p \land r)}$
when $r \geq 2 \lor (k + \gamma)$,
while $\|X_i\|_2 \leq M$ almost surely implies
$\bar\beta_{k,\gamma} \lesssim M^{k+\gamma}$, for $p \in [2, \infty]$.
For Gaussian $Y_i$ with matching means and covariances,
$\|X_i\|_\infty \leq M$ almost surely implies
$\bar\beta_{k,\gamma} \lesssim M^{k+\gamma} d^{k/2} \psi^\gamma$;
bounds under $\|X_i\|_2 \leq M$ and
$\vvvert e_j^\T X_i \vvvert_r \leq M$ remain valid
for $r < \infty$.

\begin{table}[!htbp]
  \centering
  \begin{minipage}{\textwidth}
    \centering
    \footnotesize
    \setlength{\tabcolsep}{4pt}
    \renewcommand{\arraystretch}{1.5}
    \caption{Yurinskii coupling rates in terms of average moment bounds.}
    \label{tab:rates_moments}
    \begin{tabular}{@{}lccc@{}}
      \hline
      Source & Matching & $p$ & Coupling rate \\
      \hline
      \citet{yurinskii1978error}
      & $2$ & $2$
      & $n^{-1/6} d^{1/3} \bar\beta_{3,0}^{1/3}$ \\
      \citet{le1988prokhorov}
      & $2$ & $[2,\infty]$
      & $n^{-1/6} d^{1/3} \bar\beta_{3,0}^{1/3}$ \\
      \citet{belloni2019conditional}
      & $2$ & $[2,\infty]$
      & $n^{-1/6} \bar\beta_{2,1}^{1/3} \psi^{2/3}$ \\
      \citet{cattaneo2025yurinskii}
      & $3$ & $[2,\infty]$
      & $n^{-1/4} \bar\beta_{3,1}^{1/4} \psi^{3/4}$ \\
      \hline
      Theorem~\ref{thm:yurinskii}
      & $k$ & $[2,\infty]$
      & $n^{-1/2+1/(k+\gamma)}
      \bar\beta_{k,\gamma}^{1/(k+\gamma)} \psi^{k/(k+\gamma)}$ \\
      \hline
    \end{tabular}

    \medskip
    \begin{minipage}{0.95\linewidth}
      \centering
      \footnotesize
      \emph{Notes:}
      Rates use the conventions of Table~\ref{tab:literature},
      with arbitrarily slowly diverging factors omitted.
    \end{minipage}
  \end{minipage}
\end{table}

\subsubsection{KMT-type couplings}
\label{sec:literature_details_hungarian}

\citet[Theorem~1]{komlos1976approximation} achieved
a Gaussian approximation for i.i.d.\ real-valued random variables
with finite exponential moments
at a rate of $n^{-1/2} \log n$;
their result also establishes this rate for all partial sums.
For i.i.d.\ vectors, for each $\eta > 0$,
\citet[Theorem~4]{einmahl1989extensions} derived a coupling with
\begin{equation*}
  \P\bigl(\|S_n - T_n\|_2 \geq \eta\bigr)
  \leq \frac{C_1 n}{H(C_2 \eta)},
\end{equation*}
where $C_1$, $C_2$ depend on the law of $X_1$,
$T_n$ is Gaussian and $H$ satisfies
$\E\{H(\|X_1\|_2)\} < \infty$ and some other conditions.
Taking $H(t) = t^r$ with $r > 3$ gives a rate of
$n^{-1/2 + 1/r}$ while
$H(t) = e^{\sqrt t}$ gives $n^{-1/2} (\log n)^2$.
Under exponential moments and a further smoothness condition,
his Theorem~3 gives a rate of $(\log n)/\sqrt n$.
We omit these results from Table~\ref{tab:literature}
because of their unspecified dependence on dimension.

\citet{zaitsev1998multidimensional} provided the following result.
Suppose $\Var(X_i) = \Sigma$ for each $i \in [n]$,
where $\Sigma$ is positive definite with bounded condition number,
and assume that $\|X_i\|_2\leq M$ almost surely.
Then the laws of $\Sigma^{-1/2} X_i$ belong to
$\cA_d(\tau)$ with $\tau = C M \lambda_{\min}(\Sigma)^{-1/2} \geq 1$,
where $C \geq 1$ is universal.
For each $\vartheta > 0$,
\citet[Theorem~1.3]{zaitsev1998multidimensional} gives
a Gaussian coupling satisfying
\begin{equation*}
  \E \biggl\{
    \exp
    \biggl(
      \frac{\|\Sigma^{-1/2}(S_n'-T_n')\|_\infty}
      {C_\vartheta \tau d^3 \log 2d}
    \biggr)
  \biggr\}
  \leq
  \exp\Bigl(
    C_\vartheta
    d^{9/4+\vartheta}
    \log 2n
  \Bigr),
\end{equation*}
where $C_\vartheta > 0$ depends only on $\vartheta$.
By Markov's inequality and since
$\|\Sigma^{1/2} v\|_p \leq \sqrt{d \|\Sigma\|_{2 \to 2}} \|v\|_\infty$
for all $v \in \R^d$ and $p \in [2, \infty]$, we have
\begin{equation*}
  \frac{\|S_n'-T_n'\|_p}{\sqrt n}
  \lesssim_\P
  \frac{M d^{23/4+\vartheta} (\log 2d)(\log 2n)}{\sqrt n}.
\end{equation*}
This yields a sample size factor of $n^{-1/2} \log 2n$ and
a dimension factor of $d^{23/4 + \vartheta} \log 2d$.
If instead $\|X_i\|_\infty \leq M$ almost surely,
we have $\|X_i\|_2 \leq \sqrt d M$ and the dimension factor
is instead $d^{25/4 + \vartheta} \log 2d$.
If all third moments of each $X_i$ vanish,
\citet[Theorem~1.6]{zaitsev1998multidimensional}
removes the factor $\log 2d$ from both rates.

\citet{zaitsev2001multidimensional,zaitsev2002multidimensionalII,
zaitsev2002multidimensionalIII} allowed for differing covariance matrices.
Suppose $\lambda I_d \preceq \Sigma_i \preceq c \lambda I_d$
for each $i \in [n]$, with $c \geq 1$ fixed and $\lambda > 0$,
and assume $\|X_i\|_2 \leq M$ almost surely.
Then the laws of $\lambda^{-1/2} X_i$ lie in $\cA_d(\tau)$ with
$\tau = C M \lambda^{-1/2} \geq 1$, where $C \geq 1$ is universal.
Applying \citet[Theorem~1.2, Equation~1.8]{zaitsev2001multidimensional},
whose proof uses \citet[Theorem~8.2]{zaitsev2002multidimensionalII}
and is concluded in \citet[Section~11]{zaitsev2002multidimensionalIII},
gives a Gaussian coupling satisfying
\begin{equation*}
  \E \exp\biggl(
    \frac{\|\lambda^{-1/2} (S_n'-T_n')\|_2}
    {C_\vartheta \tau d^{9/2} \log 2d}
  \biggr)
  \leq
  \exp\bigl(C_\vartheta d^{3+\vartheta} \log 2n\bigr),
\end{equation*}
where $C_\vartheta>0$ depends only on $\vartheta$ and $c$.
By Markov's inequality and as $\|v\|_p \leq \|v\|_2$ for $p\in[2,\infty]$,
\begin{equation*}
  \frac{\|S_n'-T_n'\|_p}{\sqrt n}
  \lesssim_\P
  \frac{M d^{15/2+\vartheta}(\log 2d)(\log 2n)}{\sqrt n}.
\end{equation*}
If instead $\|X_i\|_\infty\leq M$ almost surely, then
$\|X_i\|_2\leq\sqrt d M$ and the dimension factor is
$d^{8+\vartheta}\log 2d$.
If all third moments of each $X_i$ vanish,
\citet[Theorem~1.2]{zaitsev2001multidimensional}
removes the factor $\log 2d$ from both rates.

\citet{zaitsev2007estimates} allowed settings with
summands having finite polynomial moments.
Suppose $L_r \coloneqq \sum_{i=1}^n \E\bigl(\|X_i\|_2^r\bigr) \in (0,\infty)$
for $r > 2$.
Assume that $[n]$ can be partitioned into consecutive blocks
$\cB$ such that for each block $B \in \cB$,
\begin{equation*}
  4e^2 L_r^{2/r} I_d \preceq \sum_{i\in B}\Sigma_i
  \preceq 4e^2 c L_r^{2/r} I_d,
\end{equation*}
where $c \geq 1$ is fixed.
By \citet[Corollary~3]{zaitsev2007estimates},
for any $\vartheta > 0$,
there is a Gaussian coupling satisfying
\begin{equation*}
  \P\bigl(\|S_n'-T_n'\|_2 > 5\eta\bigr)
  \leq \frac{2 L_r}{\eta^{r}}
  + \exp\biggl(-\frac{\eta}{C_\vartheta L_r^{1/r} d^{9/2} \log 2d}\biggr),
\end{equation*}
for each
$\eta > C_\vartheta L_r^{1/r} d^{15/2 + \vartheta}\log 2n$,
where $C_\vartheta$ depends only on $\vartheta$ and $c$.
It follows that for $p \in [2,\infty]$,
\begin{equation*}
  \frac{\|S_n'-T_n'\|_p}{\sqrt n}
  \lesssim_\P
  \frac{L_r^{1/r}d^{15/2+\vartheta}\log 2n}{\sqrt n}.
\end{equation*}
If $\|X_i\|_2\leq M$ almost surely, then $L_r \leq n M^r$,
while if $\|X_i\|_\infty \leq M$ almost surely or
$\vvvert e_j^\T X_i\vvvert_r\leq M$, then
$L_r \leq n d^{r/2} M^r$.

\citet{gotze2009bounds} obtained the following bound for finite
polynomial moments. Let $\cB$ be a partition of $[n]$ into consecutive
blocks, and set $w \coloneqq c_1 L_r^{1/r} / \log 2|\cB|$ for fixed $c_1 > 0$.
Assume that, for fixed $c_2 \geq 1$ and $\varepsilon \in (0,1)$,
we have $w^2 I_d \preceq \sum_{i\in B} \Sigma_i \preceq c_2 w^2 I_d$
for every $B \in \cB$ and, for fixed $c_3 > 0$,
\begin{equation*}
  d^{r/2} |\cB|^\varepsilon (\log 2|\cB|)^{r+3}
  \max_{B\in\cB}
  \sum_{i\in B}
  \E\bigl(\|X_i\|_2^r\bigr)
  \leq c_3 L_r.
\end{equation*}
For $\vartheta > 0$, apply
\citet[Theorem~4, Equation~1.19]{gotze2009bounds} to obtain
a coupling with
\begin{equation*}
  \E\bigl(\|S_n'-T_n'\|_2^r\bigr)
  \leq
  C_\vartheta
  L_r d^{r(21/2+\vartheta/2)} (\log 2d)^r
  \leq
  C_\vartheta^2
  L_r d^{r(21/2+\vartheta)},
\end{equation*}
where $C_\vartheta > 0$ depends only on
$\vartheta$, $r$, $\varepsilon$, $c_1$, $c_2$ and $c_3$.
By Markov's inequality and as $\|v\|_p \leq\|v\|_2$ for $p \in [2,\infty]$,
we deduce that
\begin{equation*}
  \frac{\|S_n'-T_n'\|_p}{\sqrt n}
  \lesssim_\P
  \frac{\bigl\{\E\bigl(\|S_n'-T_n'\|_2^r\bigr)\bigr\}^{1/r}}{\sqrt n}
  \lesssim
  \frac{L_r^{1/r} d^{21/2+\vartheta}}{\sqrt n}.
\end{equation*}
The previous bounds on $L_r$ give the sample size factor
$n^{-1/2+1/r}$ and dimension factors $d^{21/2+\vartheta}$ under
$\|X_i\|_2\leq M$
and $d^{11+\vartheta}$ under
$\|X_i\|_\infty \leq M$ or
$\vvvert e_j^\T X_i\vvvert_r\leq M$.

\citet{gotze2011estimates} obtained the following bound for i.i.d.\
summands. Suppose $\Var(X_i) = \Sigma$ is positive
definite with bounded condition number.
For each $\vartheta>0$, set $\psi = (21/2+\vartheta) \land 11$.
With $A(r,\psi,d)$ as defined in \citet[Equation~8]{gotze2011estimates}, we have
\begin{equation*}
  A(r, \psi, d)^{1/r}
  = \max\bigl\{d^\psi,d^{(r+2)/4}(1\lor\log d)^{(r+1)/2}\bigr\}
  \leq (\log 2)^{-(r+1)/2}\varrho_r(d),
\end{equation*}
since $\psi \leq 21/2 + \vartheta$ and
$1 \lor \log d \leq \log(2d)/\log 2$.
Therefore, \citet[Theorem~2]{gotze2011estimates}
gives a Gaussian coupling satisfying
\begin{equation*}
  \E\bigl(\|S_n'-T_n'\|_2^r\bigr)
  \leq C_\vartheta \varrho_r(d)^r L_r,
\end{equation*}
where $C_\vartheta>0$ depends only on $\vartheta$, $r$,
and an upper bound on the condition number of $\Sigma$.
By Markov's inequality and as $\|v\|_p \leq \|v\|_2$ for $p \in [2,\infty]$,
we deduce that
\begin{equation*}
  \P\biggl(\frac{\|S_n'-T_n'\|_p}{\sqrt n}>\eta\biggr)
  \leq \frac{C_\vartheta \varrho_r(d)^r L_r}{n^{r/2}\eta^r}
\end{equation*}
for each $\eta > 0$.
The previous bounds on $L_r$ give the sample size factor
$n^{-1/2+1/r}$ and dimension factors $\varrho_r(d)$ under
$\|X_i\|_2\leq M$
and $\sqrt d\,\varrho_r(d)$ under
$\|X_i\|_\infty \leq M$ or
$\vvvert e_j^\T X_i\vvvert_r\leq M$.

\subsubsection{Wasserstein-based couplings}
\label{sec:literature_details_wasserstein}

\citet[Theorem~1.1]{zhai2018high} showed that if
$X_i$ are independent with $\E(X_i) = 0_d$
and $\Var(X_i) = \Sigma$ and satisfy $\|X_i\|_2 \leq M$ almost surely,
\begin{equation*}
  W_{2,2}\biggl(
    \frac{S_n}{\sqrt n}, N(0_d, \Sigma)
  \biggr)
  \leq \frac{5 \sqrt d M (1 + \log n)}{\sqrt n}.
\end{equation*}
For $p \in [2,\infty]$, taking an optimal coupling for $W_{2,2}$
with $T_n' \sim N(0_d,n\Sigma)$, using $\|v\|_p \leq \|v\|_2$ and
applying Markov's inequality yields, for every $\eta>0$,
\begin{equation*}
  \P\biggl(
    \frac{\|S_n'-T_n'\|_p}{\sqrt n} > \eta
  \biggr)
  \leq
  \frac{25 M^2 d (1 + \log n)^2}{n \eta^2}.
\end{equation*}
This gives a sample size factor of $n^{-1/2} \log 2n$ and a dimension
factor of $\sqrt d$. If instead $\|X_i\|_\infty \leq M$ almost surely,
then $\|X_i\|_2 \leq \sqrt d M$ and the dimension factor is $d$.

The author also gave a lattice-based lower bound.
For a fixed summand law supported on $M \Z^d$
with $\Sigma$ nonsingular and satisfying $\|X_i\|_2 \leq M$,
\citet[Proposition~1.2]{zhai2018high} gives
\begin{equation*}
  \liminf_{n\to\infty}
  \sqrt n\,
  W_{2,2}\biggl(\frac{S_n}{\sqrt n}, N(0_d, \Sigma)\biggr)
  \geq
  \frac{\sqrt d M}{4}.
\end{equation*}

\citet[Theorem~1]{bonis2020stein} gives,
for i.i.d.\ vectors with $\E(X_i) = 0_d$ and $\Var(X_i) = I_d$,
under finite fourth moments,
\begin{equation*}
  W_{2,2}\biggl(
    \frac{S_n}{\sqrt n}, N(0_d, I_d)
  \biggr)
  \leq
  \frac{C d^{1/4}}{\sqrt n}
  \bigl\|
  \E\bigl(X_1 X_1^\T \|X_1\|_2^2\bigr)
  \bigr\|_{\Fr, 2}^{1/2}.
\end{equation*}
More generally, if $\Var(X_i) = \Sigma$ with
$\lambda_{\max}(\Sigma) / \lambda_{\min}(\Sigma) \leq c^2$,
where $c \geq 1$ is fixed,
applying this bound to $\Sigma^{-1/2} X_i$ gives
\begin{equation*}
  W_{2,2}\biggl(
    \frac{S_n}{\sqrt n}, N(0_d, \Sigma)
  \biggr)
  \leq
  \frac{C c d^{1/4}}{\sqrt{n \lambda_{\min}(\Sigma)}}
  \bigl\|
  \E\bigl(X_1 X_1^\T \|X_1\|_2^2\bigr)
  \bigr\|_{\Fr, 2}^{1/2}.
\end{equation*}
For $p \in [2,\infty]$, taking an optimal coupling for $W_{2,2}$
with $T_n' \sim N(0_d, n \Sigma)$ and
applying Markov's inequality yields, for every $\eta>0$,
\begin{equation*}
  \P\biggl(
    \frac{\|S_n'-T_n'\|_p}{\sqrt n} > \eta
  \biggr)
  \leq
  \frac{C^2 c^2 d^{1/2}}{n \lambda_{\min}(\Sigma) \eta^2}
  \bigl\|
  \E\bigl(X_1 X_1^\T \|X_1\|_2^2\bigr)
  \bigr\|_{\Fr, 2}.
\end{equation*}
The sample size factor is therefore $n^{-1/2}$.
If $\|X_1\|_2 \leq M$, we have
$\E\bigl(X_1 X_1^\T \|X_1\|_2^2\bigr) \preceq M^2 \Sigma$,
giving a dimension factor of $d^{1/2}$.
If $\|X_1\|_\infty \leq M$, then
$\|X_1\|_2 \leq \sqrt d M$, so the dimension factor is $d$.
If $\vvvert e_j^\T X_i\vvvert_4 \leq M$
and $M^2 \lesssim \lambda_{\min}(\Sigma)$, then
$\bigl\|\E\bigl(X_1 X_1^\T \|X_1\|_2^2\bigr)\bigr\|_{\Fr, 2} \leq d^2 M^4$,
and the factor is $d^{5/4}$.

\citet[Theorem~1]{eldan2020clt} established that
for i.i.d.\ $X_i$ with $\E(X_i) = 0_d$, $\Var(X_i) = \Sigma$
and $\|X_i\|_2 \leq M$ almost surely,
\begin{equation*}
  W_{2,2}\biggl(\frac{S_n}{\sqrt n}, N(0_d, \Sigma)\biggr)
  \leq
  \frac{M \sqrt d \sqrt{32 + 2\log_2 n}}{\sqrt n},
\end{equation*}
thereby removing a factor of $\sqrt{\log 2n}$ from
the corresponding bound due to \citet{zhai2018high}.
The same calculations then yield an improved sample size factor of
$n^{-1/2} \sqrt{\log 2n}$.

\citet[Theorem~2.1]{mies2023sequential}
derived a Gaussian coupling assuming $X_i$ independent
with $\E(X_i) = 0_d$ and
$\bigl\{\E\bigl(\|X_i\|_2^r\bigr)\bigr\}{}^{1/r} \leq M_i$,
where $r > 2$, with
\begin{equation*}
  \biggl\{
    \E\biggl(
      \frac{\|S_n' - T_n'\|_2^2}{n}
    \biggr)
  \biggr\}^{1/2}
  \leq
  \frac{C}{\sqrt{r-2}\land1}
  \biggl(\frac d n\biggr)^{1/2-1/r}\sqrt{\log n}
  \Biggl(\frac 1 n \sum_{i=1}^n M_i^2\Biggr)^{1/2},
\end{equation*}
where $C > 0$ is universal.
For $p \in [2, \infty]$, apply $\|v\|_p \leq \|v\|_2$.
The sample size factor is therefore $n^{-1/2 + 1/r} \sqrt{\log 2n}$.
If $\|X_i\|_2 \leq M$ almost surely,
take $M_i = M$ so the dimension factor is $d^{1/2 - 1/r}$.
If $\|X_i\|_\infty \leq M$ almost surely,
then $\|X_i\|_2 \leq \sqrt d M$, yielding a factor of $d^{1 - 1/r}$.
For $\vvvert e_j^\T X_i\vvvert_r \leq M$,
we take $M_i = \sqrt d M$ also.
At $r = \infty$, we may obtain rates by
taking $r \to \infty$ in the above display.

\citet[Corollary~2.1]{fang2026martingale} showed that
for independent $X_i$ with $\E(X_i) = 0_d$,
writing $\Sigma_i \coloneqq \Var(X_i)$ and
$\Sigma \coloneqq \Var(S_n) / n$ and
assuming $\Sigma \succ 0$ and $\Sigma_i \succeq c \Sigma$
for some fixed $c > 0$,
\begin{equation*}
  W_{s, p}
  \biggl(\frac{S_n}{\sqrt n}, N(0_d, \Sigma)\biggr)
  \leq
  \frac{C s}{c \sqrt n}
  \sum_{i=1}^n
  \frac{
    \bigl\{
      \E\bigl(\|X_i\|_p^{3s}\bigr)
    \bigr\}^{1/(3s)}
    \bigl\{
      \E\bigl(\|\Sigma^{-1/2} X_i\|_2^{3s}\bigr)
    \bigr\}^{2/(3s)}
  }{n - i + 1},
\end{equation*}
where $s \geq 1$ and $C > 0$ is universal.
Taking $s = 1$ and using
$\|\Sigma^{-1/2} X_i\|_2 \leq \lambda^{-1/2} \|X_i\|_2$,
we have for $p \in [2, \infty]$ that
\begin{align*}
  &
  \bigl\{\E\bigl(\|X_i\|_p^3\bigr)\bigr\}^{1/3}
  \bigl\{\E\bigl(\|\Sigma^{-1/2} X_i\|_2^3\bigr)\bigr\}^{2/3}
  \iftoggle{journal}{\\ &\qquad}{}
  \leq
  \frac{M^3}{\lambda}
  \times
  \begin{cases}
    d^{1+1/p}, & \text{if } \|X_i\|_\infty \leq M \text{ almost surely},\\
    1, & \text{if } \|X_i\|_2 \leq M \text{ almost surely},\\
    d^{1+1/(p \land r)}, &
    \text{if }
    \vvvert e_j^\T X_i\vvvert_r \leq M
    \text{ for all } j \in [d],
  \end{cases}
\end{align*}
where $r \geq 3$.
This gives the dimension factors presented in
the corresponding row of Table~\ref{tab:literature}.
Since $\sum_{i=1}^n (n-i+1)^{-1} \leq 1 + \log n$,
the sample size factor is $n^{-1/2} \log 2n$.

In the same work,
\citet[Corollary~2.2]{fang2026martingale} show that
if the conditions of $\Sigma \succ 0$ and $\Sigma_i \succeq c \Sigma$
are dropped, then one instead obtains
\begin{align*}
  &W_{s,p}\biggl(\frac{S_n}{\sqrt n}, N(0_d,\Sigma)\biggr)
  \iftoggle{journal}{\\ &\qquad}{}
  \leq
  C\Biggl\{
    \frac{s}{n^{3/2}}
    \sum_{i=1}^n
    \bigl\{\E\bigl(\|X_i\|_p^{3s}\bigr)\bigr\}^{1/(3s)}
    \bigl\{\E\bigl(\|X_i\|_2^{3s}\bigr)\bigr\}^{2/(3s)}
  \Biggr\}^{1/3}
  \bigl\{\E\bigl(\|Z\|_p^s\bigr)\bigr\}^{2/(3s)},
\end{align*}
where $Z \sim N(0_d, I_d)$, $s \geq 1$ and
$C > 0$ is universal.
Taking $s = 1$ and $p \in [2,\infty]$,
the same moment bounds as above give
\begin{equation*}
  \bigl\{\E\bigl(\|X_i\|_p^3\bigr)\bigr\}^{1/3}
  \bigl\{\E\bigl(\|X_i\|_2^3\bigr)\bigr\}^{2/3}
  \leq M^3
  \begin{cases}
    d^{1+1/p}, & \text{if } \|X_i\|_\infty \leq M \text{ almost surely},\\
    1, & \text{if } \|X_i\|_2 \leq M \text{ almost surely},\\
    d^{1+1/(p \land r)}, & \text{if }
    \vvvert e_j^\T X_i\vvvert_r \leq M
    \text{ for all } j \in [d],
  \end{cases}
\end{equation*}
where $r \geq 3$.
We obtain a sample size factor of $n^{-1/6}$.
By the proof of Lemma~\ref{lem:gaussian_norm_tail},
$\E(\|Z\|_p) \leq \sqrt p\, d^{1/p}$ for $p \in [2, \infty)$ and
$\E(\|Z\|_\infty) \leq \sqrt{2\log 2d}$,
yielding the dimension factors presented in the
relevant row of Table~\ref{tab:literature}.

For independent summands with $\E(X_i) = 0_d$
and $\|X_i\|_2 \leq M$ almost surely,
applying \citet[Corollary~2.1]{fang2026sharp}
on the range of $\Sigma$ gives
\begin{equation*}
  W_{2,2}\biggl(\frac{S_n}{\sqrt n}, N(0_d,\Sigma)\biggr)
  \leq 4 M \sqrt{\frac d n}.
\end{equation*}
Taking an optimal coupling yields a sample size factor of $n^{-1/2}$
and a dimension factor of $d^{1/2}$.
Under $\|X_i\|_\infty \leq M$, we have
$\|X_i\|_2 \leq \sqrt d M$ and the dimension factor is $d$.
If $\vvvert e_j^\T X_i\vvvert_r \leq M$, suppose
$\lambda = \lambda_{\min}\bigl(\sum_{i=1}^n \Var(X_i) / n\bigr) > 0$.
By \citet[Theorem~2.1]{fang2026sharp},
with $X_i'$ an independent copy of $X_i$,
\begin{equation*}
  W_{2,2}\biggl(\frac{S_n}{\sqrt n}, N(0_d,\Sigma)\biggr)^2
  \leq \frac{2}{n^2\lambda}
  \sum_{i=1}^n \E\bigl(\|X_i'-X_i\|_2^4\bigr)
  \leq \frac{32}{n^2\lambda}
  \sum_{i=1}^n \E\bigl(\|X_i\|_2^4\bigr).
\end{equation*}
For $r \geq 4$, we have $\E(\|X_i\|_2^4) \leq d^2 M^4$,
and hence the dimension factor is $\lambda^{-1/2} d$.

In the same paper,
for i.i.d.\ summands with $\E(X_i) = 0_d$ and finite fourth moments,
\citet[Theorem~2.3]{fang2026sharp} gives a Gaussian coupling satisfying
\begin{equation*}
  \biggl\{
    \E\biggl(
      \frac{\|S_n' - T_n'\|_2^4}{n^2}
    \biggr)
  \biggr\}^{1/4}
  \leq
  \frac{C\log 2n}{n^{1/4}}
  \bigl\{\tr\bigl(Q^{1/2}\bigr)\bigr\}^{1/2},
\end{equation*}
where $C > 0$ is universal and
$Q = \E\bigl\{\|X_1' - X_1\|_2^2 (X_1' - X_1)(X_1' - X_1)^\T\bigr\}$,
with $X_1'$ an independent copy of $X_1$.
By Cauchy--Schwarz,
\begin{equation*}
  \bigl\{\tr\bigl(Q^{1/2}\bigr)\bigr\}^{1/2}
  \leq (d \tr Q)^{1/4}
  = \bigl\{d \, \E\bigl(\|X_1' - X_1\|_2^4\bigr)\bigr\}^{1/4}
  \leq 2 \bigl\{d \, \E\bigl(\|X_1\|_2^4\bigr)\bigr\}^{1/4}.
\end{equation*}
For $p \in [2,\infty]$, the sample size factor is $n^{-1/4}\log 2n$.
If $\|X_i\|_2 \leq M$ almost surely, the dimension factor is $d^{1/4}$.
If $\|X_i\|_\infty \leq M$ almost surely or
$\vvvert e_j^\T X_i\vvvert_r \leq M$ for all $j \in [d]$ with $r \geq 4$,
the dimension factor is $d^{3/4}$.

\subsection{Proofs for Section~\ref{sec:empirical_process}}

We provide a useful maximal inequality
for random variables with finite sub-Weibull Orlicz norm.

\begin{lemma}[Moments and tail probabilities under Orlicz norm bounds]%
  \label{lem:moment_tail_orlicz}
  Let $X_1, \ldots, X_n$ be non-negative random variables.
  Suppose that there exist $r \in (0, \infty)$ and $M > 0$
  such that $\vvvert X_i \vvvert_{\omega_r} \leq M$
  for each $i \in [n]$.
  Then with $M_n \coloneqq \max_{1 \leq i \leq n} X_i$,
  for $t > 0$, we have
  \begin{equation*}
    \P(M_n > t)
    \leq
    2 n \exp\biggl(
      - \frac{t^r}{M^r}
    \biggr).
  \end{equation*}
  Further, if $m > 0$ and $\log 2n > 2m/r$, then
  \begin{equation*}
    \E\bigl(M_n^m\bigr)
    \leq
    2 M^m (\log 2n)^{m/r}.
  \end{equation*}
  It follows that there exists $C > 0$
  depending only on $r$ and $m$ such that for all $n \in \N$,
  \begin{equation*}
    \E\bigl(M_n^m\bigr)
    \leq
    C M^m (\log 2n)^{m/r}.
  \end{equation*}
  If instead $\vvvert X_i \vvvert_\infty \leq M$ for each $i \in [n]$,
  then $M_n \leq M$ almost surely and $\E\bigl(M_n^m\bigr) \leq M^m$.
\end{lemma}

\begin{proof}[Proof of Lemma~\ref{lem:moment_tail_orlicz}]
  By Markov's inequality, a union bound and definition of
  the Orlicz norm,
  \begin{align*}
    \P(M_n > t)
    &\leq
    \sum_{i=1}^n
    \P\bigl(X_i > t\bigr)
    =
    \sum_{i=1}^n
    \P\biggl\{
      \exp\biggl(
        \frac{X_i^r}{M^r}
      \biggr)
      >
      \exp\biggl(
        \frac{t^r}{M^r}
      \biggr)
    \biggr\} \\
    &\leq
    \exp\biggl(
      - \frac{t^r}{M^r}
    \biggr)
    \sum_{i=1}^n
    \E\biggl\{
      \exp\biggl(
        \frac{X_i^r}{M^r}
      \biggr)
    \biggr\}
    \leq 2 n
    \exp\biggl(
      - \frac{t^r}{M^r}
    \biggr).
  \end{align*}
  As $\Gamma(a, x) \coloneqq \int_x^\infty t^{a-1} e^{-t} \diffi t
  \leq 2 x^{a-1} e^{-x}$
  for $a \geq 1$ and $x > 2(a - 1)$ by
  \citet[Proposition~2.7]{pinelis2020exact},
  we may deduce that if $\log 2n > 2m/r$, then
  by substituting $u \coloneqq t^{r/m} / M^r$,
  \begin{align*}
    \E\bigl(M_n^m\bigr)
    &=
    \int_0^\infty
    \P\bigl(M_n^m > t\bigr)
    \diffi t
    \leq
    \int_0^\infty
    \biggl\{
      2 n
      \exp\biggl(
        - \frac{t^{r/m}}{M^r}
      \biggr)
    \biggr\}
    \land 1
    \diffi t \\
    &=
    M^m (\log 2n)^{m/r}
    + \frac{2 n m M^m}{r}
    \int_{\log 2n}^\infty
    u^{m/r - 1} e^{-u}
    \diffi u \\
    &\leq
    M^m (\log 2n)^{m/r}
    + \frac{2 n m M^m}{r \log 2n}
    \int_{\log 2n}^\infty
    u^{m/r} e^{-u}
    \diffi u \\
    &=
    M^m (\log 2n)^{m/r}
    + \frac{2 n m M^m}{r \log 2n}
    \Gamma\biggl(\frac{m}{r} + 1, \log 2n\biggr) \\
    &\leq
    M^m (\log 2n)^{m/r}
    + \frac{4 n m M^m}{2 n r \log 2n}
    (\log 2n)^{m/r}
    \leq
    2 M^m (\log 2n)^{m/r}.
  \end{align*}
  For the last result, if $\log 2n \leq 2 m / r$ then
  $\E(M_n^m) \leq \sum_{i=1}^n \E(X_i^m) \leq C M^m$.
\end{proof}

\begin{proof}[Proof of Theorem~\ref{thm:empirical_process}]
  Define the random vectors
  \begin{align*}
    \tilde S_n
    &\coloneqq \sum_{i=1}^{n}
    \bigl\{ \tilde\cF(X_i) - P_{X,i}(\tilde\cF) \bigr\},
    &\tilde T_n
    &\coloneqq \sum_{i=1}^{n}
    \bigl\{ \tilde\cF(Y_i) - P_{Y,i}(\tilde\cF) \bigr\}.
  \end{align*}
  By Theorem~\ref{thm:yurinskii},
  for each $\eta, t > 0$,
  there is a probability space
  supporting $\tilde S_n'$ and $\tilde T_n'$
  where $\tilde S_n' \equallaw \tilde S_n$
  and $\tilde T_n' \equallaw \tilde T_n$ with
  \begin{equation}
    \label{eq:empirical_process_yurinskii}
    \P\bigl(\|\tilde S_n' - \tilde T_n'\|_\infty > 3\eta\bigr)
    \leq
    2 e^{-t^2 / 2}
    + \frac{n \bar\beta_{k,\gamma} \phi(t)^k}{\eta^{k+\gamma}}
    + \Delta\biggl(\frac{\phi(t)}{\eta}\biggr)
    + \Gamma_\eta\biggl(\frac{\phi(t)}{\eta}\biggr).
  \end{equation}
  By two applications of the Vorob'ev--Berkes--Philipp theorem
  \citep[Theorem~1.31]{dudley2014uniform},
  there exist random variables
  $\bX_n' = (X_1', \ldots, X_n')$ and
  $\bY_n' = (Y_1', \ldots, Y_n')$ such that
  $(\tilde S_n', \bX_n') \equallaw (\tilde S_n, \bX_n)$ and
  $(\tilde T_n', \bY_n') \equallaw (\tilde T_n, \bY_n)$,
  where $\bX_n \coloneqq (X_1, \ldots, X_n)$
  and $\bY_n \coloneqq (Y_1, \ldots, Y_n)$.
  Define $S_n'(f) \coloneqq \sum_{i=1}^n \{f(X_i') - P_{X,i}(f)\}$
  and $T_n'(f) \coloneqq \sum_{i=1}^n \{f(Y_i') - P_{Y,i}(f)\}$
  for each $f \in \cF$, so
  $S_n' \equallaw S_n$ and $T_n' \equallaw T_n$.
  Moreover, $\tilde S_n' = S_n'(\tilde\cF)$ and
  $\tilde T_n' = T_n'(\tilde\cF)$ almost surely.
  By the triangle inequality,
  \begin{equation}
    \label{eq:empirical_process_decomp}
    \sup_{f \in \cF} |S_n'(f) - T_n'(f)|
    \leq
    \|\tilde S_n' - \tilde T_n'\|_\infty
    + \sup_{f \in \check \cF} |S_n'(f)|
    + \sup_{f \in \check \cF} |T_n'(f)|,
  \end{equation}
  observing that the first term on the right-hand side of
  \eqref{eq:empirical_process_decomp}
  has already been controlled
  in \eqref{eq:empirical_process_yurinskii}.
  For the second term, since $\check\cF$ admits a
  countable dense subclass under pointwise convergence by
  Assumption~\ref{ass:function_class},
  applying Theorem~\ref{thm:empirical_tail} with
  $\kappa$ replaced with $2\kappa$ yields
  \begin{equation*}
    \sup_{f \in \check\cF} |S_n'(f)|
    \leq C
    \biggl\{
      \sigma \sqrt n +
      \check M (\log 2n)^{1/r}
      \biggl(
        \frac{\cJ(\nu)}{\nu} + \kappa^{1/r} t
      \biggr)
    \biggr\}
    \biggl( \frac{\cJ(\nu)}{\nu} + t \biggr)
    \leq C R(t),
  \end{equation*}
  with probability at least
  $1 - 1/n^{2\kappa} - e^{-t^2/2}$,
  where $C > 0$ depends only on $r$.
  The process $T_n'$ satisfies the same assumptions as $S_n'$,
  so admits the same bound.
  Therefore, as $1 \land (2/n^{2\kappa}) \leq 1/n^\kappa$,
  by a union bound,
  increasing $C$ if necessary, we deduce that
  \begin{align*}
    &
    \P\biggl\{
      \sup_{f \in \cF} |S_n'(f) - T_n'(f)|
      > 3 \eta
      + C R(t)
    \biggr\}
    \iftoggle{journal}{\\&\qquad}{}
    \leq
    4 e^{-t^2 / 2}
    + \frac{n \bar\beta_{k,\gamma} \phi(t)^k}{\eta^{k+\gamma}}
    + \Delta\biggl(\frac{\phi(t)}{\eta}\biggr)
    + \Gamma_\eta\biggl(\frac{\phi(t)}{\eta}\biggr)
    + \frac{1}{n^\kappa},
  \end{align*}
  as required.
\end{proof}

\begin{lemma}[Entropy of VC classes]%
  \label{lem:entropy_vc}
  Suppose that a measurable function
  $N: (0, 1] \to [1, \infty)$ satisfies
  $N(\varepsilon) \leq a \varepsilon^{-u}$
  for all $\varepsilon \in (0, 1]$, where
  $a \geq 1$ and $u \geq 0$. Then,
  for all $\varepsilon \in (0, 1]$,
  \begin{equation*}
    \cJ_N(\varepsilon)
    \coloneqq
    \int_0^\varepsilon
    \sqrt{1 + \log N(\xi)}
    \diffi \xi
    \leq
    \varepsilon
    \sqrt{1 + \log a  + u + u \log(1/\varepsilon)}.
  \end{equation*}
\end{lemma}

\begin{proof}[Proof of Lemma~\ref{lem:entropy_vc}]
  By Jensen's inequality, if $u > 0$,
  \begin{align*}
    \cJ_N(\varepsilon)
    &\leq
    \int_0^\varepsilon
    \sqrt{1 + \log a - u \log \xi}
    \diffi \xi
    =
    \sqrt u
    \int_0^\varepsilon
    \sqrt{(1 + \log a) / u - \log \xi}
    \diffi \xi \\
    &=
    \varepsilon \sqrt u
    \int_0^1
    \sqrt{(1 + \log a) / u - \log \varepsilon - \log \eta}
    \diffi \eta
    \iftoggle{journal}{\\&}{}
    \leq
    \varepsilon \sqrt u
    \biggl\{
      (1 + \log a) / u - \log \varepsilon -
      \int_0^1
      \log \eta
      \diffi \eta
    \biggr\}^{1/2} \\
    &=
    \varepsilon \sqrt u
    \sqrt{ (1 + \log a) / u - \log \varepsilon + 1 }
    =
    \varepsilon
    \sqrt{1 + \log a  + u + u \log(1/\varepsilon)}.
  \end{align*}
  Otherwise, if $u = 0$ then
  $\cJ_N(\varepsilon) \leq \varepsilon \sqrt{1 + \log a }$,
  as required.
\end{proof}

\begin{proof}[Proof of Corollary~\ref{cor:empirical_process_vc}]
  Let $C > 0$ depend only on the stated quantities;
  its value may vary from line to line.
  Set $\check F \coloneqq 2 F$,
  let $\nu \coloneqq n^{(2 - k - \gamma)/(2k + 2\gamma + u k)}$,
  and define
  \begin{equation*}
    \sigma
    \coloneqq
    \nu
    \bigl\{
      \|\check F\|_{\bar P_X,2}
      \lor
      \|\check F\|_{\bar P_Y,2}
    \bigr\}
    \leq
    2 \nu
    \bigl\{
      \|F\|_{\bar P_X,2}
      \lor
      \|F\|_{\bar P_Y,2}
    \bigr\}
    \leq
    C \nu M,
  \end{equation*}
  where the last inequality follows by
  Lemma~\ref{lem:moment_tail_orlicz}.
  By Assumption~\ref{ass:vc},
  \begin{align*}
    &N\bigl(\sigma, \cF,
      \|\cdot\|_{\bar P_X,2} \lor \|\cdot\|_{\bar P_Y,2}
    \bigr)
    \leq
    N\bigl(\sigma, \cF,
      \|\cdot\|_{\bar P_X + \bar P_Y,2}
    \bigr)
    =
    N\bigl(\sigma / \sqrt 2, \cF,
      \|\cdot\|_{\bar P_X/2 + \bar P_Y/2,2}
    \bigr) \\
    &\qquad\leq
    a \biggl(
      1 \lor
      \frac{\sqrt 2 \|F\|_{\bar P_X/2 + \bar P_Y/2,2}}
      {\sigma}
    \biggr)^u
    =
    a \biggl(
      1 \lor
      \frac{\|F\|_{\bar P_X + \bar P_Y,2}}
      {\sigma}
    \biggr)^u
    =
    a \biggl(
      1 \lor
      \frac{\|\check F\|_{\bar P_X + \bar P_Y,2}}
      {2 \sigma}
    \biggr)^u \\
    &\qquad\leq
    a \biggl(
      1 \lor
      \frac{\|\check F\|_{\bar P_X, 2}
      \lor \|\check F\|_{\bar P_Y, 2}}
      {\sqrt 2 \sigma}
    \biggr)^u
    \leq
    a \nu^{-u},
  \end{align*}
  as $\nu \leq 1$.
  Therefore, we may choose $\tilde\cF$ and
  $\pi: \cF \to \tilde\cF$ so that
  $d \leq a \nu^{-u}$ and $\check\cF$ is pointwise measurable.
  Here, pointwise measurability can be established
  by first defining $\pi$ on a countable pointwise dense class,
  then extending it to $\cF$ along subsequences on which
  $\pi$ is constant.
  For $\varepsilon \in (0, 1]$ and $Q \in \cQ$,
  by the internal--external covering inequality
  \citep[Exercise~4.2.9]{vershynin2019high}
  and Assumption~\ref{ass:vc},
  \begin{align*}
    N\bigl(
      \varepsilon \|\check F\|_{Q,2},
      \check \cF,
      \|\cdot\|_{Q,2}
    \bigr)
    &=
    N\bigl(
      2 \varepsilon \|F\|_{Q,2},
      \check \cF,
      \|\cdot\|_{Q,2}
    \bigr) \\
    &\leq
    N\bigl(
      \varepsilon \|F\|_{Q,2} / 2,
      \cF,
      \|\cdot\|_{Q,2}
    \bigr)^2
    \leq
    a^2 (\varepsilon / 2)^{-2u}.
  \end{align*}
  Since $\check F = 2 F$, we may set $\check M \coloneqq 2 M$.
  By Lemma~\ref{lem:entropy_vc},
  with $\kappa \coloneqq 2k + 2\gamma \leq 4 k$, we have
  \begin{align*}
    R(t)
    &\leq
    C \Bigl\{
      \sigma \sqrt n
      + M (\log 2n)^{1/r}
      \Bigl(\sqrt{\log(2/\nu)} + t\Bigr)
    \Bigr\}
    \Bigl(\sqrt{\log(2/\nu)} + t\Bigr) \\
    &\leq
    C M
    \Bigl\{
      \nu \sqrt n + (\log 2n)^{1/r}
      \Bigl(\sqrt{\log 2n} + t\Bigr)
    \Bigr\}
    \Bigl(\sqrt{\log 2n} + t\Bigr) \\
    &\leq
    C M
    \Bigl\{
      \nu \sqrt{n \log 2n}
      + t \nu \sqrt n
      + t^2 (\log 2n)^{1/r}
    \Bigr\},
  \end{align*}
  where we used
  $\nu \sqrt {n} = n^{(4 + u k)/(4k + 4\gamma + 2 u k)} \geq n^{1/(k+1)}$
  and $(\log 2n)^{1/r+1/2} \leq C \nu \sqrt n$.
  By H{\"o}lder's inequality and
  Lemma~\ref{lem:moment_tail_orlicz},
  using $\vvvert F(X_i)\vvvert_{k+\gamma}
  \leq\vvvert F(X_i)\vvvert_{k+1}$ and likewise for $Y_i$,
  we have
  \begin{equation}
    \label{eq:vc_beta_bound}
    \bar\beta_{k,\gamma}
    \leq
    \frac{2^{k + \gamma} d^{k/2}}{\sqrt{k!} n}
    \sum_{i=1}^n
    \E \Bigl\{
      F(X_i)^{k + \gamma} + F(Y_i)^{k + \gamma}
    \Bigr\}
    \leq
    \frac{C M^{k+\gamma}}{\nu^{u k/2}}.
  \end{equation}
  Therefore, in the notation of
  Theorem~\ref{thm:empirical_process},
  \begin{equation*}
    \frac{n \bar\beta_{k,\gamma} \phi(t)^k}{\eta^{k+\gamma}}
    \leq
    \frac{C n M^{k+\gamma}}
    {\nu^{u k/2} \eta^{k+\gamma}}
    \bigl(\log 2n + t^2\bigr)^{k/2}.
  \end{equation*}
  Since $\delta_{\alpha,i} = 0$ for each
  $i \in [n]$ and
  $\alpha \in \N_0^d$ with $|\alpha| \in [k]$,
  we have $\Delta(s) = \Gamma_\eta(s) = 0$ for each $s > 0$.
  By Theorem~\ref{thm:empirical_process}
  with $m = 1$, we obtain
  \begin{align*}
    &\P\biggl\{
      \sup_{f \in \cF} \bigl|S_n'(f) - T_n'(f)\bigr|
      > 3 \eta
      + C M
      \Bigl\{
        \nu \sqrt{n \log 2n}
        + t \nu \sqrt n
        + t^2 (\log 2n)^{1/r}
      \Bigr\}
    \biggr\} \\
    &\qquad\leq
    4 e^{-t^2 / 2}
    + \frac{C n M^{k+\gamma}}
    {\nu^{u k/2} \eta^{k + \gamma}}
    \bigl(\log 2n + t^2\bigr)^{k/2}
    + \frac{1}{n^\kappa}.
  \end{align*}
  Replacing $\eta$ with
  $C n^{\frac{4 + u k}{4k + 4\gamma + 2u k}} \sqrt{\log 2n} M \eta$
  and taking $t \coloneqq \sqrt{2 \kappa \log 2n}$,
  we have
  \begin{equation*}
    \P\biggl\{
      \sup_{f \in \cF} |S_n'(f) - T_n'(f)|
      >
      C M \eta
      n^{\frac{4 + u k}{4k + 4\gamma + 2u k}}
      \sqrt{\log 2n}
    \biggr\}
    \leq
    \frac{1}{\eta^{k + \gamma}}
    + \frac{2}{n^\kappa},
  \end{equation*}
  noting that this inequality holds trivially if $\eta \leq 1$.

  Finally, we remove the $2/n^\kappa$ term.
  By Lemma~\ref{lem:moment_tail_orlicz},
  if $r < \infty$, then for $t > 0$,
  with probability at least $1 - 2n \exp(-t^r)$,
  \begin{equation*}
    \sup_{f \in \cF} |S_n'(f)|
    \leq
    \sum_{i=1}^n F(X_i')
    + \sum_{i=1}^n P_{X,i}(F)
    \leq
    (t + C) n M.
  \end{equation*}
  The same holds for $T_n'$.
  Therefore, if $t \geq 1$,
  \begin{equation*}
    \P\biggl\{
      \sup_{f \in \cF} |S_n'(f) - T_n'(f)|
      >
      C n^2 M t
    \biggr\}
    \leq
    4 n \exp\bigl\{-(n t)^r\bigr\}
    \leq
    \frac{C}{(n t)^{2k + 2}}
    \leq
    \frac{C}{(n^2 t)^{k + 1}}.
  \end{equation*}
  This bound also holds for $r = \infty$.
  It follows that for $\eta \geq n^2$,
  \begin{equation*}
    \P\biggl\{
      \sup_{f \in \cF} |S_n'(f) - T_n'(f)|
      >
      C M \eta
      n^{\frac{4 + u k}{4k + 4\gamma + 2u k}}
      \sqrt{\log 2n}
    \biggr\}
    \leq
    \frac{C}{\eta^{k + \gamma}}.
  \end{equation*}
  On the other hand, for $\eta \leq n^2$,
  since $\kappa = 2k + 2\gamma$,
  we have $1/\eta^{k+\gamma} \geq 1/n^{2k+2\gamma} = 1/n^\kappa$.
  Increasing $C$, we have that for all $\eta > 0$
  and $n \in \N$,
  \begin{equation*}
    \P\biggl\{
      \sup_{f \in \cF} |S_n'(f) - T_n'(f)|
      >
      C M \eta
      n^{\frac{4 + u k}{4k + 4\gamma + 2u k}}
      \sqrt{\log 2n}
    \biggr\}
    \leq
    \frac{1}{\eta^{k + \gamma}}.
    \iftoggle{journal}{}{\qedhere}
  \end{equation*}
\end{proof}

\begin{proof}[Proof of Theorem~\ref{thm:bootstrap_vc}]
  Let $C > 0$ depend only on the stated quantities;
  its value may vary from line to line.
  Set $\check F \coloneqq 2F$ and $\check M \coloneqq 2M$.
  Let $\nu \coloneqq n^{-\epsilon}$, where $\epsilon \in (0, 1/2)$
  depends only on $k$ and $u$, and define
  $\sigma \coloneqq \nu \|\check F\|_{\bar P_X,2}$.
  Construct $\tilde \cF \subseteq \cF$ and $\pi$ as
  in the proof of Corollary~\ref{cor:empirical_process_vc},
  using only $\bar P_X$,
  so that $d \coloneqq |\tilde\cF| \leq C \nu^{-u}$ and $\check\cF$ is
  pointwise measurable, with
  $\sup_{f \in \cF} \|f - \pi(f)\|_{\bar P_X,2} \leq \sigma$.
  By the same calculations as
  in the proof of Corollary~\ref{cor:empirical_process_vc},
  taking $\kappa = 4$, we have
  \begin{equation*}
    R^*(t)
    \leq
    C M
    \Bigl\{
      \nu \sqrt{n \log 2n}
      + t \nu \sqrt n
      + t^2 (\log 2n)^{1/r + 1/r'}
    \Bigr\}.
  \end{equation*}
  Let $\gamma = 1$ and $s, \eta > 0$.
  In the notation of Theorem~\ref{thm:empirical_process_bootstrap},
  by Lemma~\ref{lem:moment_tail_orlicz},
  \begin{equation*}
    \frac{n s^k}{\eta}\E\bigl\{\hat\beta_{k,1}\bigr\}
    \leq
    \frac{C n d^{k/2} M^{k+1} s^k}
    {\eta}
    \leq
    \frac{C n M^{k+1} s^k}
    {\nu^{u k/2} \eta}.
  \end{equation*}
  For $a+b+c \leq 4$ and $i \in [n]$,
  in the notation of Theorem~\ref{thm:empirical_process_bootstrap},
  by Jensen's inequality, H{\"o}lder's inequality
  and Lemma~\ref{lem:moment_tail_orlicz}, we have
  \begin{equation*}
    \E\bigl(\hat\beta_{a,b,i}\bigr) \leq C d^{a/2} M^{a+b},
    \qquad
    \E\bigl(\hat\beta_{a,b,i}\hat\beta_{c,0,i}\bigr)
    \leq C d^{(a+c)/2} M^{a+b+c},
  \end{equation*}
  with corresponding bounds holding for
  $\hat\beta_{a,b}$ and $\hat\beta_{a,b}\hat\beta_{c,0}$.
  For $\ell \in \{2, k\}$, define $k_\ell \coloneqq 1 \land (k - \ell)$.
  Thus, if $n^{1/4} s \sqrt d M \leq 1$,
  then arguing as in the proof of Theorem~\ref{thm:bootstrap},
  under either condition~\ref{it:thm_ep_bootstrap_vc_nonparametric}
  or~\ref{it:thm_ep_bootstrap_vc_multiplier}, we have
  \begin{align*}
    \E\bigl\{ \hat\Delta(s) \bigr\}
    &\leq
    C s^2 d \sqrt{n} M^2
    \leq
    C s^2 \nu^{-u} \sqrt{n} M^2,
    &\E\bigl\{ \hat\Gamma_\eta(s) \bigr\}
    &\leq
    \frac{C n M^{k+1} s^k}{\nu^{u k/2} \eta}.
  \end{align*}
  Thus, taking $s \coloneqq \phi_{\infty,d}(t) / \eta$,
  by Theorem~\ref{thm:empirical_process_bootstrap} with $m = 1$,
  \begin{align*}
    &\P\biggl(
      \sup_{f \in \cF} |S_n'(f) - S_n^{*\prime}(f)|
      > 3 \eta
      + C M
      \Bigl\{
        \nu \sqrt{n \log 2n}
        + t \nu \sqrt n
        + t^2 (\log 2n)^{1/r + 1/r'}
      \Bigr\}
    \biggr) \\
    &\qquad\leq
    4 e^{-t^2 / 2}
    + \frac{C n M^{k+1} \bigl(\sqrt{\log 2n} + t \bigr)^k}
    {\nu^{u k/2} \eta^{k+1}}
    +
    \frac{C M^2
    \sqrt{n} \bigl(\sqrt{\log 2n} + t\bigr)^2}
    {\nu^{u} \eta^2}
    + \frac{1}{n^\kappa}.
  \end{align*}
  If $n^{1/4} s \sqrt d M > 1$, then this bound holds trivially,
  increasing $C$ if necessary.
  Let $t = \sqrt{2 \kappa \log 2 n}$, so that
  \begin{align*}
    &\P\biggl(
      \sup_{f \in \cF} |S_n'(f) - S_n^{*\prime}(f)|
      > 3 \eta
      + C M \sqrt{\log 2n}
      \Bigl\{
        n^{1/2-\epsilon}
        + (\log 2n)^{1/2 + 1/r + 1/r'}
      \Bigr\}
    \biggr) \\
    &\qquad\leq
    \frac{C n^{1 + \epsilon u k/2} M^{k+1} (\log 2n)^{k/2}}
    {\eta^{k+1}}
    +
    \frac{C M^2 n^{1/2 + \epsilon u} \log 2n}
    {\eta^2}
    + \frac{2}{n^\kappa}.
  \end{align*}
  Replace $\eta$ with
  $\bigl(n^{(1 + \epsilon u k / 2) / (k+1)}
  + n^{(1/2 + \epsilon u) / 2} \bigr) M \eta$
  and set $\epsilon \coloneqq (k - 1) / (2ku + 2k - 2u + 2)$,
  recalling that $k \in \{2, 3\}$ to see
  \begin{align*}
    &\P\biggl(
      \sup_{f \in \cF} |S_n'(f) - S_n^{*\prime}(f)|
      > C M \eta
      \bigl(n^{(1 + \epsilon u k / 2) / (k+1)}
      + n^{(1/2 + \epsilon u) / 2} \bigr) \\
      &\qquad\qquad+
      C M
      \Bigl\{
        n^{1/2 - \epsilon}
        + (\log 2n)^{1/2 + 1/r + 1/r'}
      \Bigr\}
      \sqrt{\log 2n}
    \biggr)
    \leq
    \frac{C \log 2n}{\eta^2}
    + \frac{2}{n^\kappa}, \\
    &\P\biggl(
      \sup_{f \in \cF} |S_n'(f) - S_n^{*\prime}(f)|
      > C M
      n^{(4-k+u) / (10 - 2k + 2u)}
      \Bigl\{ \sqrt{\log 2n} + \eta \Bigr\}
    \biggr)
    \leq
    \frac{C \log 2n}{\eta^2}
    + \frac{2}{n^\kappa},
  \end{align*}
  where we used that
  the first inequality is trivial unless
  $\eta^2 \geq \log 2 n$.
  We may assume $\eta \geq 1$, since the final bound is trivial
  otherwise. Replacing $\eta$ by $\eta\sqrt{\log 2n}$ in the
  second inequality above and increasing $C$, we obtain
  \begin{equation*}
    \P\biggl(
      \sup_{f\in\cF} |S_n'(f)-S_n^{*\prime}(f)|
      > CM\eta n^{\frac{4-k+u}{10-2k+2u}}\sqrt{\log 2n}
    \biggr)
    \leq \frac{C}{\eta^2}+\frac{2}{n^\kappa}.
  \end{equation*}

  Finally, we remove the $2/n^\kappa$ term.
  By Lemma~\ref{lem:moment_tail_orlicz},
  if $r < \infty$, then for $t > 0$,
  with probability at least $1 - 2n \exp(-t^r)$,
  \begin{equation*}
    \sup_{f \in \cF} |S_n(f)|
    \leq
    \sum_{i=1}^n F(X_i)
    + \sum_{i=1}^n P_{X,i}(F)
    \leq
    (t + C) n M.
  \end{equation*}
  Under condition~\ref{it:thm_ep_bootstrap_vc_nonparametric},
  by applying Lemma~\ref{lem:moment_tail_orlicz},
  with probability at least $1 - 2n \exp(-t^r)$,
  \begin{equation*}
    \sup_{f \in \cF} |S_n^*(f)|
    \leq
    2 n \max_{i \in [n]} F(X_i)
    \leq
    2 t n M.
  \end{equation*}
  Therefore, if $t \geq 1$,
  \begin{equation*}
    \P\biggl\{
      \sup_{f \in \cF} |S_n'(f) - S_n^{*\prime}(f)|
      >
      C n^2 M t
    \biggr\}
    \leq
    4 n \exp\bigl\{-(n t)^r\bigr\}
    \leq
    \frac{C}{(n t)^4}
    \leq
    \frac{C}{(n^2 t)^2}.
  \end{equation*}
  This bound also holds for $r = \infty$.
  Under condition~\ref{it:thm_ep_bootstrap_vc_multiplier},
  if $r' < \infty$, by Lemma~\ref{lem:moment_tail_orlicz},
  with probability at least $1 - 2n \exp(-t^r) - 2 n \exp(-t^{r'})$,
  \begin{equation*}
    \sup_{f \in \cF} |S_n^*(f)|
    \leq
    2 \max_{i \in [n]} |\chi_i|
    \sum_{j=1}^n F(X_j)
    \leq
    4 t^2 n M.
  \end{equation*}
  Therefore, if $t \geq 1$,
  \begin{align*}
    \P\biggl\{
      \sup_{f \in \cF} |S_n'(f) - S_n^{*\prime}(f)|
      >
      C n^2 M t^2
    \biggr\}
    &\leq
    4 n \exp\Bigl\{-\bigl(\sqrt n t\bigr)^r\Bigr\}
    + 2 n \exp\Bigl\{-\bigl(\sqrt n t\bigr)^{r'}\Bigr\} \\
    &\leq
    \frac{C}{(\sqrt n t)^8}
    \leq
    \frac{C}{(n^2 t^2)^2}.
  \end{align*}
  This bound also holds if $r = \infty$ or if $r' = \infty$.
  Thus, under condition~\ref{it:thm_ep_bootstrap_vc_nonparametric}
  or condition~\ref{it:thm_ep_bootstrap_vc_multiplier},
  for $\eta \geq n^2$,
  \begin{equation*}
    \P\biggl\{
      \sup_{f \in \cF} |S_n'(f) - S_n^{*\prime}(f)|
      >
      C M \eta
      n^{\frac{4 - k + u}{10 - 2k + 2u}}
      \sqrt{\log 2n}
    \biggr\}
    \leq
    \frac{C}{\eta^2}.
  \end{equation*}
  On the other hand, for $\eta \leq n^2$,
  since $\kappa = 4$,
  we have $1/\eta^2 \geq 1/n^4 = 1/n^\kappa$.
  Increasing $C$, we have that for all $\eta > 0$
  and $n \in \N$,
  \begin{equation*}
    \P\biggl\{
      \sup_{f \in \cF} |S_n'(f) - S_n^{*\prime}(f)|
      >
      C M \eta
      n^{\frac{4 - k + u}{10 - 2k + 2u}}
      \sqrt{\log 2n}
    \biggr\}
    \leq
    \frac{1}{\eta^2}.
    \iftoggle{journal}{}{\qedhere}
  \end{equation*}
\end{proof}

\subsection{Proofs for Section~\ref{sec:locpol}}

We begin by defining some useful auxiliary notation.
Set $\bW_n \coloneqq (W_1, \ldots, W_n)$
and let $\theta(\bW_n, v) \in \R^n$
with $\theta(\bW_n, v)_i \coloneqq \theta(W_i, v)$.
Define $\tilde g(v) \coloneqq \hat g(v) - \theta(\bW_n, v)$ and write
$S_{n,h}(w, v) \coloneqq \tilde N_{n,h}(w, v) \coloneqq
\hat Q_h(w)^\T \hat K_h(w) \tilde g(v) / n$
and $B_{n,h}(w, v) \coloneqq \hat Q_h(w)^\T \hat K_h(w) \theta(\bW_n, v) / n$,
so $N_{n,h}(w, v) = \tilde N_{n,h}(w, v) + B_{n,h}(w, v)$.
We have the decomposition
\begin{align}
  \nonumber
  \hat \theta_{n,h}(w, v) - \theta(w, v)
  &=
  e_1^\T \hat H_h(w)^{-1} N_{n,h}(w, v) - \theta(w, v) \\
  \label{eq:locpol_decomp_coupling}
  &=
  e_1^\T H_h(w)^{-1} S_{n,h}(w, v) \\
  \label{eq:locpol_decomp_errors}
  &\qquad+
  e_1^\T \bigl\{\hat H_h(w)^{-1} - H_h(w)^{-1}\bigr\}
  \tilde N_{n,h}(w, v) \\
  \label{eq:locpol_decomp_bias}
  &\qquad+
  e_1^\T \hat H_h(w)^{-1} B_{n,h}(w, v) - \theta(w, v).
\end{align}
The term in \eqref{eq:locpol_decomp_coupling} is the target for coupling,
while \eqref{eq:locpol_decomp_errors} is a negligible product of errors.
The term in \eqref{eq:locpol_decomp_bias} represents the bias.

For the nonparametric bootstrap, under
condition~\ref{it:thm_locpol_bootstrap_nonparametric},
let $\hat H_h^*(w)$, $N_{n,h}^*(w, v)$,
$\tilde N_{n,h}^*(w, v)$ and $B_{n,h}^*(w, v)$
be defined in the same manner as
$\hat H_h(w)$, $N_{n,h}(w, v)$,
$\tilde N_{n,h}(w, v)$ and $B_{n,h}(w, v)$ respectively,
using the bootstrapped data $(W_i^*, V_i^*)$
instead of the original data $(W_i, V_i)$. Define
$S_{n,h}^*(w, v) \coloneqq \tilde N_{n,h}^*(w, v) - \tilde N_{n,h}(w, v)$
and consider the decomposition
\begin{align}
  \nonumber
  \tilde \theta_{n,h}^*(w, v)
  &=
  e_1^\T \hat H_h^*(w)^{-1} N_{n,h}^*(w, v)
  - e_1^\T \hat H_h(w)^{-1} N_{n,h}(w, v) \\
  \label{eq:locpol_nonparametric_decomp_coupling}
  &=
  e_1^\T H_h(w)^{-1} S_{n,h}^*(w, v) \\
  \label{eq:locpol_nonparametric_decomp_errors}
  &\qquad+
  e_1^\T \bigl\{\hat H_h^*(w)^{-1} - H_h(w)^{-1}\bigr\}
  \tilde N_{n,h}^*(w, v)
  \iftoggle{journal}{\\ \nonumber &\qquad\qquad}{}
  - e_1^\T \bigl\{\hat H_h(w)^{-1} - H_h(w)^{-1}\bigr\}
  \tilde N_{n,h}(w, v) \\
  \label{eq:locpol_nonparametric_decomp_bias}
  &\qquad+
  \bigl\{e_1^\T \hat H_h^*(w)^{-1} B_{n,h}^*(w, v) - \theta(w,v)\bigr\}
  \iftoggle{journal}{\\ \nonumber &\qquad\qquad}{}
  - \bigl\{e_1^\T \hat H_h(w)^{-1} B_{n,h}(w, v) - \theta(w,v)\bigr\}.
\end{align}
The term in \eqref{eq:locpol_nonparametric_decomp_coupling}
is our coupling target,
while those in \eqref{eq:locpol_nonparametric_decomp_errors}
are negligible products of errors.
The terms in \eqref{eq:locpol_nonparametric_decomp_bias}
are related to the bias under the bootstrap and population
distributions, respectively.

For the multiplier bootstrap, under
condition~\ref{it:thm_locpol_bootstrap_multiplier},
define $\tilde g^\circ(v) \in \R^n$ with
$\tilde g^\circ(v)_i \coloneqq \chi_i \tilde g(v)_i$.
Let $S_{n,h}^\circ(w, v) \coloneqq \hat Q_h(w)^\T \hat K_h(w) \tilde
g^\circ(v) / n$
and $R_{n,h}^\circ(w, v) \coloneqq \sum_{i=1}^n
Q_h(W_i - w) K_h(W_i - w) \chi_i
\bigl\{ \hat\theta_{n,h}(W_i,v) - \theta(W_i, v) \bigr\} / n$,
so that $N_{n,h}^\circ(w, v) = S_{n,h}^\circ(w, v) - R_{n,h}^\circ(w, v)$.
For clarity of notation, write
$\tilde \theta_{n,h}^\circ(w, v)
\coloneqq \tilde \theta_{n,h}^*(w, v)$
under the multiplier bootstrap. Consider the decomposition
\begin{align}
  \nonumber
  \tilde \theta_{n,h}^\circ(w, v)
  &=
  e_1^\T \hat H_h(w)^{-1} N_{n,h}^\circ(w, v) \\
  \label{eq:locpol_multiplier_decomp_coupling}
  &=
  e_1^\T H_h(w)^{-1} S_{n,h}^\circ(w, v) \\
  \label{eq:locpol_multiplier_decomp_errors}
  &\qquad+
  e_1^\T \bigl\{\hat H_h(w)^{-1} - H_h(w)^{-1}\bigr\} S_{n,h}^\circ(w, v)
  - e_1^\T \hat H_h(w)^{-1} R_{n,h}^\circ(w, v).
\end{align}
The term in \eqref{eq:locpol_multiplier_decomp_coupling}
is the coupling target, and those
in \eqref{eq:locpol_multiplier_decomp_errors} are
negligible products of errors.

\begin{lemma}[Covering numbers of local polynomial classes]%
  \label{lem:locpol_covering}
  Let $D \in \N$ and $\alpha \in \N_0^D$.
  Suppose that $K: \R^{D} \to [0, \infty)$ is a Borel-measurable
  function supported on $[-1, 1]^{D}$ satisfying
  $K(w) \leq C$ for all $w \in \R^{D}$ where $C > 0$.
  Further, assume that $K$ satisfies
  $|K(w) - K(w')| \leq C \|w' - w\|_\infty$
  for all $w, w' \in \R^{D}$.
  For $h \in (0, 1]$,
  define a class of functions from $[0, 1]^D \to \R$ by
  \begin{equation*}
    \cF_h
    \coloneqq
    \Biggl\{
      W \mapsto
      \frac{1}{h^D}
      \biggl(\frac{W - w}{h}\biggr)^\alpha
      K\biggl(\frac{W - w}{h}\biggr):
      w \in [0, 1]^D
    \Biggr\}.
  \end{equation*}
  Then if $F_h \coloneqq C / h^D$, we have
  for $\varepsilon \in (0, 1]$ that
  \begin{equation*}
    \log
    N\bigl(\varepsilon F_h, \cF_h, \|\cdot\|_{\infty}\bigr)
    \leq
    D \log \biggl(\frac{1}{h \varepsilon} \biggr).
  \end{equation*}
\end{lemma}

\begin{proof}[Proof of Lemma~\ref{lem:locpol_covering}]
  Set $k \coloneqq |\alpha|$.
  If $k=0$, then $\alpha=0$, so
  for $u,v \in [-1, 1]^D$,
  \[
    |u^\alpha K(u)-v^\alpha K(v)|
    = |K(u)-K(v)|
    \leq C\|u-v\|_\infty.
  \]
  Now suppose $k \geq 1$ and $u,v \in [-1, 1]^D$.
  By symmetry, we may assume that
  $r \coloneqq \|u\|_\infty \geq \|v\|_\infty$.
  Since $K$ is continuous and supported on $[-1,1]^D$,
  we have $K(u)\leq C(1-r)$.
  Lemma~\ref{lem:multi_index}
  and convexity of $r \mapsto r^k$ give
  \begin{align}
    \label{eq:kernel_covering}
    |u^\alpha K(u)-v^\alpha K(v)|
    &\leq
    |u^\alpha-v^\alpha|\,|K(u)|
    + |v^\alpha|\,|K(u)-K(v)| \nonumber\\
    &\leq
    C\bigl\{k r^{k-1}(1-r)+r^k\bigr\}\|u-v\|_\infty
    \leq C\|u-v\|_\infty.
  \end{align}
  If $u \in [-1,1]^D$ and $v \in \R^D \setminus [-1,1]^D$, then
  $K(v)=0$ and $|u^\alpha|\leq 1$, so
  the last inequality in \eqref{eq:kernel_covering} holds.
  If $u, v \in \R^D \setminus [-1,1]^D$, then
  $K(u) = K(v) = 0$.
  Thus, this bound holds for all $u,v \in \R^D$.
  Writing $f_w$ for the function in $\cF_h$ indexed by $w$,
  take $u=(W-w)/h$ and $v=(W-w')/h$ to obtain
  \[
    \|f_w-f_{w'}\|_\infty
    \leq C h^{-D-1}\|w-w'\|_\infty.
  \]
  Set $m \coloneqq \lceil 1 /(2h\varepsilon)\rceil$ and
  partition $[0,1]^D$ into $m^D$ cubes of side length $1/m$
  with centers $w_1,\ldots,w_{m^D}$.
  For each $w \in [0,1]^D$, we have
  $\min_{j \in [m^D]} \|w-w_j\|_\infty
  \leq 1/(2m) \leq h\varepsilon$.
  Therefore,
  $ \min_{j \in [m^D]} \|f_w-f_{w_j}\|_\infty
  \leq \varepsilon C h^{-D} = \varepsilon F_h$.
  It follows that
  $ N(\varepsilon F_h,\cF_h,\|\cdot\|_{\infty}) \leq m^D
  \leq \bigl(1 / (h\varepsilon)\bigr)^D$,
  as $1/(h\varepsilon)\geq1$.
  Taking logarithms gives
  $\log N(\varepsilon F_h,\cF_h,\|\cdot\|_{\infty})
  \leq D\log\bigl(1/(h\varepsilon)\bigr)$.
\end{proof}

\begin{lemma}[Bounds on $H_h$]%
  \label{lem:locpol_H}
  Grant the setup of Section~\ref{sec:locpol} and Assumption~\ref{ass:locpol}.
  There exists a quantity $C > 0$ depending only on
  $\Clp$, $D_W$ and $q$ such that
  $\sup_{w \in \cW} \|H_h(w)\|_{2 \to 2} \leq C$;
  moreover, $H_h(w)$ is invertible for all $w \in \cW$,
  and satisfies
  $\sup_{w \in \cW} \|H_h(w)^{-1}\|_{2 \to 2} \leq C$.
\end{lemma}

\begin{proof}[Proof of Lemma~\ref{lem:locpol_H}]
  We extend the proof of \citet[Lemma~A.1]{ching2026efficient}
  to $h \in (0, 1]$.
  Write $\bar f_W \coloneqq \sum_{i=1}^n f_{W,i}/n$ and
  $a \coloneqq (1/\Clp) \land (1/2)$.
  For $w \in [0,1]^{D_W}$,
  \begin{equation*}
    H_h(w)
    =
    \int_{\prod_{j=1}^{D_W}[-w_j/h,\, (1-w_j)/h]}
    Q_1(z) Q_1(z)^\T K(z)\bar f_W(w+h z)
    \diffi z.
  \end{equation*}
  Since $h \leq 1$, for each $j \in [D_W]$,
  if $w_j \leq 1/2$, then
  $[0,a] \subseteq [-w_j/h, (1-w_j)/h]$.
  Otherwise, $[-a,0] \subseteq [-w_j/h, (1-w_j)/h]$.
  Thus, there exists $A = \prod_{j=1}^{D_W} [a_j, b_j]$
  such that $A \subseteq \prod_{j=1}^{D_W}\bigl[-w_j/h,\, (1-w_j)/h\bigr]$,
  where for each $j \in [D_W]$, we have either
  $[a_j, b_j] = [0, a]$ or $[a_j, b_j] = [-a, 0]$.
  On $A$, we have that $K(z)\bar f_W(w+h z) \geq 1/\Clp^2$
  and $\int_A Q_1(z) Q_1(z)^\T \diffi z$ is positive definite.
  By continuity and compactness and taking a minimum
  over the possible values of the set $A$, we obtain
  $\sup_{w \in \cW} \|H_h(w)^{-1}\|_{2 \to 2} \leq C$.
  The upper bound follows by
  boundedness and compact support of the kernel.
\end{proof}

\begin{lemma}[Bounds on $\hat H_h$]%
  \label{lem:locpol_Hhat}
  Grant the setup of Section~\ref{sec:locpol} and Assumption~\ref{ass:locpol}.
  There exists $C > 0$ depending only on
  $\Clp$, $D_W$ and $q$
  such that for all $t > 0$, with
  probability at least $1 - e^{-t^2/2}$,
  \begin{equation*}
    \sup_{w \in \cW}
    \bigl\|
    \hat H_h(w) - H_h(w)
    \bigr\|_{2 \to 2}
    \leq
    C \sqrt{\frac{1 + D_w\log 2n + t^2}{n h^{D_W}}}
    + \frac{C t^2}{n h^{D_W}}.
  \end{equation*}
  Moreover, with probability at least $1 - 2 e^{-n h^{D_W}/C}$,
  $\hat H_h(w)$ is invertible for all $w \in \cW$ and
  satisfies
  $\sup_{w \in \cW} \|\hat H_h(w)\|_{2 \to 2} \leq C$
  and
  $\sup_{w \in \cW} \|\hat H_h(w)^{-1}\|_{2 \to 2} \leq C$.
  Finally, for all $t > 0$,
  with probability at least $1 - 2 e^{-n h^{D_W}/C} - e^{-t^2/2}$,
  \begin{equation*}
    \sup_{w \in \cW}
    \bigl\|
    \hat H_h(w)^{-1} - H_h(w)^{-1}
    \bigr\|_{2 \to 2}
    \leq
    C \sqrt{\frac{1 + D_w\log 2n + t^2}{n h^{D_W}}}.
  \end{equation*}
\end{lemma}

\begin{proof}[Proof of Lemma~\ref{lem:locpol_Hhat}]
  Let $C > 0$ depend only on the stated quantities;
  its value may vary from line to line.
  For $j, k \in [D_q]$,
  $e_j, e_k \in \R^{D_q}$ and $w \in \cW$, we have
  \begin{align*}
    e_j^\T \hat H_h(w) e_k
    &=
    \frac{1}{n}
    e_j^\T \hat Q_h(w)^\T \hat K_h(w) \hat Q_h(w) e_k \\
    &=
    \frac{1}{n}
    \sum_{i=1}^n
    e_j^\T Q_h(W_i - w) K_h(W_i - w) Q_h(W_i - w)^\T e_k.
  \end{align*}
  Observe that since the monomials are bounded by one on the support
  of the kernel, as the $W_i$ variables have
  uniformly bounded densities, and as the kernel is bounded,
  \begin{equation*}
    \E\Bigl\{
      \bigl(
        e_j^\T Q_h(W_i - w) K_h(W_i - w) Q_h(W_i - w)^\T e_k
      \bigr)^2
    \Bigr\}
    \leq
    \E\bigl\{
      K_h(W_i - w)^2
    \bigr\}
    \leq C h^{-D_W}.
  \end{equation*}
  Similarly,
  \begin{equation*}
    \bigl|
    e_j^\T Q_h(W_i - w) K_h(W_i - w) Q_h(W_i - w)^\T e_k
    \bigr|
    \leq
    K_h(W_i - w)
    \leq C h^{-D_W}.
  \end{equation*}
  Consider the function class
  \begin{equation*}
    \cF_{h}
    \coloneqq
    \Bigl\{
      W \mapsto
      e_j^\T Q_h(W - w) K_h(W - w) Q_h(W - w)^\T e_k:
      w \in \cW, \, j, k \in [D_q]
    \Bigr\}.
  \end{equation*}
  By Lemma~\ref{lem:locpol_covering} and taking a union
  over $j, k \in [D_q]$, with $F_h \coloneqq \Clp / h^{D_W}$
  and $\varepsilon \in (0, 1]$,
  \begin{equation*}
    \sup_{Q \in \cQ}
    \log
    N\bigl(\varepsilon F_h, \cF_h, \|\cdot\|_{Q,2}\bigr)
    \leq
    C + C D_w \log \biggl(\frac{1}{h \varepsilon}\biggr),
  \end{equation*}
  where $\cQ$ is the set of Borel
  probability measures on $[0, 1]^{D_W}$.
  This holds because either $D_w \geq 1$, in which case
  a factor of $D_W$ may be absorbed into $C$,
  or $D_w = 0$, in which case the covering number is at most $D_q^2$.
  It follows by Lemma~\ref{lem:entropy_vc} that
  \begin{equation*}
    \cJ(\varepsilon)
    \coloneqq
    \int_0^\varepsilon
    \sup_{Q \in \cQ}
    \sqrt{1 + \log
      N\bigl(
        \xi F_h, \cF_h, \|\cdot\|_{Q,2}
    \bigr)}
    \diffi \xi
    \leq
    C \varepsilon
    \sqrt{1 + D_w \log \bigl\{1 / (h \varepsilon)\bigr\}}.
  \end{equation*}
  By Theorem~\ref{thm:empirical_tail} with $r = \infty$,
  since $n h^{D_W} \geq n^{1/\Clp}$,
  \begin{align*}
    \sup_{w \in \cW}
    \bigl\|
    \hat H_h(w) - H_h(w)
    \bigr\|_{2 \to 2}
    &\leq
    C
    \biggl\{
      \frac{1}{\sqrt{n h^{D_W}}}
      + \frac{1}{n h^{D_W}}
      \Bigl(
        \sqrt{1 + D_w \log(2/h)} + t
      \Bigr)
    \biggr\} \\
    &\qquad\quad\times
    \Bigl( \sqrt{1 + D_w \log(2/h)} + t \Bigr) \\
    &\leq
    C \sqrt{\frac{1 + D_w\log 2n + t^2}{n h^{D_W}}}
    + \frac{C t^2}{n h^{D_W}},
  \end{align*}
  with probability at least
  $1 - e^{-t^2/2}$.
  The bounds on
  $\|\hat H_h(w)\|_{2 \to 2}$
  and $\|\hat H_h(w)^{-1}\|_{2 \to 2}$
  follow by taking $t^2 = n h^{D_W} / C$ and applying
  Lemma~\ref{lem:locpol_H} and Weyl's inequality
  \citep[see, for example,][Corollary~III.2.6]{bhatia1997matrix}
  for $n \geq C$; when $n < C$ the result is vacuous
  as $2 e^{-n h^{D_W}/C} \geq 1$, increasing $C$ if necessary.
  For the final result,
  with probability at least $1 - 2 e^{-n h^{D_W}/C} - e^{-t^2/2}$,
  \begin{align}
    \nonumber
    \bigl\|
    \hat H_h(w)^{-1} - H_h(w)^{-1}
    \bigr\|_{2 \to 2}
    &=
    \bigl\|
    \hat H_h(w)^{-1}
    \bigl(
      H_h(w) - \hat H_h(w)
    \bigr)
    H_h(w)^{-1}
    \bigr\|_{2 \to 2} \\
    \nonumber
    &\leq
    \bigl\|
    H_h(w)^{-1}
    \bigr\|_{2 \to 2}
    \bigl\|
    \hat H_h(w)^{-1}
    \bigr\|_{2 \to 2}
    \bigl\|
    \hat H_h(w) - H_h(w)
    \bigr\|_{2 \to 2} \\
    \label{eq:locpol_Hhat_tail}
    &\leq
    C \sqrt{\frac{1 + D_w\log 2n + t^2}{n h^{D_W}}}
    + \frac{C t^2}{n h^{D_W}}.
  \end{align}
  On the event where $\sup_{w \in \cW}
  \|\hat H_h(w)^{-1}\|_{2 \to 2} \leq C$,
  the triangle inequality and Lemma~\ref{lem:locpol_H}
  give $\sup_{w \in \cW}
  \|\hat H_h(w)^{-1} - H_h(w)^{-1}\|_{2 \to 2} \leq 2 C$.
  Thus, the desired bound is immediate when $t \geq \sqrt{n h^{D_W}}$.
  For $t \leq \sqrt{n h^{D_W}}$,
  the second term in \eqref{eq:locpol_Hhat_tail}
  may be bounded by the first,
  again yielding the result.
\end{proof}

\begin{lemma}[Bounds on $\hat H_h^*$]%
  \label{lem:locpol_Hhatstar}
  Grant the setup of Theorem~\ref{thm:locpol_bootstrap}.
  The results of Lemma~\ref{lem:locpol_Hhat} hold
  with $\hat H_h(w)$ replaced with $\hat H_h^*(w)$.
\end{lemma}

\begin{proof}[Proof of Lemma~\ref{lem:locpol_Hhatstar}]
  By the triangle inequality,
  \begin{equation*}
    \sup_{w \in \cW}
    \bigl\|\hat H_h^*(w) - H_h(w)\bigr\|_{2 \to 2}
    \leq
    \sup_{w \in \cW}
    \bigl\|\hat H_h^*(w) - \hat H_h(w)\bigr\|_{2 \to 2}
    + \sup_{w \in \cW}
    \bigl\|\hat H_h(w) - H_h(w)\bigr\|_{2 \to 2}.
  \end{equation*}
  For the first term,
  apply Lemma~\ref{lem:bootstrap_tail} with $r = r' = \infty$,
  using the function class and entropy bounds
  from the proof of Lemma~\ref{lem:locpol_Hhat}.
  For the second term,
  apply Lemma~\ref{lem:locpol_Hhat} directly.
  Replace $t$ with $t + 2$ and take a union bound;
  the remaining assertions follow as in the proof of
  Lemma~\ref{lem:locpol_Hhat}.
\end{proof}

\begin{lemma}[Bound on $\tilde N_{n,h}$]%
  \label{lem:locpol_linear}
  Grant the setup of Section~\ref{sec:locpol}
  and Assumption~\ref{ass:locpol}, and take $\kappa \geq 1$.
  There exists $C > 0$ depending only on
  $\Clp$, $D_W$, $D_V$, $q$, $r$ and $\kappa$ such that
  for all $t > 0$,
  with probability at least $1 - 1/n^{\kappa} - e^{-t^2/2}$,
  \begin{equation*}
    \sup_{(w, v) \in \cW \times \cV}
    \bigl\|
    \tilde N_{n,h}(w, v)
    \bigr\|_2
    \leq
    C \sqrt{\frac{1 + D_w \log 2n + t^2}{n h^{D_W}}}.
  \end{equation*}
\end{lemma}

\begin{proof}[Proof of Lemma~\ref{lem:locpol_linear}]
  Throughout the proof,
  $C > 0$ depends only on the stated quantities
  and may change from line to line.
  For $j \in [D_q]$, $w \in \cW$ and $v \in \cV$, we have
  \begin{equation*}
    e_j^\T
    \tilde N_{n,h}(w, v)
    =
    \frac{1}{n}
    \sum_{i=1}^n
    e_j^\T
    Q_h(W_i - w)
    K_h(W_i - w)
    \tilde g(v)_i.
  \end{equation*}
  Since $\vvvert \tilde g(v)_i \vvvert_{\omega_r \mid W_i} \leq C$
  and $W_i$ has bounded Lebesgue density,
  Lemma~\ref{lem:moment_tail_orlicz} gives
  \begin{align*}
    \E\Bigl\{
      \bigl(
        e_j^\T Q_h(W_i - w) K_h(W_i - w) \tilde g(v)_i
      \bigr)^2
    \Bigr\}
    &\leq
    \E\bigl\{
      K_h(W_i - w)^2
      \tilde g(v)_i^2
    \bigr\}
    \leq C h^{-D_W}, \\
    \bigvvvert
    e_j^\T Q_h(W_i - w) K_h(W_i - w) \tilde g(v)_i
    \bigvvvert_{\omega_r}
    &\leq
    \Clp h^{-D_W}
    \vvvert \tilde g(v)_i \vvvert_{\omega_r}
    \leq C h^{-D_W}.
  \end{align*}

  Firstly, consider the case where $g(V, v) = V$.
  Define the function class
  \begin{equation*}
    \cF_{h}
    =
    \Bigl\{
      (W, V) \mapsto
      e_j^\T Q_h(W - w) K_h(W - w)
      \{V - \theta(W)\}:
      w \in \cW, \, j \in [D_q]
    \Bigr\},
  \end{equation*}
  which is pointwise measurable by continuity in $w$.
  Define $F_h(V) \coloneqq \Clp (|V| + \Clp) / h^{D_W}$,
  so $\vvvert F_h(V_i) \vvvert_{\omega_r} \leq C h^{-D_W}$ for $i \in [n]$.
  By Lemma~\ref{lem:locpol_covering} and taking a union
  over $j \in [D_q]$, for $\varepsilon \in (0, 1]$,
  \begin{equation*}
    \sup_{Q \in \cQ}
    \log
    N\bigl(\varepsilon \|F_h\|_{Q,2}, \cF_h, \|\cdot\|_{Q,2}\bigr)
    \leq
    C + C D_w \log \biggl(\frac{1}{h \varepsilon}\biggr),
  \end{equation*}
  where $\cQ$ is the set of Borel
  probability measures on $[0, 1]^{D_W} \times \R$
  for which $\|F_h\|_{Q,2} > 0$.
  This holds because either $D_w \geq 1$, in which case
  a factor of $D_W$ may be absorbed into $C$,
  or $D_w = 0$, in which case the covering number is at most $D_q$.
  It follows by Lemma~\ref{lem:entropy_vc} that
  \begin{equation*}
    \cJ(\varepsilon)
    \coloneqq
    \int_0^\varepsilon
    \sup_{Q \in \cQ}
    \sqrt{1 + \log
      N\bigl(
        \xi \|F_h\|_{Q,2}, \cF_h, \|\cdot\|_{Q,2}
    \bigr)}
    \diffi \xi
    \leq
    C
    \varepsilon
    \sqrt{1 + D_w \log \bigl\{2 / (h \varepsilon)\bigr\}}.
  \end{equation*}
  By Theorem~\ref{thm:empirical_tail},
  since $n h^{D_W} \geq n^{1/\Clp}$,
  \begin{align*}
    \sup_{(w, v) \in \cW \times \cV}
    \!
    \bigl\|
    \tilde N_{n,h}(w, v)
    \bigr\|_2
    &\leq
    C
    \biggl\{
      \!
      \frac{1}{\sqrt{n h^{D_W}}} +
      \frac{(\kappa \log 2n)^{1/r}}{n h^{D_W}}
      \Bigl(
        \!
        \sqrt{1 + D_w \log 2n} + t
      \Bigr)
      \!
    \biggr\}
    \iftoggle{journal}{\\&\qquad\qquad\times}{}
    \Bigl( \! \sqrt{1 + D_w \log 2n} + t \Bigr) \\
    &\leq
    C \sqrt{\frac{1 + D_w \log 2n + t^2}{n h^{D_W}}}
    + \frac{C (\kappa \log 2n)^{1/r} (1 + D_w \log 2n + t^2)}{n h^{D_W}},
  \end{align*}
  with probability at least
  $1 - 1/n^{\kappa} - e^{-t^2/2}$.
  Since $n h^{D_W} \geq n^{1/\Clp}$, the second term
  above is at most $C$ times the
  first whenever $t \leq C \sqrt{\log 2n}$.
  For larger $t$, the claim follows by replacing
  $t$ with $\sqrt{2(\kappa + 1) \log 2 n}$
  and $\kappa$ with $\kappa + 1$.

  Secondly, we handle the case where
  $g(V, v) = \I\{V \leq v\}$.
  Consider the function class
  \begin{align*}
    \cF_{h}'
    &=
    \Bigl\{
      (W, V) \mapsto
      e_j^\T Q_h(W - w) K_h(W - w)
      \bigl\{\I\{V \leq v\} - \theta(W, v)\bigr\}:
      \iftoggle{journal}{\\&\qquad\qquad}{\,}
      w \in \cW, \, v \in \cV, \, j \in [D_q]
    \Bigr\},
  \end{align*}
  which is pointwise measurable by continuity in $w$
  and coordinatewise right-continuity in $v$.
  Note that $\cF_h' = \bigl\{(W, V) \mapsto
  f_1(W) f_2(W,V): f_1 \in \cG_1, f_2 \in \cG_2\bigr\}$,
  where
  \begin{align*}
    \cG_1
    &=
    \Bigl\{
      W \mapsto
      e_j^\T Q_h(W - w) K_h(W - w):\,
      w \in \cW, \, j \in [D_q]
    \Bigr\}, \\
    \cG_2
    &=
    \Bigl\{
      (W, V) \mapsto
      \I\{V \leq v\} - \theta(W, v):\,
      v \in \cV
    \Bigr\}.
  \end{align*}
  Take a Borel probability measure $Q$ on
  $[0,1]^{D_W} \times \R^{D_V}$ and $\varepsilon \in (0,1]$,
  and set $F_h' \coloneqq \Clp h^{-D_W}$.
  If $D_w = 0$, then $\cG_1$ contains at most $D_q$ functions.
  For $D_w \geq 1$, apply Lemma~\ref{lem:locpol_covering}
  and take a union over $j \in [D_q]$.
  In either case,
  \begin{equation*}
    \log
    N\biggl(
      \frac{\Clp \varepsilon}{2 h^{D_W}},
      \cG_1, \|\cdot\|_{Q,2}
    \biggr)
    \leq
    C + C D_w
    \log
    \biggl(
      \frac{2}{h \varepsilon}
    \biggr).
  \end{equation*}
  For $\cG_2$,
  write $P_{V \mid W}$ for the common regular conditional distribution
  of $V$ given $W$, so that $\theta(w,v)=P_{V \mid W=w}\bigl((-\infty,v]\bigr)$.
  Let $Q_W$ be the marginal of $Q$ on the $W$ coordinate
  and define $Q' \coloneqq Q_W \otimes P_{V \mid W}$
  and $\bar Q \coloneqq (Q+Q')/2$.
  By the triangle inequality, Cauchy--Schwarz
  and the conditional Jensen inequality,
  for all $v, v' \in \cV$,
  \begin{align*}
    &\bigl\|\I\{V\leq v\}-\theta(W,v)
    -\I\{V\leq v'\}+\theta(W,v')\bigr\|_{Q,2} \\
    &\qquad\leq
    \bigl\|\I\{V\leq v\}-\I\{V\leq v'\}\bigr\|_{Q,2}
    +\bigl\|\I\{V\leq v\}-\I\{V\leq v'\}\bigr\|_{Q',2} \\
    &\qquad\leq
    2\bigl\|\I\{V\leq v\}-\I\{V\leq v'\}\bigr\|_{\bar Q,2}.
  \end{align*}
  An $\varepsilon/4$-cover of
  $\bigl\{\I\{V\leq v\}:v\in\cV\bigr\}$ in $L_2(\bar Q)$
  therefore induces an $\varepsilon/2$-cover of $\cG_2$
  in $L_2(Q)$.
  By the VC dimension bounds in
  \citet[Example~2.6.1 and Theorem~2.6.4]{van1996weak},
  \begin{equation*}
    \log N\biggl(\frac{\varepsilon}{2},\cG_2,\|\cdot\|_{Q,2}\biggr)
    \leq
    \log N\biggl(
      \frac{\varepsilon}{4},
      \bigl\{\I\{V\leq v\}:v\in\cV\bigr\},
      \|\cdot\|_{\bar Q,2}
    \biggr)
    \leq
    C\log\biggl(\frac{2}{\varepsilon}\biggr).
  \end{equation*}
  We deduce that
  \begin{align*}
    \log
    N\bigl(\varepsilon F_h', \cF_h', \|\cdot\|_{Q,2}\bigr)
    &\leq
    \log
    N\biggl(
      \frac{\Clp \varepsilon}{2 h^{D_W}},
      \cG_1, \|\cdot\|_{Q,2}
    \biggr)
    + \log
    N\biggl(
      \frac{\varepsilon}{2},
      \cG_2, \|\cdot\|_{Q,2}
    \biggr) \\
    &\leq
    C \log
    \biggl(
      \frac{2}{\varepsilon}
    \biggr)
    + C D_w \log
    \biggl(
      \frac{2}{h \varepsilon}
    \biggr).
  \end{align*}
  When $D_w \geq 1$,
  the rest of the proof proceeds
  as for the case $g(V, v) = V$,
  applying Theorem~\ref{thm:empirical_tail}.

  It remains to ensure that no spurious logarithmic
  factors remain when $D_w = 0$.
  In this case, write $\cW = \{w_0\}$
  and use the envelope
  $ F_h''(W,V) \coloneqq
  \Clp h^{-D_W}\I\{\|W-w_0\|_\infty\leq h\}$.
  By the same VC dimension argument,
  \begin{equation*}
    \sup_Q \log N\bigl(
      \varepsilon\|F_h''\|_{Q,2},\cF_h',\|\cdot\|_{Q,2}
    \bigr)
    \leq C\log \biggl(\frac{2}{\varepsilon}\biggr).
  \end{equation*}
  We may apply
  Theorem~\ref{thm:empirical_tail} with $r=\infty$
  as before.
\end{proof}

\begin{lemma}[Bound on $\tilde N_{n,h}^*$]%
  \label{lem:locpol_linear_bootstrap}
  Lemma~\ref{lem:locpol_linear} holds
  with $\tilde N_{n,h}(w, v)$ replaced with $\tilde N_{n,h}^*(w, v)$.
\end{lemma}

\begin{proof}[Proof of Lemma~\ref{lem:locpol_linear_bootstrap}]
  By the triangle inequality,
  \begin{align*}
    \sup_{(w,v) \in \cW \times \cV}
    \bigl\|\tilde N_{n,h}^*(w,v)\bigr\|_2
    &\leq
    \sup_{(w,v) \in \cW \times \cV}
    \bigl\|\tilde N_{n,h}^*(w,v)-\tilde N_{n,h}(w,v)\bigr\|_2
    \iftoggle{journal}{\\ &\qquad}{}
    + \sup_{(w,v) \in \cW \times \cV}
    \bigl\|\tilde N_{n,h}(w,v)\bigr\|_2.
  \end{align*}
  For the first term,
  we follow the proof of Lemma~\ref{lem:locpol_linear},
  using Lemma~\ref{lem:bootstrap_tail} with $r' = \infty$
  instead of Theorem~\ref{thm:empirical_tail} to handle the bootstrap process.
  The second term is bounded by applying Lemma~\ref{lem:locpol_linear}
  directly.
  Replace $t$ with $t + 2$ and $\kappa$ with $\kappa + 1$,
  then take a union bound to conclude.
\end{proof}

\begin{lemma}[Bound on $B_{n,h}$]%
  \label{lem:locpol_bias}
  Grant the setup of Section~\ref{sec:locpol}
  and Assumption~\ref{ass:locpol}.
  There exists $C > 0$ depending only on
  $\Clp$, $D_W$ and $q$ such that
  with probability at least $1 - 2 e^{-n h^{D_W}/C}$,
  \begin{equation*}
    \sup_{(w, v) \in \cW \times \cV}
    \bigl|
    e_1^\T \hat H_h(w)^{-1} B_{n,h}(w, v)
    - \theta(w, v)
    \bigr|
    \leq
    C h^\zeta.
  \end{equation*}
\end{lemma}

\begin{proof}[Proof of Lemma~\ref{lem:locpol_bias}]
  Throughout the proof, $C > 0$ depends only on the stated quantities
  and may change from line to line.
  By Taylor's theorem,
  for all $w, w' \in [0, 1]^{D_W}$
  and $v \in \cV$, there exists $\tilde w$
  on the line segment between $w$ and $w'$ such that
  \begin{align*}
    \theta(w',v) - \theta(w,v)
    &=
    \sum_{\alpha \in \N_0^{D_W}: |\alpha| \in [\lceil\zeta\rceil-1]}
    \frac{(w' - w)^\alpha}{\alpha!}
    \frac{\partial^{|\alpha|} \theta(w, v)}{\partial w^\alpha} \\
    &\qquad\qquad+
    \sum_{\alpha \in \N_0^{D_W}: |\alpha| = \lceil\zeta\rceil-1}
    \frac{(w' - w)^\alpha}{\alpha!}
    \biggl\{
      \frac{\partial^{|\alpha|} \theta(\tilde w, v)}{\partial w^\alpha}
      - \frac{\partial^{|\alpha|} \theta(w, v)}{\partial w^\alpha}
    \biggr\}.
  \end{align*}
  Since $q \geq \lceil\zeta\rceil - 1$,
  writing $e_\alpha$ for the standard basis vector corresponding
  to the monomial in $Q_h$ indexed by $\alpha \in \N_0^{D_W}$,
  the first term on the right-hand side satisfies
  \begin{align*}
    &\sum_{\alpha \in \N_0^{D_W}: |\alpha| \in [\lceil\zeta\rceil-1]}
    \frac{(w' - w)^\alpha}{\alpha!}
    \frac{\partial^{|\alpha|} \theta(w, v)}{\partial w^\alpha}
    \iftoggle{journal}{\\&\qquad}{}
    =
    \sum_{\alpha \in \N_0^{D_W}: |\alpha| \in [\lceil\zeta\rceil-1]}
    h^{|\alpha|} Q_h(w' - w)^\T e_\alpha
    \frac{\partial^{|\alpha|} \theta(w, v)}{\partial w^\alpha}.
  \end{align*}
  For the second term on the right-hand side,
  for $\|w' - w\|_\infty \leq h$, we have
  \begin{align*}
    &\left|
    \sum_{\alpha \in \N_0^{D_W}: |\alpha| = \lceil\zeta\rceil-1}
    \frac{(w' - w)^\alpha}{\alpha!}
    \biggl\{
      \frac{\partial^{|\alpha|} \theta(\tilde w, v)}{\partial w^\alpha}
      - \frac{\partial^{|\alpha|} \theta(w, v)}{\partial w^\alpha}
    \biggr\}
    \right| \\
    &\qquad\qquad\leq
    \sum_{\alpha \in \N_0^{D_W}:\, |\alpha| = \lceil\zeta\rceil-1}
    h^{|\alpha|}
    \Clp \|\tilde w - w\|_2^{\zeta - |\alpha|}
    \leq
    C h^{\zeta}.
  \end{align*}
  As $1_n \coloneqq (1, \ldots, 1)^\T = \hat Q_h(w) e_1 \in \R^n$,
  \begin{align*}
    &\Bigl|
    e_1^\T \hat H_h(w)^{-1} B_{n,h}(w, v)
    - \theta(w, v)
    \Bigr|
    \iftoggle{journal}{\\&\quad}{}
    =
    \Bigl|
    e_1^\T \hat H_h(w)^{-1}
    \hat Q_h(w)^\T \hat K_h(w)
    \theta(\bW_n, v) / n
    - \theta(w, v)
    \Bigr| \\
    &\quad=
    \frac{1}{n}
    \Bigl|
    e_1^\T \hat H_h(w)^{-1}
    \hat Q_h(w)^\T \hat K_h(w)
    \bigl\{
      \theta(\bW_n, v)
      - \theta(w, v) 1_n
    \bigr\}
    \Bigr| \\
    &\quad=
    \frac{1}{n}
    \biggl|
    e_1^\T \hat H_h(w)^{-1}
    \sum_{i=1}^n
    Q_h(W_i - w) K_h(W_i - w)
    \bigl\{
      \theta(W_i, v)
      - \theta(w, v)
    \bigr\}
    \biggr| \\
    &\quad\leq
    \sum_{\alpha \in \N_0^{D_W}:\, |\alpha| \in [\lceil\zeta\rceil-1]}
    \biggl|
    h^{|\alpha|}
    e_1^\T
    e_\alpha
    \frac{\partial^{|\alpha|} \theta(w, v)}{\partial w^\alpha}
    \biggr| \\
    &\qquad\qquad+
    \frac{C h^\zeta}{n}
    \bigl\|\hat H_h(w)^{-1}\bigr\|_{2 \to 2}
    \biggl\|\sum_{i=1}^n \|Q_h(W_i - w)\|_\infty 1_{D_q}
    K_h(W_i - w) \biggr\|_2 \\
    &\quad\leq
    \frac{C h^\zeta}{n}
    \bigl\|\hat H_h(w)^{-1}\bigr\|_{2 \to 2}
    \sum_{i=1}^n
    K_h(W_i - w)
    = C h^\zeta
    \bigl\|\hat H_h(w)^{-1}\bigr\|_{2 \to 2}
    e_1^\T \hat H_h(w) e_1 \\
    &\quad\leq
    C h^\zeta
    \bigl\|\hat H_h(w)^{-1}\bigr\|_{2 \to 2}
    \bigl\|\hat H_h(w)\bigr\|_{2 \to 2}.
  \end{align*}
  Thus, by Lemma~\ref{lem:locpol_Hhat},
  with probability at least $1 - 2 e^{-n h^{D_W}/C}$,
  \begin{equation*}
    \sup_{(w, v) \in \cW \times \cV}
    \Bigl|
    e_1^\T \hat H_h(w)^{-1} B_{n,h}(w, v)
    - \theta(w, v)
    \Bigr|
    \leq
    C h^\zeta.
    \iftoggle{journal}{}{\qedhere}
  \end{equation*}
\end{proof}

\begin{lemma}[Bound on $B_{n,h}^*$]%
  \label{lem:locpol_bias_bootstrap}
  Grant the setup of Theorem~\ref{thm:locpol_bootstrap}.
  The result of Lemma~\ref{lem:locpol_bias} holds with
  $\hat H_h(w)^{-1} B_{n,h}(w, v)$ replaced with
  $\hat H_h^*(w)^{-1} B_{n,h}^*(w, v)$.
\end{lemma}

\begin{proof}[Proof of Lemma~\ref{lem:locpol_bias_bootstrap}]
  The proof is identical to that of Lemma~\ref{lem:locpol_bias},
  bounding $\bigl\|\hat H_h^*(w)^{-1}\bigr\|_{2 \to 2}$
  using Lemma~\ref{lem:locpol_Hhatstar}
  rather than Lemma~\ref{lem:locpol_Hhat}.
\end{proof}

\begin{lemma}[Uniform convergence of local polynomial estimators]%
  \label{lem:locpol_convergence}
  Grant the setup of Section~\ref{sec:locpol} and Assumption~\ref{ass:locpol},
  and take $\kappa \geq 1$.
  There exists $C > 0$ depending only on
  $\Clp$, $D_W$, $D_V$, $q$, $r$ and $\kappa$
  such that for all $t > 0$,
  with probability at least
  $1 - 1/n^{\kappa} - 2 e^{-n h^{D_W}/C} - e^{-t^2/2}$,
  \begin{equation*}
    \sup_{(w, v) \in \cW \times \cV}
    \bigl|
    \hat\theta_{n,h}(w, v) - \theta(w, v)
    \bigr|
    \leq
    C \sqrt{\frac{1 + D_w \log 2n + t^2}{n h^{D_W}}}
    + C h^\zeta.
  \end{equation*}
\end{lemma}

\begin{proof}[Proof of Lemma~\ref{lem:locpol_convergence}]
  Apply Lemmas~\ref{lem:locpol_Hhat},
  \ref{lem:locpol_linear} and~\ref{lem:locpol_bias}
  to the decomposition
  \begin{equation*}
    \hat\theta_{n,h}(w, v) - \theta(w, v)
    =
    e_1^\T \hat H_h(w)^{-1}
    \tilde N_{n,h}(w, v)
    + e_1^\T \hat H_h(w)^{-1} B_{n,h}(w, v)
    - \theta(w, v).
    \iftoggle{journal}{}{\qedhere}
  \end{equation*}
\end{proof}

\begin{lemma}[Bound on $S_{n,h}^\circ$]%
  \label{lem:locpol_Scirc}
  Grant the setup of Theorem~\ref{thm:locpol_bootstrap}
  and take $\kappa \geq 1$.
  There exists $C > 0$ depending only on
  $\Clp$, $D_W$, $D_V$, $q$, $r$, $r'$ and $\kappa$ such that
  for all $t > 0$, with probability at least
  $1 - 1/n^\kappa - e^{-t^2/2}$,
  \begin{equation*}
    \sup_{(w, v) \in \cW \times \cV}
    \bigl\| S_{n,h}^\circ(w, v) \bigr\|_2
    \leq
    C \sqrt{\frac{1 + D_w \log 2n + t^2}{n h^{D_W}}}.
  \end{equation*}
\end{lemma}

\begin{proof}[Proof of Lemma~\ref{lem:locpol_Scirc}]
  By the triangle inequality,
  \begin{align*}
    \sup_{(w,v)\in\cW\times\cV}
    \|S_{n,h}^\circ(w,v)\|_2
    &\leq
    \sup_{(w,v)\in\cW\times\cV}
    \biggl\|
    S_{n,h}^\circ(w,v)
    -\frac1n\sum_{i=1}^n\chi_i
    \tilde N_{n,h}(w,v)
    \biggr\|_2 \\
    &\qquad+
    \biggl|\frac1n\sum_{i=1}^n\chi_i\biggr|
    \sup_{(w,v)\in\cW\times\cV}
    \|\tilde N_{n,h}(w,v)\|_2.
  \end{align*}
  The first term may be bounded as in the proof of
  Lemma~\ref{lem:locpol_linear},
  using Lemma~\ref{lem:bootstrap_tail} instead of
  Theorem~\ref{thm:empirical_tail}.
  For the second term, Theorem~\ref{thm:empirical_tail} gives
  $\bigl|\sum_{i=1}^n\chi_i / n\bigr| \leq C$
  with probability at least $1-n^{-\kappa}$,
  while $\sup_{(w,v)\in\cW\times\cV} \|\tilde N_{n,h}(w,v)\|_2$
  is controlled using Lemma~\ref{lem:locpol_linear}.
  Replace $\kappa$ with $\kappa+2$ and $t$ with $t+2$,
  then take a union bound to conclude.
\end{proof}

\begin{lemma}[Bound on $R_{n,h}^\circ$]%
  \label{lem:locpol_R}
  Grant the setup of Theorem~\ref{thm:locpol_bootstrap}
  and take $\kappa \geq 1$.
  There exists $C > 0$ depending only on
  $\Clp$, $D_W$, $D_V$, $q$, $r$, $r'$ and $\kappa$ such that
  with probability at least $1 - 1/n^\kappa - 2 e^{-n h^{D_W}/C}$,
  \begin{equation*}
    \sup_{(w, v) \in \cW \times \cV}
    \bigl\| R_{n,h}^\circ(w, v) \bigr\|_2
    \leq
    C \sqrt{\frac{\log 2n}{n h^{D_W}}}
    \Biggl\{
      \sqrt{\frac{\log 2n}{n h^{D_W}}}
      + h^\zeta
    \Biggr\}.
  \end{equation*}
\end{lemma}

\begin{proof}[Proof of Lemma~\ref{lem:locpol_R}]
  Throughout the proof, $C > 0$ depends only on the stated quantities
  and may change from line to line.
  Fix a realization of the data
  $(W_1, V_1), \ldots, (W_n, V_n)$ and define
  \begin{equation*}
    f_{w,v,j}(W, \chi)
    \coloneqq
    e_j^\T Q_h(W - w) K_h(W - w) \chi
    \bigl\{ \hat\theta_{n,h}(W, v) - \theta(W, v) \bigr\}
  \end{equation*}
  and the function class
  $\hat\cF \coloneqq \bigl\{(W, \chi) \mapsto f_{w,v,j}(W, \chi):
  w \in \cW, v \in \cV, j \in [D_q]\bigr\}$. Therefore,
  $e_j^\T R_{n,h}^\circ(w, v) =
  \frac{1}{n} \sum_{i=1}^n f_{w, v, j}(W_i, \chi_i)$.
  For each fixed realization of the data, $\hat\cF$ is pointwise
  measurable by continuity in $w$ and,
  when $\cV \neq \emptyset$,
  by coordinatewise right-continuity in $v$
  of both $\hat\theta_{n,h}(W,v)$ and $\theta(W,v)$.
  Define the envelope function
  \begin{equation*}
    \hat F(W, \chi)
    \coloneqq
    C h^{-D_W} |\chi|
    \sup_{(w, v) \in [0,1]^{D_W} \times \cV}
    \bigl| \hat\theta_{n,h}(w, v) - \theta(w, v) \bigr|
    + \frac{1 + |\chi|}{n},
  \end{equation*}
  so that $\hat F(W, \chi) \geq |f_{w,v,j}(W, \chi)|$
  for all $w, v$.
  By the conditional Orlicz quasi-triangle inequality,
  \begin{align*}
    \vvvert \hat F\vvvert_{\omega_{r'} \mid \bX_n}
    &\leq
    C h^{-D_W}
    \sup_{(w, v) \in [0,1]^{D_W} \times \cV}
    \bigl| \hat\theta_{n,h}(w, v) - \theta(w, v) \bigr|
    + \frac{C}{n}
    \eqqcolon \hat M, \\
    \frac{1}{n}
    \sum_{i=1}^n
    \E \bigl\{
      f_{w,v,j}(W_i, \chi_i)^2
      \mid \bX_n
    \bigr\}
    &\leq
    \frac{C}{n h^{D_W}}
    \sup_{w' \in [0,1]^{D_W}}
    \sum_{i=1}^n
    K_h(W_i - w') \\
    &\qquad\times
    \sup_{(w, v) \in [0,1]^{D_W} \times \cV}
    \bigl| \hat\theta_{n,h}(w, v) - \theta(w, v) \bigr|^2
    +
    \frac{\hat M^2}{n^2}
    \eqqcolon \hat \sigma^2.
  \end{align*}
  Take $Q$ a probability measure on $[0,1]^{D_W} \times \R$
  with $\|\hat F\|_{Q,2} < \infty$ and let $\xi\in(0,1]$.
  If $\|\chi\|_{Q,2}=0$, all functions in $\hat\cF$
  vanish $Q$-almost surely.
  Otherwise, define a probability measure $Q'$ on $[0,1]^{D_W}$ by
  \begin{equation*}
    Q'(A)\coloneqq
    \frac{\int \I\{W\in A\}\chi^2 \diffi Q(W,\chi)}
    {\int\chi^2 \diffi Q(W,\chi)}.
  \end{equation*}
  By the Lipschitz bound in the proof of Lemma~\ref{lem:locpol_covering},
  \begin{equation*}
    \|f_{w,v,j}-f_{w',v,j}\|_{Q,2}
    \leq \frac{C\|w-w'\|_\infty}{h}\|\hat F\|_{Q,2}.
  \end{equation*}
  Thus, $N\bigl(\xi\|\hat F\|_{Q,2}, \{f_{w,v,j}: w\in\cW\},
  \|\cdot\|_{Q,2}\bigr) \leq C / (h \xi)^{D_w}$.
  If $\cV\neq\emptyset$, partition $\cV$ according to the
  value of $\hat\theta_{n,h}(\cdot,v)$;
  there are at most $C(2n)^{D_V}$ such values.
  For $v,v'$ in the same part,
  \begin{equation*}
    \|f_{w,v,j}-f_{w,v',j}\|_{Q,2}
    \leq Ch^{-D_W}\|\chi\|_{Q,2}
    \|\theta(\cdot,v)-\theta(\cdot,v')\|_{Q',2}.
  \end{equation*}
  Since $\|\hat F\|_{Q,2}\geq\|\chi\|_{Q,2}/n$,
  by the conditional Jensen and VC argument in the proof of
  Lemma~\ref{lem:locpol_linear},
  $\log N\bigl(\xi h^{D_W} / n, \{\theta(\cdot,v): v\in\cV\},
  \|\cdot\|_{Q',2} \bigr) \leq C\log\bigl(n / (\xi h^{D_W})\bigr)$.
  Taking a union over $j\in[D_q]$ gives
  $\log N\bigl( \xi\|\hat F\|_{Q,2},\hat\cF,\|\cdot\|_{Q,2} \bigr)
  \leq C\log(2n/\xi)$. Therefore,
  \begin{equation*}
    \cJ(\varepsilon)
    \coloneqq
    \int_0^\varepsilon
    \sup_{Q \in \cQ}
    \sqrt{1 + \log N\bigl(
        \xi \|\hat F\|_{Q,2}, \hat\cF, \|\cdot\|_{Q,2}
    \bigr)}
    \diffi \xi
    \leq
    C \varepsilon \sqrt{\log (2 n / \varepsilon)}.
  \end{equation*}
  Note that
  $\sum_{i=1}^n \E\bigl\{\hat F(W_i,\chi_i)^2\mid\bX_n\bigr\} / n
  \leq C\hat M^2$.
  Since $\hat\sigma\geq\hat M/n$, we may take
  $\nu=1\land(C n)^{-1}$, giving
  $\cJ(\nu)/\nu\leq C\sqrt{\log 2n}$.
  Applying Theorem~\ref{thm:empirical_tail} conditional on the data gives
  that with probability at least $1 - e^{-t^2/2} - 1/n^\kappa$,
  \begin{equation*}
    \sup_{(w, v) \in \cW \times \cV}
    \bigl\| R_{n,h}^\circ(w, v) \bigr\|_2
    \leq
    C \biggl\{
      \frac{\hat \sigma}{\sqrt n} +
      \frac{\hat M (\log 2n)^{1/r'}}{n}
      \Bigl(
        \sqrt{\log 2n}
        + t
      \Bigr)
    \biggr\}
    \Bigl(
      \sqrt{\log 2n}
      + t
    \Bigr).
  \end{equation*}
  Next, $\hat\sigma$ and $\hat M$ are bounded using
  Lemma~\ref{lem:locpol_Hhat} and
  Lemma~\ref{lem:locpol_convergence},
  applied on $[0, 1]^{D_W}$ and
  taking $t = C \sqrt{\log 2n}$.
  Replacing $\kappa$ with $\kappa + 2$,
  with probability at least
  $1 - 1/n^{\kappa} - 2 e^{-n h^{D_W}/C}$, we have
  \begin{equation*}
    \sup_{(w, v) \in \cW \times \cV}
    \bigl\| R_{n,h}^\circ(w, v) \bigr\|_2
    \leq
    C \sqrt{\frac{\log 2n}{n h^{D_W}}}
    \Biggl\{
      \sqrt{\frac{\log 2n}{n h^{D_W}}}
      + h^\zeta
    \Biggr\},
  \end{equation*}
  as required.
\end{proof}

\begin{lemma}[Bootstrap approximation of the linear term]%
  \label{lem:locpol_bootstrap}
  Grant the setup of Theorem~\ref{thm:locpol_bootstrap}.
  Under condition~\ref{it:thm_locpol_bootstrap_nonparametric},
  for each $\eta > 0$, there exists a probability space supporting
  $\bX_n'$, $\bX_n''$, $\bX_n^{*\prime}$,
  $S_{n,h}'$ and $S_{n,h}^{*\prime}$, where
  $\bX_n''$ is independent of $\bX_n'$,
  $(\bX_n'', S_{n,h}') \equallaw (\bX_n, S_{n,h})$
  and $(\bX_n',\bX_n^{*\prime},S_{n,h}^{*\prime})
  \equallaw (\bX_n, \bX_n^*, S_{n,h}^*)$, with
  \begin{align*}
    &\P\Biggl\{
      \sup_{(w, v) \in \cW \times \cV}
      \sqrt{n h^{D_W}}
      \bigl\| S_{n,h}^{*\prime}(w, v) - S_{n,h}'(w, v) \bigr\|_2
      \iftoggle{journal}{\\ &\qquad\qquad}{}
      > \frac{C \eta \sqrt{\log 2n} }
      {(n h^{D_W + (4-k)D_w})^{1/(10-2k+2D_w+4D_v)}}
    \Biggr\}
    \leq
    \frac{1}{\eta^2},
  \end{align*}
  where $C > 0$ depends only on $\Clp$, $D_W$, $D_V$, $r$, $r'$ and $q$.
  Under condition~\ref{it:thm_locpol_bootstrap_multiplier},
  the same result holds with
  $\bX_n^{*}$ replaced by $(\chi_i: i\in[n])$,
  $\bX_n^{*\prime}$ replaced by $(\chi_i': i\in[n])$,
  $S_{n,h}^*$ replaced by $S_{n,h}^\circ$
  and $S_{n,h}^{*\prime}$ replaced by $S_{n,h}^{\circ\prime}$.
\end{lemma}

\begin{proof}[Proof of Lemma~\ref{lem:locpol_bootstrap}]
  Throughout the proof, $C > 0$ depends only on the stated quantities.
  Fix $\kappa = 4$ and assume for now that $n h^{D_W+(4-k)D_w} > 1$.
  Take $\delta_w \in (0, h]$ and let
  $\cW_\delta = \{w_l: l \in [d_w]\}$
  be a regular Cartesian $\delta_w$-cover of $\cW$
  under $\|\cdot\|_\infty$ of cardinality
  $d_w \leq \lceil 1/\delta_w \rceil^{D_w}$.
  If $\cV = \emptyset$, take $\delta_v \in (0, 1]$,
  set $d_v \coloneqq 1$, let
  $\cM \coloneqq \cW \times [D_q]$, define
  $\cM_\delta \coloneqq \cW_\delta \times [D_q]$,
  and, for $(w, j) \in \cM$, set
  $\pi(w, j) \coloneqq (w', j) \in \cM_\delta$,
  where $w' \in \argmin_{w'' \in \cW_\delta} \|w - w''\|_1$,
  resolving ties by taking the lexicographically largest element.
  Thus, $\|w'-w\|_\infty \leq \delta_w$.
  Otherwise, if $\cV = [0, 1]^{D_V}$, let
  $\cM \coloneqq \cW \times \cV \times [D_q]$,
  take $\delta_v \in (0, 1]$ and let
  $d_v \coloneqq \lceil 1 / \delta_v \rceil^{D_V}$.
  For each $l \in [d_w]$ and $s \in [D_V]$, let
  $\{t_{l,s,k}: k \in [\lceil 1 / \delta_v \rceil]\}$ satisfy
  $t_{l,s,k} =
  \inf\bigl\{t \in \R: \P\bigl(e_s^\T V_1 \leq t \mid
    \|W_1 - w_l\|_\infty \leq 3 h\bigr)
  \geq k \delta_v\bigr\}$ for
  $k \in [\lceil 1 / \delta_v \rceil - 1]$
  and $t_{l,s,\lceil 1 / \delta_v \rceil} = 0$.
  Set $\cV_\delta(w_l) \coloneqq
  \bigl\{(t_{l,1,k_1}, \ldots, t_{l,D_V,k_{D_V}}):
  k_1, \ldots, k_{D_V} \in [\lceil 1 / \delta_v \rceil]\bigr\}$
  and $\cM_\delta \coloneqq \bigl\{(w, v, j): w \in \cW_\delta,
  v \in \cV_\delta(w) \cap \cV, j \in [D_q]\bigr\}$.
  For $(w, v, j) \in \cM$, define
  $\pi(w,v,j) \coloneqq (w',v',j) \in \cM_\delta$,
  choosing $w'$ as when $\cV = \emptyset$ and taking
  $v'_s \coloneqq \max\{t_{l,s,k}: t_{l,s,k} \leq v_s\}$
  for each $s \in [D_V]$, where $w' = w_l$.
  This maximum is well-defined because $v_s \geq 0$
  and $t_{l,s,\lceil 1/\delta_v \rceil}=0$ for each $s$.
  Thus, $\|w' - w\|_\infty \leq \delta_w$,
  $v' \leq v$ and $\P(V_1 \leq v,\, V_1 \not\leq v'
  \mid \|W_1 - w'\|_\infty \leq 3h)
  \leq D_V \delta_v$ by a union bound.
  We have $d \coloneqq |\cM_\delta| \leq
  d_w d_v D_q \leq C \delta_w^{-D_w} \delta_v^{-D_v}$.
  For $(w, v, j) \in \cM$, define the functions
  \begin{equation*}
    f_{w,v,j}(W, V)
    \coloneqq
    \frac{1}{n}
    e_j^\T Q_h(W - w) K_h(W - w)
    \bigl\{g(V, v) - \theta(W, v)\bigr\},
  \end{equation*}
  and consider the function classes
  $\cF \coloneqq \{f_{w,v,j}: (w,v,j) \in \cM\}$,
  $\tilde\cF \coloneqq \{f_{w,v,j}: (w,v,j) \in \cM_\delta\}$ and
  $\check\cF \coloneqq \{f_{w,v,j} - f_{\pi(w,v,j)}: (w,v,j) \in \cM\}$.
  Observe that $\cF$ is pointwise measurable by continuity in $w$
  and, if $\cV \neq \emptyset$, by coordinatewise right-continuity in $v$.
  The class $\check\cF$ is also pointwise measurable;
  to see this, restrict the first $D_w$ coordinates of
  $w \in \cW$ and all coordinates of $v \in \cV$ to rational values
  (when $\cV \neq \emptyset$),
  leaving the remaining coordinates of $w$ fixed.
  Let $(w_m, v_m)$ be a sequence in this rational class with $w_m \to w$,
  $v_m \to v$, $w_m \geq w$ and $v_m \geq v$,
  so that eventually $\pi(w_m, v_m, j) = \pi(w, v, j)$.
  By continuity in $w$ and right-continuity in $v$,
  we deduce that
  $f_{w_m, v_m, j} - f_{\pi(w_m, v_m, j)} \to f_{w, v, j} - f_{\pi(w, v, j)}$
  pointwise, as required.
  If $\cV = \emptyset$, omit the $v$ coordinates.
  Note that $P_X(\tilde \cF) = 0$ and
  \begin{align*}
    \bigl|f_{w, v, j}(W_i, V_i)\bigr|
    &\leq
    \frac{1}{n}
    K_h(W_i - w)
    \bigl|g(V_i, v) - \theta(W_i, v)\bigr| \\
    &\leq
    \frac{C}{n h^{D_W}}
    \bigl(|g(V_i, v)| + 1\bigr)
    \I\bigl\{\|W_i - w\|_\infty \leq h\bigr\}.
  \end{align*}
  If $g(V, v) = V$, take $F(W, V) \coloneqq C n^{-1} h^{-D_W}
  \bigl(|g(V)| + 1\bigr)$,
  while if $g(V, v) = \I\{V \leq v\}$, set $F(W, V) \coloneqq C
  n^{-1} h^{-D_W}$.
  In either case, we have
  $F(W, V) \geq |f_{w,v,j}(W, V)|$ for all $(w, v, j) \in \cM$ and
  $\vvvert F(W_i, V_i) \vvvert_{\omega_r} \leq C n^{-1} h^{-D_W}$
  for each $i \in [n]$.
  For $a \in \{0, \ldots, k\}$ and $i \in [n]$,
  by Lemma~\ref{lem:moment_tail_orlicz},
  \begin{align*}
    \beta_{a,1,i}
    &\coloneqq
    \frac{1}{\sqrt{a!}}
    \E \Bigl\{
      \bigl\|\tilde\cF(W_i, V_i)\bigr\|_2^a
      \bigl\|\tilde\cF(W_i, V_i)\bigr\|_\infty
    \Bigr\} \\
    &\,\leq
    \frac{C}{(n h^{D_W})^{a+1}}
    \E \Biggl\{
      \Biggl(
        \sum_{(w, v, j) \in \cM_\delta}
        \I\bigl\{\|W_i - w\|_\infty \leq h\bigr\}
      \Biggr)^{a/2}
      \iftoggle{journal}{\\ &\qquad\qquad\qquad\qquad\qquad\times}{}
      \I\bigl\{\exists w\in\cW_\delta:
      \|W_i-w\|_\infty\leq h\bigr\}
    \Biggr\} \\
    &\,\leq
    \frac{C}{(n h^{D_W})^{a + 1}}
    \biggl( \frac{h}{\delta_w} \biggr)^{a D_w/2}
    \biggl( \frac{1}{\delta_v} \biggr)^{a D_v/2}
    \P\bigl(
      \exists w \in \cW_\delta:
      \|W_i - w\|_\infty \leq h
    \bigr) \\
    &\,\leq
    \frac{C}{(n h^{D_W})^{a + 1}}
    \biggl( \frac{h}{\delta_w} \biggr)^{a D_w/2}
    \biggl( \frac{1}{\delta_v} \biggr)^{a D_v/2}
    h^{D_W - D_w}.
  \end{align*}
  Similarly, with $P_n(f) \coloneqq \sum_{i=1}^n f(X_i) / n$,
  by Jensen's inequality,
  under condition~\ref{it:thm_locpol_bootstrap_nonparametric},
  \begin{align*}
    \hat\beta^*_{a,1,i}
    &\coloneqq
    \frac{1}{\sqrt{a!}}
    \E^* \Bigl\{
      \bigl\|\tilde\cF(W_i^*, V_i^*) - P_n(\tilde\cF)\bigr\|_2^a
      \bigl\|\tilde\cF(W_i^*, V_i^*) - P_n(\tilde\cF)\bigr\|_\infty
    \Bigr\} \\
    &\,\leq
    \frac{C}{n}
    \sum_{i=1}^n
    \bigl\|\tilde\cF(W_i, V_i)\bigr\|_2^a
    \bigl\|\tilde\cF(W_i, V_i)\bigr\|_\infty
    + \frac{C}{n^2}
    \sum_{i=1}^n
    \sum_{j=1}^n
    \bigl\|\tilde\cF(W_i, V_i)\bigr\|_2^a
    \bigl\|\tilde\cF(W_j, V_j)\bigr\|_\infty.
  \end{align*}
  Therefore, $\E\bigl(\hat\beta^*_{a,1,i}\bigr)$
  admits the same upper bound as $\beta_{a,1,i}$.
  Similarly, for $c \in \{0, \ldots, k\}$,
  \begin{equation*}
    \E\bigl(
      \hat\beta^*_{a,1,i}
      \hat\beta^*_{c,0,i}
    \bigr)
    \leq
    \frac{C}{(n h^{D_W})^{a + c + 1}}
    \biggl( \frac{h}{\delta_w} \biggr)^{(a+c) D_w/2}
    \biggl( \frac{1}{\delta_v} \biggr)^{(a+c) D_v/2}
    h^{D_W - D_w}.
  \end{equation*}
  Likewise, under condition~\ref{it:thm_locpol_bootstrap_multiplier},
  since $a+1 \leq 4$ and $\vvvert\chi_i\vvvert_4 \leq 2$,
  \begin{equation*}
    \hat\beta^\circ_{a,1,i}
    \coloneqq
    \frac{1}{\sqrt{a!}}
    \E^* \Bigl\{
      \bigl\|\chi_i \tilde\cF(W_i, V_i)\bigr\|_2^a
      \bigl\|\chi_i \tilde\cF(W_i, V_i)\bigr\|_\infty
    \Bigr\}
    \leq
    C \bigl\|\tilde\cF(W_i, V_i)\bigr\|_2^a
    \bigl\|\tilde\cF(W_i, V_i)\bigr\|_\infty.
  \end{equation*}
  Thus, $\E\bigl(\hat\beta^\circ_{a,1,i}\bigr)$
  and $\E\bigl(\hat\beta^\circ_{a,1,i}\hat\beta^\circ_{c,0,i}\bigr)$
  admit the same upper bounds as $\beta_{a,1,i}$ and
  $\E\bigl(\hat\beta^*_{a,1,i}\hat\beta^*_{c,0,i}\bigr)$,
  respectively.
  Write $\hat\beta_{a,b,i} \coloneqq \beta_{a,b,i} + \hat\beta^*_{a,b,i}$
  under condition~\ref{it:thm_locpol_bootstrap_nonparametric}
  and $\hat\beta_{a,b,i} \coloneqq \beta_{a,b,i} + \hat\beta^\circ_{a,b,i}$
  under condition~\ref{it:thm_locpol_bootstrap_multiplier}.
  We apply Theorem~\ref{thm:empirical_process_bootstrap} with $\gamma = 1$,
  using part (ii) of the theorem
  under condition~\ref{it:thm_locpol_bootstrap_multiplier},
  to see that for $s, \eta > 0$, we have
  \begin{equation*}
    \frac{s^k}{\eta}
    \sum_{i=1}^n
    \E\bigl\{\hat\beta_{k,1,i}\bigr\}
    \leq
    \frac{C s^k h^{- D_w}}{(n h^{D_W})^k \eta}
    \biggl( \frac{h}{\delta_w} \biggr)^{k D_w/2}
    \biggl( \frac{1}{\delta_v} \biggr)^{k D_v/2}.
  \end{equation*}
  By Lemma~\ref{lem:bound_Gamma}, with
  $\hat\beta_{a,b} \coloneqq \sum_{i=1}^n \hat\beta_{a,b,i}/n$,
  taking $k_\ell = \I\{k=3, \ell=2\}$ so $\ell + k_\ell = k$,
  \begin{align*}
    \E\bigl\{\hat\Gamma_\eta(s)\bigr\}
    &\leq
    \sum_{i=1}^n
    \sum_{\ell = 2}^k
    \frac{(2 s)^{\ell + k_\ell}}
    {2 \eta}
    \E\Bigl\{
      \bigl(\hat \beta_{k_\ell,1}
      + \hat \beta_{k_\ell,1,i}\bigr)
      \hat\beta_{\ell,0,i}
    \Bigr\} \\
    &\leq
    \sum_{i=1}^n
    \sum_{\ell = 2}^k
    \frac{(2 s)^{\ell + k_\ell}}
    {2 \eta}
    \frac{C}{(n h^{D_W})^{k_\ell + \ell + 1}}
    \biggl( \frac{h}{\delta_w} \biggr)^{(k_\ell+\ell) D_w/2}
    \biggl( \frac{1}{\delta_v} \biggr)^{(k_\ell+\ell) D_v/2}
    h^{D_W - D_w} \\
    &\leq
    \frac{C s^k h^{- D_w} }{(n h^{D_W})^k \eta}
    \biggl( \frac{h}{\delta_w} \biggr)^{k D_w/2}
    \biggl( \frac{1}{\delta_v} \biggr)^{k D_v/2}.
  \end{align*}
  For $\alpha \in \N_0^d$ with $2 \leq |\alpha| \leq k$, write
  \begin{align*}
    \hat \delta^*_{\alpha}
    &=
    \frac{1}{n}
    \sum_{i=1}^n
    \E\bigl\{
      \tilde \cF(X_i)^\alpha
    \bigr\}
    - \frac{1}{n}
    \sum_{i=1}^{n}
    \bigl\{
      \tilde \cF(X_i) - P_{n}(\tilde\cF)
    \bigr\}^\alpha, \\
    \hat \delta^\circ_{\alpha}
    &=
    \frac{1}{n}
    \sum_{i=1}^n
    \E\bigl\{
      \tilde \cF(X_i)^\alpha
    \bigr\}
    - \frac{1}{n}
    \sum_{i=1}^{n}
    \tilde \cF(X_i)^\alpha.
  \end{align*}
  Under condition~\ref{it:thm_locpol_bootstrap_nonparametric},
  \begin{align}
    \label{eq:boot_empirical_1}
    \hat \delta^*_{\alpha}
    &=
    \frac{1}{n}
    \sum_{i=1}^n
    \E\bigl\{
      \tilde \cF(X_i)^\alpha
    \bigr\}
    -
    \frac{1}{n}
    \sum_{i=1}^{n}
    \tilde \cF(X_i)^\alpha \\
    \label{eq:boot_empirical_2}
    &\qquad+
    \frac{1}{n}
    \sum_{i=1}^{n}
    \Bigl\{
      \tilde \cF(X_i)^\alpha
      - \bigl(
        \tilde \cF(X_{i}) - P_{n}(\tilde \cF)
      \bigr)^\alpha
    \Bigr\}.
  \end{align}
  Let $2 \leq \ell \leq k$.
  For \eqref{eq:boot_empirical_1}, by the multinomial theorem and
  Lemma~\ref{lem:moment_tail_orlicz},
  \begin{align*}
    &\sum_{|\alpha| = \ell}
    \frac{1}{\alpha!}
    \E \Biggl\{
      \biggl(
        \frac{1}{n}
        \sum_{i=1}^n
        \E\bigl\{
          \tilde \cF(X_i)^\alpha
        \bigr\}
        - \frac{1}{n}
        \sum_{i=1}^{n}
        \tilde \cF(X_i)^\alpha
      \biggr)^2
    \Biggr\} \\
    &\quad\leq
    \frac{1}{n^2}
    \sum_{|\alpha| = \ell}
    \frac{1}{\alpha!}
    \sum_{i=1}^n
    \E \bigl\{
      \tilde \cF(X_i)^{2\alpha}
    \bigr\}
    \leq
    \frac{C}{n^2}
    \sum_{i=1}^n
    \E \Bigl\{
      \bigl\|
      \tilde \cF(X_i)
      \bigr\|_2^{2 \ell}
    \Bigr\} \\
    &\quad\leq
    \frac{C}{n (n h^{D_W})^{2 \ell}}
    \biggl( \frac{h}{\delta_w} \biggr)^{\ell D_w}
    \biggl( \frac{1}{\delta_v} \biggr)^{\ell D_v}
    h^{D_W - D_w}.
  \end{align*}
  For \eqref{eq:boot_empirical_2}, by Jensen's inequality,
  Lemma~\ref{lem:multi_index_sum} and Cauchy--Schwarz,
  \begin{align*}
    &\sum_{|\alpha| = \ell}
    \frac{1}{\alpha!}
    \E\Biggl\{
      \biggl(
        \frac{1}{n}
        \sum_{i=1}^{n}
        \Bigl\{
          \tilde \cF(X_i)^\alpha
          - \bigl(
            \tilde \cF(X_{i}) - P_{n}(\tilde \cF)
          \bigr)^\alpha
        \Bigr\}
      \biggr)^2
    \Biggr\} \\
    &\quad\leq
    \frac{1}{n}
    \sum_{i=1}^n
    \sum_{|\alpha| = \ell}
    \frac{1}{\alpha!}
    \E\biggl(
      \Bigl\{
        \tilde \cF(X_i)^\alpha
        - \bigl(
          \tilde\cF(X_{i})
          - P_{n}(\tilde\cF)
        \bigr)^\alpha
      \Bigr\}^2
    \biggr) \\
    &\quad\leq
    \frac{C}{n}
    \sum_{i=1}^n
    \E\biggl(
      \bigl\|
      P_{n}(\tilde\cF)
      \bigr\|_2^2
      \Bigl\{
        \bigl\|
        \tilde \cF(X_i)
        \bigr\|_2^{2\ell-2}
        + \bigl\|
        \tilde\cF(X_{i}) - P_{n}(\tilde\cF)
        \bigr\|_2^{2\ell-2}
      \Bigr\}
    \biggr) \\
    &\quad\leq
    \frac{C}{n}
    \sum_{i=1}^n
    \Bigl\{\E\Bigl(
        \bigl\|
        P_{n}(\tilde\cF)
        \bigr\|_2^4
    \Bigr)\Bigr\}^{1/2}
    \Bigl\{\E\Bigl(
        \bigl\|\tilde\cF(X_{i})\bigr\|_2^{4\ell-4}
    \Bigr)\Bigr\}^{1/2}.
  \end{align*}
  By Lemma~\ref{lem:moment_tail_orlicz} again,
  \begin{equation*}
    \E\Bigl(
      \bigl\|\tilde\cF(X_{i})\bigr\|_2^{4\ell-4}
    \Bigr)
    \leq
    \frac{C}{(n h^{D_W})^{4\ell-4}}
    \biggl(\frac{h}{\delta_w}\biggr)^{(2\ell-2) D_w}
    \biggl(\frac{1}{\delta_v}\biggr)^{(2\ell-2) D_v}
    h^{D_W - D_w}.
  \end{equation*}
  For $f, f' \in \tilde\cF$, writing $g_i \coloneqq f(X_i)/n$
  and $g_i' \coloneqq f'(X_i)/n$,
  noting that these variables
  are independent and zero-mean across $i \in [n]$,
  \begin{align*}
    &\E \Bigl\{
      \bigl\|
      P_{n}(\tilde\cF)
      \bigr\|_2^4
    \Bigr\}
    =
    \E \Biggl\{
      \biggl(
        \sum_{f \in \tilde\cF}
        P_{n}(f)^2
      \biggr)^2
    \Biggr\} \\
    &\quad=
    \sum_{f \in \tilde\cF}
    \sum_{f' \in \tilde\cF}
    \sum_{i=1}^{n}
    \sum_{j=1}^{n}
    \sum_{i'=1}^{n}
    \sum_{j'=1}^{n}
    \E \bigl\{ g_{i} g'_{j} g_{i'} g'_{j'} \bigr\}
    \leq
    C
    \sum_{f \in \tilde\cF}
    \sum_{f' \in \tilde\cF}
    \sum_{i=1}^{n}
    \sum_{j=1}^{n}
    \E \bigl\{
      g_i g_i' g_j g_j'
      + g_{i}^2 g_{j}^{\prime 2}
    \bigr\} \\
    &\quad\leq
    C
    \sum_{i=1}^{n}
    \sum_{j=1}^{n}
    \E \Biggl\{
      \biggl(
        \sum_{f \in \tilde\cF} g_i g_j
      \biggr)^2
    \Biggr\}
    + C
    \Biggl\{
      \sum_{i=1}^n
      \sum_{f \in \tilde\cF}
      \E \bigl( g_{i}^2 \bigr)
    \Biggr\}^2.
  \end{align*}
  By Lemma~\ref{lem:moment_tail_orlicz},
  \begin{align*}
    \E \Biggl\{
      \biggl(
        \sum_{f \in \tilde\cF} g_i g_j
      \biggr)^2
    \Biggr\}
    &=
    \frac{1}{n^4}
    \E \Biggl\{
      \biggl(
        \sum_{f \in \tilde\cF}
        f(X_{i})
        f(X_{j})
      \biggr)^2
    \Biggr\}
    \iftoggle{journal}{\\ &}{}
    \leq
    \frac{C}{n^4 (n h^{D_W})^4}
    \biggl(\frac{h}{\delta_w}\biggr)^{2 D_w}
    \biggl(\frac{1}{\delta_v}\biggr)^{2 D_v}
    h^{D_W - D_w}.
  \end{align*}
  Similarly,
  \begin{equation*}
    \sum_{f \in \tilde\cF}
    \E \bigl( g_{i}^2 \bigr)
    =
    \frac{1}{n^2}
    \sum_{f \in \tilde\cF}
    \E \bigl\{
      f(X_{i})^2
    \bigr\}
    \leq
    \frac{C}{n^2 (n h^{D_W})^2}
    \biggl(\frac{h}{\delta_w}\biggr)^{D_w}
    \biggl(\frac{1}{\delta_v}\biggr)^{D_v}
    h^{D_W - D_w}.
  \end{equation*}
  It follows that
  \begin{equation*}
    \E \Bigl\{
      \bigl\|
      P_{n}(\tilde\cF)
      \bigr\|_2^4
    \Bigr\}
    \leq
    \frac{C}{n^2 (n h^{D_W})^4}
    \biggl(\frac{h}{\delta_w}\biggr)^{2 D_w}
    \biggl(\frac{1}{\delta_v}\biggr)^{2 D_v}
    h^{D_W - D_w}.
  \end{equation*}
  Therefore,
  \begin{align*}
    &\sum_{|\alpha| = \ell}
    \frac{1}{\alpha!}
    \E\Biggl\{
      \biggl(
        \frac{1}{n}
        \sum_{i=1}^{n}
        \Bigl\{
          \tilde \cF(X_i)^\alpha
          - \bigl(
            \tilde \cF(X_{i}) - P_{n}(\tilde \cF)
          \bigr)^\alpha
        \Bigr\}
      \biggr)^2
    \Biggr\}
    \iftoggle{journal}{\\ &\qquad}{}
    \leq
    \frac{C}{n (n h^{D_W})^{2\ell}}
    \biggl(\frac{h}{\delta_w}\biggr)^{\ell D_w}
    \biggl(\frac{1}{\delta_v}\biggr)^{\ell D_v}
    h^{D_W - D_w}.
  \end{align*}
  By combining the bounds on
  \eqref{eq:boot_empirical_1}
  and \eqref{eq:boot_empirical_2}, we deduce that
  \begin{equation*}
    \sum_{|\alpha| = \ell}
    \frac{\E\bigl(\hat \delta^{*2}_{\alpha}\bigr)}{\alpha!}
    \leq
    \frac{C}{n (n h^{D_W})^{2\ell}}
    \biggl(\frac{h}{\delta_w}\biggr)^{\ell D_w}
    \biggl(\frac{1}{\delta_v}\biggr)^{\ell D_v}
    h^{D_W - D_w}.
  \end{equation*}

  Under condition~\ref{it:thm_locpol_bootstrap_multiplier},
  the corresponding term involving
  $\hat\delta_\alpha^\circ$ may be bounded in the same manner
  as \eqref{eq:boot_empirical_1}.
  Thus, writing $\hat\delta_\alpha \coloneqq \hat\delta_\alpha^*$
  under condition~\ref{it:thm_locpol_bootstrap_nonparametric} and
  $\hat\delta_\alpha \coloneqq \hat\delta_\alpha^\circ$
  under condition~\ref{it:thm_locpol_bootstrap_multiplier},
  by Jensen's inequality,
  \begin{align*}
    1 \land
    \E\bigl\{
      \hat\Delta(s)
    \bigr\}
    &\leq
    1 \land
    \frac{n}{2}
    \Biggl(
      \sum_{\ell = 2}^k
      s^{2 \ell}
      \sum_{|\alpha| = \ell}
      \frac{\E\bigl(\hat \delta_{\alpha}^2\bigr)}{\alpha!}
    \Biggr)^{1/2} \\
    &\leq
    1 \land
    C h^{D_W/2 - D_w/2}
    \sqrt{n}
    \Biggl(
      \sum_{\ell = 2}^k
      \frac{s^{2 \ell}}{(n h^{D_W})^{2\ell}}
      \biggl(\frac{h}{\delta_w}\biggr)^{\ell D_w}
      \biggl(\frac{1}{\delta_v}\biggr)^{\ell D_v}
    \Biggr)^{1/2} \\
    &\leq
    \frac{C s^2 h^{D_W/2 - D_w/2}}{(n h^{D_W})^2}
    \biggl(\frac{h}{\delta_w}\biggr)^{D_w}
    \biggl(\frac{1}{\delta_v}\biggr)^{D_v}
    \sqrt{n}.
  \end{align*}

  Recalling that $\pi(w, v, j) = (w', v', j)$, we have
  \begin{align*}
    &\bigl|
    f_{w,v,j}(W_i, V_i)
    - f_{w',v',j}(W_i, V_i)
    \bigr| \\
    &\quad=
    \frac{1}{n}
    \Bigl|
    e_j^\T
    Q_h(W_i - w) K_h(W_i - w)
    \bigl\{g(V_i, v) - \theta(W_i, v)\bigr\} \\
    &\quad\qquad-
    e_j^\T
    Q_h(W_i - w') K_h(W_i - w')
    \bigl\{g(V_i, v') - \theta(W_i, v')\bigr\}
    \Bigr| \\
    &\quad\leq
    \frac{1}{n}
    \bigl\| Q_h(W_i - w) - Q_h(W_i - w') \bigr\|_\infty \\
    &\quad\qquad\qquad\times
    \, K_h(W_i - w)
    \, \bigl|g(V_i, v) - \theta(W_i, v)\bigr| \\
    &\quad\qquad+
    \frac{1}{n}
    \bigl\| Q_h(W_i - w') \bigr\|_\infty
    \, \I\{\|W_i - w'\|_\infty \leq 2 h\} \\
    &\quad\qquad\qquad\times
    \bigl|K_h(W_i - w) - K_h(W_i - w')\bigr|
    \, \bigl|g(V_i, v) - \theta(W_i, v)\bigr| \\
    &\quad\qquad+
    \frac{1}{n}
    \bigl\| Q_h(W_i - w') \bigr\|_\infty
    K_h(W_i - w') \\
    &\quad\qquad\qquad\times
    \bigl|
    g(V_i, v) - \theta(W_i, v)
    - g(V_i, v') + \theta(W_i, v')
    \bigr| \\
    &\quad\leq
    \frac{C \delta_w}{n h^{D_W+1}}
    \, \I\{\|W_i - w\|_\infty \leq 2 h\}
    \, \bigl|g(V_i, v) - \theta(W_i, v)\bigr| \\
    &\quad\qquad+
    \frac{C}{n h^{D_W}}
    \, \I\{\|W_i - w\|_\infty \leq 2 h\}
    \, \bigl|
    g(V_i, v) - \theta(W_i, v)
    - g(V_i, v') + \theta(W_i, v')
    \bigr|.
  \end{align*}
  If $g(V) = V$, define
  $\check F(W, V) \coloneqq C (|g(V)| + 1) n^{-1} h^{-D_W}$,
  and if $g(V, v) = \I\{V \leq v\}$, let
  $\check F(W, V) \coloneqq C n^{-1} h^{-D_W}$.
  Thus, $\vvvert \check F(W_i, V_i) \vvvert_{\omega_r} \leq
  \check M \coloneqq C n^{-1} h^{-D_W}$
  for all $i \in [n]$.
  Further, as $f_{W,i} \leq \Clp^2 f_{W,1}$,
  \begin{equation*}
    \E\Bigl(\bigl|
      f_{w,v,j}(W_i, V_i)
      - f_{w',v',j}(W_i, V_i)
    \bigr|^2 \Bigr)
    \leq
    \frac{C \delta_w^2}{n^2 h^{D_W+2}}
    + \frac{C \delta_v}{n^2 h^{D_W}}
    \eqqcolon
    \sigma^2,
  \end{equation*}
  and $\nu \coloneqq 1 \land h^{D_W/2}
  \bigl(\delta_w / h + \delta_v^{1/2}\bigr)$.
  The constant defining $\sigma^2$ is chosen sufficiently large that
  $\nu \leq \sigma / \|\check F\|_{\bar P_X,2}$.
  Let $Q$ be a Borel probability
  measure on $[0,1]^{D_W} \times \R^{D_V}$.
  Since $\check\cF \subseteq \cF-\cF$, the covering argument
  in the proof of Lemma~\ref{lem:locpol_linear} gives,
  for all $\varepsilon \in (0,1]$,
  \begin{equation*}
    \log N\bigl(
      \varepsilon\|\check F\|_{Q,2},
      \check\cF,\|\cdot\|_{Q,2}
    \bigr)
    \leq
    2\log N\bigl(
      \varepsilon\|\check F\|_{Q,2}/4,
      \cF,\|\cdot\|_{Q,2}
    \bigr)
    \leq C\log\biggl(\frac{2}{h\varepsilon}\biggr).
  \end{equation*}

  Choose $\delta_w$ so that
  $h/\delta_w = \bigl(n h^{D_W+(4-k)D_w}\bigr)^{1/(10-2k+2D_w+4D_v)}$
  and set $\delta_v \coloneqq \delta_w^2/h^2$,
  so that $\log\{2/(h\nu)\} \leq C\log 2n$.
  By Lemma~\ref{lem:entropy_vc}, we therefore have
  $\cJ(\nu) \leq C\nu\sqrt{\log 2n}$.
  Thus, for $t \leq C\sqrt{\log 2n}$,
  \begin{align*}
    R^*(t)
    &\coloneqq
    \biggl\{
      \sigma \sqrt n
      + \check M (\kappa \log 2n)^{1/r + 1/r'}
      \biggl(\frac{\cJ( \nu)}{ \nu} + t\biggr)
    \biggr\}
    \biggl(\frac{\cJ( \nu)}{ \nu} + t\biggr) \\
    &\,\leq
    C \Biggl\{
      \sqrt{\frac{\log 2n + t^2}{n h^{D_W}}}
      \biggl(
        \frac{\delta_w}{h}
        + \sqrt{\delta_v}
      \biggr)
      + \frac{(\log 2n)^{1 + 1/r + 1/r'}}{n h^{D_W}}
    \Biggr\}.
  \end{align*}
  For $(w,v,j) \in \cM$, let
  $S_n(w,v,j) \coloneqq e_j^\T S_{n,h}(w, v)$
  and $S_n^*(w, v, j) \coloneqq
  e_j^\T S_{n,h}^*(w,v)$.
  Under condition~\ref{it:thm_locpol_bootstrap_nonparametric},
  using the construction in the proof of
  Theorem~\ref{thm:empirical_process_bootstrap},
  retaining $\bX_n'$, $\bX_n''$ and $\bX_n^{*\prime}$,
  we have
  \begin{align*}
    &
    \P\biggl\{
      \sup_{(w, v, j) \in \cM}
      |S_n'(w, v, j) - S_n^{*\prime}(w, v, j)|
      > 3 \eta \\
      &\qquad\qquad\qquad+
      C
      \sqrt{\frac{\log 2n + t^2}{n h^{D_W}}}
      \biggl(
        \frac{\delta_w}{h}
        + \sqrt{\delta_v}
      \biggr)
      + \frac{C (\log 2n)^{1 + 1/r + 1/r'}}{n h^{D_W}}
    \biggr\} \\
    &\qquad\leq
    4 e^{-t^2 / 2}
    + \ \frac{C \phi(t)^2 h^{D_W/2 - D_w/2}}{(n h^{D_W})^2 \eta^2}
    \biggl(\frac{h}{\delta_w}\biggr)^{D_w}
    \biggl(\frac{1}{\delta_v}\biggr)^{D_v}
    \sqrt{n} \\
    &\qquad\qquad+
    \frac{C \phi(t)^k h^{- D_w}}{(n h^{D_W})^k \eta^{k+1}}
    \biggl( \frac{h}{\delta_w} \biggr)^{k D_w/2}
    \biggl( \frac{1}{\delta_v} \biggr)^{k D_v/2}
    + \frac{1}{n^\kappa}.
  \end{align*}
  Take $t \coloneqq \sqrt{2\kappa \log 2n}$
  and replace $\eta$ with $\eta/3$.
  Since $\log 2d \leq C\log 2n$, increasing $C$ if necessary,
  \begin{align*}
    &
    \P\biggl\{
      \sup_{(w, v, j) \in \cM}
      |S_n'(w, v, j) - S_n^{*\prime}(w, v, j)|
      > \eta
      + \frac{C \delta_w}{h}
      \sqrt{\frac{\log 2n}{n h^{D_W}}}
      + \frac{C (\log 2n)^{1 + 1/r + 1/r'}}{n h^{D_W}}
    \biggr\} \\
    &\qquad\leq
    \frac{C h^{- D_w/2} \log 2n}{(n h^{D_W})^{3/2} \eta^2}
    \biggl(\frac{h}{\delta_w}\biggr)^{D_w + 2 D_v}
    + \frac{C h^{- D_w} (\log 2n)^{k/2}}{(n h^{D_W})^k \eta^{k+1}}
    \biggl( \frac{h}{\delta_w} \biggr)^{k D_w/2 + k D_v}
    + \frac{C}{n^\kappa}.
  \end{align*}
  Replacing $\eta$ with
  \begin{equation*}
    \frac{h^{-D_w/4} \eta}{(n h^{D_W})^{3/4}}
    \biggl(\frac{h}{\delta_w}\biggr)^{D_w/2 + D_v}
    +
    \frac{h^{-D_w/(k+1)} \eta}{(n h^{D_W})^{k/(k+1)}}
    \biggl(\frac{h}{\delta_w}\biggr)^{k D_w/(2k+2) + k D_v/(k+1)}
  \end{equation*}
  gives
  \begin{align*}
    &
    \P\Biggl\{
      \sup_{(w, v, j) \in \cM}
      |S_n'(w, v, j) - S_n^{*\prime}(w, v, j)|
      >
      \frac{h^{-D_w/4} \eta}{(n h^{D_W})^{3/4}}
      \biggl(\frac{h}{\delta_w}\biggr)^{\frac{D_w}{2} + D_v} \\
      &\qquad\qquad+
      \frac{h^{-D_w/(k+1)} \eta}{(n h^{D_W})^{k/(k+1)}}
      \biggl(\frac{h}{\delta_w}\biggr)^{\frac{k D_w}{2k+2} + \frac{k D_v}{k+1}}
      + \frac{C \delta_w}{h}
      \sqrt{\frac{\log 2n}{n h^{D_W}}}
      + \frac{C (\log 2n)^{1 + 1/r + 1/r'}}{n h^{D_W}}
    \Biggr\} \\
    &\qquad\leq
    \frac{C \log 2n}{\eta^2}
    + \frac{C (\log 2n)^{k/2}}{\eta^{k+1}}
    + \frac{C}{n^\kappa}.
  \end{align*}
  Substituting the value of $\delta_w$ makes the term
  $C (\log 2n)^{1 + 1/r + 1/r'} / (n h^{D_W})$ negligible.
  Replace $\eta$ with $C \eta \sqrt{\log 2n}$ to obtain
  \begin{align*}
    &\P\Biggl\{
      \sup_{(w, v) \in \cW \times \cV}
      \sqrt{n h^{D_W}}
      \bigl\| S_{n,h}^{*\prime}(w, v) - S_{n,h}'(w, v) \bigr\|_2
      \iftoggle{journal}{\\ &\qquad\qquad}{}
      > \frac{C \eta \sqrt{\log 2n} }
      {(n h^{D_W + (4-k)D_w})^{1/(10-2k+2D_w+4D_v)}}
    \Biggr\}
    \leq
    \frac{1}{\eta^2},
  \end{align*}
  as required.
  Here, the term $C/n^\kappa$ was removed as follows.
  If $\eta \leq n^2$, then $n^{-\kappa} \leq \eta^{-2}$.
  Alternatively, for $\eta > n^2$, Jensen's inequality and
  Lemma~\ref{lem:moment_tail_orlicz} give
  \begin{equation*}
    \E\Biggl\{
      \biggl(\sqrt{n h^{D_W}}
        \sup_{(w,v)\in\cW\times\cV}
        \|S_{n,h}^{*\prime}(w,v)-S_{n,h}'(w,v)\|_2
      \biggr)^8
    \Biggr\}
    \leq
    C(n h^{D_W})^4 h^{-8D_W}
    \leq C n^8,
  \end{equation*}
  as $n h^{D_W} \geq 1$.
  Since $n h^{D_W+(4-k)D_w} \leq n$ and
  $10-2k+2D_w+4D_v \geq 4$, Markov's inequality
  gives a probability bound of $C n^{10}/\eta^8 \leq C/\eta^2$.
  Increasing $C$ yields the result.

  It remains to handle the case where $n h^{D_W+(4-k)D_w} \leq 1$.
  Choosing an independent coupling, by the triangle inequality and
  Lemmas~\ref{lem:locpol_linear},
  \ref{lem:locpol_linear_bootstrap} and~\ref{lem:locpol_Scirc},
  taking $t = \sqrt{2\kappa\log 2n}$ gives
  \begin{equation*}
    \P\Biggl\{
      \sqrt{n h^{D_W}}
      \sup_{(w,v)\in\cW\times\cV}
      \|S_{n,h}^{*\prime}(w,v)-S_{n,h}'(w,v)\|_2
      > C\sqrt{\log 2n}
    \Biggr\}
    \leq \frac{C}{n^\kappa}.
  \end{equation*}
  For $1 \leq \eta \leq n^2$, this gives the stated bound.
  For $\eta > n^2$, use the eighth-moment argument above;
  if $\eta < 1$, the bound is trivial.

  Under condition~\ref{it:thm_locpol_bootstrap_multiplier},
  the same argument gives the stated bound
  for $S_{n,h}^\circ$,
  using condition~\ref{it:thm_ep_bootstrap_multiplier}
  part (ii) of Theorem~\ref{thm:empirical_process_bootstrap}.
\end{proof}

\begin{proof}[Proof of Theorem~\ref{thm:locpol_bootstrap}]
  Throughout the proof, $C > 0$ depends only on the stated quantities.
  Fix $\kappa = 4$.
  We begin with the nonparametric bootstrap
  of condition~\ref{it:thm_locpol_bootstrap_nonparametric}.
  By~\eqref{eq:locpol_decomp_coupling},
  \eqref{eq:locpol_decomp_errors},
  \eqref{eq:locpol_decomp_bias} and
  Lemmas~\ref{lem:locpol_Hhat},
  \ref{lem:locpol_linear} and~\ref{lem:locpol_bias},
  with probability at least
  $1 - 1/n^{\kappa} - 2 e^{-n h^{D_W}/C}$,
  \begin{align*}
    &\sup_{(w, v) \in \cW \times \cV}
    \bigl|
    \hat \theta_{n,h}(w, v) - \theta(w, v)
    - e_1^\T H_h(w)^{-1} S_{n,h}(w, v)
    \bigr| \\
    &\qquad\leq
    \sup_{(w, v) \in \cW \times \cV}
    \biggl\{
      \bigl\|\hat H_h(w)^{-1} - H_h(w)^{-1}\bigr\|_{2 \to 2}
      \bigl\|\tilde N_{n,h}(w, v)\bigr\|_2
      \iftoggle{journal}{\\ &\qquad\qquad\qquad\qquad\qquad}{}
      + \bigl| e_1^\T \hat H_h(w)^{-1} B_{n,h}(w, v) - \theta(w, v) \bigr|
    \biggr\} \\
    &\qquad\leq
    \frac{C \log 2 n}{n h^{D_W}} + C h^\zeta.
  \end{align*}
  Similarly, by~\eqref{eq:locpol_nonparametric_decomp_coupling},
  \eqref{eq:locpol_nonparametric_decomp_errors},
  \eqref{eq:locpol_nonparametric_decomp_bias}
  and Lemmas~\ref{lem:locpol_Hhat},
  \ref{lem:locpol_Hhatstar},
  \ref{lem:locpol_linear},
  \ref{lem:locpol_linear_bootstrap},
  \ref{lem:locpol_bias} and
  \ref{lem:locpol_bias_bootstrap},
  with probability at least
  $1 - 1/n^{\kappa} - 2 e^{-n h^{D_W}/C}$,
  \begin{align*}
    &\sup_{(w, v) \in \cW \times \cV}
    \bigl| \tilde \theta_{n,h}^*(w, v)
    - e_1^\T H_h(w)^{-1} S_{n,h}^*(w, v) \bigr| \\
    &\qquad\leq
    \sup_{(w, v) \in \cW \times \cV}
    \biggl\{
      \bigl\|\hat H_h^*(w)^{-1} - H_h(w)^{-1}\bigr\|_{2 \to 2}
      \bigl\|\tilde N_{n,h}^*(w, v)\bigr\|_2 \\
      &\qquad\qquad\qquad\qquad+
      \bigl\|\hat H_h(w)^{-1} - H_h(w)^{-1}\bigr\|_{2 \to 2}
      \bigl\|\tilde N_{n,h}(w, v)\bigr\|_2 \\
      &\qquad\qquad+
      \bigl|e_1^\T \hat H_h^*(w)^{-1} B_{n,h}^*(w, v) - \theta(w,v)\bigr|
      + \bigl|e_1^\T \hat H_h(w)^{-1} B_{n,h}(w, v) - \theta(w,v)\bigr|
    \biggr\} \\
    &\qquad\leq
    \frac{C \log 2 n}{n h^{D_W}} + C h^\zeta.
  \end{align*}
  Take $\bX_n'$, $\bX_n''$, $\bX_n^{*\prime}$,
  $S_{n,h}'$ and $S_{n,h}^{*\prime}$
  from Lemma~\ref{lem:locpol_bootstrap}.
  Compute $\hat\theta_{n,h}'$ from $\bX_n''$ and
  $\tilde\theta_{n,h}^{*\prime}$ from
  $(\bX_n',\bX_n^{*\prime})$.
  By the preceding bounds and Lemma~\ref{lem:locpol_H},
  with probability at least
  $1 - 1/n^{\kappa} - 2 e^{-n h^{D_W}/C}$,
  \begin{align*}
    &\sup_{(w, v) \in \cW \times \cV}
    \bigl|
    \hat \theta_{n,h}'(w, v) - \theta(w, v)
    - \tilde \theta_{n,h}^{*\prime}(w, v)
    \bigr| \\
    &\qquad\leq
    \sup_{(w, v) \in \cW \times \cV}
    \Bigl\{
      \bigl\|H_h(w)^{-1}\bigr\|_{2 \to 2}
      \bigl\| S_{n,h}^{*\prime}(w, v) - S_{n,h}'(w, v) \bigr\|_2
    \Bigr\}
    + \frac{C \log 2 n}{n h^{D_W}} + C h^\zeta \\
    &\qquad\leq
    C \sup_{(w, v) \in \cW \times \cV}
    \bigl\| S_{n,h}^{*\prime}(w, v) - S_{n,h}'(w, v) \bigr\|_2
    + \frac{C \log 2 n}{n h^{D_W}} + C h^\zeta.
  \end{align*}
  By Lemma~\ref{lem:locpol_bootstrap},
  \begin{align*}
    &\P\Biggl\{
      \sup_{(w, v) \in \cW \times \cV}
      \sqrt{n h^{D_W}}
      \bigl|
      \hat \theta_{n,h}'(w, v) - \theta(w, v)
      - \tilde \theta_{n,h}^{*\prime}(w, v)
      \bigr| \\
      &\qquad> \frac{C \eta \sqrt{\log 2n} }
      {(n h^{D_W + (4-k)D_w})^{1/(10-2k+2D_w+4D_v)}}
      + \frac{C \log 2 n}{\sqrt{n h^{D_W}}} + C \sqrt{n h^{D_W}} h^\zeta
    \Biggr\}
    \iftoggle{journal}{\\ &\quad}{}
    \leq
    \frac{1}{n^\kappa} + 2 e^{-n h^{D_W}/C} + \frac{1}{\eta^2}.
  \end{align*}
  We may assume $\eta \geq 1$;
  since $n h^{D_W + (4-k)D_w} \leq n h^{D_W}$,
  $10 - 2k + 2D_w + 4D_v \geq 4$ and $n h^{D_W} \geq n^{1/\Clp}$,
  \begin{equation*}
    \frac{\log 2n}{\sqrt{n h^{D_W}}}
    \leq
    \frac{C \eta \sqrt{\log 2n}}
    {(n h^{D_W + (4-k)D_w})^{1/(10-2k+2D_w+4D_v)}},
  \end{equation*}
  and so the $\log 2 n / \sqrt{n h^{D_W}}$ term above may be removed
  by increasing $C$.
  To remove the $1/n^\kappa$ term we argue as follows.
  If $\eta \leq n^2$, we have $n^{-\kappa} \leq \eta^{-2}$.
  For $\eta > n^2$, restrict to the event of
  Lemmas~\ref{lem:locpol_Hhat} and~\ref{lem:locpol_Hhatstar}
  on which
  \begin{align*}
    \sup_{w \in [0, 1]^{D_W}}
    \bigl\|\hat H_h'(w)\bigr\|_{2 \to 2}
    &\leq C,
    &\sup_{w \in [0, 1]^{D_W}}
    \bigl\|\hat H_h'(w)^{-1}\bigr\|_{2 \to 2}
    &\leq C, \\
    \sup_{w \in [0, 1]^{D_W}}
    \bigl\|\hat H_h''(w)\bigr\|_{2 \to 2}
    &\leq C,
    &\sup_{w \in [0, 1]^{D_W}}
    \bigl\|\hat H_h''(w)^{-1}\bigr\|_{2 \to 2}
    &\leq C, \\
    \sup_{w \in [0, 1]^{D_W}}
    \bigl\|\hat H_h^{*\prime}(w)\bigr\|_{2 \to 2}
    &\leq C,
    &\sup_{w \in [0, 1]^{D_W}}
    \bigl\|\hat H_h^{*\prime}(w)^{-1}\bigr\|_{2 \to 2}
    &\leq C,
  \end{align*}
  where $\hat H_h'$ is computed using $\bX_n'$,
  $\hat H_h''$ is computed using $\bX_n''$,
  and $\hat H_h^{*\prime}$ uses $\bX_n^{*\prime}$;
  this event has probability at least $1 - 2 e^{- n h^{D_W} / C}$,
  increasing $C$ if necessary.
  On this event, the argument given in
  the proof of Lemma~\ref{lem:locpol_bootstrap} gives
  a probability bound of $C n^{10}/\eta^8 \leq C/\eta^2$.

  Now we turn to the multiplier bootstrap
  of condition~\ref{it:thm_locpol_bootstrap_multiplier}.
  By~\eqref{eq:locpol_multiplier_decomp_coupling},
  \eqref{eq:locpol_multiplier_decomp_errors}
  and Lemmas~\ref{lem:locpol_Hhat},
  \ref{lem:locpol_Scirc} and \ref{lem:locpol_R},
  with probability at least
  $1 - 1/n^{\kappa} - 2 e^{-n h^{D_W}/C}$,
  \begin{align*}
    &\sup_{(w, v) \in \cW \times \cV}
    \bigl|
    \tilde \theta_{n,h}^\circ(w, v)
    - e_1^\T H_h(w)^{-1} S_{n,h}^\circ(w, v)
    \bigr| \\
    &\qquad\leq
    \sup_{(w, v) \in \cW \times \cV}
    \biggl\{
      \bigl\|\hat H_h(w)^{-1} - H_h(w)^{-1}\bigr\|_{2 \to 2}
      \bigl\|S_{n,h}^\circ(w, v)\bigr\|_2
      \iftoggle{journal}{\\ &\qquad\qquad\qquad\qquad\qquad}{}
      + \bigl\|\hat H_h(w)^{-1}\bigr\|_{2 \to 2}
      \bigl\|R_{n,h}^\circ(w, v)\bigr\|_2
    \biggr\} \\
    &\qquad\leq
    \frac{C \log 2n}{n h^{D_W}}
    + C h^\zeta.
  \end{align*}
  The rest of the proof proceeds as for the nonparametric bootstrap,
  applying Lemma~\ref{lem:locpol_bootstrap}.
\end{proof}

\begin{lemma}[Linear term variance estimator]%
  \label{lem:Omega_hat}
  Grant the setup of Section~\ref{sec:locpol}
  and Assumption~\ref{ass:locpol}, and suppose that
  $\Var\bigl\{g(V_i, v) \mid W_i\bigr\} \geq 1 /\Clp$
  almost surely for each $i \in [n]$ and $v \in \cV$.
  Define $\Omega_h(w, v, w', v') \coloneqq
  \Cov\bigl\{S_{n,h}(w, v), S_{n,h}(w', v')\bigr\}$
  for $(w, v), (w', v') \in \cW \times \cV$,
  and take $\kappa \geq 1$.
  Then $\Omega_h(w, v, w, v)$
  is invertible for all $(w, v) \in \cW \times \cV$.
  Further, there exists $C > 0$ depending only on
  $\Clp$, $D_W$, $D_V$, $r$, $q$, $\kappa$ and $\zeta$ such that
  \begin{align*}
    \sup_{(w, v), (w', v') \in \cW \times \cV}
    \bigl\|\Omega_h(w, v, w', v')\bigr\|_{2 \to 2}
    &\leq
    C n^{-1} h^{-D_W}, \\
    \sup_{(w, v) \in \cW \times \cV}
    \bigl\|\Omega_h(w, v, w, v)^{-1}\bigr\|_{2 \to 2}
    &\leq
    C n h^{D_W}.
  \end{align*}
  Moreover, defining
  $\hat\Omega_h(w, v, w', v') \coloneqq
  \hat Q_h(w)^\T \hat K_h(w) \hat D_h(v, v') \hat K_h(w') \hat Q_h(w') / n^2$,
  for all $t > 0$, with probability at least
  $1 - 1/n^{\kappa} - 2 e^{-n h^{D_W} / C} - e^{-t^2/2}$,
  \begin{equation*}
    \sup_{(w, v), (w', v') \in \cW \times \cV}
    \bigl\|
    \hat\Omega_h(w, v, w', v') - \Omega_h(w, v, w', v')
    \bigr\|_{2 \to 2}
    \leq
    \frac{C}{n h^{D_W}}
    \Biggl\{
      \sqrt{\frac{\log 2n + t^2}{n h^{D_W}}}
      + h^\zeta
    \Biggr\}.
  \end{equation*}
  Finally, with probability at least
  $1 - 1/n^{\kappa} - 2 e^{-n h^{D_W}/C}$, we have that
  $\hat\Omega_h(w, v, w, v)$ is invertible for all
  $(w, v) \in \cW \times \cV$ and both
  \begin{align*}
    \sup_{(w, v), (w', v') \in \cW \times \cV}
    \bigl\|\hat\Omega_h(w, v, w', v')\bigr\|_{2 \to 2}
    &\leq C n^{-1} h^{-D_W}, \\
    \sup_{(w, v) \in \cW \times \cV}
    \bigl\|\hat\Omega_h(w, v, w, v)^{-1}\bigr\|_{2 \to 2}
    &\leq C n h^{D_W}.
  \end{align*}
\end{lemma}

\begin{proof}[Proof of Lemma~\ref{lem:Omega_hat}]
  Throughout the proof,
  $C > 0$ depends only on the stated quantities.
  Since $S_{n,h}(w, v)
  = \hat Q_h(w)^\T \hat K_h(w) \tilde g(v) / n$,
  using the lower bound on the conditional variance of $g(V_i, v)$
  and the argument in the proof of Lemma~\ref{lem:locpol_H}
  with $K(w)$ replaced by $K(w)^2$, we obtain that
  \begin{align*}
    \Omega_h(w,v,w,v)
    &=
    \Var\bigl\{S_{n,h}(w, v)\bigr\}
    \iftoggle{journal}{\\ &}{}
    =
    \frac{1}{n^2}
    \sum_{i=1}^{n}
    \E\bigl\{
      Q_h(W_i - w)
      Q_h(W_i - w)^\T
      K_h(W_i - w)^2
      \tilde g(v)_i^2
    \bigr\} \\
    &\succeq
    \frac{1}{\Clp n^2}
    \sum_{i=1}^{n}
    \E\bigl\{
      Q_h(W_i - w)
      Q_h(W_i - w)^\T
      K_h(W_i - w)^2
    \bigr\}
    \succeq
    \frac{I_{D_q}}{C n h^{D_W}},
  \end{align*}
  where $\succeq$ is the positive semidefinite partial order.
  Similarly, $\|\Omega_h(w,v,w',v')\|_{2 \to 2} \leq C / (n h^{D_W})$
  by Lemma~\ref{lem:moment_tail_orlicz}.
  Next, define
  \begin{align*}
    &\tilde\Omega_h(w, v, w', v')
    \iftoggle{journal}{\\ &\quad}{}
    \coloneqq
    \frac{1}{n^2}
    \sum_{i=1}^{n}
    Q_h(W_i - w)
    Q_h(W_i - w')^\T
    K_h(W_i - w)
    K_h(W_i - w')
    \tilde g(v)_i
    \tilde g(v')_i,
  \end{align*}
  so that we have
  \begin{align*}
    &\hat\Omega_h(w, v, w', v')
    - \tilde\Omega_h(w, v, w', v') \\
    &\quad=
    \frac{1}{n^2}
    \sum_{i=1}^{n}
    Q_h(W_i - w)
    Q_h(W_i - w')^\T
    K_h(W_i - w)
    K_h(W_i - w') \\
    &\qquad\qquad\times
    \Bigl\{
      \tilde g(v)_i \bigl(\theta(W_i, v') - \hat\theta_{n,h}(W_i, v')\bigr)
      + \tilde g(v')_i \bigl(\theta(W_i, v) - \hat\theta_{n,h}(W_i, v)\bigr) \\
      &\qquad\qquad\qquad\quad+
      \bigl(\theta(W_i, v) - \hat\theta_{n,h}(W_i, v)\bigr)
      \bigl(\theta(W_i, v') - \hat\theta_{n,h}(W_i, v')\bigr)
    \Bigr\}.
  \end{align*}
  As in the proof of Lemma~\ref{lem:locpol_linear},
  it suffices to consider only
  $t\leq C\sqrt{\log 2n}$;
  larger values are handled by replacing $t$ with
  $\sqrt{2(\kappa+1)\log 2n}$ and $\kappa$ with $\kappa+1$.
  By Lemma~\ref{lem:locpol_convergence} applied on $[0, 1]^{D_W}$,
  with probability at least
  $1 - 1/n^{\kappa} - 2 e^{-n h^{D_W}/C} - e^{-t^2/2}$,
  \begin{align}
    \nonumber
    &\sup_{(w, v), (w', v') \in \cW \times \cV}
    \bigl\|\hat\Omega_h(w, v, w', v')
    - \tilde\Omega_h(w, v, w', v')\bigr\|_{2 \to 2} \\
    &\quad\leq
    \nonumber
    \sup_{(w, v), (w', v') \in \cW \times \cV}
    \frac{C}{n^2}
    \sum_{i=1}^{n}
    \|Q_h(W_i - w)\|_2
    \|Q_h(W_i - w')\|_2
    K_h(W_i - w)
    K_h(W_i - w') \\
    &\qquad\times
    \label{eq:Omega_hat_sup}
    \Biggl\{
      \sqrt{\frac{\log 2n + t^2}{n h^{D_W}}}
      + h^\zeta
    \Biggr\}
    \Biggl\{
      |\tilde g(v)_i| + |\tilde g(v')_i|
      + \sqrt{\frac{\log 2n + t^2}{n h^{D_W}}}
      + h^\zeta
    \Biggr\}.
  \end{align}
  Apply Theorem~\ref{thm:empirical_tail} to the class
  \begin{equation*}
    \bigl\{
      (W,V)\mapsto K_h(W-w)
      \bigl(1+|g(V,v)-\theta(W,v)|\bigr):
      (w,v)\in\cW\times\cV
    \bigr\},
  \end{equation*}
  with envelope $Ch^{-D_W}(|V|+\Clp)$ in the regression case
  and $Ch^{-D_W}$ in the indicator case.
  The entropy argument in the proof of
  Lemma~\ref{lem:locpol_linear} gives
  \begin{equation*}
    \sup_{(w,v)\in\cW\times\cV}
    \frac{1}{n}\sum_{i=1}^n K_h(W_i-w)
    \bigl(1+|\tilde g(v)_i|\bigr)
    \leq C,
  \end{equation*}
  with probability at least $1-n^{-\kappa}-e^{-t^2/2}$.
  Combining this with \eqref{eq:Omega_hat_sup},
  replacing $\kappa$ with $\kappa+1$ and
  $t$ with $t+2$ in both estimates,
  with probability at least
  $1 - 1/n^{\kappa} - 2 e^{-n h^{D_W}/C} - e^{-t^2/2}$,
  we have
  \begin{equation}
    \label{eq:Omega_hat_convergence}
    \sup_{(w, v), (w', v') \in \cW \times \cV}
    \bigl\|\hat\Omega_h(w, v, w', v')
    - \tilde\Omega_h(w, v, w', v')\bigr\|_{2 \to 2}
    \leq
    \frac{C}{n h^{D_W}}
    \Biggl\{
      \sqrt{\frac{\log 2n + t^2}{n h^{D_W}}}
      + h^\zeta
    \Biggr\}.
  \end{equation}
  Since
  $\E\bigl\{\tilde\Omega_h(w,v,w',v')\bigr\} = \Omega_h(w,v,w',v')$,
  we apply Theorem~\ref{thm:empirical_tail} with $r'\coloneqq r/2$
  to the class
  \begin{align*}
    \Bigl\{
      &(W,V)\mapsto
      e_j^\T Q_h(W-w)Q_h(W-w')^\T e_k
      K_h(W-w)K_h(W-w') \\
      \iftoggle{journal}{&\qquad}{&\qquad\qquad\qquad}
      \times
      \bigl(g(V,v)-\theta(W,v)\bigr)
      \bigl(g(V,v')-\theta(W,v')\bigr):
      (w,v),(w',v')\in\cW\times\cV,\ j,k\in[D_q]
    \Bigr\},
  \end{align*}
  with envelope $Ch^{-2D_W}(|V|+1)^2$ in the regression
  case and $Ch^{-2D_W}$ in the indicator case, each
  having $\omega_{r/2}$-Orlicz quasi-norm bounded by $Ch^{-2D_W}$.
  Take $\sigma=Ch^{-3D_W/2}$ and $\nu=h^{D_W/2}/C$.
  The covering number argument from the proof of
  Lemma~\ref{lem:locpol_linear} gives
  $\cJ(\nu)/\nu\leq C\sqrt{\log 2n}$.
  Thus, with probability at least $1-e^{-t^2/2}-n^{-\kappa}$,
  we have
  \begin{equation*}
    \sup_{(w, v), (w', v') \in \cW \times \cV}
    \bigl\|\tilde\Omega_h(w, v, w', v')
    - \Omega_h(w, v, w', v')\bigr\|_{2 \to 2}
    \leq
    \frac{C}{n h^{D_W}}
    \sqrt{\frac{\log 2n + t^2}{n h^{D_W}}}.
  \end{equation*}
  Combining this with \eqref{eq:Omega_hat_convergence},
  using $\kappa + 1$ and $t + 2$ in place of $\kappa$ and $t$,
  taking a union bound, and increasing $C$ gives the desired
  convergence of $\hat\Omega_h(w,v,w',v')$.

  Finally, we show the bounds on the singular values.
  Set $t^2 \coloneqq 2 n h^{D_W}/C'$, with $C'$ sufficiently large.
  Since $n h^{D_W} \geq n^{1/\Clp}$ and $h \leq n^{-1/\Clp}$,
  we have for $n \geq C$ that
  \begin{equation*}
    \sup_{(w, v) \in \cW \times \cV}
    \bigl\|
    \hat \Omega_h(w, v, w, v) - \Omega_h(w, v, w, v)
    \bigr\|_{2 \to 2}
    \leq
    \frac{1}{2}
    \inf_{(w, v) \in \cW \times \cV}
    \lambda_{\min}\bigl\{\Omega_h(w, v, w, v)\bigr\},
  \end{equation*}
  with probability at least
  $1 - 1/n^\kappa - 2e^{-n h^{D_W}/C} - e^{-n h^{D_W}/C'}$,
  where $\lambda_{\min}$ denotes the smallest eigenvalue.
  Weyl's inequality and the triangle inequality give the
  required bounds, increasing $C$ if necessary.
  For $n \leq C$, the bound is trivial as $h \leq 1$.
\end{proof}

\begin{lemma}[Local polynomial variance estimator]%
  \label{lem:Sigma_hat}
  Grant the setup of Section~\ref{sec:locpol}
  and Assumption~\ref{ass:locpol}, and suppose that
  $\Var\bigl\{g(V_i, v) \mid W_i\bigr\} \geq 1 /\Clp$
  almost surely for each $i \in [n]$ and $v \in \cV$.
  Take $\kappa \geq 1$.
  There exists $C > 0$ depending only on
  $\Clp$, $D_W$, $D_V$, $r$, $q$, $\kappa$ and $\zeta$
  such that
  $\sup_{(w,v), (w',v') \in \cW \times \cV} |\Sigma_h(w,v,w',v')|
  \leq C n^{-1} h^{-D_W}$
  and $\sup_{(w,v) \in \cW \times \cV} |\Sigma_h(w,v,w,v)^{-1}|
  \leq C n h^{D_W}$.
  For all $t > 0$, with probability at least
  $1 - 1/n^{\kappa} - 2 e^{-n h^{D_W} / C} - e^{-t^2/2}$, we have
  \begin{equation*}
    \sup_{(w, v), (w', v') \in \cW \times \cV}
    \bigl|
    \hat\Sigma_{n,h}(w, v, w', v')
    - \Sigma_h(w, v, w', v')
    \bigr|
    \leq
    \frac{C}{n h^{D_W}}
    \Biggl\{
      \sqrt{\frac{\log 2n + t^2}{n h^{D_W}}}
      + h^\zeta
    \Biggr\}.
  \end{equation*}
  With probability at least
  $1 - 1/n^{\kappa} - 2 e^{-n h^{D_W}/C}$, we have
  \begin{align*}
    \sup_{(w,v), (w',v') \in \cW \times \cV}
    |\hat\Sigma_{n,h}(w,v,w',v')|
    &\leq
    C n^{-1} h^{-D_W}, \\
    \sup_{(w,v) \in \cW \times \cV}
    |\hat \Sigma_{n,h}(w,v,w,v)^{-1}|
    &\leq
    C n h^{D_W}.
  \end{align*}
\end{lemma}

\begin{proof}[Proof of Lemma~\ref{lem:Sigma_hat}]
  Throughout the proof,
  $C > 0$ depends only on the stated quantities.
  Recall the definitions of
  $\Omega_h(w, v, w', v')$
  and $\hat\Omega_h(w, v, w', v')$
  from Lemma~\ref{lem:Omega_hat}.
  The upper and lower bounds for $\Sigma_h(w, v, w', v')$
  follow from Lemmas~\ref{lem:locpol_H}
  and~\ref{lem:Omega_hat}, since
  \begin{equation*}
    \Sigma_h(w, v, w', v')
    =
    e_1^\T H_h(w)^{-1}
    \Omega_h(w, v, w', v')
    H_h(w')^{-1} e_1.
  \end{equation*}
  Next, observe that
  \begin{align*}
    &\hat\Sigma_{n,h}(w,v,w',v') - \Sigma_h(w,v,w',v') \\
    &\quad=
    e_1^\T \hat H_h(w)^{-1}
    \hat\Omega_h(w, v, w', v')
    \hat H_h(w')^{-1} e_1
    - e_1^\T H_h(w)^{-1}
    \Omega_h(w, v, w', v')
    H_h(w')^{-1} e_1 \\
    &\quad=
    e_1^\T
    \bigl\{ \hat H_h(w)^{-1} -  H_h(w)^{-1} \bigr\}
    \hat\Omega_h(w, v, w', v')
    \hat H_h(w')^{-1} e_1 \\
    &\qquad+
    e_1^\T H_h(w)^{-1}
    \bigl\{ \hat\Omega_h(w, v, w', v') - \Omega_h(w, v, w', v') \bigr\}
    \hat H_h(w')^{-1} e_1 \\
    &\qquad+
    e_1^\T
    H_h(w)^{-1}
    \Omega_h(w, v, w', v')
    \bigl\{ \hat H_h(w')^{-1} -  H_h(w')^{-1} \bigr\}
    e_1.
  \end{align*}
  Applying Lemmas~\ref{lem:locpol_H},
  \ref{lem:locpol_Hhat} and \ref{lem:Omega_hat}
  with $\kappa + 1$ and $t + 2$ in place of $\kappa$ and $t$,
  taking a union bound, and increasing $C$ gives the required
  bound for $\hat\Sigma_{n,h}(w,v,w',v') - \Sigma_h(w,v,w',v')$.
  The remaining results follow by Weyl's inequality,
  as in the proof of Lemma~\ref{lem:Omega_hat}.
\end{proof}

\begin{lemma}[Uniform convergence of bootstrapped local polynomial estimators]%
  \label{lem:locpol_bootstrap_convergence}
  Grant the setup of Theorem~\ref{thm:locpol_bootstrap}
  and take $\kappa \geq 1$.
  Under either
  condition~\ref{it:thm_locpol_bootstrap_nonparametric}
  or condition~\ref{it:thm_locpol_bootstrap_multiplier},
  there exists $C > 0$ depending only on
  $\Clp$, $D_W$, $D_V$, $q$, $r$, $r'$ and $\kappa$
  such that for all $t > 0$,
  with probability at least
  $1 - 1/n^{\kappa} - 2 e^{-n h^{D_W}/C} - e^{-t^2/2}$,
  \begin{equation*}
    \sup_{(w, v) \in \cW \times \cV}
    \bigl| \tilde\theta_{n,h}^*(w, v) \bigr|
    \leq
    C \sqrt{\frac{1 + D_w \log 2n + t^2}{n h^{D_W}}}
    + C h^\zeta.
  \end{equation*}
\end{lemma}

\begin{proof}[Proof of Lemma~\ref{lem:locpol_bootstrap_convergence}]
  By the proof of Theorem~\ref{thm:locpol_bootstrap} and
  Lemmas~\ref{lem:locpol_H}, \ref{lem:locpol_linear} and
  \ref{lem:locpol_linear_bootstrap},
  with probability at least
  $1 - 1/n^{\kappa} - 2 e^{-n h^{D_W}/C} - e^{-t^2/2}$,
  \begin{align*}
    \sup_{(w, v) \in \cW \times \cV}
    \bigl|\tilde \theta_{n,h}^*(w, v)\bigr|
    &\leq
    \sup_{(w, v) \in \cW \times \cV}
    \bigl|e_1^\T H_h(w)^{-1} S_{n,h}^*(w, v)\bigr|
    + \frac{C \log 2n}{n h^{D_W}}
    + C h^\zeta \\
    &\leq
    C \sqrt{\frac{1 + D_w \log 2n + t^2}{n h^{D_W}}}
    + C h^\zeta.
  \end{align*}
  Similarly, by Lemmas~\ref{lem:locpol_H} and \ref{lem:locpol_Scirc},
  with probability at least
  $1 - 1/n^{\kappa} - 2 e^{-n h^{D_W}/C} - e^{-t^2/2}$,
  \begin{align*}
    \sup_{(w, v) \in \cW \times \cV}
    \bigl|
    \tilde \theta_{n,h}^\circ(w, v)
    \bigr|
    &\leq
    \sup_{(w, v) \in \cW \times \cV}
    \bigl|
    e_1^\T H_h(w)^{-1} S_{n,h}^\circ(w, v)
    \bigr|
    + \frac{C \log 2n}{n h^{D_W}}
    + C h^\zeta \\
    &\leq
    C \sqrt{\frac{1 + D_w \log 2n + t^2}{n h^{D_W}}}
    + C h^\zeta,
  \end{align*}
  as required.
\end{proof}

\begin{lemma}[Anti-concentration construction]%
  \label{lem:anti-concentration_construction}
  Let $d \in \N$ and $X$ be a random vector
  taking values in $\R^d$.
  Suppose that $\E(\|X\|_2^2) < \infty$ and
  write $\mu \coloneqq \E(X)$.
  Take $\zeta$ to be a real-valued random variable independent
  of $X$ satisfying
  $\zeta = 1 + \sqrt{6}/2$ with probability $1/2 - 1/\sqrt{6}$
  and $\zeta = 1 - \sqrt{6}/2$ with probability $1/2 + 1/\sqrt{6}$.
  Let $\Sigma \coloneqq \Var(X) \in \R^{d \times d}$
  and take $Z \sim N(0_d, \Sigma/2)$
  independent of $(X, \zeta)$.
  Define $\tilde X \coloneqq \mu + \zeta (X - \mu) + Z$.
  Then $\E(\tilde X) = \mu$ and
  $\Var(\tilde X) = \Sigma$,
  so $\E(\tilde X^\alpha) = \E(X^\alpha)$
  for all $\alpha \in \N_0^d$ with $|\alpha| \in \{1,2\}$.
  If additionally $\E(\|X\|_2^3) < \infty$, this identity
  also holds when $|\alpha| = 3$.
  Further, $\tilde X \sim N\bigl(\mu + \zeta(X - \mu), \Sigma/2\bigr)$
  conditional on $(\zeta, X)$.
  Finally, for $p \in [1, \infty]$, we have
  $\|\tilde X\|_p \leq \|\mu\|_p
  + \bigl(1 + \sqrt{6}/2\bigr) \|X - \mu\|_p + \|Z\|_p$.
\end{lemma}

\begin{proof}[Proof of Lemma~\ref{lem:anti-concentration_construction}]
  By independence and since $\E(\zeta) = 0$ and $\E(\zeta^2) = 1/2$,
  we have
  \begin{align*}
    \E(\tilde X)
    &=
    \mu + \E(\zeta) \E(X - \mu) + \E(Z) = \mu, \\
    \Var(\tilde X)
    &=
    \E(\zeta^2) \Var(X)
    + \Var(Z)
    = \Sigma/2 + \Sigma/2 = \Sigma.
  \end{align*}
  If $\E(\|X\|_2^3)<\infty$, then for
  $\alpha\in\N_0^d$ with $|\alpha|=3$,
  the multi-binomial theorem and independence give
  \begin{equation*}
    \E\bigl\{(\tilde X-\mu)^\alpha\bigr\}
    =
    \sum_{\nu\leq\alpha}\binom{\alpha}{\nu}
    \E(\zeta^{|\nu|})
    \E\bigl\{(X-\mu)^\nu\bigr\}
    \E(Z^{\alpha-\nu})
    = \E\bigl\{(X-\mu)^\alpha\bigr\},
  \end{equation*}
  by symmetry of the Gaussian distribution, and as
  $\E(\zeta)=0$ and $\E(\zeta^3)=1$.
  Combining this with the mean and covariance identities,
  another application of the multi-binomial theorem gives
  \begin{equation*}
    \E(\tilde X^\alpha)
    =
    \sum_{\nu\leq\alpha}\binom{\alpha}{\nu}
    \mu^{\alpha-\nu}\E\bigl\{(\tilde X-\mu)^\nu\bigr\}
    =
    \E(X^\alpha),
  \end{equation*}
  for $|\alpha|\in\{1,2\}$, and also for $|\alpha|=3$
  when $\E(\|X\|_2^3)<\infty$.
  Conditional Gaussianity follows by definition of $\tilde X$,
  and the final claim is by the triangle inequality.
\end{proof}

\begin{lemma}[Bootstrap inference for local polynomial estimators]%
  \label{lem:locpol_bootstrap_ks}
  Grant the setup of Theorem~\ref{thm:locpol_bootstrap}
  and assume that $h^\zeta \sqrt{n h^{D_W} \log 2n} \to 0$ and
  $\Var\bigl\{g(V_i, v) \mid W_i\bigr\} \geq 1 /\Clp$
  almost surely for each $i \in [n]$ and $v \in \cV$.
  Suppose that
  $n h^{D_W + (4-k)D_w} / (\log 2n)^{10-2k+2D_w+4D_v} \to \infty$.
  Then under condition~\ref{it:thm_locpol_bootstrap_nonparametric},
  \begin{align*}
    &\E\Biggl\{
      \sup_{\xi \in \R}
      \,
      \Biggl|
      \P\Biggl(
        \sup_{(w, v) \in \cW \times \cV}
        \,
        \Biggl|
        \frac{\hat\theta_{n,h}(w,v)-\theta(w,v)}
        {\sqrt{\smash[b]{\hat \Sigma_{n,h}(w,v,w,v)}}}
        \Biggr| \leq \xi
      \Biggr) \\
      &\qquad\qquad\qquad-
      \P^*\Biggl(
        \sup_{(w, v) \in \cW \times \cV}
        \,
        \Biggl|
        \frac{\tilde\theta^*_{n,h}(w, v)}
        {\sqrt{\smash[b]{\hat \Sigma_{n,h}(w,v,w,v)}}}
        \Biggr| \leq \xi
      \Biggr)
      \Biggr|
    \Biggr\}
    \to 0.
  \end{align*}
  The same result holds under
  condition~\ref{it:thm_locpol_bootstrap_multiplier}.
  We interpret both ratios as zero when
  $\hat\Sigma_{n,h}(w,v,w,v) = 0$.
\end{lemma}

\begin{proof}[Proof of Lemma~\ref{lem:locpol_bootstrap_ks}]
  We combine ideas from the proofs of
  Theorems~\ref{thm:bootstrap_ks}
  and~\ref{thm:locpol_gaussian_ks},
  beginning by constructing an auxiliary coupling to
  a conditionally Gaussian process.
  To this end,
  for $i \in [n]$, let $\zeta_i$ be independent
  variables following the two-point distribution
  defined in Lemma~\ref{lem:anti-concentration_construction},
  independent of $\bX_n$.
  For $i \in [n]$, let $T_i$ be independent
  zero-mean Gaussian processes
  indexed on $\cW \times \cV \times [D_q]$,
  independent of both $\bX_n$ and $(\zeta_i: i \in [n])$, with
  \begin{align*}
    &\Cov\bigl\{
      T_i(w, v, j),
      T_i(w', v', j')
    \bigr\} \\
    &\qquad=
    \frac{1}{2}
    \Cov\Bigl\{
      e_j^\T Q_h(W_i - w) K_h(W_i - w) \tilde g(v)_i, \,
      \iftoggle{journal}{\\ &\qquad\qquad\qquad\qquad}{}
      e_{j'}^\T Q_h(W_i - w') K_h(W_i - w') \tilde g(v')_i
    \Bigr\}.
  \end{align*}
  Define $\tilde S_{n,h}(w, v) \coloneqq
  \bigl(\tilde S_{n}(w, v, j): j \in [D_q]\bigr)$, where
  \begin{equation*}
    \tilde S_{n}(w, v, j)
    \coloneqq
    \frac{1}{n}
    \sum_{i=1}^n
    \Bigl\{
      e_j^\T
      Q_h(W_i - w) K_h(W_i - w) \tilde g(v)_i \zeta_i
      + T_i(w, v, j)
    \Bigr\}.
  \end{equation*}

  Use the discretization $\cM_\delta$,
  and associated notation, from the proof of
  Lemma~\ref{lem:locpol_bootstrap}.
  Apply Theorem~\ref{thm:yurinskii} with
  $p=\infty$ and $\gamma=1$ to the random vectors
  $S_{n,h}(\cM_\delta)$ and $\tilde S_{n,h}(\cM_\delta)$.
  By Lemma~\ref{lem:anti-concentration_construction},
  their summands have exactly matching moments up to order three,
  so $\Delta=\Gamma_\eta=0$.
  For $i \in [n]$ and $a \in \{2, k\}$,
  boundedness of $\zeta_i$, independence
  and Lemma~\ref{lem:gaussian_norm} give
  \begin{equation*}
    \beta_{a,1,i}
    \leq
    \frac{C}{(n h^{D_W})^{a+1}}
    \biggl(\frac{h}{\delta_w}\biggr)^{a D_w/2}
    \biggl(\frac{1}{\delta_v}\biggr)^{a D_v/2}
    h^{D_W - D_w},
  \end{equation*}
  as in the proof of Lemma~\ref{lem:locpol_gaussian}.
  Extend the coupling to $\cM$ using the
  Vorob'ev--Berkes--Philipp theorem.

  Define $f_{w,v,j}$ and $\pi$ as in
  the proof of Lemma~\ref{lem:locpol_bootstrap},
  and let $(w',v',j') \coloneqq \pi(w, v, j)$.
  Since $\zeta_i$ is bounded, the envelope and entropy bounds
  from that proof remain valid.
  By Theorem~\ref{thm:empirical_tail} applied to $(W_i,V_i,\zeta_i)$,
  with probability at least $1-e^{-t^2/2}-n^{-\kappa}$, we therefore have
  \begin{align*}
    &\sup_{(w,v,j)\in\cM}
    \biggl|\sum_{i=1}^n\zeta_i
    \bigl\{f_{w,v,j}(W_i,V_i)-f_{(w',v',j')}(W_i,V_i)\bigr\}\biggr| \\
    &\qquad\leq
    C\sqrt{\frac{\log 2n+t^2}{n h^{D_W}}}
    \biggl(\frac{\delta_w}{h}+\sqrt{\delta_v}\biggr)
    +\frac{C(\log 2n)^{1+1/r}}{n h^{D_W}}.
  \end{align*}
  By Dudley's theorem \citep[Corollary~2.2.8]{van1996weak}
  and Borell's inequality \citep[Proposition~A.2.1]{van1996weak},
  as in the proof of Lemma~\ref{lem:locpol_gaussian},
  with probability at least $1-e^{-t^2/2}$,
  \begin{equation*}
    \sup_{(w,v,j)\in\cM}
    \biggl|\frac1n\sum_{i=1}^n
    \bigl\{T_i(w,v,j)-T_i(w',v',j')\bigr\}\biggr|
    \leq
    C\sqrt{\frac{\log 2n+t^2}{n h^{D_W}}}
    \biggl(\frac{\delta_w}{h}+\sqrt{\delta_v}\biggr).
  \end{equation*}
  The fluctuations of $S_{n,h}$ are bounded as in the proof of
  Lemma~\ref{lem:locpol_bootstrap}.
  Choosing $\delta_w$, $\delta_v$, $\kappa$, $t$ and $\eta$
  as in the proof of Lemma~\ref{lem:locpol_gaussian},
  we deduce that there exist processes
  $S_{n,h}'$ and $\tilde S_{n,h}'$, where
  $S_{n,h}' \equallaw S_{n,h}$ and
  $\tilde S_{n,h}' \equallaw \tilde S_{n,h}$, with
  \begin{align*}
    &\P\Biggl\{
      \sup_{(w, v) \in \cW \times \cV}
      \!
      \sqrt{n h^{D_W}}
      \bigl\| S'_{n,h}(w,v) - \tilde S'_{n,h}(w, v) \bigr\|_2
      \iftoggle{journal}{\\ &\qquad\qquad}{}
      >
      \frac{C \eta \sqrt{\log 2n}}
      {\bigl\{(n h^{D_W})^{k-1} h^{2D_w}\bigr\}^{1/(2k+2 + k D_w + 2k D_v)}}
    \Biggr\}
    \leq
    \frac{1}{\eta^{k+1}}.
  \end{align*}
  As in the proof of Theorem~\ref{thm:locpol_gaussian},
  by the Vorob'ev--Berkes--Philipp theorem we may extend
  the probability space to carry $\bX_n'$ with
  $(\bX_n', S_{n,h}') \equallaw (\bX_n, S_{n,h})$,
  and compute $\hat\theta'_{n,h}$ from $\bX_n'$.
  Define $\tilde Z'_{n,h}(w, v) \coloneqq
  e_1^\T H_h(w)^{-1} \tilde S_{n,h}'(w, v)$ so that,
  as in Theorem~\ref{thm:locpol_gaussian},
  \begin{align}
    \nonumber
    &\P\Biggl\{
      \sup_{(w, v) \in \cW \times \cV}
      \sqrt{n h^{D_W}}
      \bigl|
      \hat\theta'_{n,h}(w, v) - \theta(w, v) - \tilde Z'_{n,h}(w, v)
      \bigr| \\
      &\qquad>
      \frac{C \eta \sqrt{\log 2n}}
      {\bigl\{(n h^{D_W})^{k-1} h^{2D_w}\bigr\}^{1/(2k+2 + k D_w + 2k D_v)}}
      + C \sqrt{n h^{D_W}} h^\zeta
    \Biggr\}
    \iftoggle{journal}{\nonumber \\ &\quad}{}
    \leq
    2 e^{-{n h^{D_W} / C}}
    + \frac{1}{\eta^{k+1}}.
    \label{eq:locpol_bootstrap_ks_auxiliary_coupling}
  \end{align}
  Let $\bigl(\hat\theta''_{n,h},\tilde\theta^{*\dprime}_{n,h},\bX_n''\bigr)$
  be the coupling constructed in Theorem~\ref{thm:locpol_bootstrap}.
  Apply the Vorob'ev--Berkes--Philipp theorem to
  $\bigl(\hat\theta'_{n,h},\tilde Z'_{n,h}\bigr)$ and
  $\bigl(\hat\theta''_{n,h},\tilde\theta^{*\dprime}_{n,h},\bX_n''\bigr)$,
  yielding
  $\bigl(\hat\theta'''_{n,h},\tilde Z'''_{n,h},
  \tilde\theta^{*\tprime}_{n,h},\bX_n'''\bigr)$;
  in particular, $\hat\theta'''_{n,h}$ is independent of $\bX_n'''$.

  Next, by an argument similar to that used to show \eqref{eq:ks_decomp},
  defining $\hat\Sigma_{n,h}'''$ using the data $\bX_n'''$,
  under condition~\ref{it:thm_locpol_bootstrap_nonparametric},
  we have
  \begin{align}
    \nonumber
    &\sup_{\xi \in \R}
    \Biggl|
    \P\Biggl(
      \sup_{(w, v) \in \cW \times \cV}
      \Biggl|
      \frac{\hat\theta_{n,h}(w, v) - \theta(w, v)}
      {\sqrt{\smash[b]{\hat \Sigma_{n,h}(w,v,w,v)}}}
      \Biggr| \leq \xi
    \Biggr)
    \iftoggle{journal}{\\ \nonumber &\qquad\qquad}{}
    -
    \P\Biggl(
      \sup_{(w, v) \in \cW \times \cV}
      \Biggl|
      \frac{\tilde\theta^{*\tprime}_{n,h}(w, v)}
      {\sqrt{\smash[b]{\hat \Sigma_{n,h}'''(w,v,w,v)}}}
      \Biggr| \leq \xi
      \biggm| \bX_n'''
    \Biggr)
    \Biggr| \\
    \label{eq:locpol_bootstrap_ks_coupling}
    &\qquad\leq
    \P\Biggl(
      \sup_{(w, v) \in \cW \times \cV}
      \Biggl|
      \frac{\hat\theta_{n,h}'''(w, v) - \theta(w, v)}
      {\sqrt{\Sigma_h(w, v, w, v)}}
      - \frac{\tilde\theta_{n,h}^{*\tprime}(w, v)}
      {\sqrt{\Sigma_h(w, v, w, v)}}
      \Biggr| > \eta
      \biggm| \bX_n'''
    \Biggr) \\
    \label{eq:locpol_bootstrap_ks_gaussian_coupling}
    &\qquad\qquad+
    \P\Biggl(
      \sup_{(w, v) \in \cW \times \cV}
      \Biggl|
      \frac{\hat\theta_{n,h}'''(w, v) - \theta(w, v)}
      {\sqrt{\Sigma_h(w, v, w, v)}}
      - \frac{\tilde Z'''_{n,h}(w, v)}
      {\sqrt{\Sigma_h(w, v, w, v)}}
      \Biggr| > \eta
    \Biggr) \\
    \label{eq:locpol_bootstrap_ks_anti-concentration}
    &\qquad\qquad+
    \sup_{\xi \in \R}
    \P\Biggl(
      \xi <
      \sup_{(w, v) \in \cW \times \cV}
      \Biggl|
      \frac{\tilde Z'''_{n,h}(w, v)}
      {\sqrt{\Sigma_h(w, v, w, v)}}
      \Biggr| \leq \xi + 5\eta
    \Biggr) \\
    \label{eq:locpol_bootstrap_ks_convergence_thetahat}
    &\qquad\qquad+
    \P\Biggl(
      \sup_{(w, v) \in \cW \times \cV}
      \Biggl|
      \frac{\hat\theta_{n,h}(w, v) - \theta(w, v)}
      {\sqrt{\smash[b]{\hat \Sigma_{n,h}(w,v,w,v)}}}
      - \frac{\hat\theta_{n,h}(w, v) - \theta(w, v)}
      {\sqrt{\Sigma_h(w, v, w, v)}}
      \Biggr| > \eta
    \Biggr) \\
    \label{eq:locpol_bootstrap_ks_convergence_thetahatstar}
    &\qquad\qquad+
    \P\Biggl(
      \sup_{(w, v) \in \cW \times \cV}
      \Biggl|
      \frac{\tilde\theta^{*\tprime}_{n,h}(w, v)}
      {\sqrt{\smash[b]{\hat \Sigma_{n,h}'''(w,v,w,v)}}}
      - \frac{\tilde\theta^{*\tprime}_{n,h}(w, v)}
      {\sqrt{\Sigma_h(w, v, w, v)}}
      \Biggr| > \eta
      \biggm| \bX_n'''
    \Biggr).
  \end{align}
  We control each term in turn.
  Firstly, Lemma~\ref{lem:Sigma_hat} and
  Theorem~\ref{thm:locpol_bootstrap} ensure that
  the expectation of \eqref{eq:locpol_bootstrap_ks_coupling}
  converges to zero whenever
  \begin{equation*}
    \eta \gg
    \frac{\sqrt{\log 2n} }
    {(n h^{D_W + (4-k)D_w})^{1/(10-2k+2D_w+4D_v)}}
    + \sqrt{n h^{D_W}} h^\zeta.
  \end{equation*}
  By \eqref{eq:locpol_bootstrap_ks_auxiliary_coupling},
  we have that \eqref{eq:locpol_bootstrap_ks_gaussian_coupling}
  converges to zero if
  \begin{equation*}
    \eta \gg
    \frac{\sqrt{\log 2n}}
    {\bigl\{(n h^{D_W})^{k-1} h^{2D_w}\bigr\}^{1/(2k+2 + k D_w + 2k D_v)}}
    + \sqrt{n h^{D_W}} h^\zeta.
  \end{equation*}
  By Lemma~\ref{lem:anti-concentration_gp} conditional
  on $\bX_n$ and $(\zeta_i: i \in [n])$,
  bounding the expected supremum by Dudley's theorem
  \citep[Corollary~2.2.8]{van1996weak}
  and the entropy bounds in the proof of
  Lemma~\ref{lem:locpol_gaussian},
  and using Lemma~\ref{lem:Sigma_hat} to obtain
  a lower bound on the variance function, we have that
  \begin{equation*}
    \sup_{\xi \in \R}
    \P\Biggl(
      \xi <
      \sup_{(w, v) \in \cW \times \cV}
      \Biggl|
      \frac{\tilde Z'''_{n,h}(w, v)}
      {\sqrt{\Sigma_h(w, v, w, v)}}
      \Biggr| \leq \xi + \eta
    \Biggr)
    \leq C \eta \sqrt{\log 2n},
  \end{equation*}
  and so the expectation of
  \eqref{eq:locpol_bootstrap_ks_anti-concentration}
  is negligible whenever $\eta \sqrt{\log 2n} \to 0$.
  The terms in \eqref{eq:locpol_bootstrap_ks_convergence_thetahat}
  and \eqref{eq:locpol_bootstrap_ks_convergence_thetahatstar}
  are handled exactly as in the proof of
  Theorem~\ref{thm:locpol_gaussian_ks},
  applying Lemmas~\ref{lem:locpol_convergence}, \ref{lem:Sigma_hat},
  and~\ref{lem:locpol_bootstrap_convergence}.
  The same arguments handle
  condition~\ref{it:thm_locpol_bootstrap_multiplier}.
  We conclude that a suitable sequence of $\eta$ values exists if
  \begin{equation*}
    \frac{\log 2n}
    {(n h^{D_W + (4-k)D_w})^{1/(10-2k+2D_w+4D_v)}}
    + h^\zeta \sqrt{n h^{D_W} \log 2n}
    \to 0,
  \end{equation*}
  which is implied by the conditions of the lemma.
\end{proof}

\begin{proof}[Proof of Theorem~\ref{thm:locpol_bootstrap_cb}]
  Under condition~\ref{it:thm_locpol_bootstrap_nonparametric},
  by applying the Dvoretzky--Kiefer--Wolfowitz inequality
  \citep{dvoretzky1956asymptotic} conditional on $\bX_n$
  and by Lemma~\ref{lem:locpol_bootstrap_ks}, as $n, B \to \infty$,
  \begin{align*}
    &\sup_{\xi \in \R}
    \,
    \Biggl|
    \frac{1}{B}\sum_{b=1}^B
    \I\Biggl\{
      \sup_{(w,v) \in \cW \times \cV}
      \Biggl|
      \frac{\tilde\theta_{n,h}^{*(b)}(w,v)}
      {\sqrt{\smash[b]{\hat\Sigma_{n,h}(w,v,w,v)}}}
      \Biggr| \leq \xi
    \Biggr\} \\
    &\qquad\qquad-
    \P\Biggl(
      \sup_{(w,v) \in \cW \times \cV}
      \Biggl|
      \frac{\hat\theta_{n,h}(w,v)-\theta(w,v)}
      {\sqrt{\smash[b]{\hat\Sigma_{n,h}(w,v,w,v)}}}
      \Biggr| \leq \xi
    \Biggr)
    \Biggr|
    \to_\P 0.
  \end{align*}
  Fix $\varepsilon > 0$.
  By the definition of $\hat\xi_{\alpha,B}$ and the previous display,
  with probability approaching one,
  \begin{align*}
    &\inf\Biggl\{
      \xi\in\R:
      \P\Biggl(
        \sup_{(w,v)\in\cW\times\cV}
        \Biggl|
        \frac{\hat\theta_{n,h}(w,v)-\theta(w,v)}
        {\sqrt{\smash[b]{\hat\Sigma_{n,h}(w,v,w,v)}}}
        \Biggr|\leq\xi
      \Biggr)
      \geq 1-\alpha-\varepsilon
    \Biggr\}
    \leq \hat\xi_{\alpha,B} \\
    &\qquad\qquad\leq
    \inf\Biggl\{
      \xi\in\R:
      \P\Biggl(
        \sup_{(w,v)\in\cW\times\cV}
        \Biggl|
        \frac{\hat\theta_{n,h}(w,v)-\theta(w,v)}
        {\sqrt{\smash[b]{\hat\Sigma_{n,h}(w,v,w,v)}}}
        \Biggr|\leq\xi
      \Biggr)
      \geq 1-\alpha+\varepsilon
    \Biggr\}.
  \end{align*}
  Further, by the anti-concentration argument in the
  proof of Lemma~\ref{lem:locpol_bootstrap_ks}, we have
  \begin{equation*}
    \sup_{\xi \in \R}
    \P\Biggl(
      \sup_{(w,v)\in\cW\times\cV}
      \Biggl|
      \frac{\hat\theta_{n,h}(w,v)-\theta(w,v)}
      {\sqrt{\smash[b]{\hat\Sigma_{n,h}(w,v,w,v)}}}
      \Biggr| = \xi
    \Biggr)
    \to 0.
  \end{equation*}
  It follows that for sufficiently large $n$ and $B$,
  \begin{align*}
    1-\alpha-2\varepsilon
    &\leq
    \P\Biggl(
      \theta(w,v) \in
      \biggl[
        \hat\theta_{n,h}(w,v)
        \pm \hat\xi_{\alpha,B}
        \sqrt{\smash[b]{\hat\Sigma_{n,h}(w,v,w,v)}}
      \biggr]
      \text{ for all }(w,v) \in \cW \times \cV
    \Biggr) \\
    &\leq 1-\alpha+2\varepsilon.
  \end{align*}
  Taking $\varepsilon \to 0$ yields the result.
  Under condition~\ref{it:thm_locpol_bootstrap_multiplier},
  repeat the argument with draws from
  the multiplier bootstrap process.
\end{proof}

\subsection{Proofs for Appendix~\ref{sec:consequences}}

\begin{proof}[Proof of Theorem~\ref{thm:prokhorov}
  and \eqref{eq:prokhorov_matching}]
  The theorem follows from Lemma~\ref{lem:borel}.
  For \eqref{eq:prokhorov_matching},
  if $n\bar\beta_{k,\gamma}\psi^k \geq 1$, then
  the result is immediate.
  Otherwise, take
  \begin{align*}
    t
    &\coloneqq
    \sqrt{\frac{2}{k+\gamma+1}
    \log\frac{1}{n\bar\beta_{k,\gamma}\psi^k}},
    &\eta
    &\coloneqq
    \bigl\{n\bar\beta_{k,\gamma}\phi(t)^k\bigr\}^{1/(k+\gamma+1)}.
  \end{align*}
  Since $\phi(t) \geq \psi$, we have
  \begin{align*}
    e^{-t^2/2}
    &= \bigl(n\bar\beta_{k,\gamma}\psi^k\bigr)^{1/(k+\gamma+1)}
    \leq \eta,
    &\frac{n\bar\beta_{k,\gamma}\phi(t)^k}{\eta^{k+\gamma}}
    &= \eta.
  \end{align*}
  Thus, Theorem~\ref{thm:prokhorov} gives
  $\rho_p(S_n,T_n) \leq 3\eta$.
\end{proof}

\begin{proof}[Details for \eqref{eq:prokhorov_rate_p}
  and \eqref{eq:prokhorov_rate_infty}]
  The derivation of \eqref{eq:prokhorov_matching}
  works with $n\bar\beta_{k,\gamma}$ replaced by any larger value.
  By the proof of \eqref{eq:yurinskii_rate_general},
  rescaling $S_n$ and $T_n$ by $1/\sqrt n$, we have
  $n \bar\beta_{k,\gamma} \leq
  2 n^{1 - k/2 - \gamma/2} d^{k/2 + \gamma / (p \land r)} M^{k+\gamma}$.
  Since $\psi \geq 1$, we therefore have
  \begin{equation*}
    \rho_p\biggl(\frac{S_n}{\sqrt n}, \frac{T_n}{\sqrt n}\biggr)
    \lesssim
    n^{\frac{2 - k - \gamma}{2k + 2\gamma + 2}}
    d^{\frac{k/2 + \gamma / (p \land r)}{k + \gamma + 1}}
    M^{\frac{k+\gamma}{k+\gamma+1}}
    \phi
    \Bigl( \sqrt{ \log 2n } \Bigr)^{\frac{k}{k + \gamma + 1}}.
  \end{equation*}
  We now proceed by cases on $p$.
  If $p \in [1, 2)$, we have
  $\phi(t) = \sqrt{p} d^{1/p} + d^{1/p - 1/2} t$, so
  \begin{align*}
    \rho_p\biggl(\frac{S_n}{\sqrt n}, \frac{T_n}{\sqrt n}\biggr)
    &\lesssim
    n^{\frac{2 - k - \gamma}{2k + 2\gamma + 2}}
    d^{\frac{k/2 + \gamma / (p \land r)}{k + \gamma + 1}}
    M^{\frac{k+\gamma}{k+\gamma+1}}
    \Bigl(
      \sqrt{p} d^{1/p} + d^{1/p - 1/2} \sqrt{\log 2n}
    \Bigr)^{\frac{k}{k + \gamma + 1}} \\
    &\lesssim
    n^{\frac{2 - k - \gamma}{2k + 2\gamma + 2}}
    d^{\frac{k/2 + \gamma/(p \land r) + k/p}{k + \gamma + 1}}
    M^{\frac{k+\gamma}{k+\gamma+1}}
    \biggl(
      p + \frac{\log 2n}{d}
    \biggr)^{\frac{k}{2k + 2\gamma + 2}}.
  \end{align*}
  For $p \in [2, \infty)$, we have
  $\phi(t) = \sqrt{p} d^{1/p} + t$, so
  \begin{align*}
    \rho_p\biggl(\frac{S_n}{\sqrt n}, \frac{T_n}{\sqrt n}\biggr)
    &\lesssim
    n^{\frac{2 - k - \gamma}{2k + 2\gamma + 2}}
    d^{\frac{k/2 + \gamma / (p \land r)}{k + \gamma + 1}}
    M^{\frac{k+\gamma}{k+\gamma+1}}
    \Bigl(
      \sqrt{p} d^{1/p} + \sqrt{\log 2n}
    \Bigr)^{\frac{k}{k + \gamma + 1}} \\
    &\lesssim
    n^{\frac{2 - k - \gamma}{2k + 2\gamma + 2}}
    d^{\frac{k/2 + \gamma/(p \land r) + k/p}{k + \gamma + 1}}
    M^{\frac{k+\gamma}{k+\gamma+1}}
    \biggl(
      p + \frac{\log 2n}{d^{2/p}}
    \biggr)^{\frac{k}{2k + 2\gamma + 2}}.
  \end{align*}
  Finally, if $p = \infty$, then
  $\phi(t) = \sqrt{2 \log 2 d} + t$, giving
  \begin{align*}
    \rho_p\biggl(\frac{S_n}{\sqrt n}, \frac{T_n}{\sqrt n}\biggr)
    &\lesssim
    n^{\frac{2 - k - \gamma}{2k + 2\gamma + 2}}
    d^{\frac{k/2 + \gamma / r}{k + \gamma + 1}}
    M^{\frac{k+\gamma}{k+\gamma+1}}
    \Bigl(
      \sqrt{2 \log 2d} + \sqrt{\log 2n}
    \Bigr)^{\frac{k}{k + \gamma + 1}} \\
    &\lesssim
    n^{\frac{2 - k - \gamma}{2k + 2\gamma + 2}}
    d^{\frac{k/2 + \gamma/r}{k + \gamma + 1}}
    M^{\frac{k+\gamma}{k+\gamma+1}}
    \bigl(
      \log 2 n d
    \bigr)^{\frac{k}{2k + 2\gamma + 2}}.
    \iftoggle{journal}{}{\qedhere}
  \end{align*}
\end{proof}

\begin{lemma}[Gaussian moment matching on a lattice]%
  \label{lem:gaussian_moment_matching}
  Let $Z \sim N(0, 1)$.
  For each $k \in \N_0$,
  there exist $\varepsilon_k > 0$
  and a distribution $P_k$
  with at most $k + 1$ atoms, all contained in
  $\varepsilon_k \Z \coloneqq \{\varepsilon_k t: t \in \Z\}$,
  such that if $X \sim P_k$ then
  $\E(Z^r) = \E(X^r)$ for $0 \leq r \leq k$.
\end{lemma}

\begin{proof}[Proof of Lemma~\ref{lem:gaussian_moment_matching}]
  The result is trivial for $k = 0$, so fix $k \in \N$.
  Define $\mu_r \coloneqq \E(Z^r)$ for $0 \leq r \leq 2k$.
  Write $\cP_k$ for the space of polynomials on
  $\R$ with degree at most $k$, and define the linear functional
  $L: \cP_{2k} \to \R$ by $L(f) \coloneqq \E\{f(Z)\}$.
  If $f \in \cP_k$ with $L(f^2) = 0$, then as $Z$ has a
  strictly positive Lebesgue density, we must have $f \equiv 0$.
  It follows that the sequence $(\mu_r)_{0 \leq r \leq 2k}$
  is positive definite in the sense of
  \citet[Equation~9.1]{schmudgen2017moment}.
  By \citet[Corollary~9.8]{schmudgen2017moment},
  there exist distinct points $x_0, \ldots, x_k \in \R$ and
  quantities $p_0, \ldots, p_k > 0$ such that
  $\sum_{j=0}^k x_j^r p_j = \mu_r$ whenever $0 \leq r \leq 2k$.
  Define $p \coloneqq (p_0, \ldots, p_k)^\T$
  and $\mu \coloneqq (\mu_0, \ldots, \mu_k)^\T$.
  Consider the matrix $A \in \R^{(k+1) \times (k+1)}$
  given by $A_{r j} \coloneqq x_j^r$
  for $0 \leq r \leq k$ and $0 \leq j \leq k$,
  so that $A p = \mu$.
  Then $A$ is a transposed square Vandermonde matrix, and
  $x_0, \ldots, x_k$ are distinct, so $A$ is invertible
  \citep[Section~4.3, Exercise~22]{friedberg2003linear}
  and $p = A^{-1} \mu$.
  For $0 \leq j \leq k$,
  let $y_j \in \varepsilon_k \Z$
  satisfy $|y_j - x_j| \leq \varepsilon_k / 2$, where
  $0 < \varepsilon_k < \inf_{i \neq j} |x_i - x_j|$,
  so it follows that $y_0, \ldots, y_k$ are distinct.
  Define the matrix $B \in \R^{(k+1) \times (k+1)}$ by $B_{r j} = y_j^r$
  for $0 \leq r \leq k$ and $0 \leq j \leq k$,
  so that $B$ is an invertible transposed Vandermonde matrix.
  Set $q \coloneqq B^{-1} \mu$
  so that $\sum_{j=0}^k y_j^r q_j = \mu_r$ whenever $0 \leq r \leq k$.
  The difference satisfies
  \begin{equation*}
    \|q - p\|_2
    =
    \bigl\|\bigl(B^{-1} - A^{-1}\bigr) \mu\bigr\|_2
    \leq
    \|B^{-1}\|_{2 \to 2}
    \|B - A\|_{2 \to 2}
    \|A^{-1}\|_{2 \to 2}
    \|\mu\|_2,
  \end{equation*}
  which can be made arbitrarily small since
  $\|B - A\|_{2 \to 2} \to 0$ as $\varepsilon_k \to 0$.
  Therefore, we may assume that
  $q_j > 0$ for each $0 \leq j \leq k$.
  To conclude, let $P_k$ be the probability measure on
  $\{y_0, \ldots, y_k\}$ with
  $P_k(\{y_j\}) = q_j$ for $0 \leq j \leq k$.
\end{proof}

\begin{proof}[Proof of Lemma~\ref{lem:lower_bound}]
  Let $P_k$ and $\varepsilon_k$ be as in
  Lemma~\ref{lem:gaussian_moment_matching};
  we may assume $\varepsilon_k \leq 1$ without loss of generality.
  Let $A_n \coloneqq (\varepsilon_k / \sqrt n) \Z \in \cB(\R)$,
  so that $S_n/\sqrt n \in A_n$ almost surely.
  For $\eta > 0$, let $A_n^\eta$ be the $\eta$-enlargement of $A_n$
  as in \eqref{eq:enlargement}.
  For $\eta < \varepsilon_k / (2 \sqrt n)$
  we have the disjoint union of open intervals
  \begin{equation*}
    A_n^\eta
    = \bigsqcup_{j = -\infty}^{\infty}
    \biggl(
      \frac{j \varepsilon_k}{\sqrt n} - \eta,
      \,
      \frac{j \varepsilon_k}{\sqrt n} + \eta
    \biggr).
  \end{equation*}
  Since $Z \sim N(0, 1)$, by a Riemann sum,
  \begin{align*}
    \P(Z \in A_n^\eta)
    &=
    \sum_{j = -\infty}^{\infty}
    \P \biggl(
      \frac{j \varepsilon_k}{\sqrt n} - \eta
      \leq Z \leq
      \frac{j \varepsilon_k}{\sqrt n} + \eta
    \biggr) \\
    &=
    \P \bigl(-\eta \leq Z \leq \eta\bigr)
    + 2 \sum_{j=1}^{\infty}
    \P \biggl(
      \frac{j \varepsilon_k}{\sqrt n} - \eta
      \leq Z \leq
      \frac{j \varepsilon_k}{\sqrt n} + \eta
    \biggr) \\
    &\leq
    \frac{2 \eta}{\sqrt{2 \pi}}
    + \frac{4 \eta}{\sqrt{2 \pi}}
    \sum_{j=1}^{\infty}
    \exp \Biggl\{- \frac{1}{2}
      \biggl(\frac{j \varepsilon_k}{\sqrt n} - \eta\biggr)^2
    \Biggr\}
    \leq
    \frac{2 \eta}{\sqrt{2 \pi}}
    + \frac{4 \eta}{\sqrt{2 \pi}}
    \sum_{j=1}^{\infty}
    \exp \biggl( -\frac{j^2 \varepsilon_k^2}{8 n} \biggr) \\
    &\leq
    \frac{2 \eta}{\sqrt{2 \pi}}
    + \frac{4 \eta}{\sqrt{2 \pi}}
    \int_0^\infty
    \exp \biggl( -\frac{\varepsilon_k^2 t^2}{8 n} \biggr)
    \diffi t
    =
    \frac{2 \eta}{\sqrt{2 \pi}}
    + \frac{4 \eta \sqrt{n}}{\varepsilon_k}
    \leq
    \frac{5 \eta \sqrt{n}}{\varepsilon_k}.
  \end{align*}
  Thus, if
  $\P\bigl(S_n/\sqrt n \in A_n\bigr) \leq
  \P\bigl(Z \in A_n^\eta\bigr) + \eta$,
  we must have that either
  $\eta \geq \varepsilon_k / (2 \sqrt n)$ or
  $1 \leq 5 \eta \sqrt n / \varepsilon_k + \eta$,
  so that $\eta \geq \varepsilon_k / (6 \sqrt n)$.
  It follows that
  $\rho\bigl(S_n/\sqrt n, Z\bigr) \geq \varepsilon_k / (6 \sqrt n)$,
  as required.
\end{proof}

\begin{lemma}[Anti-concentration property]%
  \label{lem:anti-concentration}
  Let $n, d \in \N$. For $i \in [n]$,
  let $X_i$ be independent random vectors taking values in $\R^d$.
  Suppose that $\E(\|X_i\|_2^2) < \infty$ and define
  $\mu_i \coloneqq \E(X_i)$ for each $i \in [n]$.
  Let $\Sigma_i \coloneqq \Var(X_i) \in \R^{d \times d}$
  and $\Sigma \coloneqq \sum_{i=1}^n \Sigma_i$,
  and suppose that $\sigma^2 \coloneqq \min_{j \in [d]} \Sigma_{j j} > 0$.
  For $i \in [n]$, construct $\tilde X_i$ independently
  using $X_i$ as in Lemma~\ref{lem:anti-concentration_construction}.
  With $\tilde S_n \coloneqq \sum_{i=1}^n \tilde X_i$,
  for any $\eta > 0$,
  \begin{equation*}
    \sup_{\xi \in \R}
    \P\bigl(\xi < \|\tilde S_n\|_\infty \leq \xi + \eta\bigr)
    \leq
    \frac{2 \eta \sqrt{\log 2d}}{\sigma}.
  \end{equation*}
\end{lemma}

\begin{proof}[Proof of Lemma~\ref{lem:anti-concentration}]
  By Lemma~\ref{lem:anti-concentration_construction},
  conditional on $(\zeta_i, X_i)$, we have
  $\tilde X_i \sim N(\mu_i + \zeta_i (X_i - \mu_i), \Sigma_i/2)$.
  Therefore, conditional on
  $\cH \coloneqq (\zeta_i, X_i: i \in [n])$, we have
  $\tilde S_n \sim N\bigl(\sum_{i=1}^n \{\mu_i + \zeta_i (X_i - \mu_i)\},
  \Sigma/2\bigr)$, by independence.
  By \citet[Lemma~SA.4]{cattaneo2025yurinskii},
  \begin{equation*}
    E \coloneqq
    \E\Biggl(
      \max_{j \in [d]}
      \frac{
        \bigl|e_j^\T\bigl\{\tilde S_n-\E(\tilde S_n\mid\cH)\bigr\}\bigr|
      }{\sqrt{\Sigma_{j j}/2}}
      \biggm| \cH
    \Biggr)
    \leq \sqrt{2\log 2d}.
  \end{equation*}
  Applying Lemma~\ref{lem:anti-concentration_gp}
  conditionally on $\cH$ and taking an expectation yields the result.
\end{proof}

\begin{lemma}[Gaussian norm inequality]%
  \label{lem:gaussian_norm}
  Take $d \in \N$ and let $\Sigma \in \R^{d \times d}$
  be positive semidefinite.
  Suppose that $\max_{j \in [d]}\Sigma_{j j} \leq \sigma^2$
  for some $\sigma^2 > 0$.
  Then, if $X \sim N(0_d, \Sigma)$,
  $k \in \N_0$ and $\gamma \in [0, 1]$,
  \begin{equation*}
    \E\Bigl\{
      \|X\|_2^k
      \|X\|_\infty^\gamma
    \Bigr\}
    \leq
    2^\gamma \sqrt{(2k - 1)!!} \, \sigma^{k+\gamma}
    d^{k/2} (\log 2d)^{\gamma/2},
  \end{equation*}
  where $(2k-1)!! \coloneqq (2k-1) \times (2k-3)
  \times \cdots \times 3 \times 1$
  is the double factorial, defining $(-1)!! \coloneqq 1$.
\end{lemma}

\begin{proof}[Proof of Lemma~\ref{lem:gaussian_norm}]
  By \citet[supplementary material, proof of Lemma~SA.3]{%
  cattaneo2025yurinskii},
  \begin{equation*}
    \E\bigl\{
      \|X\|_\infty^2
    \bigr\}
    \leq
    4 \sigma^2 \log 2d.
  \end{equation*}
  By Jensen's inequality and the moments
  of a Gaussian distribution,
  \begin{equation*}
    \E\bigl\{
      \|X\|_2^{2k}
    \bigr\}
    =
    \E\Biggl\{
      \biggl(
        \sum_{j=1}^d
        X_j^2
      \biggr)^k
    \Biggr\}
    \leq
    d^{k-1}
    \sum_{j=1}^d
    \E\bigl\{
      X_j^{2 k}
    \bigr\}
    \leq
    d^k (2k - 1)!! \,\sigma^{2k}.
  \end{equation*}
  Therefore, by Cauchy--Schwarz and Jensen's inequality,
  \begin{align*}
    \E\Bigl\{
      \|X\|_2^k
      \|X\|_\infty^\gamma
    \Bigr\}
    &\leq
    \Bigl\{\E\bigl(
        \|X\|_2^{2k}
    \bigr)\Bigr\}^{1/2}
    \, \Bigl\{\E\bigl(
        \|X\|_\infty^2
    \bigr)\Bigr\}^{\gamma/2} \\
    &\leq
    2^{\gamma} \sqrt{(2k - 1)!!} \, \sigma^{k+\gamma}
    d^{k/2} (\log 2d)^{\gamma/2}.
    \iftoggle{journal}{}{\qedhere}
  \end{align*}
\end{proof}

\begin{proof}[Proof of Theorem~\ref{thm:ks}]
  The proof proceeds by constructing couplings for
  $(S_n, T_n)$ and $(S_n, \tilde S_n)$,
  where $\tilde S_n$ is defined as in
  Lemma~\ref{lem:anti-concentration}.
  We write $C > 0$ for a universal constant;
  the value of $C$ may vary from line to line.
  For $\xi \in \R$ and $\eta > 0$, we have
  \begin{align*}
    &\P\bigl(\|T_n\|_\infty \leq \xi\bigr)
    - \P\bigl(\|S_n\|_\infty \leq \xi\bigr) \\
    &\qquad\leq
    \P\bigl(\xi < \|S_n\|_\infty \leq \xi + \eta\bigr)
    + \P\bigl(\|S_n' - T_n'\|_\infty > \eta\bigr) \\
    &\qquad\leq
    \P\bigl(\xi - \eta < \|\tilde S_n\|_\infty \leq \xi + 2\eta\bigr)
    + \P\bigl(\|S_n' - T_n'\|_\infty > \eta\bigr)
    + \P\bigl(\|S_n'' - \tilde S_n''\|_\infty > \eta\bigr),
  \end{align*}
  where $S_n \equallaw S_n' \equallaw S_n''$,
  $T_n \equallaw T_n'$ and
  $\tilde S_n \equallaw \tilde S_n''$.
  Similarly,
  \begin{align*}
    &\P\bigl(\|S_n\|_\infty \leq \xi\bigr)
    - \P\bigl(\|T_n\|_\infty \leq \xi\bigr)
    \leq
    \P\bigl(\xi - \eta < \|S_n\|_\infty \leq \xi\bigr)
    + \P\bigl(\|S_n' - T_n'\|_\infty > \eta\bigr) \\
    &\qquad\leq
    \P\bigl(\xi - 2\eta < \|\tilde S_n\|_\infty \leq \xi + \eta\bigr)
    + \P\bigl(\|S_n' - T_n'\|_\infty > \eta\bigr)
    + \P\bigl(\|S_n'' - \tilde S_n''\|_\infty > \eta\bigr).
  \end{align*}
  It follows that
  \begin{align}
    \nonumber
    &\sup_{\xi \in \R} \,
    \Bigl|
    \P\bigl(\|S_n\|_\infty \leq \xi\bigr)
    - \P\bigl(\|T_n\|_\infty \leq \xi\bigr)
    \Bigr|
    \leq
    \sup_{\xi \in \R}
    \P\bigl(\xi < \|\tilde S_n\|_\infty \leq \xi + 3\eta\bigr) \\
    \label{eq:ks_decomp}
    &\qquad+
    \P\bigl(\|S_n' - T_n'\|_\infty > \eta\bigr)
    + \P\bigl(\|S_n'' - \tilde S_n''\|_\infty > \eta\bigr).
  \end{align}
  The terms in \eqref{eq:ks_decomp} involving $\|S_n' - T_n'\|_\infty$
  and $\|S_n'' - \tilde S_n''\|_\infty$
  are controlled by two applications of
  Theorem~\ref{thm:yurinskii}, as follows.
  As $r \in [k + \gamma, \infty]$ and $k \geq 2$,
  H{\"o}lder's inequality gives
  $\|x\|_\infty \leq \|x\|_r$
  and $\|x\|_2 \leq d^{1/2 - 1/r} \|x\|_r$
  for all $x \in \R^d$.
  Thus, by Jensen's inequality,
  \begin{align*}
    \E \bigl(
      \|X_i\|_2^k
      \| X_i\|_\infty^{\gamma}
    \bigr)
    &\leq
    d^{k/2 - k/r}
    \E \bigl(
      \| X_i\|_r^{k+\gamma}
    \bigr)
    \leq
    d^{k/2 - k/r}
    \bigvvvert
    \| X_i\|_r
    \bigvvvert_r^{k+\gamma} \\
    &\leq
    d^{k/2 - k/r}
    d^{k/r + \gamma/r}
    M^{k+\gamma}
    =
    d^{k/2 + \gamma / r} M^{k+\gamma}.
  \end{align*}
  The same bound holds with $X_i$ replaced by $Y_i$.
  Write $c \coloneqq 1 + \sqrt{6}/2$ and let
  $\mu_i \coloneqq \E(X_i)$
  and $U_i \coloneqq \mu_i + \zeta_i (X_i - \mu_i)$, so that
  $\tilde X_i = U_i + Z_i$.
  Since $|\zeta_i| \leq c$, for each $j \in [d]$,
  by Jensen's inequality,
  \begin{equation*}
    \vvvert e_j^\T U_i \vvvert_r
    \leq
    |e_j^\T \mu_i|
    + c \vvvert e_j^\T (X_i - \mu_i) \vvvert_r
    \leq (1 + 2c) M,
  \end{equation*}
  and $(\Sigma_i)_{j j} \leq M^2$.
  By independence and Lemma~\ref{lem:gaussian_norm},
  \begin{align*}
    \E \bigl(
      \|\tilde X_i\|_2^k
      \|\tilde X_i\|_\infty^\gamma
    \bigr)
    &\leq
    C \, \E \Bigl\{
      \bigl(
        \|U_i\|_2^k + \|Z_i\|_2^k
      \bigr)
      \bigl(
        \|U_i\|_\infty^\gamma + \|Z_i\|_\infty^\gamma
      \bigr)
    \Bigr\} \\
    &\leq
    C d^{k/2}
    \bigl(
      d^{\gamma/r} + (\log 2d)^{\gamma/2}
    \bigr)
    M^{k + \gamma}.
  \end{align*}
  It follows that
  \begin{equation*}
    \frac{1}{\sqrt{k!}}
    \E \Bigl(
      2 \|X_i\|_2^k \|X_i\|_\infty^{\gamma}
      + \|Y_i\|_2^k \|Y_i\|_\infty^{\gamma}
      + \|\tilde X_i\|_2^k \|\tilde X_i\|_\infty^{\gamma}
    \Bigr)
    \leq
    C d^{k/2}
    \bigl( d^{1 / r} + \sqrt{\log 2d} \bigr)^{\gamma}
    M^{k + \gamma}.
  \end{equation*}
  By Lemma~\ref{lem:anti-concentration_construction},
  $\E(X_i^\alpha) = \E(\tilde X_i^\alpha)$
  for all $\alpha \in \N_0^d$ with $|\alpha| \in [k]$.
  Apply Theorem~\ref{thm:yurinskii} with $p=\infty$ first to
  the summands $(X_i, Y_i)$ and then to
  $(X_i, \tilde X_i)$. In the second
  application, set $m=0$; exact moment matching gives
  $\Delta=\Gamma_\eta=0$.
  Replacing $\eta$ by $\sqrt{2}\eta$ and using
  $\phi(t)/\sqrt{2}\leq\sqrt{\log 2d}+t$ and
  $\Gamma_{\sqrt{2}\eta}(s)\leq\Gamma_\eta(s)$
  yields, for $t > 0$,
  \begin{align*}
    &\P\bigl(\|S_n' - T_n'\|_\infty > 3\sqrt{2}\eta\bigr)
    + \P\bigl(\|S_n'' - \tilde S_n''\|_\infty > 3\sqrt{2}\eta\bigr)
    \leq
    4 e^{-t^2 / 2}
    + \Delta\biggl(\frac{\sqrt{\log 2d} + t}{\eta}\biggr) \\
    &\qquad+
    \frac{C d^{k/2} \bigl( d^{1 / r} + \sqrt{\log 2d} \bigr)^{\gamma}
    \bigl(\sqrt{\log 2 d} + t\bigr)^{k} n M^{k + \gamma}}
    {\eta^{k + \gamma}}
    + \Gamma_\eta\biggl(\frac{\sqrt{\log 2d} + t}{\eta}\biggr).
  \end{align*}
  By Lemma~\ref{lem:anti-concentration} and
  \eqref{eq:ks_decomp} with $3\sqrt{2}\eta$ in place of $\eta$, we deduce that
  for any $\eta, t > 0$,
  \begin{align*}
    &\sup_{\xi \in \R} \,
    \Bigl|
    \P\bigl(\|S_n\|_\infty \leq \xi\bigr)
    - \P\bigl(\|T_n\|_\infty \leq \xi\bigr)
    \Bigr|
    \iftoggle{journal}{\\ &\quad}{}
    \leq
    4 e^{-t^2 / 2}
    + \Delta\biggl(\frac{\sqrt{\log 2d} + t}{\eta}\biggr)
    + \Gamma_\eta\biggl(\frac{\sqrt{\log 2d} + t}{\eta}\biggr) \\
    &\qquad+
    \frac{C \eta \sqrt{\log 2 d}}{\sigma}
    + \frac{C d^{k/2 + \gamma / r}
    \bigl(\sqrt{\log 2 d} + t\bigr)^{k + \gamma} n M^{k + \gamma}}
    {\eta^{k + \gamma}}.
    \iftoggle{journal}{}{\qedhere}
  \end{align*}
\end{proof}

\begin{proof}[Proof of Theorem~\ref{thm:bootstrap_ks}]
  The proof proceeds similarly to that of
  Theorem~\ref{thm:ks}.
  We write $C > 0$ for a universal constant;
  the value of $C$ may vary from line to line.
  For $\xi \in \R$ and $\eta > 0$, by
  the same argument used to show \eqref{eq:ks_decomp},
  with $\tilde S_n$ constructed as in
  Lemma~\ref{lem:anti-concentration}
  (using $X_i - \mu$ in place of $X_i$),
  with $S_n \equallaw S_n' \equallaw S_n''$,
  $(\bX_n, S_n^*) \equallaw (\bX_n', S_n^{*\prime})$,
  $\tilde S_n \equallaw \tilde S_n''$, and
  $S_n'$ independent of $\bX_n'$, we have
  \begin{align*}
    &\sup_{\xi \in \R} \,
    \Bigl|
    \P\bigl(\|S_n\|_\infty \leq \xi\bigr)
    - \P\bigl(\|S_n^{*\prime}\|_\infty \leq \xi \mid \bX_n'\bigr)
    \Bigr|
    \leq
    \sup_{\xi \in \R}
    \P\bigl(\xi < \|\tilde S_n\|_\infty \leq \xi + 3\eta\bigr) \\
    &\qquad+
    \P\bigl(\|S_n' - S_n^{*\prime}\|_\infty > \eta \mid \bX_n'\bigr)
    + \P\bigl(\|S_n'' - \tilde S_n''\|_\infty > \eta\bigr).
  \end{align*}
  Minkowski's and H{\"o}lder's inequalities give
  $M_2\leq\sqrt d M$ and $M_\infty\leq d^{1/r}M$.
  Thus, applying Theorem~\ref{thm:bootstrap}
  with $p = \infty$, $k = 3$ and $\gamma = 1$,
  and by the proof of Theorem~\ref{thm:ks}, since $r \geq 4$,
  \begin{align*}
    &\P\bigl(\|S_n' - S_n^{*\prime}\|_\infty > \eta\bigr)
    + \P\bigl(\|S_n'' - \tilde S_n''\|_\infty > \eta\bigr) \\
    &\qquad\leq
    4 e^{-t^2 / 2}
    + \frac{C \sqrt{n} d (\log 2d + t^2) M^2}{\eta^2}
    + \frac{C n d^{3/2}
      \bigl(d^{1/r} + \sqrt{\log 2 d}\bigr)
      \bigl(\sqrt{\log 2 d} + t\bigr)^3
    M^4}{\eta^4}.
  \end{align*}
  The last term above is at most
  $C^2 \bigl\{\sqrt n d \bigl(\log 2d + t^2\bigr)M^2 / \eta^2\bigr\}^2$
  since $d^{1/r} \leq \sqrt d$, so is bounded by
  $C \sqrt n d \bigl(\log 2d + t^2\bigr)M^2 / \eta^2$
  whenever this quantity is at most one.
  Otherwise, the right-hand side exceeds one and the bound is trivial.
  By Lemma~\ref{lem:anti-concentration},
  noting that the supremum is measurable by right continuity
  of the distribution functions and taking an expectation,
  we deduce that for $\eta, t > 0$,
  \begin{align*}
    &\E \biggl\{
      \sup_{\xi \in \R} \,
      \Bigl|
      \P\bigl(\|S_n\|_\infty \leq \xi\bigr)
      - \P^*\bigl(\|S_n^*\|_\infty \leq \xi\bigr)
      \Bigr|
    \biggr\} \\
    &\qquad\leq
    4 e^{-t^2 / 2}
    + \frac{6 \eta \sqrt{\log 2 d}}{\sigma}
    + \frac{C \sqrt{n} d (\log 2d + t^2) M^2}{\eta^2}.
    \iftoggle{journal}{}{\qedhere}
  \end{align*}
\end{proof}

\begin{proof}[Details for Example~\ref{ex:prony}]
  The polynomial $p(x)$ satisfies
  \begin{equation}
    \label{eq:prony_polynomial}
    \E\bigl\{X^i p(X)\bigr\}
    = \mu_{k+i} + \sum_{j=0}^{k-1} v_j \mu_{j+i} = 0
  \end{equation}
  for $0 \leq i \leq k-1$.
  Suppose that $z = a + b \imath$ is a root of $p(x)$,
  where $b \neq 0$ and $\imath$ is the imaginary unit.
  Then since $a - b \imath$ is also a root, we deduce that
  $p(x)$ is divisible by $(x-a)^2 + b^2 \geq b^2 > 0$.
  Write $q(x) \coloneqq p(x) / \bigl((x-a)^2 + b^2\bigr)$,
  which is a polynomial of degree $k - 2$.
  By \eqref{eq:prony_polynomial}, we have
  $0 = \E\{q(X) p(X)\} = \E\bigl\{q(X)^2 \bigl((X-a)^2 + b^2\bigr)\bigr\}$
  and so $\E\{q(X)^2\} = 0$.
  However, writing $q(x) = \sum_{j=0}^{k-1} u_j x^j$, we have
  $0 = \E\{q(X)^2\} = \sum_{i=0}^{k-1} \sum_{j=0}^{k-1} u_i u_j H_{i j}$
  and so $q(x) \equiv 0$ since $H$ is positive definite.
  This is a contradiction, so we conclude that every root of $p(x)$ is real.
  Now suppose that $y$ is a repeated real root of $p(x)$.
  Write $q(x) \coloneqq p(x) / (x - y)$, so that by \eqref{eq:prony_polynomial}
  we have $0 = \E\{q(X) p(X) / (X - y)\} = \E\{q(X)^2\}$. By the same argument,
  we conclude that $q(x) \equiv 0$ and reach a contradiction.
  Therefore, every root of $p(x)$ is real and non-repeated.
  For $0 \leq j \leq k-1$, we have
  $\sum_{i=1}^k x_i^{k+j} w_i
  = -\sum_{r=0}^{k-1} v_r \sum_{i=1}^k x_i^{r+j} w_i$.
  By induction on $j$ and using the definition of $v_0, \ldots, v_{k-1}$,
  we obtain $\sum_{i=1}^k x_i^{k+j} w_i
  = -\sum_{r=0}^{k-1} v_r \mu_{r+j} = \mu_{k+j}$.
  It remains to check that $w_i \geq 0$ for each $i \in [k]$.
  Let $V \in \R^{k \times k}$ satisfy
  $V_{i j}= x_j^i$ for $0 \leq i < k$ and $1 \leq j \leq k$.
  We have $H_{i j}=\mu_{i + j}=\sum_{r=1}^k w_r x_r^i x_r^j$,
  so $H=V\diag(w_1,\ldots,w_k)V^\T$.
  Since $x_1, \ldots, x_k$ are distinct,
  $\det V=\prod_{1\leq i<j\leq k}(x_j-x_i)\neq0$,
  so $V$ is invertible. Thus
  $\diag(w_1,\ldots,w_k)$
  is positive definite, giving $w_i > 0$ for each $i \in [k]$.
\end{proof}

\subsection{Proofs for Appendix~\ref{sec:further_empirical}}

\begin{proof}[Proof of Theorem~\ref{thm:gaussian_process}]
  The proof is similar
  to that of Theorem~\ref{thm:empirical_process}.
  Firstly, note that $\E\bigl\{\tilde\cF(Y_i)\bigr\} = 0_d$ and,
  if $m \geq 2$ then further
  $\E\bigl\{\tilde\cF(Y_i) \tilde \cF(Y_i)^\T\bigr\}$
  does not depend on~$i$.
  By Gaussianity, it follows that for each
  $\alpha \in \N_0^d$ with $|\alpha| \in [m]$,
  there exists $\mu_{Y,\alpha} \in \R$ with
  $\E\bigl\{\tilde\cF(Y_i)^\alpha\bigr\} = \mu_{Y,\alpha}$
  for all $i \in [n]$.

  Next, the conclusion is immediate if $\cJ(\nu) = \infty$,
  so suppose that $\cJ(\nu) < \infty$.
  Set $d(f,g) \coloneqq \|f-g\|_{\bar P_X,2}$.
  Since $N(2\varepsilon,\cF,d) \leq |\tilde\cF| N(\varepsilon,\check\cF,d)$,
  we deduce that $(\cF, d)$ is totally bounded and
  \begin{equation*}
    \int_0^\varepsilon
    \sqrt{1+\log N(\xi,\cF,d)}
    \diffi \xi
    \to 0,
  \end{equation*}
  as $\varepsilon \to 0$.
  Since
  \begin{equation*}
    \frac{1}{n}\E\bigl(\{T_n(f)-T_n(g)\}^2\bigr)
    = \frac{1}{n}
    \sum_{i=1}^n
    \Var\{f(X_i) - g(X_i)\}
    \leq d(f,g)^2,
  \end{equation*}
  for $f, g \in \cF$,
  Dudley's continuity criterion
  \citep[Corollary~2.2.8 and Problem~2.2.17]{van1996weak}
  ensures that there exists a version of
  $T_n$ with uniformly $d$-continuous paths.
  This version is a Borel random element of the Polish space
  $\cC(\cF, d)$ of uniformly $d$-continuous functions
  from $\cF$ to $\R$, equipped with the supremum norm.
  Extensions from $\tilde\cF$ to $\cF$ then proceed via the
  Vorob'ev--Berkes--Philipp theorem
  \citep[Theorem~1.31]{dudley2014uniform},
  using $\cX^n$ for the empirical variables and
  $\cC(\cF,d)$ for the Gaussian process.

  Finally, we must bound $\sup_{f \in \check \cF} |T_n(f)|$,
  where we define $T_n\bigl(f - \pi(f)\bigr)
  \coloneqq T_n(f) - T_n\bigl(\pi(f)\bigr)$.
  Since $\sup_{f \in \check \cF} \|f\|_{\bar P_X, 2} \leq \sigma$,
  by Dudley's theorem \citep[Corollary~2.2.8]{van1996weak}
  and a change of variables,
  since $\varepsilon \mapsto \cJ(\varepsilon) / \varepsilon$
  is decreasing, for a universal constant $C > 0$, we have
  \begin{align*}
    \E\biggl(
      \sup_{f \in \check \cF}
      |T_n(f)|
    \biggr)
    &\leq
    C \sqrt n
    \int_0^{\sigma}
    \sqrt{1 + \log N\bigl(\xi, \check \cF, \|\cdot\|_{\bar P_X, 2}\bigr)}
    \diffi \xi \\
    &\leq
    C \sqrt n \|\check F\|_{\bar P_X, 2}
    \int_0^{\sigma / \|\check F\|_{\bar P_X, 2}}
    \sqrt{1 + \log N\bigl(\xi \|\check F\|_{\bar P_X, 2}, \check \cF,
    \|\cdot\|_{\bar P_X, 2}\bigr)}
    \diffi \xi \\
    &\leq
    C \sqrt n \|\check F\|_{\bar P_X, 2}
    \cJ\biggl(\frac{\sigma}{\|\check F\|_{\bar P_X, 2}}\biggr)
    \leq
    \frac{C \sigma \sqrt n \cJ(\nu)}{\nu}.
  \end{align*}
  Further, by Borell's inequality
  \citep[Proposition~A.2.1]{van1996weak}, for all $t > 0$,
  increasing $C$ if necessary,
  with probability at least $1 - e^{-t^2/2}$,
  \begin{equation*}
    \sup_{g \in \check \cF}
    |T_n(g)|
    \leq
    \frac{C \sigma \sqrt n \cJ(\nu)}{\nu}
    + C t \sigma \sqrt n
    \leq C R(t).
    \iftoggle{journal}{}{\qedhere}
  \end{equation*}
\end{proof}

\begin{proof}[Proof of Corollary~\ref{cor:gaussian_vc}]
  We write $C > 0$ for a quantity depending
  only on the stated quantities; its value may vary from line to line.
  Define $\nu \coloneqq n^{(2 - k - \gamma)/(2k + 2\gamma + u k)}$.
  As in the proof of Corollary~\ref{cor:empirical_process_vc},
  choose $\tilde\cF$ so that $\check\cF$ is pointwise measurable
  with $d \leq a\nu^{-u}$ and
  $\|f-\pi(f)\|_{\bar P_X,2} \leq \sigma$ for all $f\in\cF$,
  taking $\check F \coloneqq 2F$, $\check M \coloneqq 2M$ and
  $\sigma \coloneqq \nu\|\check F\|_{\bar P_X,2}$.
  The same entropy calculations, with $\kappa \coloneqq 2k+4$, give
  $R(t) \leq C M \nu \sqrt{n \log 2n}$
  for $t \coloneqq \sqrt{2\kappa\log 2n}$.
  By Lemma~\ref{lem:moment_tail_orlicz},
  $Y_i(f)$ has variance at most $C M^2$ for each $f \in \cF$.
  Thus, Lemma~\ref{lem:gaussian_norm} and
  \eqref{eq:vc_beta_bound} yield
  \begin{equation*}
    \bar\beta_{k,\gamma}
    \leq C M^{k + \gamma} d^{k/2} (\log 2d)^{\gamma/2}
    \leq C M^{k + \gamma} \nu^{-u k/2} (\log 2n)^{\gamma/2}.
  \end{equation*}
  Take $m=1$; the moment matching assumption gives
  $\Delta(s)=\Gamma_\eta(s)=0$.
  As $n\nu^{-u k/2} = n^{(k + \gamma)/2} \nu^{k + \gamma}$ and
  $\phi(t) \leq C \sqrt{\log 2n}$,
  replacing $\eta$ with $M \eta \nu \sqrt{n \log 2n}$ in
  Theorem~\ref{thm:gaussian_process} gives, for $\eta \geq 1$,
  \begin{equation*}
    \P\biggl\{
      \sup_{f\in\cF}
      |S_n'(f)-T_n'(f)|
      > C M \eta \nu \sqrt{n \log 2n}
    \biggr\}
    \leq \frac{C}{\eta^{k + \gamma}}+\frac{2}{n^\kappa}.
  \end{equation*}
  For $1 \leq \eta \leq n^2$, the last term above is at most
  $2/\eta^{k + \gamma}$.
  To handle $\eta\geq n^2$, Dudley's theorem
  \citep[Corollary~2.2.8]{van1996weak}
  and Borell's inequality \citep[Proposition~A.2.1]{van1996weak}
  along with Lemma~\ref{lem:entropy_vc} give
  $\bigvvvert\sup_{f\in\cF} |T_n'(f)|\bigvvvert_{2k+2} \leq C M \sqrt n$,
  as in the proof of Theorem~\ref{thm:gaussian_process}.
  By Minkowski's inequality and Lemma~\ref{lem:moment_tail_orlicz},
  $\bigvvvert\sup_{f\in\cF} |S_n'(f)|\bigvvvert_{2k+2} \leq C M n$.
  Since $\nu \sqrt{n \log 2n}\geq 1/C$, Markov's inequality therefore bounds
  the same probability by $C (n/\eta)^{2k+2} \leq C / \eta^{k + \gamma}$.
  This also handles $n=1$, while $\eta < 1$ is trivial.
  Increasing $C$ gives the stated result.
\end{proof}

\begin{proof}[Proof of
  Theorem~\ref{thm:empirical_process_bootstrap}]
  We write $C > 0$ for a quantity depending only on $r$ and $r'$,
  which may vary from line to line.
  Let $\tilde S_n \coloneqq \sum_{i=1}^{n}
  \bigl(\tilde\cF(X_{i}) - P_{X,i}(\tilde\cF)\bigr)$
  and $\tilde S^*_n \coloneqq \sum_{i=1}^{n} Y_i^*(\tilde\cF)$,
  and write $\P^*$ for conditional probability given the data $\bX_n$.
  Observe that, under condition~\ref{it:thm_ep_bootstrap_nonparametric}
  and under condition~\ref{it:thm_ep_bootstrap_multiplier}
  parts (i) and (ii), the summands
  $Y_i^*(\tilde\cF)$ have conditional mean zero and are
  conditionally independent, given the data.
  Applying Lemma~\ref{lem:borel}
  with $p = \infty$ and $m = 1$
  to the law of $\tilde S_n$ and the conditional law of
  $\tilde S^*_n$ given $\bX_n$,
  then taking an outer expectation, we obtain that for $t, \eta > 0$,
  \begin{align*}
    &\E^\star \Biggl(
      \sup_{A \in \cB(\R^d)}
      \Bigl\{
        \P(\tilde S_n \in A)
        - \P^*\big( \tilde S^*_n \in A_\infty^{3\eta} \big)
      \Bigr\}
    \Biggr) \\
    &\qquad\leq
    2 e^{-t^2 / 2}
    + \E\Biggl\{
      \frac{n \hat\beta_{k,\gamma} \phi(t)^k}{\eta^{k+\gamma}}
      + \hat\Delta\biggl(\frac{\phi(t)}{\eta}\biggr)
      + \hat\Gamma_\eta\biggl(\frac{\phi(t)}{\eta}\biggr)
    \Biggr\}.
  \end{align*}
  As in the proof of Theorem~\ref{thm:bootstrap},
  under condition~\ref{it:thm_ep_bootstrap_nonparametric}
  or~\ref{it:thm_ep_bootstrap_multiplier},
  applying a conditional
  version of Strassen's theorem
  \citep[Lemma~SA.1]{cattaneo2025yurinskii},
  we deduce that there is a probability space
  supporting $\bX_n'$, $\tilde S_n'$ and $\tilde S_n^{*\prime}$, where
  $\tilde S_n' \equallaw \tilde S_n$,
  $\tilde S_n'$ is independent of $\bX_n'$,
  and $(\bX_n', \tilde S_n^{*\prime})
  \equallaw (\bX_n, \tilde S_n^*)$, with
  \begin{equation}
    \label{eq:empirical_process_bootstrap_yurinskii}
    \P\Bigl(
      \bigl\|\tilde S_n' - \tilde S_n^{*\prime}\bigr\|_\infty
      > 3\eta
    \Bigr)
    \leq
    2 e^{-t^2 / 2}
    + \E\Biggl\{
      \frac{n \hat\beta_{k,\gamma} \phi(t)^k}{\eta^{k+\gamma}}
      + \hat\Delta\biggl(\frac{\phi(t)}{\eta}\biggr)
      + \hat\Gamma_\eta\biggl(\frac{\phi(t)}{\eta}\biggr)
    \Biggr\}.
  \end{equation}

  By two applications of the Vorob'ev--Berkes--Philipp theorem
  \citep[Theorem~1.31]{dudley2014uniform},
  firstly with common coordinate $\tilde S_n'$ and secondly with common
  coordinate $(\bX_n', \tilde S_n^{*\prime})$, there exist random variables
  $\bX_n'' = (X_i'': i \in [n])$ and $\bX_n^{*\prime}$ such that
  $(\tilde S_n', \bX_n'') \equallaw (\tilde S_n, \bX_n)$,
  $(\tilde S_n', \bX_n'')$ is independent of $\bX_n'$, and
  $(\bX_n', \bX_n^{*\prime}, \tilde S_n^{*\prime}) \equallaw
  (\bX_n, \bX_n^*, \tilde S_n^*)$, with substitutions for multiplier weights
  under condition~\ref{it:thm_ep_bootstrap_multiplier} as in the
  theorem statement.
  Here, independence of $(\tilde S_n', \bX_n'')$ and $\bX_n'$ holds
  because $\tilde S_n'$ is independent of $\bX_n'$, and the first application
  makes $\bX_n''$ conditionally independent of $\bX_n'$ given $\tilde S_n'$.
  For $f \in \cF$, let $S_n'(f) \coloneqq \sum_{i=1}^{n}
  \bigl(f(X_{i}'') - P_{X,i}(f)\bigr)$, and compute $S_n^{*\prime}(f)$ from
  $(\bX_n', \bX_n^{*\prime})$ in the same fashion
  as $S_n^*(f)$ is computed from
  $(\bX_n, \bX_n^*)$. Then $\tilde S_n' = S_n'(\tilde\cF)$ and
  $\tilde S_n^{*\prime} = S_n^{*\prime}(\tilde\cF)$ almost surely.
  We have
  \begin{equation}
    \label{eq:empirical_process_bootstrap_decomp}
    \sup_{f \in \cF}
    \bigl|S_n'(f) - S_n^{*\prime}(f)\bigr|
    \leq
    \bigl\|\tilde S_n' - \tilde S_n^{*\prime}\bigr\|_\infty
    + \sup_{f \in \check\cF}
    |S_n'(f)|
    + \sup_{f \in \check\cF}
    |S_n^{*\prime}(f)|,
  \end{equation}
  observing that the first term on the right-hand side of
  \eqref{eq:empirical_process_bootstrap_decomp}
  has already been bounded (in probability)
  in \eqref{eq:empirical_process_bootstrap_yurinskii}.
  We may assume $n \geq 2$, since the result is trivial for $n = 1$.
  For the second and third terms, by
  Theorem~\ref{thm:empirical_tail}
  and Lemma~\ref{lem:bootstrap_tail} replacing $\kappa$ with $2 \kappa$
  in each, then taking a union bound, we have
  \begin{align*}
    \sup_{f \in \check\cF} |S_n'(f)|
    + \sup_{f \in \check\cF} |S_n^{*\prime}(f)|
    &\leq C
    \biggl\{
      \sigma \sqrt n
      + \check M (\kappa \log 2n)^{1/r + 1/r'}
      \biggl(
        \frac{\cJ(\nu)}{\nu} + t
      \biggr)
    \biggr\}
    \biggl( \frac{\cJ(\nu)}{\nu} + t \biggr) \\
    &= C R^*(t),
  \end{align*}
  with probability at least
  $1 - 2/n^{2\kappa} - 2 e^{-t^2/2} \geq 1 - 1/n^\kappa - 2 e^{-t^2/2}$.
  By another union bound, we conclude that
  \begin{align*}
    &
    \P\biggl\{
      \sup_{f \in \cF}
      \bigl|S_n'(f) - S_n^{*\prime}(f)\bigr|
      > 3 \eta
      + C R^*(t)
    \biggr\} \\
    &\qquad\leq
    4 e^{-t^2 / 2}
    + \E\Biggl\{
      \frac{n \hat\beta_{k,\gamma} \phi(t)^k}{\eta^{k+\gamma}}
      + \hat\Delta\biggl(\frac{\phi(t)}{\eta}\biggr)
      + \hat\Gamma_\eta\biggl(\frac{\phi(t)}{\eta}\biggr)
    \Biggr\}
    + \frac{1}{n^\kappa},
  \end{align*}
  as required.
\end{proof}

\begin{proof}[Proof of Lemma~\ref{lem:multiplier_weights}]
  We consider each case in turn.
  \begin{enumerate}

    \item Rademacher.
      Clearly, $\E(\chi) = 0$, $\E(\chi^2) = 1$,
      $\vvvert \chi \vvvert_4 = 1 \leq 2$ and
      $\vvvert \chi \vvvert_\infty = 1 \leq 2$.

    \item Gaussian.
      Clearly, $\E(\chi) = 0$ and $\E(\chi^2) = 1$.
      Further, $\E(\chi^4) = 3$ so
      $\vvvert \chi \vvvert_4 = 3^{1/4} \leq 1.32 \leq 2$.
      Finally, $\vvvert \chi \vvvert_{\omega_2} = \sqrt{8/3} \leq 1.64 \leq 2$
      by \citet{leskela2026sharp}.

    \item Centered exponential.
      Clearly, $\E(\chi) = 0$ and $\E(\chi^2) = 1$.
      We have $\E(\chi^4) = \int_0^\infty (x - 1)^4 e^{-x} \diffi x = 9$,
      $\vvvert \chi \vvvert_4 = 9^{1/4} \leq 1.74 \leq 2$.
      Moreover, for $t > 1$,
      \begin{align*}
        \E\bigl\{\exp(|\chi| / t) - 1\bigr\}
        &=
        \int_0^\infty
        \bigl\{
          \exp(|x - 1| / t) - 1
        \bigr\}
        e^{-x}
        \diffi x \\
        &= \int_0^\infty
        \exp(|x - 1| / t - x)
        \diffi x
        - 1 \\
        &=
        \int_0^1
        \exp\{(1 - x) / t - x\}
        \diffi x
        + \int_1^\infty
        \exp\{(x - 1) / t - x\}
        \diffi x
        - 1 \\
        &=
        \frac{t \bigl(e^{1+1/t} - 1\bigr)}{e (t+1)}
        + \frac{t}{e (t-1)}
        - 1.
      \end{align*}
      Setting $t \coloneqq 2$ yields
      $\E\bigl\{\exp(|\chi| / t) - 1\bigr\}
      \leq 0.59 \leq 1$,
      so $\vvvert \chi \vvvert_{\omega_1} \leq 2$.

    \item Mammen.
      By calculating the expectations,
      $\E(\chi) = 0$, $\E(\chi^2) = 1$ and $\E(\chi^3) = 1$.
      Further, $\E(\chi^4) = 2$ and
      $\vvvert \chi \vvvert_\infty = 1/2 + \sqrt{5}/2 \leq 1.62 \leq 2$.

    \item Centered Poisson.
      Let $N \sim \Poisson(1)$ and $\chi = N-1$.
      Then $\E(N)=1$, $\E(N^2)=2$, $\E(N^3)=5$ and $\E(N^4)=15$.
      Thus, $\E(\chi)=0$, $\E(\chi^2)=\E(\chi^3)=1$ and
      $\E(\chi^4)=4$, so $\vvvert \chi \vvvert_4=\sqrt{2}\leq2$.
      Further, since $|N-1|=N-1$ except when $N=0$,
      \begin{align*}
        \E\bigl\{\exp(|\chi|/2)\bigr\}
        &=
        \E\biggl\{\exp\biggl(\frac{N-1}{2}\biggr)\biggr\}
        + \P(N=0)\bigl(e^{1/2}-e^{-1/2}\bigr) \\
        &= \sum_{n=0}^{\infty}
        \exp\biggl(\frac{n-1}{2}\biggr)
        \frac{e^{-1}}{n!}
        + e^{-1}\bigl(e^{1/2}-e^{-1/2}\bigr) \\
        &= \exp\bigl(e^{1/2}-3/2\bigr)
        + e^{-1/2} - e^{-3/2}
        < 1.55 < 2.
      \end{align*}
      Hence, $\vvvert \chi \vvvert_{\omega_1}\leq2$.
  \end{enumerate}
  This establishes the claims.
\end{proof}

\begin{proof}[Proof of Theorem~\ref{thm:empirical_tail}]
  We write $C > 0$ for a quantity depending only on $r$;
  the value of $C$ may vary from line to line.
  Without loss of generality, since replacing $\sigma$ with
  $\sigma \land \|F\|_{\bar P_X, 2}$ preserves the assumptions
  and decreases the claimed bound, we assume that
  $\sigma \leq \|F\|_{\bar P_X, 2}$.
  Define $\tau \coloneqq M (2 \kappa \log 2n)^{1/r}$
  and let $\cF_\tau \coloneqq
  \bigl\{x \mapsto f(x) \I\{F(x) \leq \tau\}: f \in \cF\bigr\}$.
  Observe that $\cF_\tau$ is pointwise measurable
  because if $(f_m)$ is a sequence in $\cF'$
  with $f_m(x) \to f(x)$ for each $x \in \cX$,
  then with $f_{m,\tau} \coloneqq f_m \I\{F \leq \tau\}$,
  we have that $(f_{m,\tau})$ forms a sequence in
  the countable class $\cF'_\tau \coloneqq
  \bigl\{x \mapsto f(x) \I\{F(x) \leq \tau\}: f \in \cF'\bigr\}$
  with $f_{m,\tau}(x) \to f_\tau(x) \coloneqq f(x) \I\{F(x) \leq \tau\}$
  for each $x \in \cX$. By the triangle inequality, write
  \begin{align}
    \sup_{f \in \cF} |S_n(f)|
    &\leq
    \sup_{f \in \cF_\tau}
    |S_n(f)|
    + \sup_{f \in \cF}
    \Biggl|
    \sum_{i=1}^n
    f(X_i) \I\{F(X_i) > \tau\}
    \Biggr|
    \iftoggle{journal}{\nonumber\\ &\qquad}{}
    + \sup_{f \in \cF}
    \Biggl|
    \sum_{i=1}^n
    P_{X,i}\bigl(f \I\{F > \tau\}\bigr)
    \Biggr|.
    \label{eq:tail_bound_truncation}
  \end{align}
  Let $\varepsilon_1, \ldots, \varepsilon_n$ be independent Rademacher
  variables, independent of the data.
  For the first term on the right-hand side of \eqref{eq:tail_bound_truncation},
  by a symmetrization inequality \citep[Lemma~2.3.6]{van1996weak}
  and a contraction principle \citep[Proposition A.1.10]{van1996weak},
  \begin{align*}
    \E \biggl( \sup_{f \in \cF_\tau} |S_n(f)| \biggr)
    &=
    \E \Biggl(
      \sup_{f \in \cF_\tau}
      \Biggl|
      \sum_{i=1}^n
      \bigl\{
        f(X_i) - P_{X,i}(f)
      \bigr\}
      \Biggr|
    \Biggr) \\
    &\leq
    2 \, \E \Biggl(
      \sup_{f \in \cF_\tau}
      \Biggl|
      \sum_{i=1}^n
      \varepsilon_i f(X_i)
      \Biggr|
    \Biggr) \\
    &=
    2 \, \E \Biggl(
      \sup_{f \in \cF}
      \Biggl|
      \sum_{i=1}^n
      \varepsilon_i f(X_i)
      \I\bigl\{F(X_i) \leq \tau\bigr\}
      \Biggr|
    \Biggr) \\
    &\leq
    2 \, \E \Biggl(
      \sup_{f \in \cF}
      \Biggl|
      \sum_{i=1}^n
      \varepsilon_i f(X_i)
      \Biggr|
    \Biggr).
  \end{align*}
  Let $g_f(x, \varepsilon) \coloneqq \varepsilon f(x)$ and
  note that $(x, \varepsilon) \mapsto F(x)$ forms an envelope
  for $\bigl\{(x, \varepsilon) \mapsto g_f(x, \varepsilon): f \in \cF\bigr\}$.
  Further, for any finitely supported probability measure $Q$ on
  $\cX \times \{-1, 1\}$, we have
  $\|g_{f_1} - g_{f_2}\|_{Q,2} \leq \|f_1 - f_2\|_{Q_\cX, 2}$,
  where $Q_\cX$ is the marginalization of $Q$ onto $\cX$.
  Thus, writing $M_n \coloneqq \max_{1 \leq i \leq n} F(X_i)$,
  by Lemma~SA24 in \citet[supplementary material]{cattaneo2024uniform}
  applied with the strictly positive envelope
  $F + \delta$, then taking $\delta \to 0$,
  since $\varepsilon \mapsto \cJ(\varepsilon) / \varepsilon$ is decreasing,
  and using Lemma~\ref{lem:moment_tail_orlicz},
  \begin{align*}
    \E \biggl( \sup_{f \in \cF_\tau} |S_n(f)| \biggr)
    &\leq
    C \|F\|_{\bar P_X,2}
    \cJ\biggl(\frac{\sigma}{\|F\|_{\bar P_X,2}}\biggr)
    \sqrt n
    + \frac{C \vvvert M_n \vvvert_2 \|F\|_{\bar P_X,2}^2}{\sigma^2}
    \cJ\biggl(\frac{\sigma}{\|F\|_{\bar P_X,2}}\biggr)^2 \\
    &\leq
    \frac{C \sigma \cJ(\nu) \sqrt n}{\nu}
    + \frac{C M (\log 2n)^{1/r} \cJ(\nu)^2}{\nu^2}.
  \end{align*}

  We now control the fluctuations of $\sup_{f \in \cF_\tau} |S_n(f)|$
  above its expected value.
  For $f \in \cF_\tau$ and $i \in [n]$, define the random variable
  \begin{equation*}
    U_i(f)
    \coloneqq
    \frac{1}{2 \tau}
    \bigl\{
      f(X_i) - P_{X,i}(f)
    \bigr\},
  \end{equation*}
  so that $S_n(f) = 2 \tau \sum_{i=1}^n U_i(f)$.
  Further, $\E\{U_i(f)\} = 0$ and $|U_i(f)| \leq 1$ almost surely.
  Define
  \begin{equation*}
    \tilde\sigma^2
    \coloneqq
    \sup_{f \in \cF_\tau}
    \sum_{i=1}^n \E\bigl\{U_i(f)^2\bigr\}
    \leq
    \frac{1}{4 \tau^2}
    \sup_{f \in \cF_\tau} \sum_{i=1}^n
    \|f\|_{P_{X,i},2}^2
    \leq
    \frac{n \sigma^2}{4 \tau^2}.
  \end{equation*}
  By \citet[Theorem~11.8]{boucheron2013concentration}
  applied to $\cF_\tau \cup (- \cF_\tau)$,
  increasing $C$ if necessary, define
  \begin{align*}
    \Sigma^2
    &\coloneqq
    \E\biggl(
      \sup_{f \in \cF_\tau}
      \sum_{i=1}^n
      U_i(f)^2
    \biggr)
    \leq
    8 \, \E\biggl(
      \sup_{f \in \cF_\tau}
      \biggl|
      \sum_{i=1}^n
      U_i(f)
      \biggr|
    \biggr)
    + \tilde\sigma^2 \\
    &\,=
    \frac{4}{\tau} \,
    \E\biggl(
      \sup_{f \in \cF_\tau}
      |S_n(f)|
    \biggr)
    + \tilde\sigma^2
    \leq
    \frac{C \sqrt n \sigma \cJ(\nu)}{\tau \nu}
    + \frac{C \cJ(\nu)^2}{\kappa^{1/r} \nu^2}
    + \frac{C n \sigma^2}{\tau^2}.
  \end{align*}
  Hence, applying \citet[Theorem~12.2]{boucheron2013concentration}
  to $\cF_\tau \cup (- \cF_\tau)$,
  for $t > 0$, increasing $C$ again,
  \begin{align*}
    &\P \Biggl(
      \sup_{f \in \cF_\tau}
      \biggl|
      \sum_{i=1}^n
      U_i(f)
      \biggr|
      \geq
      \E \biggl(
        \sup_{f \in \cF_\tau}
        \biggl|
        \sum_{i=1}^n
        U_i(f)
        \biggr|
      \biggr)
      + t
    \Biggr) \\
    &\qquad\leq
    \exp \biggl(
      - \frac{t^2}{
        C \cJ(\nu)^2 \nu^{-2} \kappa^{-1/r}
        + C n^{1/2} \sigma \cJ(\nu) \tau^{-1} \nu^{-1}
        + C n \sigma^2 \tau^{-2}
        + C t
      }
    \biggr).
  \end{align*}
  Therefore, with probability at least $1 - e^{-t^2 / 2}$,
  \begin{align*}
    \sup_{f \in \cF_\tau}
    \bigl| S_n(f) \bigr|
    &\leq
    \E \biggl\{
      \sup_{f \in \cF_\tau}
      \bigl| S_n(f) \bigr|
    \biggr\}
    + \frac{C \tau \cJ(\nu) t}{\nu \kappa^{1/(2 r)}}
    + \frac{C n^{1/4} \sqrt {\sigma \tau \cJ(\nu)} t}{\sqrt{\nu}}
    + C \sigma \sqrt n t
    + C \tau t^2 \\
    &\leq
    C \Biggl(
      \frac{\sigma \cJ(\nu) \sqrt n}{\nu}
      + \frac{\tau \cJ(\nu)^2}{\kappa^{1/r} \nu^2}
      + \sqrt{
        \frac{\tau \cJ(\nu)^2}{\kappa^{1/r} \nu^2}
        \cdot \tau t^2
      }
      \iftoggle{journal}{\\ &\qquad\qquad}{}
      + \sqrt{
        \frac{\sigma \cJ(\nu) \sqrt n}{\nu}
        \cdot \tau t^2
      }
      + \sigma \sqrt n t
      + \tau t^2
    \Biggr) \\
    &\leq
    C \biggl\{
      \sigma \sqrt n
      + \tau
      \biggl(
        \frac{\cJ(\nu)}{\kappa^{1/r} \nu}
        + t
      \biggr)
    \biggr\}
    \biggl(\frac{\cJ(\nu)}{\nu} + t\biggr) \\
    &\leq
    C \biggl\{
      \sigma \sqrt n
      + M (\log 2n)^{1/r}
      \biggl(
        \frac{\cJ(\nu)}{\nu}
        + \kappa^{1/r} t
      \biggr)
    \biggr\}
    \biggl(\frac{\cJ(\nu)}{\nu} + t\biggr).
  \end{align*}

  If $r = \infty$ then the second and third terms in
  \eqref{eq:tail_bound_truncation} vanish and the proof is complete.
  Otherwise, for the second term in \eqref{eq:tail_bound_truncation},
  by Lemma~\ref{lem:moment_tail_orlicz},
  \begin{align*}
    \P\Biggl(
      \sup_{f \in \cF}
      \biggl|
      \sum_{i=1}^n
      f(X_i) \I\{F(X_i) > \tau\}
      \biggr|
      > 0
    \Biggr)
    &\leq
    \P(M_n > \tau)
    \leq
    2 n \exp\biggl(-\frac{\tau^r}{M^r}\biggr) \\
    &= 2 n e^{- 2 \kappa \log 2 n}
    = \frac{1}{(2 n)^{2\kappa -1}}
    \leq \frac{1}{n^\kappa}.
  \end{align*}
  For the third term in \eqref{eq:tail_bound_truncation},
  by Cauchy--Schwarz and Lemma~\ref{lem:moment_tail_orlicz},
  since $\kappa \geq 1$,
  \begin{align*}
    \sup_{f \in \cF}
    \Biggl|
    \sum_{i=1}^n P_{X,i}\bigl(f \I\{F > \tau\}\bigr)
    \Biggr|
    &\leq
    \sum_{i=1}^n
    P_{X,i}\bigl(F \I\{F > \tau\}\bigr)
    \leq
    \sum_{i=1}^n
    \vvvert F(X_i) \vvvert_2
    \P\bigl(F(X_i) > \tau\bigr)^{1/2} \\
    &\leq
    C n M \exp\biggl(-\frac{\tau^r}{2M^r}\biggr)
    = C n M (2n)^{-\kappa}
    \leq C M.
  \end{align*}
  We conclude that for all $n \geq 1$,
  with probability at least
  $1 - e^{-t^2/2} - 1/n^\kappa$,
  \begin{equation*}
    \sup_{f \in \cF}
    \bigl| S_n(f) \bigr|
    \leq
    C \biggl\{
      \sigma \sqrt n
      + M (\log 2n)^{1/r}
      \biggl(
        \frac{\cJ(\nu)}{\nu}
        + \kappa^{1/r} t
      \biggr)
    \biggr\}
    \biggl(\frac{\cJ(\nu)}{\nu} + t\biggr).
    \iftoggle{journal}{}{\qedhere}
  \end{equation*}
\end{proof}

\begin{proof}[Proof of Lemma~\ref{lem:empirical_norm}]
  Write $C > 0$ for a quantity depending only on $r$,
  which may vary from line to line.
  Define $\tau \coloneqq M (4 \kappa \log 2n)^{1/r}$
  and let $\cF_{\tau}$ be as in the proof of
  Theorem~\ref{thm:empirical_tail}.
  For $f \in \cF_{\tau}$,
  \begin{equation*}
    \bigl\|f^2\bigr\|_{\bar P_X,2}^2
    \leq
    \|f\|_{\infty}^2
    \|f\|_{\bar P_X,2}^2
    \leq  \tau^2 \sigma^2.
  \end{equation*}
  The function $\tau F$
  is an envelope for
  $\cF^2_{\tau} \coloneqq
  \bigl\{x \mapsto f(x)^2: f \in \cF_{\tau}\bigr\}$, and
  $\vvvert \tau F(X_i) \vvvert_{\omega_{r}} \leq \tau M$
  for each $i \in [n]$.
  Further, $\cF^2_{\tau}$ is pointwise measurable.
  Note that
  $\bigl\|\tau F\bigr\|_{\bar P_X, 2}^2
  =  \tau^2 \|F\|_{\bar P_X,2}^2$
  and $\tau \sigma /
  \bigl(\tau \|F\|_{\bar P_X,2}\bigr) \geq \nu$.
  For $f_1, f_2 \in \cF_{\tau}$ and $Q$ a
  probability measure on $(\cX, \cH)$ we have
  \begin{equation*}
    \bigl\|f_1^2 - f_2^2\bigr\|_{Q,2}
    \leq
    \|f_1 + f_2\|_{\infty}
    \|f_1 - f_2\|_{Q,2}
    \leq
    2  \tau
    \|f_1 - f_2\|_{Q,2}.
  \end{equation*}
  It follows that
  $N\bigl(2 \tau \varepsilon, \cF_{\tau}^2, \|\cdot\|_{Q,2}\bigr)
  \leq N\bigl(\varepsilon, \cF_{\tau}, \|\cdot\|_{Q,2}\bigr)
  \leq N\bigl(\varepsilon, \cF, \|\cdot\|_{Q,2}\bigr)$
  for all $\varepsilon > 0$. So, for $\xi \in (0, 1]$, we have
  $N\bigl(\xi \| \tau F\|_{Q,2},
  \cF_{\tau}^2, \|\cdot\|_{Q,2}\bigr)
  \leq N\bigl(\xi \|F\|_{Q,2} / 2, \cF, \|\cdot\|_{Q,2}\bigr)$.
  Thus, by applying Theorem~\ref{thm:empirical_tail}
  to $\cF_{\tau}^2$ with envelope $\tau F$,
  for all $t > 0$, with probability at least
  $1 - 1/n^{2\kappa} - e^{-t^2/2}$, we have
  \begin{equation*}
    \sup_{f \in \cF_{ \tau}}
    \Bigl\{
      \|f\|_{P_n,2}^2
      - \|f\|_{\bar P_X,2}^2
    \Bigr\}
    \leq
    C  \tau
    \biggl\{
      \frac{\sigma}{\sqrt n}
      + \frac{M (\log 2n)^{1/r}}{n}
      \biggl(
        \frac{\cJ(\nu)}{\nu} + \kappa^{1/r} t
      \biggr)
    \biggr\}
    \biggl( \frac{\cJ(\nu)}{\nu} + t \biggr).
  \end{equation*}

  Next, we handle the truncation errors.
  If $r < \infty$, by Lemma~\ref{lem:moment_tail_orlicz},
  \begin{equation*}
    \P\Bigl(
      \max_{1 \leq i \leq n} F(X_i)
      > \tau
    \Bigr)
    \leq
    2 n \exp\biggl(-\frac{\tau^r}{M^r}\biggr)
    \leq \frac{1}{n^{2\kappa}}.
  \end{equation*}
  If $r = \infty$, the left-hand side above is zero.
  Noting that the lemma is vacuous if $n = 1$,
  suppose $n \geq 2$.
  Increasing $C$ if necessary,
  as $2/n^{2\kappa} \leq 1/n^\kappa$,
  with probability at least $1 - 1/n^\kappa - e^{-t^2/2}$,
  \begin{align*}
    \sup_{f \in \cF}
    \|f\|_{P_n,2}^2
    &\leq
    \sigma^2
    + C \tau
    \biggl\{
      \frac{\sigma}{\sqrt n}
      + \frac{M (\log 2n)^{1/r}}{n}
      \biggl(
        \frac{\cJ(\nu)}{\nu} + \kappa^{1/r} t
      \biggr)
    \biggr\}
    \biggl( \frac{\cJ(\nu)}{\nu} + t \biggr) \\
    &\leq
    2 \sigma^2
    + \frac{C \tau^2}{n}
    \biggl( \frac{\cJ(\nu)}{\nu} + t \biggr)^2
    \leq
    2 \sigma^2
    + \frac{C M^2 (\kappa \log 2n)^{2/r}}{n}
    \biggl( \frac{\cJ(\nu)}{\nu} + t \biggr)^2.
    \iftoggle{journal}{}{\qedhere}
  \end{align*}
\end{proof}

\begin{proof}[Proof of Lemma~\ref{lem:bootstrap_tail}]
  Under condition \ref{it:bootstrap_tail_nonparametric},
  let $T_n^*(f) \coloneqq S_n^*(f)$,
  while under condition \ref{it:bootstrap_tail_multiplier},
  define $T_n^*(f) \coloneqq \sum_{i=1}^n \chi_i f(X_i)$.
  For $t > 0$ and with $C > 0$
  a quantity to be determined which depends only on $r$ and $r'$,
  define the events
  \begin{align*}
    \cA_{1,n,t}
    &\coloneqq
    \Biggl\{
      \sup_{f \in \cF}
      \|f\|_{P_n, 2}
      \leq
      \sigma'
      \coloneqq
      2 \sigma
      + \frac{C  M (\kappa \log 2n)^{1/r}}{\sqrt n}
      \biggl(
        \frac{\cJ(\nu)}{\nu} + t
      \biggr)
    \Biggr\}, \\
    \cA_{2,n,t}
    &\coloneqq
    \biggl\{
      \max_{1 \leq i \leq n}
      F(X_{i})
      \leq
      M'
      \coloneqq
      M (5 \kappa \log 2n)^{1/r}
    \biggr\}, \\
    \cA_{3,n,t}
    &\coloneqq
    \biggl\{
      \| F\|_{P_n,2}
      \leq
      2 \| F\|_{\bar P_X,2}
      + \frac{C  M (\kappa \log 2n)^{1/r} (1+t)}{\sqrt n}
    \biggr\}.
  \end{align*}
  By Lemma~\ref{lem:empirical_norm}, with $C$ sufficiently large,
  we have $\P(\cA_{1,n,t}) \geq 1 - 1/n^{4\kappa} - e^{-t^2/2}$.
  If $r < \infty$ then $\P(\cA_{2,n,t}) \geq 1 - 1/n^{4\kappa}$
  by Lemma~\ref{lem:moment_tail_orlicz},
  while $r = \infty$ implies $\P(\cA_{2,n,t}) = 1$.
  Applying Lemma~\ref{lem:empirical_norm} to the singleton
  function class $\{ F\}$, we have
  $\P(\cA_{3,n,t}) \geq 1 - 1/n^{4\kappa} - e^{-t^2/2}$.
  Defining $\cA_{n,t} \coloneqq \cA_{1,n,t} \cap \cA_{2,n,t}
  \cap \cA_{3,n,t}$, we have by a union bound that
  $\P(\cA_{n,t}) \geq 1 - 3/n^{4\kappa} - 2 e^{-t^2/2}$.
  We now bound the process $T_n^*(f)$ on
  $\cA_{n,t}$, under either bootstrap scheme.
  Under condition~\ref{it:bootstrap_tail_nonparametric},
  define $\cF' \coloneqq \cF$
  and let $F'(x) \coloneqq F(x)$ for $x \in \cX$.

  Under condition~\ref{it:bootstrap_tail_multiplier},
  let $\cF' \coloneqq \bigl\{(x, \xi) \mapsto \xi f(x): f \in \cF\bigr\}$
  and take $F'(x, \xi) \coloneqq |\xi| F(x)$ for
  $x \in \cX$ and $\xi \in \R$.
  For a probability measure $Q'$ on $\cX \times \R$,
  set $a \coloneqq \int \xi^2 \diffi Q'(x,\xi)$
  and suppose for now that $0 < a < \infty$.
  For $A \in \cH$, define
  $Q(A) \coloneqq \int \xi^2 \I(A)(x) \diffi Q'(x, \xi) / a$.
  It follows that
  $N\bigl(\varepsilon\|F'\|_{Q',2},\cF',\|\cdot\|_{Q',2}\bigr)
  = N\bigl(\varepsilon\|F\|_{Q,2},\cF,\|\cdot\|_{Q,2}\bigr)$.
  Define
  \begin{equation*}
    \cJ'(\varepsilon)
    \coloneqq
    \int_0^\varepsilon
    \sup_{Q' \in \cQ'}
    \sqrt{1 + \log N\bigl(
        t \|F'\|_{Q',2}, \cF', \|\cdot\|_{Q',2}
    \bigr)}
    \diffi t,
  \end{equation*}
  where $\cQ'$ is the set of all
  probability measures $Q'$ on $\cX \times \R$
  with $\|F'\|_{Q',2} > 0$.
  If $a = \infty$ and $\|F'\|_{Q',2} \in (0, \infty)$,
  fix $0 < \varepsilon < \varepsilon'$ and take $L > 0$ so that
  \begin{align*}
    \|F' \I\{|\xi| \leq L\}\|_{Q',2}
    &> 0,
    &2 \|F' \I\{|\xi| > L\}\|_{Q',2}
    &\leq
    (\varepsilon'-\varepsilon)\|F'\|_{Q',2}.
  \end{align*}
  For $A \in \cH$, define
  \begin{equation*}
    Q_L(A)
    \coloneqq
    \frac{\int \xi^2 \I(A)(x) \I\{|\xi| \leq L\}
    \diffi Q'(x,\xi)}
    {\int \xi^2 \I\{|\xi| \leq L\} \diffi Q'(x,\xi)}.
  \end{equation*}
  Let $f_1,\ldots,f_N$ be an
  $\varepsilon\|F\|_{Q_L,2}$-cover of $\cF$
  under $\|\cdot\|_{Q_L,2}$.
  For each $f \in \cF$, choose $j \in [N]$ such that
  $\|f-f_j\|_{Q_L,2} \leq \varepsilon\|F\|_{Q_L,2}$.
  By the definition of $Q_L$ and since $|f-f_j| \leq 2F$,
  writing $f'(x, \xi) \coloneqq \xi f(x)$,
  \begin{align*}
    \|f' - f_j'\|_{Q',2}
    &\leq
    \bigl\|\xi(f - f_j) \I\{|\xi| \leq L\}\bigr\|_{Q',2}
    + \bigl\|\xi(f - f_j) \I\{|\xi| > L\}\bigr\|_{Q',2} \\
    &\leq
    \bigl\|\xi \I\{|\xi| \leq L\}\bigr\|_{Q',2}
    \|f - f_j\|_{Q_L,2}
    + 2\|F' \I\{|\xi| > L\}\|_{Q',2} \\
    &\leq
    \varepsilon \bigl\|F' \I\{|\xi| \leq L\}\bigr\|_{Q',2}
    + 2\|F' \I\{|\xi| > L\}\|_{Q',2}
    \leq \varepsilon' \|F'\|_{Q',2}.
  \end{align*}
  Thus, $\{f_1',\ldots,f_N'\}$ forms an
  $\varepsilon'\|F'\|_{Q',2}$-cover of $\cF'$
  under $\|\cdot\|_{Q',2}$, so
  \begin{equation*}
    N\bigl(\varepsilon'\|F'\|_{Q',2},
    \cF', \|\cdot\|_{Q',2}\bigr)
    \leq
    \sup_{Q \in \cQ}
    N\bigl(\varepsilon\|F\|_{Q,2},
    \cF, \|\cdot\|_{Q,2}\bigr).
  \end{equation*}
  Setting $\varepsilon' \coloneqq \varepsilon/c$
  for $c \in (0,1)$ and integrating
  gives $\cJ'(\varepsilon/c) \leq \cJ(\varepsilon)/c$,
  and so $\cJ'(\varepsilon) \leq \cJ(\varepsilon)/c$
  since $\cJ'$ is increasing.
  Taking $c \to 1$ yields $\cJ'(\varepsilon) \leq \cJ(\varepsilon)$.

  Now, on the event $\cA_{n,t}$ and
  conditional on $\bX_n$, we have that
  $T_n^*(f)$ satisfies Assumption~\ref{ass:general_process}
  for $f \in \cF$ with envelope $F'$,
  Orlicz exponent $r'$,
  Orlicz bound $2 M'$,
  variance bound $\sigma^{\prime 2} \land \| F\|_{P_n,2}^2$
  and ratio parameter
  $ \nu' \coloneqq
  \bigl( \sigma' / \| F\|_{P_n, 2}\bigr) \land 1$.
  We may assume that $\| F\|_{P_n, 2} > 0$
  as otherwise $\sup_{f \in \cF} |T_n^*(f)| = 0$.
  Since $\cJ(\varepsilon) / \varepsilon \geq 1$ and,
  for $a, b > 0$ and $x \geq 0$,
  $\{(a+x)/(b+x)\} \land 1 \geq (a/b) \land 1$,
  we have that on $\cA_{n,t}$,
  \begin{equation*}
    \nu'
    =
    \frac{ \sigma'}{\| F\|_{P_n, 2}}
    \land 1
    \geq
    \frac{2 \sigma + C  M
    (\kappa \log 2n)^{1/r} (1 + t) / \sqrt n}
    {2 \| F\|_{\bar P_X, 2} + C  M
    (\kappa \log 2n)^{1/r} (1 + t) / \sqrt n}
    \land 1
    \geq
    \frac{\sigma}{\| F\|_{\bar P_X, 2}}
    \land 1
    \geq \nu.
  \end{equation*}
  Since $\varepsilon \mapsto \cJ(\varepsilon) / \varepsilon$
  is decreasing, by Theorem~\ref{thm:empirical_tail}
  with $\kappa$ replaced by $4\kappa$, we have
  \begin{align*}
    \sup_{f \in \cF}
    |T_n^*(f)|
    &\leq
    C \biggl\{
      \sigma' \sqrt n +
      M' (4 \kappa \log 2n)^{1/r'}
      \biggl(
        \frac{\cJ(\nu')}{\nu'} + t
      \biggr)
    \biggr\}
    \biggl( \frac{\cJ(\nu')}{\nu'} + t \biggr) \\
    &\leq
    C \biggl\{
      \sigma
      \sqrt n
      + M (\kappa \log 2n)^{1/r + 1/r'}
      \biggl(
        \frac{\cJ(\nu)}{\nu} + t
      \biggr)
    \biggr\}
    \biggl( \frac{\cJ(\nu)}{\nu} + t \biggr),
  \end{align*}
  with probability at least $1 - e^{-t^2/2} - 1/n^{4\kappa}$
  conditional on $\bX_n$, on $\cA_{n,t}$.
  Thus, this inequality holds marginally with probability at least
  $0 \lor (1 - 4/n^{4\kappa} - 3 e^{-t^2/2})
  \geq 1 - 1/n^\kappa - 3 e^{-t^2/2}$.
  Replacing $t$ with $t + 2$ and increasing $C$
  concludes the proof under condition
  \ref{it:bootstrap_tail_nonparametric}
  and under condition
  \ref{it:bootstrap_tail_multiplier} part (ii),
  using the fact that $\cJ(\nu) / \nu \geq 1$.

  Under condition
  \ref{it:bootstrap_tail_multiplier} part (i), it remains to bound
  $\sup_{f \in \cF} |S_n^*(f) - T_n^*(f)|$.
  Applying Theorem~\ref{thm:empirical_tail} conditionally on
  $\bX_n$ again, this time to the singleton identity class,
  note that
  \begin{align*}
    \sup_{f \in \cF}
    \bigl|
    S_n^*(f)
    - T_n^*(f)
    \bigr|
    &=
    \sup_{f \in \cF}
    \Biggl|
    \sum_{i=1}^n
    \chi_i P_n(f)
    \Biggr|
    =
    \Biggl|
    \sum_{i=1}^n
    \chi_i
    \Biggr|
    \sup_{f \in \cF} |P_n(f)| \\
    &\leq
    C
    \sup_{f \in \cF} |P_n(f)|
    \Bigl\{
      \sqrt n
      + (4 \kappa \log 2n)^{1/r'}
      (1 + t)
    \Bigr\}
    (1 + t),
  \end{align*}
  with probability at least $1 - e^{-t^2/2} - 1/n^{4\kappa}$,
  conditional on $\bX_n$.
  On $\cA_{n,t}$, Cauchy--Schwarz yields
  $\sup_{f \in \cF} |P_n(f)|
  \leq \sup_{f \in \cF} \|f\|_{P_n,2} \leq \sigma'$,
  and we also have
  $\sup_{f \in \cF} |P_n(f)| \leq \max_{1 \leq i \leq n} F(X_i) \leq M'$.
  Therefore, on $\cA_{n,t}$,
  \begin{equation*}
    \sup_{f \in \cF}
    \bigl|
    S_n^*(f)
    - T_n^*(f)
    \bigr|
    \leq
    C
    \Bigl\{
      \sigma' \sqrt n
      + M' (\kappa \log 2n)^{1/r'} (1 + t)
    \Bigr\}
    (1 + t),
  \end{equation*}
  with probability at least $1 - e^{-t^2/2} - 1/n^{4\kappa}$,
  conditional on $\bX_n$.
  Therefore, as $\cJ(\nu)/\nu \geq 1$,
  \begin{equation*}
    \sup_{f \in \cF}
    |S_n^*(f)|
    \leq
    C \biggl\{
      \sigma
      \sqrt n
      + M (\kappa \log 2n)^{1/r + 1/r'}
      \biggl(
        \frac{\cJ(\nu)}{\nu} + t
      \biggr)
    \biggr\}
    \biggl( \frac{\cJ(\nu)}{\nu} + t \biggr),
  \end{equation*}
  with probability at least
  $1 - 4 e^{-t^2/2} - 5/n^{4\kappa}$.
  Replace $t$ with $t + 2$ and
  increase $C$ again to conclude.
  If $r = r' = \infty$ and $n \geq 2$,
  letting $\kappa \to \infty$ removes $1/n^\kappa$.
  For $n=1$, the result follows as
  $\sup_{f\in\cF} |S_1^*(f)| \leq 2M$, increasing $C$ if necessary.
\end{proof}

\begin{proof}[Proof of Lemma~\ref{lem:anti-concentration_gp}]
  We combine ideas from
  \citet[Theorem~2.1 and Corollary~2.1]{chernozhukov2014anti}
  and \citet[Corollary~7]{belloni2024anti}.
  Note that $\mu_{t_n(t)} \to \mu_t$
  and $\sigma_{t_n(t)} \to \sigma_t$ for each $t \in T$
  because convergence of Gaussian variables in law implies
  convergence of the associated means and variances.
  Therefore, $\sup_{t \in T} |X_t - \mu_t| / \sigma_t
  = \sup_{t \in T_0} |X_t - \mu_t| / \sigma_t$ almost surely.
  Finiteness of $E$ then follows from the
  Landau--Shepp--Fernique theorem
  \citep[Lemma~2.10]{dudley2014uniform}.
  The lower bound on $E$ follows as
  $\sup_t |X_t - \mu_t| / \sigma_t
  \geq |X_{t_0} - \mu_{t_0}| / \sigma_{t_0}$
  for any fixed $t_0 \in T$.
  Write $Z \coloneqq \sup_{t \in T} X_t$ and
  let $T_n$ be a sequence of finite subsets of $T$ with
  $Z_n \coloneqq \max_{t \in T_n} X_t \to Z$ almost surely
  as $n \to \infty$. Fix $\xi \in \R$.
  Since $|Z_n - \xi| \to |Z - \xi|$ almost surely,
  there exists a countable subset $\cN_\xi$ of $[0, \infty)$
  such that for all $\eta \in [0, \infty) \setminus \cN_\xi$, we have
  $\lim_{n \to \infty} \P(|Z_n - \xi| \leq \eta)
  = \P(|Z - \xi| \leq \eta)$.
  Now, by \citet[Corollary~7]{belloni2024anti},
  \begin{equation*}
    \P(|Z_n - \xi| \leq \eta)
    = \P(\xi - \eta \leq Z_n \leq \xi + \eta)
    \leq
    \frac{2 \eta}{\sigma}
    \E \biggl(
      \max_{t \in T_n}
      \frac{|X_t - \mu_t|}{\sigma_t}
    \biggr)
    \leq
    \frac{2 \eta E}{\sigma},
  \end{equation*}
  for all $\eta \geq 0$.
  Therefore, $\P(|Z - \xi| \leq \eta) \leq 2 \eta E / \sigma$
  for all $\eta \in [0, \infty) \setminus \cN_\xi$.
  By right continuity of $\eta \mapsto \P(|Z - \xi| \leq \eta)$,
  this bound in fact holds for all $\eta \geq 0$.
  Since $\xi$ was arbitrary, we conclude that
  \begin{equation*}
    \sup_{\xi \in \R}
    \P\Bigl(\xi <
      \sup_{t \in T} X_t
      \leq \xi + \eta
    \Bigr)
    \leq
    \frac{\eta E}{\sigma}.
  \end{equation*}
  The bound for $\sup_{t \in T} |X_t|$ follows by applying
  the preceding argument to the process $(\pm X_t : t \in T)$,
  which satisfies the hypotheses with the same $E$ and $\sigma$.
\end{proof}

\subsection{Proofs for Appendix~\ref{sec:gaussian_approx_locpol}}

\begin{lemma}[Gaussian approximation of the linear term]%
  \label{lem:locpol_gaussian}
  Grant the setup of Theorem~\ref{thm:locpol_gaussian}.
  There exists $C > 0$ depending only on
  $\Clp, D_W, D_V, r$ and $q$
  such that for each $\eta > 0$, there is
  a probability space supporting
  stochastic processes $S_{n,h}'$ and $T_{n,h}'$, where
  $S_{n,h}' \equallaw S_{n,h}$ and
  $T_{n,h}' \equallaw T_{n,h}$, with
  \begin{align*}
    &\P\biggl\{
      \sup_{(w, v, j) \in \cW \times \cV \times [D_q]}
      \sqrt{n h^{D_W}}
      \bigl| e_j^\T S'_{n,h}(w,v) - e_j^\T T'_{n,h}(w, v) \bigr| \\
      &\qquad\qquad>
      C \eta
      \Bigl\{
        \bigl(n h^{D_W}\bigr)^{1 - k}
        h^{- 2 D_w}
      \Bigr\}^{\frac{1}{2k + 2 + k D_w + 2 k D_v}}
      \sqrt{\log 2n}
    \biggr\}
    \leq
    \frac{1}{\eta^{k + 1}}.
  \end{align*}
\end{lemma}

\begin{proof}[Proof of Lemma~\ref{lem:locpol_gaussian}]
  Throughout the proof,
  $C > 0$ depends only on the stated quantities.
  Fix $\kappa = 2(k+1)$ so $n^{-\kappa} \leq
  \eta^{-(k+1)}$ when $\eta \leq n^2$.
  Let $\delta_w$, $\delta_v$, $\pi$, $\cW_\delta$,
  $\cM_\delta$, $\cM$, $d_w$, $d_v$, $d$,
  $f_{w,v,j}$, $\cF$, $\tilde\cF$ and $\check\cF$
  be as in the proof of Lemma~\ref{lem:locpol_bootstrap}.
  Assume for now that $(n h^{D_W})^{k-1} h^{2D_w}>1$.
  We apply Theorem~\ref{thm:gaussian_process} with $\gamma = 1$.
  For $i \in [n]$, let $Y_i$
  be independent zero-mean Gaussian processes indexed on $\cM$ with
  \begin{equation*}
    \E\bigl\{Y_i(w, v, j) Y_i(w', v', j')\bigr\}
    = \E\bigl\{f_{w,v,j}(W_i, V_i) f_{w',v',j'}(W_i, V_i)\bigr\}
  \end{equation*}
  for each $(w,v,j), (w',v',j') \in \cM$.
  Using $\E\bigl\{f_{w,v,j}(W_i,V_i)\bigr\}=0$ and
  the covariance identity above, we have
  $\delta_{\alpha,i}=0$ for $|\alpha|\in\{1,2\}$.
  For $k=3$, the condition
  $\E\bigl\{ \tilde g(v)_i \tilde g(v')_i
  \tilde g(v'')_i \mid W_i\bigr\} = 0$ gives
  $\E\bigl\{\tilde\cF(W_i,V_i)^\alpha\bigr\}=0$ for $|\alpha|=3$ also.
  Thus, $\Delta=\Gamma_\eta=0$.
  By the proof of Lemma~\ref{lem:locpol_bootstrap},
  \begin{align*}
    \bigl|f_{w, v, j}(W_i, V_i)\bigr|
    &\leq
    \frac{C}{n h^{D_W}}
    \bigl(|g(V_i, v)| + 1\bigr)
    \I\bigl\{\|W_i - w\|_\infty \leq h\bigr\}, \\
    \E\bigl\{
      f_{w, v, j}(W_i, V_i)^2
    \bigr\}
    &\leq
    \frac{1}{n^2}
    \E\bigl\{
      K_h(W_i - w)^2
      g(V_i, v)^2
    \bigr\}
    \leq
    \frac{C}{n^2 h^{D_W}}.
  \end{align*}
  Define $F(W, V)$ as in the proof of Lemma~\ref{lem:locpol_bootstrap}, so
  $\vvvert F(W_i, V_i) \vvvert_{\omega_r} \leq C n^{-1} h^{-D_W}$
  for each $i \in [n]$.
  For $i \in [n]$,
  as in the proof of Lemma~\ref{lem:locpol_bootstrap},
  \begin{align*}
    \beta_{X,k,1,i}
    &\coloneqq
    \frac{1}{\sqrt{k!}}
    \E \Bigl(
      \|\tilde \cF(W_i, V_i)\|_2^k
      \|\tilde \cF(W_i, V_i)\|_\infty
    \Bigr) \\
    &\,\leq
    \frac{C}{(n h^{D_W})^{k + 1}}
    \biggl( \frac{h}{\delta_w} \biggr)^{k D_w/2}
    \biggl( \frac{1}{\delta_v} \biggr)^{k D_v/2}
    h^{D_W - D_w}.
  \end{align*}
  Write $\tilde \cF(Y_i) \coloneqq
  \bigl(Y_i(w,v,j): (w,v,j) \in \cM_\delta\bigr)$.
  By Lemma~\ref{lem:gaussian_norm},
  if $\delta_w \land \delta_v \geq n^{-C}$,
  \begin{align*}
    \beta_{Y,k,1,i}
    &\coloneqq
    \frac{1}{\sqrt{k!}}
    \E \Bigl(
      \|\tilde \cF(Y_i)\|_2^k
      \|\tilde \cF(Y_i)\|_\infty
    \Bigr)
    \leq
    \frac{C (d_w d_v)^{k/2} (\log 2 d_w d_v)^{1/2}}
    {(n h^{D_W / 2})^{k + 1}} \\
    &\,\leq
    \frac{C \delta_w^{-k D_w/2} \delta_v^{-k D_v/2}
    (\log 2 n)^{1/2}}
    {(n h^{D_W / 2})^{k + 1}} \\
    &\,=
    \frac{C}{(n h^{D_W})^{k + 1}}
    \biggl( \frac{h}{\delta_w} \biggr)^{k D_w/2}
    \biggl( \frac{1}{\delta_v} \biggr)^{k D_v/2}
    h^{k (D_W - D_w) / 2}
    h^{D_W / 2}
    (\log 2 n)^{1/2} \\
    &\,\leq
    \frac{C}{(n h^{D_W})^{k + 1}}
    \biggl( \frac{h}{\delta_w} \biggr)^{k D_w/2}
    \biggl( \frac{1}{\delta_v} \biggr)^{k D_v/2}
    h^{D_W - D_w}
    \bigl( h \log 2 n \bigr)^{1/2} \\
    &\,\leq
    \frac{C}{(n h^{D_W})^{k + 1}}
    \biggl( \frac{h}{\delta_w} \biggr)^{k D_w/2}
    \biggl( \frac{1}{\delta_v} \biggr)^{k D_v/2}
    h^{D_W - D_w}.
  \end{align*}
  It follows that
  \begin{equation*}
    \bar\beta_{k,1}
    \coloneqq
    \frac{1}{n}
    \sum_{i=1}^n
    \bigl(\beta_{X,k,1,i} + \beta_{Y,k,1,i}\bigr)
    \leq
    \frac{C}{(n h^{D_W})^{k+1}}
    \biggl( \frac{h}{\delta_w} \biggr)^{k D_w/2}
    \biggl( \frac{1}{\delta_v} \biggr)^{k D_v/2}
    h^{D_W - D_w}.
  \end{equation*}
  The term $R(t)$ appearing in Theorem~\ref{thm:gaussian_process}
  has already been controlled in the
  proof of Lemma~\ref{lem:locpol_bootstrap};
  for $t \leq C \sqrt{\log 2n}$, we have
  \begin{equation*}
    R(t)
    \leq
    C \Biggl\{
      \sqrt{\frac{\log 2n + t^2}{n h^{D_W}}}
      \biggl(
        \frac{\delta_w}{h}
        + \sqrt{\delta_v}
      \biggr)
      + \frac{(\log 2n)^{1 + 1/r}}{n h^{D_W}}
    \Biggr\}.
  \end{equation*}
  For $(w,v,j) \in \cM$, let
  $S_n(w,v,j) \coloneqq e_j^\T S_{n,h}(w, v)
  = \sum_{i=1}^n f_{w,v,j}(W_i, V_i)$
  and, similarly, define $T_n(w, v, j) \coloneqq
  e_j^\T T_{n,h}(w,v) \coloneqq
  \sum_{i=1}^n Y_i(w, v, j)$.
  By Theorem~\ref{thm:gaussian_process},
  \begin{align*}
    &\P\biggl\{
      \sup_{(w, v, j) \in \cM}
      |S_n'(w, v, j) - T_n'(w, v, j)|
      > 3 \eta \\
      &\qquad\qquad\qquad+ C
      \sqrt{\frac{\log 2n + t^2}{n h^{D_W}}}
      \biggl(
        \frac{\delta_w}{h}
        + \sqrt{\delta_v}
      \biggr)
      + \frac{C (\log 2n)^{1 + 1/r}}{n h^{D_W}}
    \biggr\} \\
    &\qquad\leq
    4 e^{-t^2 / 2}
    + \frac{1}{\eta^{k + 1}}
    \Bigl\{\sqrt{\log 2n} + t\Bigr\}^k
    \frac{C h^{- D_w}}{(n h^{D_W})^k}
    \biggl( \frac{h}{\delta_w} \biggr)^{k D_w/2}
    \biggl( \frac{1}{\delta_v} \biggr)^{k D_v/2}
    + \frac{1}{n^\kappa}.
  \end{align*}
  We take $\delta_v \coloneqq \delta_w^2 / h^2$ so that
  \begin{align*}
    &\P\biggl\{
      \sup_{(w, v, j) \in \cM}
      |S_n'(w, v, j) - T_n'(w, v, j)|
      > 3 \eta
      \iftoggle{journal}{\\&\qquad\qquad}{}
      + C
      \sqrt{\frac{\log 2n + t^2}{n h^{D_W}}}
      \frac{\delta_w}{h}
      + \frac{C (\log 2n)^{1+1/r}}{n h^{D_W}}
    \biggr\} \\
    &\qquad\leq
    4 e^{-t^2 / 2}
    + \frac{1}{\eta^{k + 1}}
    \Bigl\{\sqrt{\log 2n} + t\Bigr\}^k
    \frac{C h^{- D_w}}{(n h^{D_W})^{k}}
    \biggl( \frac{h}{\delta_w} \biggr)^{k D_w/2 + k D_v}
    + \frac{1}{n^\kappa}.
  \end{align*}
  Next, we set
  $\delta_w / h \coloneqq
  h^{-\frac{2 D_w}{2k + 2 + k D_w + 2k D_v}}
  (n h^{D_W})^{\frac{1 - k}{2k + 2 + k D_w + 2 k D_v}}$,
  giving
  \begin{align*}
    &\P\biggl\{
      \sup_{(w, v, j) \in \cM}
      |S_n'(w, v, j) - T_n'(w, v, j)|
      > 3 \eta \\
      &+ C
      \sqrt{\log 2n + t^2}
      h^{-\frac{2 D_w}{2k + 2 + k D_w + 2k D_v}}
      (n h^{D_W})^{\frac{- 4k - k D_w - 2 k D_v}
      {4k + 4 + 2 k D_w + 4 k D_v}}
      + \frac{C (\log 2n)^{1+1/r}}{n h^{D_W}}
    \biggr\} \\
    &\qquad\leq
    4 e^{-t^2 / 2}
    + \frac{1}{n^\kappa}
    \iftoggle{journal}{\\&\qquad\qquad}{}
    + \frac{1}{\eta^{k + 1}}
    \Bigl\{\sqrt{\log 2n} + t\Bigr\}^k
    C
    h^{- \frac{2 D_w (k + 1)}
    {2k + 2 + k D_w + 2k D_v}}
    (n h^{D_W})^{\frac{(k + 1)(- 4k - k D_w - 2 k D_v)}
    {4k + 4 + 2 k D_w + 4 k D_v}}.
  \end{align*}
  Rescaling $\eta$ yields
  \begin{align*}
    &\P\biggl\{
      \sup_{(w, v, j) \in \cM}
      |S_n'(w, v, j) - T_n'(w, v, j)|
      > C
      h^{- \frac{2 D_w }
      {2k + 2 + k D_w + 2k D_v}}
      (n h^{D_W})^{\frac{- 4k - k D_w - 2 k D_v}
      {4k + 4 + 2 k D_w + 4 k D_v}} \\
      &\qquad\times
      \sqrt{\log 2n + t^2 + \eta^2}
      + \frac{C (\log 2n)^{1+1/r}}{n h^{D_W}}
    \biggr\}
    \leq
    4 e^{-t^2 / 2}
    + \frac{\bigl(\sqrt{\log 2n} + t\bigr)^k}{\eta^{k + 1}}
    + \frac{1}{n^\kappa}.
  \end{align*}
  Taking $t = \sqrt{2 \kappa \log 2n}$ and absorbing
  $C(\log 2n)^{1+1/r}/(n h^{D_W})$ into the first term, we have
  \begin{align*}
    &\P\biggl\{
      \sup_{(w, v, j) \in \cM}
      |S_n'(w, v, j) - T_n'(w, v, j)|
      > C
      h^{- \frac{2 D_w }
      {2k + 2 + k D_w + 2k D_v}}
      (n h^{D_W})^{\frac{- 4k - k D_w - 2 k D_v}
      {4k + 4 + 2 k D_w + 4 k D_v}} \\
      &\qquad\qquad\times
      \sqrt{\kappa \log 2n + \eta^2}
    \biggr\}
    \leq
    \frac{(\kappa \log 2n)^{k/2}}{\eta^{k + 1}}
    + \frac{C}{n^\kappa}.
  \end{align*}
  The $C/n^\kappa$ term is removed by arguing as
  in the proof of Lemma~\ref{lem:locpol_bootstrap},
  using sixteenth moments instead of eighth moments.
  For the Gaussian process, Dudley's theorem
  \citep[Corollary~2.2.8]{van1996weak}
  and Borell's inequality
  \citep[Proposition~A.2.1]{van1996weak} give
  \begin{equation*}
    \E\Biggl\{
      \biggl(\sqrt{n h^{D_W}}
        \sup_{(w,v)\in\cW \times \cV}\|T_{n,h}'(w,v)\|_2
      \biggr)^{16}
    \Biggr\}
    \leq C(\log 2n)^8.
  \end{equation*}
  For $\eta > n^2$, we therefore obtain
  a probability bound of $C n^{20}/\eta^{16} \leq C/\eta^{k+1}$.
  Thus,
  \begin{align*}
    &\P\biggl\{
      \sup_{(w, v, j) \in \cM}
      |S_n'(w, v, j) - T_n'(w, v, j)| \\
      &\qquad\qquad>
      C \eta
      h^{- \frac{2 D_w }
      {2k + 2 + k D_w + 2k D_v}}
      (n h^{D_W})^{\frac{- 4k - k D_w - 2 k D_v}
      {4k + 4 + 2 k D_w + 4 k D_v}}
      \sqrt{\log 2n}
    \biggr\}
    \leq
    \frac{1}{\eta^{k + 1}}.
  \end{align*}
  The result follows by rescaling the processes by
  $\sqrt{n h^{D_W}}$.

  It remains to handle the case where
  $(n h^{D_W})^{k-1} h^{2D_w} \leq 1$.
  By independent coupling, the triangle inequality,
  Lemma~\ref{lem:locpol_linear}, and Dudley's and Borell's
  inequalities for $T_{n,h}$, taking
  $t = \sqrt{2\kappa\log 2n}$ gives
  \begin{equation*}
    \P\Biggl\{
      \sqrt{n h^{D_W}}
      \sup_{(w,v)\in\cW\times\cV}
      \|S_{n,h}'(w,v)-T_{n,h}'(w,v)\|_2
      > C\sqrt{\log 2n}
    \Biggr\}
    \leq \frac{C}{n^\kappa}.
  \end{equation*}
  For $1 \leq \eta \leq n^2$, this gives the stated bound.
  For $\eta > n^2$, use the sixteenth-moment argument above;
  for $\eta < 1$, the bound is trivial.
\end{proof}

\begin{proof}[Proof of Theorem~\ref{thm:locpol_gaussian}]
  Throughout the proof, $C > 0$ depends only on the stated quantities.
  Fix $\kappa = 2(k+1)$.
  As in the proof of Theorem~\ref{thm:locpol_bootstrap},
  by~\eqref{eq:locpol_decomp_coupling},
  \eqref{eq:locpol_decomp_errors},
  \eqref{eq:locpol_decomp_bias} and
  Lemmas~\ref{lem:locpol_Hhat},
  \ref{lem:locpol_linear} and~\ref{lem:locpol_bias},
  with probability at least $1-C/n^\kappa-2e^{-n h^{D_W}/C}$,
  \begin{equation*}
    \sup_{(w,v)\in\cW\times\cV}
    \bigl|
    \hat\theta_{n,h}(w,v)-\theta(w,v)
    -e_1^\T H_h(w)^{-1}S_{n,h}(w,v)
    \bigr|
    \leq \frac{C\log 2n}{n h^{D_W}}+Ch^\zeta.
  \end{equation*}
  Take $S_{n,h}'$ and $T_{n,h}'$ as constructed in
  Lemma~\ref{lem:locpol_gaussian}.
  By the Vorob'ev--Berkes--Philipp theorem,
  extend the probability space
  to support $\bX_n'$ with
  $(\bX_n', S_{n,h}', \hat\theta_{n,h}')
  \equallaw (\bX_n, S_{n,h}, \hat\theta_{n,h})$.
  Define $Z_{n,h}'(w,v) \coloneqq e_1^\T H_h(w)^{-1}T_{n,h}'(w,v)$.
  Lemma~\ref{lem:locpol_H} and the triangle inequality
  now provide the stated bound, up to additional
  $\log 2n/\sqrt{n h^{D_W}}$ and $C/n^\kappa$ terms.
  The former is negligible for $\eta\geq1$, since
  $(k-1)/(2k+2+kD_w+2kD_v)\leq1/4$
  and $n h^{D_W}\geq n^{1/\Clp}$;
  for $\eta<1$, the desired result is trivial.
  For the latter, note that
  $\eta\leq n^2$ implies $n^{-\kappa}\leq\eta^{-(k+1)}$.
  If $\eta>n^2$, restrict to the event of
  Lemma~\ref{lem:locpol_Hhat} on which
  \begin{align*}
    \sup_{w\in\cW}
    \bigl\|\hat H_h'(w)\bigr\|_{2\to2}
    &\leq C,
    &\sup_{w\in\cW}
    \bigl\|\hat H_h'(w)^{-1}\bigr\|_{2\to2}
    &\leq C,
  \end{align*}
  where $\hat H_h'$ is computed using $\bX_n'$;
  this event has probability at least $1-2e^{-n h^{D_W}/C}$.
  The sixteenth moment argument in the proof of
  Lemma~\ref{lem:locpol_gaussian} then gives a probability
  bound of $C n^{20}/\eta^{16} \leq C/\eta^{k+1}$.
  Increasing $C$ yields the result.
\end{proof}

\begin{proof}[Proof of Theorem~\ref{thm:locpol_gaussian_ks}]
  We write
  \begin{align}
    \nonumber
    &
    \sup_{\xi \in \R}
    \Biggl|
    \P\Biggl(
      \sup_{(w, v) \in \cW \times \cV}
      \,
      \Biggl|
      \frac{\hat\theta_{n,h}(w,v)-\theta(w,v)}
      {\sqrt{\smash[b]{\hat \Sigma_{n,h}(w,v,w,v)}}}
      \Biggr| \leq \xi
    \Biggr)
    \nonumber
    \iftoggle{journal}{\\ \nonumber &\qquad\qquad}{}
    -
    \P^*\Biggl(
      \sup_{(w, v) \in \cW \times \cV}
      \,
      \Biggl|
      \frac{\hat Z_{n,h}(w,v)}
      {\sqrt{\smash[b]{\hat \Sigma_{n,h}(w,v,w,v)}}}
      \Biggr| \leq \xi
    \Biggr)
    \Biggr| \\
    \label{eq:locpol_ks_gaussian_comparison}
    &\ \leq
    \sup_{\xi \in \R}
    \Biggl|
    \P^*\Biggl(
      \sup_{(w, v) \in \cW \times \cV}
      \,
      \Biggl|
      \frac{\hat Z_{n,h}(w,v)}
      {\sqrt{\smash[b]{\hat \Sigma_{n,h}(w,v,w,v)}}}
      \Biggr| \leq \xi
    \Biggr)
    \iftoggle{journal}{\\ \nonumber &\qquad\qquad}{}
    -
    \P\Biggl(
      \sup_{(w, v) \in \cW \times \cV}
      \,
      \Biggl|
      \frac{Z_{n,h}(w,v)}
      {\sqrt{\Sigma_h(w,v,w,v)}}
      \Biggr| \leq \xi
    \Biggr)
    \Biggr|
    \\
    \label{eq:locpol_ks_anti-concentration}
    &\qquad\qquad+
    \sup_{\xi \in \R}
    \P\Biggl(
      \xi <
      \sup_{(w, v) \in \cW \times \cV}
      \,
      \Biggl|
      \frac{Z_{n,h}(w,v)}{\sqrt{\Sigma_h(w,v,w,v)}}
      \Biggr| \leq \xi + 3\eta
    \Biggr) \\
    \label{eq:locpol_ks_coupling}
    &\qquad\qquad+
    \P\Biggl(
      \sup_{(w, v) \in \cW \times \cV}
      \,
      \Biggl|
      \frac{\hat\theta_{n,h}'(w,v) - \theta(w,v) - Z_{n,h}'(w, v)}
      {\sqrt{\Sigma_h(w,v,w,v)}}
      \Biggr| > \eta
    \Biggr) \\
    \label{eq:locpol_estimator_hat}
    &\qquad\qquad+
    \P\Biggl(
      \sup_{(w, v) \in \cW \times \cV}
      \,
      \Biggl|
      \frac{\hat\theta_{n,h}(w,v) - \theta(w,v)}
      {\sqrt{\smash[b]{\hat \Sigma_{n,h}(w,v,w,v)}}}
      - \frac{\hat\theta_{n,h}(w,v) - \theta(w,v)}
      {\sqrt{\Sigma_h(w,v,w,v)}}
      \Biggr| > \eta
    \Biggr).
  \end{align}
  We begin by controlling \eqref{eq:locpol_ks_gaussian_comparison}.
  Let $\cM_\delta$ and $\pi$ be as in the proof of
  Lemma~\ref{lem:locpol_bootstrap}, with
  $\delta_w/h=n^{-B}$ and $\delta_v=n^{-2B}$ for fixed $B > 0$.
  Write $(w',v',j)\coloneqq\pi(w,v,j)$.
  By the fluctuation bounds in the proofs of
  Lemmas~\ref{lem:locpol_bootstrap} and~\ref{lem:locpol_gaussian},
  together with Lemmas~\ref{lem:locpol_Hhat}, \ref{lem:locpol_R}
  and~\ref{lem:Sigma_hat},
  \begin{align}
    \nonumber
    &\sup_{(w,v,j)\in\cM}
    \Biggl|
    \frac{Z_{n,h}(w,v)}{\sqrt{\Sigma_h(w,v,w,v)}}
    -\frac{Z_{n,h}(w',v')}{\sqrt{\Sigma_h(w',v',w',v')}}
    \Biggr|
    \lesssim_\P \frac{1}{n^c}, \\
    \label{eq:locpol_ks_gaussian_fluctuations}
    &\sup_{(w,v,j)\in\cM}
    \Biggl|
    \frac{\hat Z_{n,h}(w,v)}
    {\sqrt{\smash[b]{\hat \Sigma_{n,h}(w,v,w,v)}}}
    -\frac{\hat Z_{n,h}(w',v')}
    {\sqrt{\smash[b]{\hat \Sigma_{n,h}(w',v',w',v')}}}
    \Biggr|
    \lesssim_\P \frac{1}{n^c},
  \end{align}
  for some fixed $c > 0$.
  Next apply \citet[Lemma~3.1]{chernozhukov2013gaussian} to
  the vectors
  \begin{equation*}
    \Biggl(
      \pm
      \frac{Z_{n,h}(w,v)}{\sqrt{\Sigma_h(w,v,w,v)}}:
      (w,v,j)\in\cM_\delta
    \Biggr),
    \qquad
    \Biggl(
      \pm
      \frac{\hat Z_{n,h}(w,v)}
      {\sqrt{\smash[b]{\hat \Sigma_{n,h}(w,v,w,v)}}}:
      (w,v,j)\in\cM_\delta
    \Biggr),
  \end{equation*}
  to see that, by Lemma~\ref{lem:Sigma_hat},
  \begin{align}
    \nonumber
    &\sup_{\xi\in\R}
    \Biggl|
    \P^*\Biggl(
      \max_{(w,v,j)\in\cM_\delta}
      \frac{|\hat Z_{n,h}(w,v)|}
      {\sqrt{\smash[b]{\hat \Sigma_{n,h}(w,v,w,v)}}}
      \leq\xi
    \Biggr)
    -
    \P\Biggl(
      \max_{(w,v,j)\in\cM_\delta}
      \frac{|Z_{n,h}(w,v)|}
      {\sqrt{\Sigma_h(w,v,w,v)}}
      \leq\xi
    \Biggr)
    \Biggr| \\
    \label{eq:locpol_ks_gaussian_comparison_finite}
    &\qquad\lesssim_\P
    n^{-c/3}(\log 2n)^{2/3}
    \to 0.
  \end{align}
  By \eqref{eq:locpol_ks_gaussian_fluctuations},
  \eqref{eq:locpol_ks_gaussian_comparison_finite},
  \eqref{eq:locpol_ks_gaussian_anti_concentration}
  and dominated convergence,
  the expectation of \eqref{eq:locpol_ks_gaussian_comparison}
  converges to zero.

  For the term in \eqref{eq:locpol_ks_anti-concentration},
  we use the fact that
  $(w, v) \mapsto Z_{n,h}(w,v) \big/ \sqrt{\Sigma_h(w,v,w,v)}$
  is a separable zero-mean Gaussian process on $\cW \times \cV$
  with pointwise variance one.
  By the covering bounds in the proof of
  Lemma~\ref{lem:locpol_linear} and Dudley's theorem
  \citep[Corollary~2.2.8]{van1996weak},
  \begin{equation*}
    \E\biggl\{\sup_{(w,v)\in\cW\times\cV}
    \|T_{n,h}(w,v)\|_2\biggr\}
    \leq C\sqrt{\frac{\log 2n}{n h^{D_W}}}.
  \end{equation*}
  Since $Z_{n,h}=e_1^\T H_h^{-1}T_{n,h}$ and by
  Lemmas~\ref{lem:locpol_H} and~\ref{lem:Sigma_hat},
  we may take $E \leq C\sqrt{\log 2n}$ in
  Lemma~\ref{lem:anti-concentration_gp} to deduce that
  \begin{equation}
    \label{eq:locpol_ks_gaussian_anti_concentration}
    \sup_{\xi \in \R}
    \P\Biggl(
      \xi <
      \sup_{(w, v) \in \cW \times \cV}
      \,
      \Biggl|
      \frac{Z_{n,h}(w,v)}{\sqrt{\Sigma_h(w,v,w,v)}}
      \Biggr| \leq \xi + 3 \eta
    \Biggr)
    \leq C \eta \sqrt{\log 2n}
    \to 0,
  \end{equation}
  whenever $\eta \sqrt{\log 2n} \to 0$.
  For \eqref{eq:locpol_ks_coupling},
  by Theorem~\ref{thm:locpol_gaussian} and Lemma~\ref{lem:Sigma_hat},
  \begin{equation*}
    \P\Biggl(
      \sup_{(w, v) \in \cW \times \cV}
      \,
      \Biggl|
      \frac{\hat\theta_{n,h}'(w,v) - \theta(w,v) - Z_{n,h}'(w, v)}
      {\sqrt{\Sigma_h(w,v,w,v)}}
      \Biggr| > \eta
    \Biggr)
    \to 0,
  \end{equation*}
  whenever $\eta \gg \bigl\{ \bigl(n h^{D_W}\bigr)^{1 - k} h^{- 2 D_w}
  \bigr\}^{1/(2k + 2 + k D_w + 2 k D_v)} \sqrt{\log 2n}
  + \sqrt{n h^{D_W}} h^\zeta$.
  For \eqref{eq:locpol_estimator_hat},
  using $1/\sqrt{a} - 1/\sqrt{b} = (b - a) / (a \sqrt b + b \sqrt a)$,
  by Lemmas~\ref{lem:locpol_convergence} and~\ref{lem:Sigma_hat},
  \begin{equation*}
    \Biggl|
    \frac{\hat\theta_{n,h}(w,v) - \theta(w,v)}
    {\sqrt{\smash[b]{\hat \Sigma_{n,h}(w,v,w,v)}}}
    - \frac{\hat\theta_{n,h}(w,v) - \theta(w,v)}
    {\sqrt{\Sigma_h(w,v,w,v)}}
    \Biggr|
    \lesssim_\P
    \frac{\log 2n}{\sqrt{n h^{D_W}}}
    + \sqrt{n h^{D_W}} h^{2 \zeta}.
  \end{equation*}
  Therefore, the term in \eqref{eq:locpol_estimator_hat}
  converges to zero whenever
  $\eta \gg \log 2n / \sqrt{n h^{D_W}} + \sqrt{n h^{D_W}} h^{2\zeta}$.
  We conclude that a suitable sequence of $\eta$ values exists if
  \begin{equation*}
    \Bigl\{ \bigl(n h^{D_W}\bigr)^{1 - k} h^{- 2 D_w}
    \Bigr\}^{1/(2k + 2 + k D_w + 2 k D_v)} \log 2n
    +  h^\zeta \sqrt{n h^{D_W} \log 2n}
    \to 0,
  \end{equation*}
  which is implied by the conditions of the theorem.
\end{proof}

\subsection{Details for Appendix~\ref{sec:literature_empirical_process}}
\label{sec:literature_details_empirical_process}

We provide details for the results of
Appendix~\ref{sec:literature_empirical_process}
and Table~\ref{tab:empirical_process_literature}.

\subsubsection{Yurinskii-type couplings}

\citet{dudley1983invariance} approximated normalized empirical processes
indexed by classes of sets using Kiefer processes;
for VC classes with VC index $V$ and i.i.d.\ data,
their Theorem~7.4 gives a rate of $n^{-\lambda}$
for every $\lambda < 1/(2700 V)$.
The proof relies on Strassen's theorem and their Lemma~2.12,
itself a special case of an extension by
\citet{dehling1983limit} of a theorem due to \citet{yurinskii1978error}.

\citet[Theorem~5.10a']{massart1986rates} assumes i.i.d.\ data,
an envelope $F$ with $\E\bigl\{F(X_i)^{2+\delta}\bigr\} < \infty$
for some $\delta \in (0, 1]$, and
uniform $L_2$-packing numbers of order
$\varepsilon^{-u}$ relative to $F$.
He constructed a Brownian bridge approximation
to the normalized empirical process with a rate of
$(\log 2n)^{(1+u)/2} n^{-\delta/\{8+u(4+\delta)\}}$.
If the envelope has finite third moment,
taking $\delta = 1$ gives a rate of
$n^{-1/(8+5u)}$ up to logarithmic factors.

\citet[Theorem~2]{le1988prokhorov} constructed a
Gaussian coupling for the normalized empirical process
of independent observations
over a VC class of sets with VC exponent $v$.
His covering estimate corresponds to uniform $L_2$-covering numbers
of order $\varepsilon^{-u}$ with $u = 4v$;
this yields a rate of $n^{-1/(8+5u)}$ up to logarithmic factors.
The sharper covering bound discussed in his Section~6 replaces $u = 4v$
by $u = 2v$ and gives a correspondingly improved rate.

\citet[Theorem~1]{rio1996vitesses} constructed Gaussian approximations
to the normalized empirical process (indexed by axis-aligned boxes)
of i.i.d.\ variables in $[0, 1]^D$ with
bounded Lebesgue density, attaining an almost sure rate of
$n^{-1/12}(\log n)^{(5D+1)/6}$.
His proof uses the Gaussian approximation of
\citet[Theorem~1.1]{zaitsev1987gaussian} and Strassen's theorem.

\citet[Proposition~1]{berthet2006revisiting} assumed i.i.d.\ data and a
bounded, pointwise measurable class with uniform $L_2$ covering numbers of
order $\varepsilon^{-u}$ with respect to a constant envelope.
They obtained a normalized Gaussian approximation rate of
$n^{-1/(2+5u)} (\log n)^{(4+5u)/(4+10u)}$.
Their Theorem~1 gives a version controlling all partial sums, with a
slower rate,
while their Proposition~2 treats bracketing entropy, attaining
logarithmic rates.

\citet[Theorem~3.4]{settati2009gaussian} assumed i.i.d.\ data and a
bounded, pointwise measurable class satisfying his random entropy
condition (F.3) with $\varphi_n = \log n$ and $\beta_n = n$.
Writing $u$ for the coefficient in his population $L_2(P)$ entropy
bound (F.6), his normalized Gaussian approximation rate is
$n^{-1/(2+5u)}(\log n)^{5u/(4+10u)}$.

For i.i.d.\ data, a uniformly bounded class, and uniform
covering numbers of order
$\varepsilon^{-u}$, the discussion after Proposition~3.1 in
\citet{cattaneo2025yurinskii} gives the rate $n^{-1/(6+2u)} \sqrt{\log n}$,
up to any $R_n \to \infty$.
For i.n.i.d.\ data, we center each summand by its own expectation and
apply the proof of their Proposition~3.1 on a $\delta$-net in
$L_2(\bar P_X)$, where $\bar P_X = n^{-1}\sum_{i=1}^n P_{X,i}$.
Independence makes the covariance discrepancy on the net zero,
and boundedness gives the same finite-dimensional moment bounds.
At $\delta = n^{-1/(6+2u)}$, symmetrization and the uniform entropy bound
give a discretization error of order $\delta\sqrt{n\log n}$ in probability
for the unnormalized sum.
The Gaussian process $T$ with the covariance of the sum satisfies
$\Var\{T(f)-T(f')\} \leq n\|f-f'\|_{L_2(\bar P_X)}^2$,
so Dudley's inequality gives the same order.

Let $X$ and $Y$ be $\R^D$-valued random vectors
with laws $P$ and $Q$, respectively.
Suppose that
$\E\bigl\{(\I_A(X)-P(A))^\ell\bigr\} = \E\bigl\{(\I_A(Y)-Q(A))^\ell\bigr\}$
for every $A\in\cA$ and $\ell \in \{2,3\}$,
where $\cA = \bigl\{(-\infty, x_1] \times \cdots \times (-\infty, x_D):
x_1, \ldots, x_D \in \R\bigr\}$.
Fix $A \in \cA$ with $p \coloneqq \P(X \in A) \in (0, 1)$.
Then the second-order matching condition yields
$\E\bigl\{(\I_A(Y)-Q(A))^2\bigr\} = \E\bigl\{(\I_A(X)-P(A))^2\bigr\}
= p (1 - p)$ and so $Q(A) \in \{p, 1-p\}$.
The third-order condition similarly implies
$\E\bigl\{(\I_A(Y)-Q(A))^3\bigr\} = p (1 - p) (1 - 2p)$,
so in fact $Q(A) = P(A)$.
This also holds if $P(A) = 0$ by considering $A \cap A'$
with $P(A') \in (0, 1)$; likewise if $P(A) = 1$.
If no such $A'$ exists, both centered indicator processes are identically zero.
Thus the centered indicator processes have the same law in either case.

\subsubsection{KMT-type couplings}

The claimed rates of
\citet{komlos1975approximation},
\citet{gine2004kernel},
\citet{gine2010confidence},
\citet{cattaneo2024uniform},
\citet{cattaneo2025strong},
\citet{tusnady1977remark},
\citet{castelle1998strong} and
\citet{borisov1982approximation}
follow directly from the cited results.

For i.i.d.\ data uniform on $[0, 1]^D$ with $D \geq 2$,
\citet[Theorem~1]{massart1989strong} considered a VC class $\cA$
of sets satisfying the boundary regularity condition
\begin{equation*}
  \sup_{A \in \cA}
  \P\bigl(\|X_1 - \partial A\|_2 < \varepsilon\bigr)
  \leq C\varepsilon,
\end{equation*}
for all $\varepsilon \in (0,1]$.
His bound gives a Gaussian approximation
rate of $n^{-1/(2D)}(\log n)^{3/2}$.

Under the same boundary regularity condition,
\citet[Section~11b]{koltchinskii1994komlos}
showed that the indicator class $\{\I_A: A \in \cA\}$
satisfies his condition (11.3) with $\beta = 1/D$.
His Theorem~11.3 then yields the same rate.

For i.i.d.\ data with a strictly positive $C^1$ density on $[0, 1]^D$,
where $D \geq 2$, \citet[Theorem~1.1]{rio1994local}
gave a normalized Gaussian approximation rate of
$n^{-1/(2D)} \sqrt{\log n} + n^{-1/2}(\log n)^{3/2}$
for uniformly bounded VC-type classes with uniformly bounded variation,
by taking $t = C \log n$ in his bound.
Uniform $L_2$ covering numbers of order $\varepsilon^{-u}$
imply his uniform $L_1$ packing condition (1.1);
the proof proceeds via the Rosenblatt transformation,
reducing the problem to uniform data
while preserving uniform boundedness
and bounds on the entropy and total variation.

For i.i.d.\ data with a strictly positive $C^1$ density on $[0,1]^D$,
\citet[Corollary~2]{cattaneo2025strong} gave normalized Gaussian
approximation rates of $n^{-1/2}\log n$ for $D = 2$ and
$n^{-1/D}\sqrt{\log n}$ for $D \geq 3$, for uniformly bounded and
uniformly Lipschitz VC-type classes.
The uniform Lipschitz bound ensures total variation of order
$\ell^D$ on cubes of side length $\ell \leq 1$,
thus satisfying their local variation condition.
The proof is again via the Rosenblatt transformation.

\subsubsection{Wasserstein-based couplings}

For i.i.d.\ $\R^D$-valued data with arbitrary distribution
function $F$, \citet[Theorem~3.2]{dedecker2014strong} constructs
a Kiefer process $K$ with
\begin{equation*}
  \E\Biggl\{
    \sup_{s \in \R^D, \ k \in [n]}
    \Biggl|
    \sum_{i=1}^k \bigl(\I\{X_i \leq s\} - F(s)\bigr) - K(s, k)
    \Biggr|
  \Biggr\}
  \leq C n^{1/3} (\log 2n)^{(2D+3)/3},
\end{equation*}
where $C > 0$ depends only on $D$ and $F$.
Dividing by $\sqrt n$ gives the
approximation rate $n^{-1/6}$ up to logarithmic factors.
Their Theorem~3.1 gives a rate of $n^{1/p - 1/2}$, up to logarithmic factors,
for strictly stationary data with $\beta$-mixing coefficients of
order $n^{1-p}$ with $p \in (2, 3]$,
and their Theorem~2.1 shows that this rate cannot be
improved for continuous Gaussian approximations
in general under these conditions.

\citet[Theorem~2.1]{dedecker2013strong} assume
univariate strictly stationary sequences with dependence coefficients
$\beta_{2,X}(n) \lesssim n^{-1-\delta}$
for some $\delta > 0$. Taking $t = 1$ in their Equation (1.5) and
dividing by $\sqrt n$ gives, for the class $\cF$ of lower half-lines,
\begin{equation*}
  \sup_{f \in \cF}
  \frac{|S_n'(f)-T_n'(f)|}{\sqrt n}
  = O\bigl(n^{-\delta^2/\{22(\delta+2)^2\}}\bigr),
\end{equation*}
almost surely.

\subsubsection{Bootstrap couplings}

For i.i.d.\ data,
\citet[Theorem~3.2]{neumann1998strong} considered the kernel density
estimator $\hat f_h(x)$ with a Lipschitz second-order kernel
supported on $[-1,1]^D$, and the same estimator $\hat f_h^*(x)$ computed
from a nonparametric bootstrap sample drawn from a copy of the data.
Under further conditions on the densities,
it is possible to pair the two samples so that
\begin{align*}
  &\sup_{x\in\R^D}
  \Bigl|
  \hat f_h(x) - \E\bigl\{\hat f_h(x)\bigr\}
  - \hat f_h^*(x) + \E^*\bigl\{\hat f_h^*(x)\bigr\}
  \Bigr| \\
  &\qquad\lesssim_\P
  \frac{\log 2n}{\sqrt n}
  + \frac{\log 2n}{n h^D}
  + \sqrt{\frac{\log n}{n h^D}}
  \inf_{h' > 0} \biggl\{\sqrt{\mu_n(h')} + \frac{h'}{h}\biggr\},
\end{align*}
where $h' > 0$ is an auxiliary bandwidth and
$\mu_n(h') \coloneqq h^{\prime 2} +\sqrt{\log n / (n h^{\prime D})}$.
This yields a normalized bootstrap approximation rate of $(n h^D)^{-1/(D+4)}$
for sufficiently small $h$, up to logarithmic factors.
The rate $(n h^D)^{-1/(D+2)}$ quoted by
\citet[Remark~3.1ii]{chernozhukov2014gaussian}
follows from the bound printed in Neumann's Theorem~3.2;
the proof does not appear to justify the replacement of $D+4$ with $D+2$.

For univariate i.i.d.\ data,
\citet[Theorem~1.2]{csorgo2000approximation} approximate
the normalized bootstrap empirical process based on $m$
draws by Brownian bridges with error of order
$m^{-1/2} \log m + n^{-1/2} \log n$ and exponential tails.
Taking $m=n$ for the nonparametric bootstrap gives $n^{-1/2}\log n$.
 
\bibliographystyle{apalike}
\bibliography{references}

\end{document}